\documentclass[11pt]{amsart}

\usepackage{amsmath,yhmath}
\usepackage{amssymb}
\usepackage{amsfonts}
\usepackage{mathrsfs} 
\usepackage[mathscr]{eucal} 
\usepackage{bbm} 
\usepackage{mathtools}
\usepackage[all]{xy}
\usepackage{tikz}
\usepackage{tikz-cd}
\usepackage[inline]{enumitem}
\allowdisplaybreaks
\usepackage{stmaryrd} 
\usepackage{epigraph}

\numberwithin{equation}{section}
\usepackage{caption}
\usepackage{geometry}
\usepackage{parskip}
\usepackage{setspace}
\usepackage{tikz}
\tikzset{Equal/.style={-,double line with arrow={-,-}}} 
\usetikzlibrary{decorations.markings,arrows.meta}

\usepackage[colorlinks, linkcolor=black,citecolor=blue, urlcolor=cyan, unicode, bookmarksnumbered]{hyperref}

\usepackage{aliascnt}

\theoremstyle{plain}

\newtheorem{Thm}{Theorem}[section]

\newtheorem*{Thm*}{Theorem}

\newaliascnt{prop}{Thm}
\newtheorem{Prop}[prop]{Proposition}
\aliascntresetthe{prop}

\newaliascnt{lemma}{Thm}
\newtheorem{Lemma}[lemma]{Lemma}
\aliascntresetthe{lemma}

\newaliascnt{coro}{Thm}
\newtheorem{Coro}[coro]{Corollary}
\aliascntresetthe{coro}

\newtheorem{conjecture}[Thm]{Conjecture}

\newaliascnt{conjecture}{Thm}
\newtheorem{Conjecture}[conjecture]{Conjecture}
\aliascntresetthe{conjecture}

\newtheorem{Lemma*}{Lemma}
\newtheorem*{Prop*}{Proposition}
\newtheorem*{Coro*}{Corollary}

\theoremstyle{definition}

\newaliascnt{def}{Thm}
\newtheorem{Def}[def]{Definition}
\aliascntresetthe{def}

\newtheorem*{Def*}{Definition}

\newaliascnt{Eg}{Thm}
\newtheorem{eg}[Eg]{Example}
\aliascntresetthe{Eg}

\newtheorem*{eg*}{Example}

\theoremstyle{remark}

\newaliascnt{rmk}{Thm}
\newtheorem{RMK}[rmk]{Remark}
\aliascntresetthe{rmk}

\newtheorem*{RMK*}{Remark}

\newcommand{\bC} { {\mathbb{C}}}
\newcommand{\bF} { {\mathbb{F}}}
\newcommand{\bQ} { {\mathbb{Q}}}
\newcommand{\bZ} { {\mathbb{Z}}}
\newcommand{\bR} { {\mathbb{R}}}

\newcommand{\bN} { {\mathbb{N}}}
\newcommand{\bO} { {\mathbb{N}_{\geq 2}}}
\newcommand{\bS} { {\mathbb{S}}}

\newcommand{\cC} { {\mathcal{C}}}

\newcommand{\cE} { {\mathcal{E}}}
\newcommand{\cK} { {\mathcal{K}}}

\newcommand{\cS} { {\mathcal{S}}}

\newcommand{\sT} { {\mathscr{T}}}
\newcommand{\supp}{{\operatorname{supp}}}

\newcommand{\HOM}{{ \textnormal{Hom}  }}

\newcommand{\pt}{{\operatorname{pt}}}
\newcommand{\bq} { {\textbf{q}}}
\newcommand{\bp} { {\textbf{p}}}

\newcommand{\id}{{\operatorname{id}}}
\newcommand{\Id}{{\operatorname{Id}}}

\newcommand{\fpw}[1]{{\llbracket #1\rrbracket}}
\newcommand{\fls}[1]{{(\!(#1 )\!) }}
\newcommand{\Sh} { {\operatorname{Sh}}}

\newcommand{\colim}{\operatornamewithlimits{colim}}

\newcommand\End{\operatorname{End}}
\newcommand\Fun{\operatorname{Fun}}
\newcommand\HH{\operatorname{HH}}
\newcommand{\Thom}{\operatorname{Th}}

\newcommand\Tr{\operatorname{Tr}}

\newcommand\im{ {\operatorname{im}} }

\newcommand\Hom{\operatorname{Hom}}

\newcommand{\Map}{{\operatorname{Map}}}

\newcommand{\op}{{\operatorname{op}}}

\newcommand\Ext{\operatorname{Ext}}
\newcommand\Tor{\operatorname{Tor}}
\newcommand\cof{\operatorname{cof}}

\newcommand\fib{\operatorname{fib}}

\newcommand\Sp{ {\operatorname{Sp} } }

\newcommand\Loc{ {\operatorname{Loc}} }
\newcommand\mhom{\operatorname{\mu hom}}

\newcommand\ms{\operatorname{SS}}

\newcommand\msnz{ \dot{\operatorname{SS} }}

\newcommand\sHom{\operatorname{\mathscr{H}om}}

\newcommand\Spec{\operatorname{Spec}}
\newcommand\Perf{\operatorname{Perf}}

\newcommand\PrLV[1][\sT]{  {  \operatorname{Pr}^{\operatorname{L} }_{\sT,st} } }
\newcommand{\bLambda}{{\boldsymbol{\Lambda}}}

\newcommand\nT{ {\textnormal{T}} }

\newcommand{\larbR} { \overleftarrow{\mathbb{R}}}
\newcommand{\rarbZ} { \overrightarrow{\mathbb{Z}}}
\newcommand{\Ex}{\operatorname{Ex}}
\newcommand{\Real}{\operatorname{Real}}

\def\lf{\left\lfloor}   
\def\rf{\right\rfloor}
\makeatletter
\newcommand\boxstar{\mathrel{\boxcls@{\star}}}
\newcommand{\boxcls@}[1]{%
  \vphantom{\Box}
  \mathpalette\boxcls@@{#1}%
}
\newcommand{\boxcls@@}[2]{%
  \ooalign{$\m@th#1\Box$\cr
  \hidewidth\boxcls@fix{#1}\hbox{$\m@th#1#2\mkern 2.4mu$}\cr}}

\newcommand\boxcls@fix[1]{
  \ifx#1\displaystyle
    \raise.225ex
  \else
    \ifx#1\textstyle
      \raise.225ex
    \else
      \ifx#1\scriptstyle
        \raise.180ex
      \else
        \raise.150ex
      \fi
    \fi
  \fi
}
\makeatother

\title[Capacities from the Chiu-Tamarkin complex II]{Capacities from the Chiu-Tamarkin complex II\\ Product formula and toric domains}
\author{Wenyuan Li}
\address{Wenyuan Li,\newline School of Mathematical Sciences, Zhejiang University, 866 Yuhangtang Rd, Hangzhou, 310058, P. R. China.}
\email{wenyuan.li@zju.edu.cn}

\author{Bingyu Zhang}
\address{Bingyu Zhang,\newline Department of Mathematics, Kyiv School of Economics, 3 Mykoly Shpaka St, Kyiv, 02000, Ukraine.}
\email{bzhang@kse.org.ua}

\curraddr{School of Mathematics, Institute for Advanced Study, 1 Einstein Dr, Princeton, NJ 08540, USA.}
\email{byzhang@ias.edu}

\date{\today}

\begin{document}

\begin{abstract}
In this article, we continue our study of equivariant symplectic capacities of open sets in a cotangent bundle defined via the Chiu-Tamarkin complex constructed from the microlocal theory of sheaves. 

The paper consists of two parts: 1) We develop the theory of $S^1$-equivariant Chiu-Tamarkin complex over complex-oriented ring spectra. 2) Applications of the theory to symplectic capacities: We show that a product formula for the capacities, whose Ekeland-Hofer capacities version is folklore since 1988, is not true via an (uncountable family of) counterexample(s). The counterexample (family) is the product of a concave toric domain of dimension $6$ ($2n\geq 6$) and a $2$-dimensional disk. We also prove a Weyl law of the capacities for all toric domains in dimension $2n$. The algebraic verification is based on a chain model for the Chiu-Tamarkin complex of all toric domains in $\bR^{2n}$, which recovers Irie's conjectural chain model of $S^1$-equivariant symplectic homology for star-shaped toric domains (defined only) in $\bR^4$.

\end{abstract}
\maketitle

\vspace{1cm}

\epigraph{Where the precursor was, there the ephebe shall be.}{Harold Bloom, \emph{The Anxiety of Influence}}

\vspace{1.5cm}

\tableofcontents

\newpage
\section{Introduction}

Symplectic embedding problems have been one of the central topics since the pioneering non-squeezing theorem of Gromov \cite{gromov1985}. One of the most powerful invariants to study symplectic embedding problems is called symplectic capacities \cite{EkelandHoferII,viterbo1992symplectic,quantitivesymplectic2010,hutchings2011,McDuffSchlenk2012,gutthutchings2018capacities,SiegelHigherCapacities}.

The Ekeland-Hofer capacities are a sequence of symplectic capacities $(c_k^{EH})_{k\geq 1}$ defined for star-shaped domains in $\bR^{2n}$ via variational principles for Hamiltonian action function \cite{EkelandHoferII}. Floer, in his talk in 1988, suggested defining another sequence of capacities using equivariant Floer homology \cite[Section 11.2]{Floer_book} (notes taken by Weinstein \cite{Floer_notes}), which is sketched in \cite{viterbo-functorsI}, and Chekanov, in a talk in 2006, explained similar ideas (notes taken by Schlenk \cite{Chekanov_notes}). Later, Gutt--Hutchings finally defined a new sequence of capacities $(c_k^{GH})_{k\geq 1}$ using equivariant symplectic homology \cite{gutthutchings2018capacities}.

Discussions on the behavior of capacities under products perhaps go back to Floer \cite{Floer_notes}, where he proposed the following question on the following formula for Ekeland--Hofer capacities of product domains:
\[c_k(U \times W) \stackrel{?}{=} \min_{i+j=k,i,j\geq 0}\{c_i(U) + c_j(W)\}.\]
This was mentioned in \cite[Proposition 4.4]{Viterbo_Seminaire_Bourbaki} as an ongoing work by Floer--Hofer. Then, during a hike on \"Uetliberg on a sunny day in May 2004, Chekanov also proposed the above formula \cite[Section 3.8]{quantitivesymplectic2010} and explained more in a talk two years later \cite{Chekanov_notes}. Gutt--Hutchings conjectured a similar product formula for their capacities \cite[Conjecture 1.9]{gutthutchings2018capacities}.
Finally, a similar conjecture for symplectic $p$-product (where the original product formula corresponds to $p=\infty$), was proposed later in \cite[Conjecture 1.6]{HaimKislevOstrover2023}.

Floer and Chekanov's proposal, according to their talks \cite{Floer_notes,Chekanov_notes}, consists of two parts. First, they proposed the equivalence between the Ekeland--Hofer capacities and the Gutt--Hutchings capacities (even though the latter only came much later after their proposals), which we will discuss later; then, the product formula for the latter was expected to be an application of the Künneth formula for equivariant symplectic cohomology. However, the only documented result for the latter claim seems to be a one-sided inequality on the symplectic cohomology capacity of A.-L. Biolley \cite[Proposition 2.5]{Biolley} (who strengthened the result in \cite[Section 4.3]{Oancea-kunneth}).

The product formula for symplectic capacities, if true, would have a simple corollary, that Gutt--Hutchings (and Ekeland--Hofer) capacities are invariant under stabilizations:
\[c_k(U \times \bC^n) = c_k(U).\]
However, even the stabilization invariance does not have a complete documented proof. For generating function capacities, a one-side inequality was obtained in \cite[Lemma 2.7]{Humiliere2008}; for Gutt--Hutchings capacities, the only documented result we find is \cite[Proposition 6.4]{SiegelHigherCapacities} combined with \cite[Theorem 7.64]{Pereira}, which proved the result assuming an index bound \cite[Assumption 5.6]{SiegelHigherCapacities}.

Meanwhile, discussions on the behavior of asymptotics of Gutt--Hutchings (and Ekeland--Hofer) capacities already appeared in \cite{quantitivesymplectic2010}, where they discussed the limit
\[\lim_{k\to \infty}\frac{c_k(U)}{k}.\]
The above limit is known to exist for convex or concave toric domains \cite{quantitivesymplectic2010,gutthutchings2018capacities} and is related to the Lagrangian capacity \cite[Conjecture 1.5]{CieliebakMohnke}, proved by Faisal--Li \cite{Faisal-Li} (after the result of Pereira \cite{Pereira}). However, in general, we do not know whether the limit exists. One partial progress in this direction is by Irie \cite{Irie_ChainModel}, who conjectured a chain model for equivariant symplectic homology for all toric domains in $\bC^2$, which would imply that the limit exists.

\subsection{Main results}

In this paper, we will consider the sequence of symplectic capacities $\widehat{c}_k(U)$ using the equivariant Chiu--Tamarkin complex \cite{mc-tamarkin,chiu2017} constructed from microlocal theory of sheaves \cite{KS90}. Similar but not identical capacities already appeared in \cite{Capacities2021, CyclicZHANG}. We show that the above folklore product formula is true for $k=1$ but can be a strict inequality in general for $k \geq 2$, even for a product of two concave toric domains. Meanwhile, the stabilization invariance is true in general.

Moreover, we develop a chain model and a spectral sequence for the Chiu-Tamarkin complex of all (not necessarily star-shaped) toric domains in $\bC^n$, which recovers Irie's conjectural chain model of $S^1$-equivariant symplectic homology for star-shaped toric domains in $\bC^2$. The chain model is our main technical tool to obtain the counterexample for the product formula. We also use it to study the Weyl law for all toric domains in dimension $2n$. We can also study some interesting examples, which will be explained in the following.

\subsubsection{Results on capacities}
We construct a sequence of $ S^1$-equivariant symplectic capacities $\widehat{c}_k(U)$ (as well as a discrete version $\widehat{d}_k(U)$) for open sets $U\subset T^*M$ for oriented $M$ using the equivariant Chiu--Tamarkin complex. We will also discuss certain comparison results about $\widehat{c}_k(U)$ and some other equivariant capacities constructed in earlier work of the second named author \cite{Capacities2021,CyclicZHANG}, denoted by $c_k(U)$ and $d_k(U)$ respectively here (they are denoted by $\overline{c}_k(U)$ and $c_k(U)$ respectively in loc. cit.). We refer to \autoref{section: capacities} for more discussion on their definition and comparison.

In summary, they satisfy $\widehat{d}_k(U)=\widehat{c}_k(U)\leq c_k(U)= d_k(U)$ for all $k$. The inequality generally should not be an equality, but the equality holds when $U$ is dynamically convex (see \autoref{coro: c=hat c, dynamically convex case} for more general discussion) and for all (not necessarily convex) toric domains (see \autoref{prop: toric c = hat_c in general}).

We then prove the following results concerning a product inequality for the equivariant symplectic capacities $\widehat{c}_k(U)$, whose Ekeland--Hofer and Gutt--Hutchings capacity versions have been folklore for a few decades. In particular, we can confirm the stabilization invariance of these symplectic capacities. 

\begin{Thm}\label{thm: product inequality intro}
    Let $U \subset T^*M$, $W \subset T^*N$ be open subsets, we set $\widehat c_0(U) \coloneqq 0$. Then we have
\begin{enumerate}[fullwidth]
    \item (\autoref{prop: stabilization invariance}) For all $k, n \in \bN$, we have \[\widehat{c}_k(U \times \mathbb{R}^{2n}) = \widehat{c}_k(U).\]
    \item (\autoref{thm: product inequality}) When $U,W$ has RCT boundary, then for all $k\in \bN$,
    \[\widehat{c}_k(U \times W) \leq \min_{i+j=k,\,i,j\geq 0}\{\widehat{c}_i(U) + \widehat{c}_j(W)\}.\] 
    \item (\autoref{Product formula for c1 RCT}) When $U,W$ have non-degenerate RCT boundary, we have
    \[\widehat{c}_1(U \times W) = \min\{\widehat{c}_1(U), \widehat{c}_1(W)\}.\]
    \item (\autoref{prop: product equality 1} and \autoref{prop: product equality summand}) When $U, W$ have non-degenerate RCT boundary, and a certain homological obstruction in the Chiu--Tamarkin complex vanishes, then 
    \[\widehat{c}_k(U \times W) = \min_{i+j=k,\,i,j\geq 0}\{\widehat{c}_i(U) + \widehat{c}_j(W)\}.\] 
\end{enumerate}
\end{Thm}
\begin{RMK}We can relax the non-degenerate condition to a more algebraic condition, which is satisfied by certain non-smooth sets, for example, weakly convex toric domains containing the origin and concave toric domains in \autoref{section:convex and concave toric domains}.

Notice that if $U,W\subset \bR^{2n}$, we may also remove the non-degenerate condition using the usual perturbation argument (c.f. \cite[Proposition 3]{EkelandHoferII}).
\end{RMK}

\begin{RMK}
For the capacities $c_k$, which are variants of $\widehat{c}_k$, we also prove similar results:
\begin{enumerate}[fullwidth]
    \item (\autoref{prop: stabilization invariance} and \autoref{thm: product inequality}) For all $k, n \in \bN$, we have
    \begin{gather*}
    {c}_k(U \times \mathbb{R}^{2n}) = {c}_k(U),\\
    {c}_k(U \times W) \leq \min_{i+j=k,\,i,j\geq 0}\{{c}_i(U) + {c}_j(W)\}.
    \end{gather*}
    \item (\autoref{Product formula for c1 RCT}) When $U, W$ have non-degenerate RCT boundary, we have
    \[{c}_1(U \times W) = \min\{{c}_1(U), {c}_1(W)\}.\]
\end{enumerate}
\end{RMK}

\begin{RMK}\label{Remark: non-equivariant capcities equivalence}

By \cite{Hochschild-Kuo-Shende-Zhang}, for connected open $U$ with non-degenerate RCT boundary, we have ${c}_1(U)=c_{SH}(\overline{U})=c_{Viterbo}(U)$, the symplectic cohomology capacity (see \autoref{eq: c1=c viterbo} for more discussion). Consequently, we have 
\[{c}_{SH}(\overline{U} \times \overline{W}) = \min\{{c}_{SH}(\overline{U}), {c}_{SH}(\overline{W})\},\]
which generalizes \cite[Proposition 2.5]{Biolley}. Therefore, we also have
\[c_{Viterbo}(U\times \bR^{2n})=c_{Viterbo}(U),\]
which generalizes \cite[Lemma 2.7]{Humiliere2008}.
\end{RMK} 
\begin{RMK}\label{remark: large radius disk stablization intro}By \autoref{coro: disk stablization}, we also have the following version of stabilization invariance: For connected open $U$ with non-degenerate RCT boundary, $r>0$, and $W=D(r)\subset \bR^2$, suppose $\widehat{c}_k(U), {c}_k(U)\leq r$, then 
\[\widehat{c}_k(U\times D(r))=\widehat{c}_k(U), c_k(U\times D(r))=c_k(U).\]
This is similar to \cite[Proposition 6.3]{SiegelHigherCapacities}.
\end{RMK}

The stabilization invariance has the following application to a stable embedding problem. For $\mathbf{a}=(a_1,\dots,a_n)$ with $0<a_1\leq \dots\leq a_n$ and $t\in (0,1)$, we consider the following toric annulus \[A_t(\mathbf{a})\coloneqq E(\mathbf{a})\setminus \overline{E(t\mathbf{a})}.\]
\begin{Prop}[{\autoref{prop: ellipsoid stabilization}}]Let $n\geq 2$ and $N\geq 0 $. With the notation above, we have that if there exists an isotopy of Hamiltonian embedding $\varphi_s:A_t(\mathbf{a})\times \bC^N \to \bC^{n+N}$ such that
\[\varphi_1(A_t(\mathbf{a})\times \bC^N)\subset B^{2n}(A)\times \bC^N,\]
then we have $\lceil k/n \rceil A\geq \widehat{c}_k(E(\mathbf{a}))$ for all $\ k$.
In particular $A\geq \widehat{c}_n(E(\mathbf{a}))$.
\end{Prop}
\begin{RMK}If $a_n<2a_1$, we have $\widehat{c}_n(E(\mathbf{a}))=a_n$,and then the inclusion
\[A_t(\mathbf{a})\subset E(\mathbf{a}) \subset B^{2n}(a_n)\]shows that $A\geq a_n$ is sharp. If $N=0$, i.e., the unstable case, the result is well known, for example \cite[Theorem 1]{Schlenk2003}. In fact, it is known that the negative boundary should not contribute to embedding obstruction. However, the stable ($N\geq 1$) case seems unknown, which makes us confused.
\end{RMK}

On the other hand, we show that in certain cases of toric domains, the product inequality may become a strict inequality (We will explain in \autoref{section: product counterexample} its algebraic verification and a potential dynamics explanation):
\begin{Thm}[{\autoref{prop: couterexample capacity estimation} and \autoref{coro: smooth model of the product counterexample}}]\label{Thm: counterexample intro}There exists a bounded smooth concave toric domain $X_\Omega\subset \bR^6$ such that for the disk $D(2)\subset \bR^2$, we have
\[\widehat{c}_5(X_\Omega \times D(2)) \leq 9 ~{\color{red}<}~  \min_{i+j=5}\{\widehat{c}_i(X_\Omega) + \widehat{c}_j(D(2))\}.\]  
\end{Thm}
Compared to \autoref{remark: large radius disk stablization intro}, we see that it is crucial that $\widehat{c}_5(X_\Omega)>9>2$ here.
 
\begin{RMK}\label{remark: higher dimension example formula}There exists an uncountable family of counterexamples $X_\Omega\subset \bR^{2n+2}$ for $n\geq 2$. 
Let $E(a_1,\cdots,a_{n+1})$ be an ellipsoid in $\bR^{2n+2}$ for $n\geq 2$. Let $\sigma>0$ be real and $m\geq 2$ be an integer.  
We take $X_\Omega^0\subset \bR^{2n+2}$ to be 
\[E(R,a,\cdots,a)\cup E(b,R,\cdots,R) \subset \bR^{2n+2}.\]

Let us take $a,b,R,r$ as a solution of the following system of inequalities, i.e. \autoref{equation: inequality system for counterexample}:
\begin{equation*}
    \begin{cases}
&A_\sigma =  R\left(1-\frac{\sigma}{a}\right)
=b\left(1-\frac{\sigma}{R}\right)>0,\\
&nA_\sigma\leq r<\sigma,\\
& m\sigma+A_\sigma+r
<\min\{R,(mn+1)r\}.\\
    \end{cases}
\end{equation*}
Then we have
\[
\widehat{c}_{mn+1}(X_\Omega^0\times D(r))
\leq m\sigma+A_\sigma+r
< 
\min_{i+j=mn+1}
\{
\widehat{c}_i(X_\Omega^0)+\widehat{c}_j(D(r))
\}.
\]
We then take $X_\Omega$ to be an inner but close enough smoothing of $X_\Omega^0$ such that the above strict inequality is still true (but the equality in parentheses fails, which is harmless).

For simplicity, in \autoref{section: product counterexample}, we give the detailed computation only
for the case $n=m=2$ and $\sigma= 3$,  in which $X_\Omega$
is six-dimensional, and we take $R=20,r=2$ and $a(3)=1$ (then one can solve $a,b$) as a solution of the system. The general case follows from the same argument as explained in \autoref{remark: higher dimension example formula in the end}, where the system \autoref{equation: inequality system for counterexample} is also explained.
\end{RMK}

\begin{RMK}In \cite[Conjecture 1.6]{HaimKislevOstrover2023}, a conjecture on symplectic $p$-product is also proposed: For two star-shaped domains $\overline{U}\subset \bR^{2n}$, $\overline{W}\subset \bR^{2m}$, and $k\geq 1$ and $p\geq 1$, we have 
\[
c^{EH}_k(\overline{U}\times_p \overline{W})
=
\begin{cases}
\displaystyle
\min_{\substack{i+j=k\\ i,j\geq 0}}
\left[
c^{EH}_i(\overline{U})^{\frac{p}{p-2}}
+
c^{EH}_j(\overline{W})^{\frac{p}{p-2}}
\right]^{\frac{p-2}{p}},
& 2\leq p\leq\infty,\\[1.5em]
\displaystyle
\max_{\substack{i+j=k+1\\ i,j\geq 1}}
\left[
c^{EH}_i(\overline{U})^{\frac{p}{p-2}}
+
c^{EH}_j(\overline{W})^{\frac{p}{p-2}}
\right]^{\frac{p-2}{p}},
& 1\leq p\leq 2.
\end{cases}
\]
The $p=\infty$ case recovers the original product formula, and $p=2$ case should be understood under the limit $p\to 2$. 

We can also formulate the same question for $\widehat{c}_k$ (and also for $c_k^{GH}$ as stated therein). As \cite[Theorem 1.8]{HaimKislevOstrover2023} showed, when $1\leq p \leq 2$, the conjecture is true if both $\overline{U}\subset \bR^{2n}$, $\overline{W}\subset \bR^{2m}$ are concave toric domains; when $p\geq 2$, the conjecture is true if both $\overline{U}\subset \bR^{2n}$, $\overline{W}\subset \bR^{2m}$ are convex toric domains. As in \autoref{section:convex and concave toric domains}, by direct computation, $\widehat{c}_k=c_k^{GH}$ for both convex and concave toric domains. Then \cite[Theorem 1.8]{HaimKislevOstrover2023} also proves the corresponding result for $\widehat{c}_k$ for those cases.

Here, however, we claim that our example in \autoref{Thm: counterexample intro} also disproves this conjecture for sufficiently large $p\gg 2$. Our higher-dimensional examples \autoref{remark: higher dimension example formula} also provide more counterexamples for the equality; see \autoref{rmk: HaimKislevOstrover}. We believe that our chain model provides enough tools to disprove the equalities for all $p\neq 2$ and $k\geq 2$. However, we believe that \cite[Conjecture 1.6]{HaimKislevOstrover2023} is true for $p=2$ for $\widehat{c}_k$ for all $k$. For $k=1$ on star-shaped domains (which is clearly true for convex domains by \cite[Proposition 1.5]{HaimKislevOstrover2023}), we are less certain for its Ekeland--Hofer version.
\end{RMK}

Now we state our third main result. Following Irie's argument \cite{Irie_ChainModel}, we also prove the following Weyl law of ALL toric domains in $\bR^{2n}$. We denote by $X_\Omega \subset \bR^{2n}$ the toric domain whose moment image is $\Omega \subset \bR^n_{\geq 0}$.
\begin{Thm}[{\autoref{thm: weyl law for toric}}]Let $\Omega\subset \mathbb R_{\geq 0}^n$ be an non-empty open set, and we denote $a(\Omega)\coloneqq \sup_{x\in {\Omega}}\min_i x_i $, then we have
\[ ka(\Omega)\leq c_k(X_\Omega) \leq (k+n-1)a(\Omega). \]
In particular, 
\[\lim_{k\to \infty}\frac{\widehat{c}_k(X_\Omega)}{k}=\lim_{k\to \infty}\frac{{c}_k(X_\Omega)}{k} = a(\Omega).\]
\end{Thm}

\begin{RMK}As a corollary (see \autoref{coro: neither convex nor concave}), we find that the union of polydisks with $a_i\geq 2$
\[X_\Omega=\bigcup_i P(1,\dots,a_i,\dots,1) \subset \bC^n,\]
which is a monotone toric domain \cite{GuttHutchingsRamos}, but neither Hamiltonian diffeomorphic to a bounded weakly convex toric domain for all $n\geq 2$, nor Hamiltonian diffeomorphic to a bounded concave toric domain when $n\geq 3$.
\end{RMK}

Moreover, we have the following corollary concerning the Non-Disjoint Union of Cylinders (NDUC) capacity (compare to \cite[Definition 4]{GPR22}):
\begin{Coro}[{\autoref{coro: NDUC}}]Let $C_r=\{\zeta \in \bR^n_{>=0} : \min_i \zeta_i < r\}$, and $N_n(r)=X_{C_r}$. For $X\subset \bC^{n}$, we define the \textit{Hamiltonian NDUC capacity}
\[c^{N,h}(X)=\inf\{r>0\mid \exists\,\varphi_s:X \lhook\joinrel\xrightarrow{Ham} \bR^{2n}, \quad \varphi_1(X)\subset N_n(r)\}.\]
Then for all open $\Omega\subset \bR^n_{\geq 0}$, we have
\[\lim_{k\to \infty}\frac{\widehat{c}_k(X_\Omega)}{k}=\lim_{k\to \infty}\frac{{c}_k(X_\Omega)}{k}=c^{N,h}(X_\Omega).\]    
\end{Coro}
\begin{RMK}Here, we introduce $c^{N,h}$ since our toric domain $X_\Omega$ is not assumed to be smooth star-shaped. However, by \cite[Principle 4.3]{Sch18}, when $X_\Omega$ is smooth star-shaped, the result is true for  $c^{N}(X_\Omega)$ of \cite[Definition 4]{GPR22} instead of  $c^{N,h}(X_\Omega)$ here.
\end{RMK}

\begin{RMK}Another application of the Weyl law is a version of \cite[Theorem 1.9]{HaimKislevOstrover2023} for $\widehat{c}_k$. Notice that the proof therein assumes Conjecture 1.6, which is disproved for $\widehat{c}_k$. However, one can prove it directly by elementary methods if all factors are toric domains: Let $X_{\Omega_i}$ be finitely many toric domains such that $0<a(\Omega_i)<\infty$. Their $p$-product is given by the moment region
\[\Omega=\bigcup_{\substack{t_i\geq 0,~\sum t_i=1}} t_1^{2/p} \Omega_1\times \cdots \times t_m^{2/p} \Omega_m.\]
Therefore, by the Weyl law, it remains to show the following elementary identity
\[a(\Omega)=\sup_{\substack{t_i\geq 0,~\sum t_i=1}} \min_{1\leq i\leq m} \{t_i^{2/p} a(\Omega_i)\}=\Big(\sum_{i=1}^m a(\Omega_i)^{-p/2}\Big)^{-2/p},\]
which is clear.
\end{RMK}

\subsubsection{Equivariant cohomology package}
We summarize the main technical ingredients that go into the proofs. Since our capacities $\widehat{c}_k(U),\widehat{d}_k(U)$ (and $c_k(U),{d}_k(U)$) are defined by the $S^1$ and $C_\ell$-equivariant Chiu--Tamarkin complexes \cite{chiu2017,Capacities2021,CyclicZHANG}, the main results will be deduced from properties of the Chiu--Tamarkin complexes, which are the cyclic (co)homology of Tamarkin categories. For our future application, the equivariant cohomology package (and many related results) is proven for a complex-oriented commutative ring spectrum $E$, nevertheless, our applications only need discrete commutative rings in this article.

Recall that the Chiu--Tamarkin complex is an $\bR$-filtered chain complex. For $G = S^1$ or $C_\ell$ and $L \in \bR$, we denote by $F^G (U; E)$ the sheaf/persistence module of $G$-equivariant Chiu--Tamarkin complex, and $F^G (U; E)_L$ the stalk of the sheaf at $L\in \bR$, following \cite{Capacities2021,CyclicZHANG}. Recall that $S^1$-equivariant complexes are linear over $A= C^*(BS^1; E)$. We choose $u$ as a complex orientation of $BS^1$ and then $H^*A= E^*\fpw{u}\coloneqq H^{*}(E)\fpw{u}$. Then if $E=\Bbbk$ is a discrete commutative ring, we have $A\simeq \Bbbk\fpw{u}$ (at least as $E_2$ rings).

\begin{Thm}[Gysin sequence, \autoref{thm: Gysin}, {\cite[Theorem 1.7]{Li-Zhang26-1}}]
    Let $U \subset T^*M$ be an open subset and $E$ be a complex-oriented commutative ring spectrum. Then there is a fiber sequence
    \[\ F^{S^1}(U; E)^{hS^1} [-2]\xrightarrow{[\ell] (u)} F^{S^1}(U; E)^{hS^1}\to F^{C_\ell}(U; E)^{hC_\ell}  \xrightarrow{+1},\]
    where $[\ell](u)$ denotes the $\ell$-series of the formal group law of $E$. 
    
    In particular, when $E=\Bbbk$ is a discrete ring of characteristic $\ell$, we have $[\ell](u)=\ell u =0$ and
    \[F^{C_\ell}(U; \Bbbk)^{hC_\ell} \simeq F^{S^1}(U; \Bbbk)^{hS^1} \oplus F^{S^1}(U; \Bbbk)^{hS^1}[-1].\]
\end{Thm}
\begin{RMK}We temporarily assume $E=\Bbbk$ to be a discrete commutative ring. For general $S^1$ and $C_\ell$-complexes such that $\Bbbk$ has characteristic $\ell$, the Gysin sequence for cyclic chains of small $A_\infty$-categories, in particular Fukaya categories, is proved by Chen \cite{Chen1}\footnote{Thanks to Chen for explaining to us that his proof is actually true for all discrete commutative rings.}. For the symplectic (co)homology and $\ell = 1$, the Gysin sequence is proved by Bourgeois--Oancea \cite{Bourgeois-Oancea-17}; for $\ell \geq 2$, it is proposed by Seidel \cite[Section 2d]{Seidel-Equivariant} and proved by D.~Sen \cite{Sen2023}.
\end{RMK}

The main ingredient for \autoref{thm: product inequality intro} is the following Künneth formula of filtered Chiu--Tamarkin complexes. Our main observation is that to describe the filtered cohomology of the product completely, it is necessary to consider sheaf convolution (equivalently, tensor product over the Novikov ring plus a nice filtration, see \cite{APT-Tatsuki-Zhang}).

\begin{Thm}[Künneth formula, \autoref{thm:Künneth tamarkin}, {\cite[Theorem 1.8]{Li-Zhang26-1}}]
    Let $U \subset T^*M$, $W \subset T^*N$ be open subsets and $E$ be a complex-oriented commutative ring spectrum. Then there is a filtered Künneth isomorphism
    \[F^{S^1}(U \times W; E)^{hS^1} \simeq   F^{S^1} (U; E)^{hS^1} \widehat\star_{1_\sT^{hS^1}} F^{S^1} (W; E)^{hS^1}.\]
\end{Thm}
\begin{RMK}
    For the non-equivariant symplectic (co)homology, the Künneth formula is proved by Oancea \cite{Oancea-kunneth}. For the equivariant symplectic (co)homology, the Künneth formula (as an equivalence) has no written-down proofs even forgetting the filtration, to the best of our knowledge.
\end{RMK}

We show that when $U$ has restricted contact boundary and the Reeb dynamics of $\partial U$ is Morse--Bott non-degenerate, the Chiu--Tamarkin complex is generated by constant loops in $U$ and Morse--Bott families of parametrized Reeb orbits on $\partial U$, and the circle action is exactly reparametrization of closed Reeb orbits:

\begin{Thm}[{\autoref{prop: mhom CZ index} and \autoref{cor: mhom circle action trace}}]
    Let $U \subset T^*M$ be a domain with restricted contact boundary and $E$ be a commutative ring spectrum. Let $\gamma_j$, $j\in J$ be all the closed Reeb orbits of action $L$. Then, for any sufficiently small $\varepsilon$, we have a $S^1$-equivariant fiber sequence
    \[F^{S^1}(U; E)_{L+\epsilon} \to F^{S^1}(U; E)_{L-\epsilon} \to \bigoplus_{j\in J}\Gamma(\gamma_j, or_{\gamma_j})[-\mu_{CZ}(\gamma_j)],\]
    where $or_{\gamma}$ is a rank 1 local system on $\gamma$ and $\mu_{CZ}(\gamma)$ is the Conley--Zehnder index.
\end{Thm}

Finally, we also show in \autoref{thm: localization} the localization formula after inverting $u$, and identify the computational result with the ones in \cite{zhaoperiodic}, \cite{AlbersCieliebakFrauenfelder2016} and \cite{CoteKartal24}. We also show in \autoref{prop: ellpsoid embedding from complex-oriented ring spectrum} that naively considering complex-oriented ring spectra coefficients does not give new obstructions to the higher-dimensional ellipsoids embedding problem.

\subsubsection{Chain model for toric domains}
The novel technical tool for all toric computations is the following chain model. For convex toric domains in $\bR^{2n}$, such a chain model has been obtained by Siegel \cite{CompSiegelHigherCapacities} in his computation of their equivariant symplectic homologies. For general star-shaped toric domains in $\bR^4$, the chain model was proposed by Irie \cite{Irie_ChainModel} as a conjectural chain model of their equivariant symplectic homologies. The chain model has recently been confirmed on monotone toric domains in $\bR^4$ for cylindrical contact homology \cite{liang2026monotone}, with an application for computing capacities for those toric domains; and the same chain model is discussed in \cite{HutchingsConvexYoutube} with an application on classifying generic convex toric domains in $\bR^4$.

We generalize the chain model to all (not necessarily star-shaped) toric domains in $\bR^{2n}$ for all $n \in \bN$, and show that it computes their equivariant Chiu--Tamarkin complex (i.e.~cyclic traces of their Tamarkin categories).

We introduce the following notation. Let $m \in \bZ^n$ and $T \in \bR$, set
\begin{gather*}
\Omega(T; m) = \{\zeta \in \Omega \mid T + \left< m, \zeta\right> \geq 0\}, \quad \Omega(-\infty;m)=\varnothing,\\
I(m) = m_1 + \dots + m_n.
\end{gather*}

Let $\Bbbk$ be a commutative ring. For $\epsilon<T$, we then consider the graded $\Bbbk\fpw{u}$-module
\[I_{(\epsilon,T]}^{S^1}(X_\Omega)\coloneqq \bigoplus_{m\in \bZ^n}\left( C^{*+n+2I(m)}_c({\Omega(T;m)},\Omega(\epsilon;m) ;\Bbbk  ) \fpw{u} \otimes \bigwedge\nolimits^\bullet \Bbbk ^n\right),\]
with the differential
\begin{equation*}d|_m=d_{\mathrm{sing}}\otimes (-1)^N+\sum_{i=1}^n  (m_iu+t_i)\otimes \iota_{e_i},\end{equation*}
where $d_\mathrm{sing}$ should be understood as the usual singular differential tensored with the identity map on $\Bbbk\fpw{u}$, $t_{i,m}$ is the map on singular cochain induced by natural inclusions $\Omega(T;m-e_i)\subset \Omega(T;m)$, $(-1)^N|_{\wedge^q}=(-1)^q$ is a parity map on the exterior algebra, and $\iota_{e_i}$ the contraction map on the exterior algebra. We write $\widehat{I}_{(\epsilon,T]}^{S^1}(X_\Omega)$ the $u$-adic completion of $I_{(\epsilon,T]}^{S^1}(X_\Omega)$. 

\begin{Thm}[{\autoref{thm: chain model toric domain}}]
For $-\infty \leq \epsilon<T<\infty$. Let $\Omega\subset \bR_{\geq 0}^n$ be an open set, and $X_\Omega$ be its toric domain. Then $I^{S^1}_{(\epsilon,T]}(X_\Omega)$ and $\widehat{I}^{S^1}_{(\epsilon,T]}(X_\Omega)$ satisfying the following properties:

\begin{enumerate}[fullwidth]
\item The chain complex $\widehat{I}^{S^1}_{(\epsilon,T]}(X_\Omega)$ over $\Bbbk\fpw{u}$ is a chain model of $F_{(\epsilon,T]}^{S^1}(X_\Omega)^{hS^1}$. 

\item If $(\epsilon_2,T_2]\subset (\epsilon_1,T_1]$, then the natural morphism $F_{(\epsilon_1,T_1]}^{S^1}(X_\Omega)^{hS^1}\to F_{(\epsilon_2,T_2]}^{S^1}(X_\Omega)^{hS^1}$ is represented by a natural morphism $\widehat{I}_{(\epsilon_1,T_1]}^{S^1}(X_\Omega)\to\widehat{I}_{(\epsilon_2,T_2]}^{S^1}(X_\Omega)$. 

\item If $\Omega_1\subset \Omega_2 \subset \bR^n_{\geq 0}$ are two open sets, then the natural map $F_{(\epsilon,T]}^{S^1}(X_{\Omega_1})^{hS^1}\to F_{(\epsilon,T]}^{S^1}(X_{\Omega_2})^{hS^1}$ is represented by a natural map $\widehat{I}^{S^1}_{(\epsilon,T]}(X_{\Omega_1})\to \widehat{I}^{S^1}_{(\epsilon,T]}(X_{\Omega_2})$, and we have the following Mayer-Vietoris short exact sequence
\[0\to \widehat{I}^{S^1}_{(\epsilon,T]}(X_{\Omega_1\cap \Omega_2})\to \widehat{I}^{S^1}_{(\epsilon,T]}(X_{\Omega_1})\oplus \widehat{I}^{S^1}_{(\epsilon,T]}(X_{\Omega_2})\to \widehat{I}^{S^1}_{(\epsilon,T]}(X_{\Omega_1\cup \Omega_2})\to 0.\]

\item If $\Omega$ is a Lipschitz star-shaped domain, then $I^{S^1}_{(\epsilon,T]}(X_\Omega)$ is $u$-adic complete $\widehat{I}^{S^1}_{(\epsilon,T]}(X_\Omega)\simeq I^{S^1}_{(\epsilon,T]}(X_\Omega)$.

\item  If we assume in addition that $\Bbbk$ is a field containing $\bQ$ and if $\Omega$ is a Lipschitz star-shaped domain, there exists a filtration on $\widehat{I}^{S^1}_{(\epsilon,T]}(X_\Omega)\simeq F_{(\epsilon,T]}^{S^1}(X_\Omega)^{hS^1}$, whose associated spectral sequence has $E_1$-page as follows 
\[E_1^{p,q}=\bigoplus_{m\notin \bZ^n_{\geq 0}, I(m)= p,j} H^{q+n+3p-j}_c({\Omega(T;m)},\Omega(\epsilon;m) ;\Bbbk )   \otimes \wedge^j \ker(m) \Rightarrow H^{p+q} F_{(\epsilon,T]}^{S^1}(X_\Omega,\Bbbk)^{hS^1},\]
and the spectral sequence is $E_n$-degenerate. 
\end{enumerate} 
\end{Thm}

At first glance, it may not look like Irie's chain model. However, we will show in \autoref{section: Irie's model} that it is homotopy equivalent to Irie's chain model in $\bR^4$.

\begin{RMK}In \autoref{coro: chain model compute SH}, we explain that one can use the chain model to compute the non-equivariant symplectic cohomology for toric domains. Moreover, \autoref{conjecture:SH=HH, S1} predicts that such an equivalence can be lifted to $S^1$-equivariant symplectic cochains.    
\end{RMK}

This chain model allows us to compute equivariant symplectic capacities for all concave and weakly convex toric domains in \autoref{prop: concave capacity} and \autoref{coro: capacity weakly convex toric}. Such computation leads to the asymptotic behavior of $\widehat{c}_k(X_\Omega)/k$ following the argument of Irie as we explained before.

Our more explicit chain model for concave/convexx toric domains and the Mayer--Vietoris sequence also allows us to compute capacities of certain monotone toric domains in higher dimension, which leads to the counterexamples to the product formula stated at the beginning.

\subsection{Relation to other capacities}
We discuss the relation between the sheaf-theoretic equivariant capacities and other symplectic capacities.

\subsubsection{Gutt-Hutchings capacities} Let $C_T^{S^1}(U,\bQ)=\Hom_{S^1}(F^{S^1}(U,\bQ)_T,\bQ[-n])$ be the linear dual of $F^{S^1}(U,\bQ)_T$, which is a module over $\bQ\fpw{u}$. 
First, we will see in \autoref{prop: unfiltered CT=SH} that, forgetting the filtration, $S^1$-equivariant Chiu--Tamarkin complex $C_\infty^{S^1}(U,\bQ)$ computes (a version of) the $S^1$-equivariant symplectic cohomology when $U\subset \bR^{2n}$ is smooth star-shaped. Second, as proven in \cite{Hochschild-Kuo-Shende-Zhang}, there is an isomorphism between the filtered non-equivariant Chiu--Tamarkin complex $H^*C_T(U,\bQ)$ and the symplectic cohomology $\mathrm{SH}^{*}_{(-\infty,T]}(\overline{U}; \bQ)$.

Moreover, we also show that the Chiu--Tamarkin complex is generated by constant loops in the interior and Morse--Bott families of Reeb orbits on the boundary, and the $S^1$-action on the Chiu--Tamarkin complex can be geometrically characterised by reparameterization of Reeb orbits (see \autoref{section: Geometric circle action} for more details). The above evidence motivates the conjecture:

\begin{Conjecture}\label{conjecture:SH=HH, S1}If $M$ is Spin and $U \subset T^*M$ has non-degenerate restricted contact type boundary, for $T\notin \Spec(\partial U)$, there exists an isomorphism of $\bQ\fpw{u}$-modules
\[(\mathrm{SC}^{*}_{S^1,(-\infty,T]}(\overline{U}; \bQ),\mathbf{1}_{\mathrm{SC}}) \simeq ( C_{T}^{S^1}(U; \bQ),\eta_T),\]
where $\mathbf{1}_{\mathrm{SC}}$ is the symplectic cohomology unit and $\eta_T$ is the fundamental class (also the Hochschild cohomology unit, see \autoref{section: fundamental class}).
\end{Conjecture}
\begin{RMK}In fact, the only missing part for the conjecture is a lifting of the main result of \cite{Hochschild-Kuo-Shende-Zhang} to an $S^1$-equivariant \textit{morphism} since the forgetting $S^1$-action map is conservative.
\end{RMK}

Based on \autoref{conjecture:SH=HH, S1}, we have an equivalence between $S^1$-equivariant sheaf capacities $\widehat c_k$ and the Gutt-Hutchings capacities $c_k^{GH}$ \cite{gutthutchings2018capacities}.
\begin{Coro}[{\autoref{prop: GH=CT capacity level}}]\label{cor: GH=CT capacity}For bounded open sets $U\subset T^*M$ with RCT boundary, assuming \autoref{conjecture:SH=HH, S1}, we have for all $k\in \bN$ that 
\[\widehat{c}_k(U)=c_k^{{GH}}(\overline{U}).\] 
\end{Coro}

\subsubsection{Ekeland-Hofer capacities}
On the other hand, in a forthcoming article, we will show that
\begin{Thm}For a bounded open smooth star-shaped domain $U\subset \bR^{2n}$, we have
\[\widehat{c}_1(U)\leq c^{EH}_1(U),\qquad \widehat{c}_k(U)= c^{EH}_k(U),\,k \geq 2,\]
and the first inequality is generally strict for a certain example and becomes equality if $U$ is convex. 
\end{Thm}

Note that this result combined with \autoref{conjecture:SH=HH, S1} and \autoref{cor: GH=CT capacity} will disprove the long-standing folklore conjecture that Gutt--Hutchings capacities agree with Ekeland--Hofer capacities for star-shaped domains for $c_1$.
\begin{RMK}We thank Gutt and Ramos for informing us in November 2026 that they also found a counterexample where $c_1^{EH} < c_1^{GH}$, which is even dynamically convex (but not convex). 
\end{RMK}

\subsection{Convention}
In this article, we use the cohomological convention, even when working on ring spectra: We denote by $H^q=\pi_{-q}$ the cohomology/cohomotopy functor. For example, $H^*(M, E)$ means the homotopy of $\mathrm{Map}(\Sigma_+^\infty M, E)$ and $H_*(M, E)$ means the homotopy of $\Sigma_+^\infty M \otimes E$. For a complex-oriented ring spectrum, $H^*(BS^1, E) = H^*(E)\fpw{u}$ and $H_*(BS^1, E) = H^*(E)[u^{-1}]$.

In this article, if not mentioned otherwise, all limits/colimits, $\otimes$, and $\Hom$ are derived. We will particularly use $\Tor$ and $\Ext$ to denote cohomology of $\otimes$ and $\Hom$. For example, $M \otimes N$ means the derived tensor product of derived modules $M$ and $N$, and $H^*M \otimes H^*N$ means the underived tensor product of graded modules $H^*M$ and $H^*N$.

For the grading on symplectic cohomology, we use the convention that the unit $1$ is in degree zero of the cohomology $SH^*(U)$. For the $S^1$-equivariant symplectic cohomology (and the periodic version), we use the version where one takes homotopy fixed points (or Tate construction) in a finite action window and then take colimits with respect to action filtration; see \autoref{rem: versions of SH}.

We use $C_\ell$ to denote the finite cyclic group of order $\ell$. For a prime number $p$, we use $\bF_p$ to denote the finite field of integers modulo $p$.

\subsection*{Acknowledgment}As this work is a spiritual continuation of BZ's thesis, in which the prototype of many results was already shown, BZ would like to give Stéphane Guillermou special thanks for his encouragement. We would like to thank Chris Kuo, who is deeply involved in the discussions at the early stage of this project. We thank Kai Cieliebak, Helmut Hofer, Janko Latschev, Felix Schlenk and Claude Viterbo for explaining the historical process along the line of the product formula. We thank Laurent C\^ot\'e, Sheel Ganatra, Jean Gutt, Vincent Humilière, Michael Hutchings, Kei Irie, Alexandru Oancea, Vinicius G. B. Ramos, Kyler Siegel, Yuan Yao and Zhengyi Zhou for helpful discussions. 

WL was partially supported by the AMS Simons Travel Grant from fall 2024 to spring 2026. BZ was supported by the Novo Nordisk Foundation grant NNF20OC0066298 and VILLUM FONDEN, VILLUM Investigator grant 37814 for part of the work. Also, part of this work was completed during BZ's membership in the Institute for Advanced Study on the Special Year Program on\textit{ Conformally Symplectic Dynamics and Geometry}; BZ thanks the anonymous donor for the sponsorship.

\subsection*{AI Disclosure}
The project began in February 2024 partly as a rewriting of \cite{CyclicZHANG} after \cite{Kuo-Li-CY}. Many results were obtained, for example, most of the results in \autoref{section: Geometry of the trace of Tamarkin category},~\ref{section: capacities}, before large language models (LLMs) could help with mathematical research. However, due to the recent development of LLMs, we (mainly the second named author) used them, specifically ChatGPT 5.6 Sol, 6 Astra, and 6.1 Sol, to assist with the rest of this project. 

\begin{itemize}[fullwidth]
    \item The current form of \autoref{thm: product inequality intro}~(1)-(3) does not involve any AI-generated content, and most of it exists in 2024 or even earlier during the second named author's doctoral studies. For part (4), we worked from December 2024 till the beginning of September 2026, and thought we had obtained the ``proof'' of the product formula unconditionally. However, when we tried to simplify some of our ``uncomfortable'' arguments using ChatGPT 6 Astra, a mistake in the ``proof'' was pointed out, together with the counterexample in \autoref{Thm: counterexample intro} based on the chain model. Afterward, the algebraic obstruction of the product formula is developed by the authors, correcting our previous mistakes. Meanwhile, we prompted ChatGPT 6 Astra to figure out the potential dynamical explanation of the counterexample, which now becomes \autoref{subsection: dynamical explanation of the example} and its higher-dimensional generalizations.

    \item During the conference \textit{From Hamiltonian Dynamics to Symplectic Topology and Beyond}, June 2024 in Paris, the second named author learned Irie's conjectural chain model for $S^1$-equivariant symplectic homology of star-shaped toric domains in $\bR^4$ in \cite{Irie_ChainModel} after his talk. The second named author believes Irie's chain model could be derived from sheaf theory from then, because the sheaf $\cE^{\boxtimes n}$ (which already exists in the second named author's PhD thesis \cite{zhangthesis}) has a very similar formula. However, the second author forgot the idea until August 2026, when he recalled it in a discussion with ChatGPT 5.6 Sol on the Weyl law of capacities. ChatGPT 5.6 Sol then points out that  Irie's chain model is a \textit {Koszul}-type complex. The word immediately reminded the authors of their experience on the relation between lattice constructible sheaves and Koszul complexes, which turns out to have been developed into the current form of \autoref{thm: chain model toric domain} by the authors. Afterward, we prompted ChatGPT 5.6 Sol about the precise relation between our chain model and Irie's chain model, and then the homological perturbation argument was suggested. The spectral sequence mainly follows Irie's original argument, but ChatGPT 5.6 Sol helped us to distinguish many computations therein.

\item In \autoref{section: application of the chain model}, ChatGPT 5.6 Sol \& 6 Astra helped a lot in numerous meticulous homological algebra computations under the author's supervision, and the argument for the Weyl law belongs to Irie \cite{Irie_ChainModel}. However, the original idea of \autoref{prop: toric c = hat_c in general} was generated by ChatGPT 5.6 Sol, with only  \autoref{prop: c=hat c, direct summand criterion} and \autoref{thm: chain model toric domain} as the input, and then we simplified the proof afterward. The example in \autoref{coro: neither convex nor concave} is suggested by ChatGPT 6 Astra after we prompted it with ``can you tell me any interesting application of the Weyl law?'', and then the computation using the Mayer--Vietoris argument is figured out by the authors. 

\item After all the discussion over $\bQ$, we tried to consider complex-oriented coefficients based on our known computation, for at least ellipsoids. ChatGPT 6 Astra helped us find the reference \cite{Hausmann2022} on the equivariant Lazard isomorphism, especially \autoref{prop: decomposition of ell series}, which leads to the proof of \autoref{prop: ellpsoid embedding from complex-oriented ring spectrum}.

\item Based on our painful lesson on the fallacious proof of the product formula, we eventually use ChatGPT 6.1 Sol to thoroughly check the correctness of the mathematical content written by us, and then we correct their mistakes ourselves.

\item We use ChatGPT 6.1 Sol to produce all figures and improve the language of the article by providing a report of potential polishing.
\end{itemize}

The current article is eventually written by the authors themselves, who fully understand the mathematical arguments as they appear in the paper, and take complete responsibility for their correctness. The AI disclosure is written by the authors without help of LLMs.

\section{Some linear algebra}

\subsection{\texorpdfstring{Algebraic $S^1$-action}{}}
Here, we recall the basic notions of algebraic $S^1$-action on categories. Let $\cC$ be a stable category that admits both limits and colimits.

We recall $BS^1$ is the classifying anima of the circle group $S^1$, which is often presented as the anima of the topological space $\bC P^\infty$. Write $j: pt \hookrightarrow BS^1$ and $\pi: BS^1 \to pt$. Let  $ \Fun(BS^1,\cC)$ be the category of (Borel) $S^1$-action on $\cC$, and we write $j^*$ for the forgetful functor, and $\pi^*$ for the trivial $S^1$-action functor. Therefore, an object in $M\in \Fun(BS^1,\cC)$ an (Borel) $S^1$-action on the underlying object $j^*M$. And for $M\in \cC$, we often simply write $M=\pi^*M$ as the object equipped with the trivial $S^1$-action. 

Since $\cC$ admits enough limits/colimits, we can write $\pi_!$ for the homotopy orbit functor, and $\pi_*$ for the homotopy fixed point functor \cite[Example 6.1.6.2]{HA}:
\begin{align*}
    M_{hS^1} &= \pi_!M =  \varinjlim\nolimits_{BS^1}M, \\
    M^{hS^1} &= \pi_*M =  \varprojlim\nolimits_{BS^1}M.
\end{align*}

In practice, we often use a combinatorial model of $S^1$-action based on Connes's cyclic category $\bLambda$. Objects of the category $\boldsymbol{\Lambda}$ are $[n]=\{1,\dots,n\} $ for $n\in \bN_{\geq 1}$  and morphisms are generated by
\begin{align*}
    d_i: [n] \to [n+1], \;s_i: [n+1] \to [n], \;t_{n}: [n] \to [n]
\end{align*}
such that the face morphism $d_i$, the degeneracy morphism $s_i$ and the cyclic rotation $t_n$ satisfy the relations in the cyclic category (with the relation $t_n^{n} = \id$) \cite[Definition 6.1.1]{LodayCyclic}\footnote{For convenience in our applications, our convention has a difference with the convention of \cite[Definition 6.1.1]{LodayCyclic} for the cyclic category, and then also for the simplicial category. The difference is harmless.}.

Let $\boldsymbol{\Lambda}$ be the cyclic category \cite{Connes1983}. Then a cyclic object in $\cC$ is an object in $\Fun(\boldsymbol{\Lambda}, \cC)$. Let $\Fun_\simeq (\boldsymbol{\Lambda}, \cC)$ be the subcategory of functors that send all morphisms in $\boldsymbol{\Lambda}$ to equivalences in $\cC$, we call objects here cartesian cyclic objects. 
We write 
\[()^\#: \Fun (\boldsymbol{\Lambda}, \cC) \to \Fun_\simeq (\boldsymbol{\Lambda}, \cC)\]
for the natural localization functors.

cartesian cyclic objects in $\cC$ are known as locally constant presheaves, or equivalently, a parametrized spectrum over $BS^1$; see for example the exposition in \cite[Section 2]{Li-Zhang26-1}.
\[\Fun_\simeq (\boldsymbol{\Lambda}, \cC)\simeq \Fun(BS^1,\cC).\] 
In the following, we will then alternately use both terminologies: cartesian cyclic objects and $S^1$-objects.

In this section, we are particularly interested in the case $\cC$ admits a presentably (closed) symmetric monoidal structure $(\cC,\otimes,\Hom_\cC,1)$. In this case, the category $\Fun(BS^1,\cC)$ admits a presentably (closed) symmetric monoidal structure.
In particular, we define the $\cC$-enriched Hom $\Hom_{S^1}(M,N)\in \cC$ of $\Fun(BS^1,\cC)$, through the following (representable functor) formula
\[\Map_{\cC}(K,\Hom_{S^1}(M,N))\simeq \Map_{\Fun(BS^1,\cC)}(\pi^*K\otimes M,N),\quad \forall K\in \cC.\]

Let $A=A_\cC=\End_{S^1}(1)$, which is a commutative ring spectrum since $\cC$ is symmetric monoidal. Then we have that i) $M_{hS^1}$ is a (left) $A$-module via 
\[A\otimes M_{hS^1}=A\otimes \varinjlim\nolimits_{BS^1}M\simeq \varinjlim\nolimits_{BS^1} A\otimes M \xrightarrow{\varinjlim_{BS^1} \operatorname{ev}} \varinjlim\nolimits_{BS^1}   M  =M_{hS^1}; \]and ii) $M^{hS^1}$ is a (right) module over $A$ via
\[A\otimes M^{hS^1}=A\otimes \varprojlim\nolimits_{BS^1}M
\rightarrow\varprojlim\nolimits_{BS^1}(A\otimes M)
\xrightarrow{\varprojlim_{BS^1}\operatorname{ev}}
\varprojlim\nolimits_{BS^1}M=M^{hS^1}.
\]

The following simple lemma is clear but useful in our discussion below. Recall that $A$ is commutative, so one can identify left/right $A$-modules.
\begin{Lemma}\label{lemma: linear dual of S1-action}We have  equivalences of (left) $A$-modules:
\[\Hom_\cC(M_{hS^1},1)\simeq \Hom_{S^1}(M,\pi^*1)\simeq \Hom_\cC(j^*M,1)^{hS^1}, \]
where we equip the $A$-module structures (1) on the first term by the dual (left) $A$-module structure on the homotopy orbit; (2) on the second term by post-composition; (3) and on the third term by the (right) $A$-module structure on the homotopy fixed point.
\end{Lemma}
\begin{proof}
Then the first equivalence follows from that $\pi_!$ is the left adjoint of $\pi^*$. The second equivalence follows by $\pi_!M=\varinjlim\nolimits_{BS^1}M$. It is clear that all of three are left modules over $A$ and the equivalences identify their left $A$-module structures by the universal properties of limits/colimis.
\end{proof}

In this section, we are particularly interested in the case $E$ is a commutative ring spectrum and $D(E)$ is the stable category of $E$-modules \cite[Section 7.1.1]{HA}. Then we have
\begin{align*}
    M_{hS^1}  = (\Sigma^\infty_+ES^1  \otimes M)_{S^1}, \quad 
    M^{hS^1} = \Hom(\Sigma^\infty_+ES^1, M)^{S^1}.
\end{align*}

Here, for $A=\Hom_{S^1}(E,E)$, we have that i) the (left) $A$-module structure on $M_{hS^1}$ is induced via $M_{hS^1}\simeq E\otimes_{S^1\wedge E} M$ where $A$ acts on the $E$-factor; and ii) the (right) $A$-module structure on $M^{hS^1}$ is induced via $M^{hS^1}\simeq \Hom_E(\Sigma^\infty_+ES^1\wedge E, M)^{S^1}$ where $A$ acts on the $E$-factor.

In this case,  we have
\[A=A_{D(E)}\simeq C^*(BS^1,E)\coloneqq \Hom_E(\Sigma^\infty_+BS^1,E) = E^{hS^1}\]
as a commutative ring spectrum. In fact, we have an equivalence
\begin{Thm}[{\cite[Theorem 7.29, 7.43]{MNN_nilpotence_equivariant_homotopy}}]\label{thm: fixed point as a monoidal functor}For a ($E_1$-)complex-oriented commutative ring spectrum $E$, we have a symmetric monoidal equivalence
\[\Fun(BS^1,D(E))\to D(A)_{comp},\quad M\mapsto M^{hS^1},\]
where the subscript means that $N=M^{hS^1}$ is complete with respect to the augmentation $A\simeq C^*(BS^1,E)\to E$: I.e., for all $A$-module $P$, we have $P\otimes_A E =0$ implies $\Hom_{A}(P,N)=0$.
\end{Thm}
\begin{RMK}By the Greenlees-May Duality \cite{GreenleesMay1995}, (or \cite[Proposition 7.2.3.23~(5)]{HA}), $D(A)_{comp}$ is also equivalent to $D(A)_{tor}$, where $N$ is torsion with respect to the augmentation $A \to E$: I.e., for all $A$-module $P$, we have $P\otimes_A E =0$ implies $\Hom_{A}(N,P)=0$.  This is also equivalent to $N \otimes_A A[u^{-1}]=0$.

In the $\Fun(BS^1,D(E))$ case, we can see by the norm cofiber sequence that the equivalence is explicitly given by
\[\Fun(BS^1,D(E))\simeq D(A)_{tor},\quad M\mapsto M_{hS^1}[1].\]
\end{RMK}

\begin{RMK}\label{remark: formality for discrete ring}When $E=\Bbbk$ is a discrete commutative ring, which is particularly complex-oriented, we have the $E_2$-formality $A\simeq H^*(BS^1,\Bbbk)= \Bbbk\fpw{u}$ for $|u|=2$ \cite{Franz_E2formality,Horel_E2formality}; and if $\Bbbk$ contains $\bQ$, one can improve $A\simeq H^*(BS^1,\Bbbk)=\Bbbk\fpw{u}$ to an $E_\infty$ equivalence since $\Bbbk\fpw{u}$ is the Sullivan minimal model of $C^*(BS^1,\Bbbk)$. Therefore, the equivalence $\Fun(BS^1, D(\Bbbk))\simeq D(\Bbbk\fpw{u})_{comp}$ is a monoidal equivalence for all discrete commutative rings, and is a symmetric monoidal equivalence for $\Bbbk$ containing $\bQ$. In this case, $N$ being complete is equivalent to $N\simeq \varprojlim_n N\otimes_{\Bbbk\fpw{u}} \Bbbk\fpw{u}/u^n$. 

Notice that for all discrete commutative rings $\Bbbk$, we should be careful that $\Bbbk\fpw{u}$ is the completion in the category of graded rings, then $(\Bbbk\fpw{u})_{2d}=\prod_{i+2j=2d}(\Bbbk)_i[-2j]\simeq \Bbbk[-2d]$ and $(\Bbbk\fpw{u})_{2d+1}=0$, and then  $\Bbbk\fpw{u}=\oplus_{2d}(\Bbbk\fpw{u})_{2d}=\Bbbk[u]$. In this article, we use the notation $\Bbbk\fpw{u}$ to emphasize.
\end{RMK}

\begin{RMK}
It is beyond the scope of the paper to explain the proof of the result in great generality; however, we may explain the completeness in the cases $ E=\bQ$ and $ M=\Sigma^\infty_+ X$. In fact, a simple observation is that $ES^1\simeq S^\infty\simeq \varinjlim_n S^{2n-1}$, then we have 
\[C^*(M_{hS^1},\bQ)=C^*(X\times_{S^1} S^{\infty},\bQ)=\varprojlim_nC^*(X\times_{S^1} S^{2n-1},\bQ)\simeq \varprojlim_nC^*(X)\otimes_{\bQ\fpw{u}} \bQ\fpw{u}/(u^n).\]
\end{RMK}

In the end, we recollect some basic results on the formal group law. Let $E$ be a complex-oriented commutative ring spectrum; let $\zeta$ be the tautological line bundle over $BS^1\simeq \bC P^\infty$, then the complex orientation $u=c_1^E(\zeta)\in H^2(BS^1,E)$ determines an equivalence $H^*(BS^1,E)\simeq H^*(E)\fpw{u}$, and then the Künneth isomorphism gives $H^*(BS^1\times BS^1,E)\simeq H^*(E)\fpw{u_1,u_2}$. Let $\zeta_i$ be the tautological line bundle over the $i$-th copy of $BS^1\times BS^1$, and let $\mu:BS^1\times BS^1 \to BS^1$ the classifying map of  $\zeta_1\boxtimes\zeta_2$, we set
\[F(u_1,u_2)=c_1^E(\zeta_1\boxtimes\zeta_2)=\mu^*u\in H^*(BS^1\times BS^1,E)\simeq H^*(E)\fpw{u_1,u_2}.\]
Then $F$ is a formal group law in the sense of \cite{Quillen1969FormalGroupLaws}. We then define the $\ell$-series of $F$ as
\[[\ell]_E(u)=[\ell](u)\coloneqq c_1^E(\gamma^{\otimes \ell}) \in H^*(E)\fpw{u},\]
which can be given inductively using the formal group law by
\[[\ell+1](u)\coloneqq F([\ell](u),u)\in H^*(E)\fpw{u}.\]

For the complex-cobordism ring $MU$, we have the following result about $[\ell](u)$.

\begin{Prop}[{\cite[Proposition 5.46]{Hausmann2022}}]\label{prop: decomposition of ell series}Let $\bar{\psi}_m\in MU^*\fpw{u}$ be the ($u$-adic completion of the) image of $\psi_m$ therein under the equivariant Lazard isomorphism \cite[Theorem A]{Hausmann2022}, we have 
\[[\ell]_{MU}(u)=\prod_{m\mid \ell}\bar{\psi}_m \in MU^*\fpw{u}.\]
\end{Prop}
\begin{RMK}By the (multiplicative) Möbius inversion, where $\mu$ denotes the Möbius function, we can see 
\[\bar{\psi}_m=\prod_{d|m}[d]_{MU}(u)^{\mu(m/d)}.\]
\end{RMK}

\subsection{\texorpdfstring{Perfectness of $S^1$-modules}{}}

For $S^1$-objects over a ring spectrum $E$, we also need some technical lemmas. For a spectrum $E$, we denote its cohomotopy groups by $H^*(E)$.

The following result gives a criterion on when the cohomology groups of an $S^1$ object over a complex-oriented ring spectrum is perfect.

\begin{Lemma}\label{lem:perfectness}Let $M$ be an $S^1$-object over a complex-oriented ring spectrum $E$ where $H^*(E)$ is regular Noetherian with finite Krull dimension. Suppose $j^*M$ is a perfect $E$-module, we have $H^*(M^{hS^1})$ is a perfect module over $H^*(E^{hS^1})=H^*(E)\fpw{u}$.
\end{Lemma}
\begin{proof}
Taking the long exact sequence of cohomotopy groups associated to the Gysin sequence $M^{hS^1}[-2]\to M^{hS^1}\to j^*M$ \cite[Proposition 2.22]{Li-Zhang26-1}, we have an injection of $H^*(E)$ modules
\[H^*(M^{hS^1})/uH^*(M^{hS^1}) \hookrightarrow H^*(j^*M).\]
By \cite[Proposition 2.13]{AntieauMathewNikolaus}, when $H^*(E)$ is regular Notherian with finite Krull dimension, the $E$-module $j^*M$ is perfect if and only if its cohomology $H^*(j^*M)$ is a finitely generated $H^*(E)$-module. In particular, we have $H^*(M^{hS^1})/uH^*(M^{hS^1}) $ is a finitely generated $H^*(E)$-module.

Therefore, we pick a finite rank graded free $H^*(E)\fpw{u}$ module $\oplus_i H^*(E)\fpw{u}[a_i]$ and a homomorphism $F\to H^*(M^{hS^1})$ such that $\oplus_i H^*(E)[a_i] \to H^*(M^{hS^1})/uH^*(M^{hS^1})$ is surjective. This can be done by picking a set of homogeneous generators of $ H^*(M^{hS^1})/uH^*(M^{hS^1})$ and lifting them to $ H^*(M^{hS^1})$ (but we do not know if they are generators after lifting).

We know that both $\oplus_i H^*(E)\fpw{u}[a_i]$ and $H^*(M^{hS^1})$ are complete $H^*(E)\fpw{u}$-modules. The cokernel of $\oplus_{i}H^*(E)\fpw{u}[a_i]\to H^*(M^{hS^1})$ is complete and satisfies $C/uC=0$. Then the derived Nakayama lemma \cite[Proposition 6.5]{DwyerGreenlees} shows that $C=0$ and then $\oplus_i H^*(E)\fpw{u}[a_i]\to H^*(M^{hS^1})$ is surjective. In particular, we have that $H^*(M^{hS^1})$ is finitely generated as a $H^*(E)\fpw{u}$ module. Then we conclude that $H^*(M^{hS^1})$ is a perfect $H^*(E)\fpw{u}$ module because $H^*(E)\fpw{u}$ is also a Noetherian regular ring of finite Krull dimension.\end{proof}

\begin{RMK}
One can also consider an alternative proof of the result by the homotopy fixed point spectral sequence. For any complex-oriented spectrum $E$, we have $H^*(BS^1, H^*(j^*M))\simeq H^*(E)\fpw{u}\otimes H^*(j^*M) \Rightarrow H^{p+q}(M^{hS^1})$. Since $H^*(E)$ is a Notherian regular ring of finite Krull dimension, we know $j^*M$ is perfect over $E$ means that $H^*(j^*M)$ is a finitely generated $H^*(E)$-module. By the spectral sequence, we can conclude that $H^{p+q}(M^{hS^1})$ is also a finitely generated $H^*(E)\fpw{u}$-module, and thus $M^{hS^1}$ is perfect over $E^{hS^1}$.
\end{RMK}

In the $\bZ$-coefficient case, we can say more precisely.
\begin{Lemma}\label{lem: perfectness over Z}Let $E=\bZ$ be the discrete ring of integers. Let $M$ be an $S^1$-object over $\bZ$. Suppose $j^*M$ is a perfect $\bZ$-module.
Then there exists an integer $N > 0$ such that $H^*(M^{hS^1}) \otimes \bZ[1/N]$ is isomorphic to\[
 \bigoplus_i\bZ[1/N]\fpw{u}[a_i]
       \oplus\bigoplus_j\bigl(\bZ[1/N]\fpw{u}/(u^{b_j})\bigr)[c_j].
\]
In particular, for all $p$ prime such that $p\nmid N$, we then have
\begin{equation*}
    \begin{split}
        H^*(M^{hS^1}\otimes\bQ) &\simeq  \bigoplus_i\bQ\fpw{u}[a_i]
       \oplus\bigoplus_j\bigl(\bQ\fpw{u}/(u^{b_j})\bigr)[c_j]   ,\\
        H^*(M^{hS^1}\otimes\bF_p) &\simeq   \bigoplus_i\bF_p\fpw{u}[a_i]
       \oplus\bigoplus_j\bigl(\bF_p\fpw{u}/(u^{b_j})\bigr)[c_j]  .
    \end{split}
\end{equation*}
\end{Lemma}
\begin{proof}
We discuss the decomposition. We consider $H^*(M^{hS^1}\otimes \bQ)$. By flatness of $\bQ$ over $\bZ$, we can take the limit over finite skeleta $\bC P^k \subset \bC P^\infty \simeq BS^1$ 
to see that $H^*(M^{hS^1}\otimes \bQ)\simeq H^*((M\otimes \bQ)^{hS^1})$, as well as the entire homotopy fixed point spectral sequence that $H^p(BS^1, H^q(j^*M \otimes \bQ)) \simeq H^p(BS^1, H^q(j^*(M \otimes \bQ))$. Now, since $\bQ\fpw{u}\simeq \bZ\fpw{u}\otimes\bQ$ is a PID, then we have a decomposition of $H^*((M\otimes \bQ)^{hS^1})$ as a finite shifted direct sum 
\[
 \bigoplus_i\bQ\fpw{u}[a_i]
       \oplus\bigoplus_j\bigl(\bQ\fpw{u}/(u^{b_j})\bigr)[c_j].
\]
Since there are only finitely many $\bQ\fpw{u}$ generators; then, by carefully choosing generators, we can find $N_0>0$ and a morphism
\[f:  \bigoplus_i\bZ[1/N_0]\fpw{u}[a_i]
       \oplus\bigoplus_j\bigl(\bZ[1/N_0]\fpw{u}/(u^{b_j})\bigr)[c_j] \to H^*(M^{hS^1}) \otimes \bZ[1/N_0],\]
which is an isomorphism after tensoring with $\bQ$.
Moreover, by finite generation, we can pick a further integer $t>0$ such that both $\ker(f)$ and $\operatorname{coker}(f)$ are annihilated by $t$. 

Then we pick $N=tN_0$, and $f\otimes_{\bZ[1/N_0]} \bZ[1/N]$ induces an isomorphism
\[
 H^*(M^{hS^1}) \otimes \bZ[1/N]\simeq \bigoplus_i\bZ[1/N]\fpw{u}[a_i]
       \oplus\bigoplus_j\bigl(\bZ[1/N]\fpw{u}/(u^{b_j})\bigr)[c_j].
\]    
Now, if $p\nmid N$, direct computation verifies the last two isomorphisms.
\end{proof}
\begin{Coro}\label{Coro: finite restriction for large prime determine the free part}Let $M,N$ be two $S^1$-modules over $\bZ$ with perfect underlying $\bZ$-modules. Let $j_\ell^*$ the restriction functor $j_\ell^*:\Fun(BS^1,D(\bZ))\to \Fun(BC_\ell,D(\bZ))$. Suppose for all integer $l>0$, there exists an odd prime $p>l$, such that $j^*_p(M\otimes\bF_p)\simeq j^*_p(N\otimes\bF_p)$ over $\bF_p$. Then the free parts of $ H^*((M\otimes \bQ)^{hS^1}) $ and $ H^*((N\otimes \bQ)^{hS^1}) $ are isomorphic.
\end{Coro}
\begin{proof}We need to compare two free graded $\bQ\fpw{u}$-modules. Therefore, we only need to compare numbers and degrees of their generators.

Recall the $C_\ell$-Gysin sequence from \cite[Proposition 2.25]{Li-Zhang26-1}: For any $S^1$-module $M$, we have the fiber sequence
\[M^{hS^1}[-2] \xrightarrow{\ell u} M^{hS^1}\to (j^*_\ell M)^{C_\ell}.\]
Therefore,  for $\ell=p$  we have a splitting for $M_p\coloneqq M\otimes\bF_p$
\[H^*(j^*_p M_p)^{hC_\ell}\simeq H^*(M_p^{hS^1})\oplus H^*(M_p^{hS^1})[-1].\]

We define $p(M_p;t )$ as the counting polynomial of the free generators, weighted by their degrees, of $H^*(M_p^{hS^1})$ as a finitely generated graded module over the PID $\bF_p\fpw{u}$. Equivalently, $p(M_p;t )=(1-t^2)h(M_p;t)$ for the Hilbert series  $h(M_p;t)$ of the free part of $H^*(M_p^{hS^1})$ as a graded $\bF_p$-module. We similarly define $p(j^*_pM_p;t )$.  
Therefore, the splitting above implies 
\[p(j^*_pM_p;t )=(1+t)p(M_p;t),\quad p(j^*_pN_p;t )=(1+t)p(N_p;t).\]

Now, let $d_M,d_N$ be the nonzero integers
given by \autoref{lem: perfectness over Z} for $M,N$. Let $L>d_Md_N$ be an integer. We pick prime $p>L$ such that $j^*_p(M_p)\simeq j^*_p(N_p)$ by the condition. Therefore, we have the equality of the counting polynomial
\[(1+t)p(M_p;t)=p(j^*_pM_p;t )=p(j^*_pN_p;t )=(1+t)p(N_p;t).\]

Consequently, as $\bZ[t,t^{-1}]$ is an integral domain, we have
\[ p(M_p;t) = p(N_p;t).\]

On the other hand, since we pick $p>L>d_Md_N$, the polynomial $p(M_p;t)$ also counts the number of free generators of $ H^*((M\otimes \bQ)^{hS^1}) $ by \autoref{lem: perfectness over Z}. This concludes the proof. 
\end{proof}

\subsection{\texorpdfstring{Convolution of $S^1$-modules}{}}
We will study the convolution of diagrams augmented $S^1$-modules. Let $A = C^*(BS^1, E)$. First, the Eilenberg--Moore spectral sequence below will enable one to compute the homology of the tensor product of $A$-modules:

\begin{Lemma}\label{prop:eilenberg-moore}Let $M, N \in D(A)$ and $f: M \to A, g: N \to A$ be augmentations. Then there is a commutative diagram of spectral sequences
\[\begin{tikzcd}
\Tor^{p,q}_{H^*A}(H^*(M), H^*(N)) \ar[r,Rightarrow] \ar[d, "H^*f \otimes H^{q-*}g" left] & H^{q-p}(M \otimes_A N) \ar[d, "H^*(f \otimes g)"] \\
H^qA \ar[r, Rightarrow] & H^{q}A.
\end{tikzcd}\]

In particular, let $A=\bQ\fpw{u}$. Then there is a short exact sequence
\[\begin{tikzcd}
0 \ar[r] & H^*M \otimes_{H^*A} H^*N \ar[r] \ar[d,"H^*f \otimes H^*g"] & H^*(M \otimes_A N) \ar[r] \ar[d,"H^*(f \otimes g)"] & \Tor_{H^*A}^{1,*+1}(H^*M, H^*N) \ar[r] \ar[d] & 0\\
0 \ar[r] & A \ar[r] & A \ar[r] & 0 \ar[r] & 0.
\end{tikzcd}\]
Moreover, the horizontal morphisms in the left commutative square are isomorphisms on the torsion-free parts, and in particular, $\im(H^*f \otimes H^*g) = \im(H^*(f \otimes g))$.
\end{Lemma}
\begin{proof}
Consider the Eilengberg--Moore spectral sequence \cite[Chapter IV Theorem 6.1]{EKMM}, we can conclude the commutative diagram of spectral sequences. 

Consider the case when $A = \bQ\fpw{u}$. Since $A$ is a graded PID, the spectral sequence degenerates, and we obtain a short exact sequence of graded $A$-modules
\[0 \to H^*M \otimes_{H^*A} H^*N \to H^*(M \otimes_A N) \to \Tor_{H^*A}^{1,*}(H^*M, H^*N) \to 0.\]
Here, the torsion group is between graded $A$-modules and hence has an internal grading as graded modules. Note that over a PID, torsion-free modules are equivalent to flat modules \cite[Tag 0AUW]{StacksProject}. Now, we have a commutative diagram of torsion-free modules
\[\begin{tikzcd}
H^*M \otimes_{H^*A} H^*N \ar[r] \ar[d] & H^*(M \otimes_A N)\ar[d] \\
\overline{H^*M} \otimes_{H^*A} \overline{H^*N} \ar[r,"\sim"] & \overline{H^*(M \otimes_A N)}.
\end{tikzcd}\]
Since any augmentation factors through the torsion-free parts as observed above
\[H^*f: H^*M \to \overline{H^*M} \to A, \quad H^*g: H^*N \to \overline{H^*N} \to A,\] 
we can conclude that the diagram in the proposition commutes.
\end{proof}

\begin{Coro}\label{cor: complete eilenberg-moore}
Let $A = \bQ\fpw{u}$, $M, N \in D(A)_{comp}$ and $f: M \to A, g: N \to A$ be augmentations. Suppose $M$ is a perfect $A$-module. Then there is a commutative diagram
\[\begin{tikzcd}
{M \otimes_A N} \arrow[r,"{\simeq}"] \arrow[d,"{f\otimes g}"] & {M \,\widehat{\otimes}_A\, N} \arrow[d, "{f \,\widehat{\otimes}\,g}"] \\
A \arrow[r,"{\simeq}"] & A.
\end{tikzcd}\]
In particular, $\im(H^*f \otimes H^*g) = \im(H^*(f \otimes g))= \im(H^*(f \,\widehat\otimes\, g))$.
\end{Coro}
\begin{proof}
When $M$ is perfect, we know that it is dualizable and $M \otimes_A -$ commutes with limits. Then one can check that $M \otimes_A N$ is derived complete and $M \otimes_A N \simeq M\,\widehat\otimes_A\,N$.
\end{proof}
\begin{RMK}We are particularly interested in the case that $M,N$ are two
$S^1$-modules over $E$, and $M^{hS^1}, N^{hS^1}$ are the 
$A=E^{hS^1}$-modules. In this case, if we assume that the underlying $E$-module $M$ is perfect, then $M^{hS^1}$ is perfect by \autoref{lem:perfectness}, and \autoref{thm: fixed point as a monoidal functor} and \autoref{cor: complete eilenberg-moore} show that $(M\otimes_E N)^{hS^1}\simeq M^{hS^1} \otimes_A N^{hS^1}$.
\end{RMK}

Finally, we consider the behavior of finite limits over $A$-modules. Let $\Gamma$ be a finite diagram, $P \in \Fun(\Gamma, D(A))$ be a diagram of perfect $A$-modules and $A \in \Fun(\Gamma, D(A))$ be the constant diagram of $A$-modules. Let $p: P \to A$ be a functor. We can apply the Grothendieck spectral sequence:

\begin{Lemma}\label{lem: finite lim index}
Let $A = \bQ\fpw{u}$, $\Gamma$ be a finite diagram, $P \in \Fun(\Gamma, D(A))$ be a diagram of perfect $A$-modules and $p: P \to A$ be a functor. Suppose $\varprojlim_{\Gamma}A = A$. Then there is a commutative diagram of spectral sequences
\[\begin{tikzcd}
\mathbf{R}^p\varprojlim{}_{\Gamma}H^{q}P \ar[r,Rightarrow] \ar[d, "\mathbf{R}^p\varprojlim{}_{\Gamma}H^{q}\pi" left] & H^{q+p}(\varprojlim{}_{\Gamma}P) \ar[d, "H^{q+p}(\varprojlim{}_{\Gamma}\pi)"] \\
\mathbf{R}^p\varprojlim{}_{\Gamma}H^{q}A \ar[r, Rightarrow] & H^{q}A.
\end{tikzcd}\]
In particular, $\im(\varprojlim H^*\pi) = \im(H^*(\varprojlim \pi))$ if $\varprojlim{}H^*P$ is a direct summand of $H^*(\varprojlim{}P)$.
\end{Lemma}
\begin{proof}
Consider the Grothendieck spectral sequence \cite[Tag 015N]{StacksProject}. We obtain the commutative diagram of spectral sequences where, by assumption $\varprojlim{}_\Gamma A = A$.

For the case of equality, note that $\varprojlim{}_\Gamma H^*P$ is on the $E_1$-page of the Grothendieck spectral sequence where all entries below vanish. Then there are natural maps $H^*(\varprojlim{}_\Gamma P) \to \varprojlim{}_\Gamma H^*P$ that intertwine the maps $\pi: P \to A$. Since all complexes are perfect, it follows from the classification of finitely generated $A$-modules. 
\end{proof}

\section{Tamarkin category}

\subsection{Tamarkin category of a point}
\label{Tamarkin point}

Let us write $\Sh_-(\bR)$ for the category of sheaves with nonpositive microsupport, and correspondingly $\Sh_+(\bR)$ for the sheaves with nonnegative microsupport. Recall the Tamarkin category (of a point) is
\[\sT := \Sh(\bR) / \Sh_-(\bR).\]
Consider the left adjoint of the quotient $\sT \hookrightarrow \Sh_+(\bR)$. It realizes $\sT$ as the full subcategory on $\Sh(\bR)$ of objects with trivial compactly supported global sections \cite[Proposition 4.2]{Hochschild-Kuo-Shende-Zhang}. 

Let $s: \bR \times \bR \rightarrow \bR$ be the addition map.
Consider the map 
$\star: \Sh(\bR) \times \Sh(\bR) \rightarrow \Sh(\bR)$  which is given by the formula
\[ (G, F) \mapsto s_! (G \boxtimes F) = s_! (p_2^* G \otimes p_1^* F ).\]
This defines a symmetric monoidal structure on $\Sh(\bR)$, with unit $1_{0}$.  More generally, if we write the 
translation map on $\bR$ as $\nT_c(t)=t+c$, then 
one checks readily that $T_{c*} H \cong 1_{c} \star H$ for any $H \in \Sh(\bR)$. 

Since the symetric monoidal product $\star$ preserves colimits, we find that $\Sh(\bR)$ is a closed symmetric monoidal category. We can describe the internal Hom $\sHom^\star$ by the formula
\[\sHom^\star(F,G)=p_{2*}\sHom(p_1^{*}F,s^!G),\quad F,G\in \Sh(\bR),\]
such that, for $F,G,H\in \Sh(\bR)$, 
\[\Hom_{\Sh(\bR)}(F\star G,H)=\Hom_{\Sh(\bR)}(F,\sHom^\star(G,H)).\]

\begin{Lemma}[{\cite[Corollary 4.14]{GS2014}}] \label{lem: convolution estimates}
The following estimates hold
\begin{align*}
    \star : &\,\Sh(\bR)\times \Sh_{\pm }(\bR) \rightarrow \Sh_{\pm }(\bR),\, \text{ or } \Sh_{\pm }(\bR)\times \Sh(\bR) \rightarrow \Sh_{\pm }(\bR);\\
   \sHom^\star : &\,\Sh(\bR)^{op}\times \Sh_{\pm }(\bR) \rightarrow \Sh_{\pm }(\bR),\, \text{ or } \Sh_{\pm }(\bR)^{op}\times \Sh(\bR) \rightarrow \Sh_{\pm }(\bR).
\end{align*}
\end{Lemma}

Lemma \ref{lem: convolution estimates} implies that
the closed symmetric monoidal structure determined by $\star$
descends to $\sT$.  Note that in $\sT$, we have 
$[1_{\bR_{\geq 0}}]\cong [1_{0}]$, giving two useful expressions for
the unit.  
Let us note that the inclusion $\sT \hookrightarrow \Sh(\bR), [F] \mapsto F \star 1_{\bR_{\geq 0}}$ is a tensor ideal with respect to the monoidal structure $\star$.  Indeed: $ 1_{\bR_{\geq 0}}\star1_{\bR_{\geq 0}}=1_{\bR_{\geq 0}}$, so
\[([F]\star [G])\star 1_{\bR_{\geq 0}}=[F\star G] \star 1_{\bR_{\geq 0}}= F\star G \star 1_{\bR_{\geq 0}}=(F \star 1_{\bR_{\geq 0}}) \star (G\star 1_{\bR_{\geq 0}}).\]
\begin{Lemma}[{\cite[Lemma 4.4]{Hochschild-Kuo-Shende-Zhang}}]\label{lemma: hom computation}For $F,G\in \sT$, we have
\[\Hom_{\sT}(F,\nT_{c*}G)=\Hom_{\Sh(\bR)}(F\star 1_{\bR_{\geq 0}},\nT_{c*}G)=\Gamma_{[-c,\infty)}(\bR,\sHom^{\star}(F,G)).\]
\end{Lemma}

We mention that there is a monoidal equivalence $\Gamma_{[-,+\infty)}(\bR, -)$ from the Tamarkin category of point to torsion persistence modules \cite[Lemma 4.9]{Hochschild-Kuo-Shende-Zhang}, \cite[Proposition 4.22]{Efimov-K-theory}, or \cite{APT-Tatsuki-Zhang}; in particular, for any $F \in \sT$, the system $\Gamma_{[c,+\infty)}(\bR; F)$, ${c\in\bR}$, forms an $\bR$-filtered complex that is a persistence module.

\subsection{Tamarkin category of a domain}\label{sec:tamarkin domain}
We set
\[\sT(M)= \Sh(M \times \bR)/\Sh_{\{\tau\leq 0\}}(M \times \bR)\simeq \Sh(M;\sT),\]
where the last equivalence shows that $\sT(M)$ is $\sT$-linear.

We can further define categories of sheaves valued in $\sT$ with microsupport conditions. They are called Tamarkin categories \cite{tamarkin2013,GS2014}.

\begin{Def}[Tamarkin category]
Let $U \subseteq T^*M$ and $Z = T^*M \setminus U$. Consider the maps $ T^*(M \times \bR) \xleftarrow{i} J^1 M \xrightarrow{\pi} T^*M$ and  define the $\sT$-linear full-subcategory
\begin{align}
  \Sh_Z(M,\sT)\coloneqq \Sh_{\pi^{-1}(Z) \cup \{\tau \leq 0\}}(M \times \bR) = \{F\, |\, 
\pi(i^{-1} {\ms}(F)) \subset Z\}.
\end{align}
The notion of microsupport descends to define
a map $\underline{ss}$ from objects of $\sT(U)$ to closed coisotropic subsets of $U$, and we can write
\[\Sh_Z(M,\sT) = \{F\, |\, \underline{ss}(F) \subset Z\}.\]
We define the Tamarkin category for $U \subset T^*M$ to be the $\sT$-linear quotient
\begin{align}
\sT(U) \coloneqq \Sh(M; \sT) /\Sh_Z(M;\sT).
\end{align}
\end{Def}

Note that for $Z \subset Y$ closed conic subsets, there is a natural inclusion $\Sh_Z(M) \hookrightarrow \Sh_Y(M)$, which admits both left and right adjoints. 
For
$j: V \hookrightarrow U$, there is an evident natural $\sT$-linear surjection 
$j^*: \sT(U) \to \sT(V)$
induced by the $\sT$-linear inclusion $i_*: \Sh_{T^*M \setminus U}(M; \sT) \hookrightarrow \Sh_{T^*M \setminus V}(M; \sT)$, which preserves both limit and colimit. 
Then, both $i_*$ and $j^*$ admit left (right) adjoint $i^*(i^!)$ and $j_!(j_*)$ given by either the wrapping formula or, equivalently, the microlocal cut-off formula; see, for example \cite{Hochschild-Kuo-Shende-Zhang,zhang-cutoff}.

For a closed conic subset $Z \subset T^*M$, we introduce the notation
\[P_Z: \Sh(M) \to \Sh(M)/\Sh_Z(M), \qquad Q_Z: \Sh(M) \to \Sh_Z(M)\]
for the projector with images $\Sh(M)/\Sh_Z(M)$ and $\Sh_Z(M)$.
We also introduce the notation 
\[P_U : \sT(T^*M) \to \sT(T^*M),\qquad Q_U: \sT(T^*M) \to \sT(T^*M) \]
for the projector with image $\sT(U)$ and $\Sh_Z(M;\sT)$. 
Recalling that 
$\sT(T^*M) = \Sh(M, \sT)$, the projectors
$P_U,Q_U$ can be expressed via 
integral kernels $P_U,Q_U \in \Sh(M \times M, \sT)$ by \cite[Proposition 5.16]{Hochschild-Kuo-Shende-Zhang}. They form a fiber sequence 
$P_U \to 1_{\Delta_{M}}^\sT \to Q_U .$

\subsection{Cyclic traces of Tamarkin categories}\label{sec:Tamarkin trace}
We will review the results in \cite[Section 4.1 \& 5.3]{Li-Zhang26-1} which show that the equivariant Chiu--Tamarkin complex is the linear dual of the cyclic trace or cyclic cohomology of the Tamarkin category.

Therefore, the following results immediately follow from \cite[Theorem 1.1 \& 1.5]{Li-Zhang26-1}.

\begin{Thm}\label{thm:S1-action-sheaf}
Let $M$ be a manifold, $Z \subset T^*M$ be a closed conic subset and $U = T^*M \setminus Z$. Then the cyclic trace for $\Sh_Z(M; \cE)$ is given by
\[\Tr(\Id_{\Sh_Z(M)}, n)=\pi_{n!}\widetilde{\Delta}_n^* (Q_{Z}\boxtimes \cdots\boxtimes Q_{Z}), \]
where $\pi_{n} : X^n  \rightarrow \pt$ are projections, and $\Tilde{\Delta}_n$ are the twisted diagonals:
\begin{align*}
    \Tilde{\Delta}_n&: X^n \rightarrow X^{2n}, \qquad (x_1,\cdots,x_n)\mapsto (x_n,x_1,x_1,\cdots,x_{n-1},x_{n-1},x_n).
\end{align*}
The degeneracy maps are given by the morphism of sheaves $1_{\Delta_{M}} \to K_Z$, and the face maps are given by inverse of degeneracy maps.
\end{Thm}

\begin{Thm}\label{thm:S1-action}
    Let $U \subset T^*M$ be an open subset. Then we have a canonical functorial $S^1$-action on the Chiu--Tamarkin invariant $F^{S^1}_\bullet(U; E)$ defined by the cyclic trace $\Tr(\Id_{\sT(U)}, -): \boldsymbol{\Lambda} \to \sT$ such that
    \begin{align}\label{eq:trace-circle}
    F^{S^1}_n(U; E) := \Tr^\sT(\Id_{\sT(U)}, n) = \pi_{n!}^\sT\widetilde{\Delta}_{n,\sT}^*(P_U^{\boxtimes_\sT n}),
    \end{align}
    such that the cyclic rotation is induced by the cyclic permutation $t_{n}: (M \times M)^{n} \to (M \times M)^{n}$ of order $n$, and the face maps $d_i$ and degeneration maps $s_i$ are induced by the triangle identities of the unit and counit.
\end{Thm}
\begin{RMK}
    Using the functoriality of the Tamarkin category $\sT(U)$ with respect to inclusions and Hamiltonian homeomorphisms \cite{asanoike2022complete,guillermouviterbo_gammasupport}, we can immediately obtain the functoriality of the $S^1$-equivariant Hochschild homologies.
\end{RMK}

Similarly, considering the $C_\ell$-action on the Hochschild homology, the following results follow immediately from  \cite[Theorem 5.9]{Li-Zhang26-1}:

\begin{Thm}\label{thm:Z/ell-action}
    Let $U \subset T^*M$ be an open subset. We have a canonical functorial $C_\ell$-action on the Chiu--Tamarkin invariant $F^{C_\ell}(U; E)$ defined by the $C_\ell$-cyclic trace $\Tr(\Id_{\sT(U)}, -): \boldsymbol{\Delta}_\ell \to \sT$ such that
    \begin{align}\label{eq:trace-finite-group}
    F^{C_\ell}(U; E) := \Tr^\sT(\Id_{\sT(U)}, \ell) = \pi_{\ell!}^\sT\widetilde{\Delta}_{\ell\sT}^*(P_U^{\boxtimes_\sT \ell}),
    \end{align}
    such that the cyclic rotation is induced by the cyclic permutation $t_\ell: (M \times M)^{\ell} \to (M \times M)^{\ell}$ of order $\ell$.
\end{Thm}

\begin{RMK}
    Using the functoriality of the Tamarkin category $\sT(U)$ with respect to inclusions and Hamiltonian homeomorphisms \cite{asanoike2022complete,guillermouviterbo_gammasupport}, we can immediately obtain the functoriality of the $C_\ell$-equivariant Hochschild homologies.
\end{RMK}

Composing the stalk functor $\sT\hookrightarrow \Sh(\bR,D(E))\to D(E),F\mapsto F_T=i_T^*F$ for $i_T:\{T\}\to \bR$, we can regard both $F^{S^1}_\bullet(U; E)$ and $F^{C_\ell}(U; E)$ as equivariant spectra.

\begin{Def}[{\cite{Capacities2021,CyclicZHANG}}]
Let $T\in \bR$, and $G=S^1,C_\ell$. We equip the coefficient spectrum  $E$  with the trivial $G$-action. We define
\begin{equation}\label{eq: C_T formula}
C^{G}_T(U,E)=\HOM_{G}(F^G(U,E)_{T},E [-\dim M]) = \Gamma_{(-\infty,T]}(\bR, \Hom_G^\sT(P_U, \Delta_*\omega_M^\sT)).
\end{equation}

Since $F^G(U,E)\in\sT$, we have morphisms 
\[C^{G}_S(U,E) \to C^{G}_T(U,E),\quad S\leq T,\]
and then we set
\[C^{G}_{(S,T]}(U,E)=\cof(C^{G}_S(U,E) \to C^{G}_T(U,E)),\qquad C^{G}_\infty(U,E) =\varinjlim_{T\to \infty} C^{G}_T(U,E).\]
\end{Def}
\begin{RMK}By \autoref{lemma: linear dual of S1-action} (and its finite cyclic group version, which is true by the same proof), we can understand $C^{G}_T(U,E)$ (as well as $C^{G}_{(S,T]}(U,E)$) as
\[\Hom_E(((F^G(U,E)_{T})_{hG},E[-d])\simeq \Hom_E(F^G(U,E)_{T},E[-d])^{hG}\]
where we equip $\Hom_E(F^G(U,E)_{T},E[-d])$ the $G$-action induced from $F^G(U,E)_{T}$. All of them have the same module structure over $A=\End_G(E)$. But for $C^{G}_{\infty}(U,E)$, we should be careful that it takes the colimit over $T$ after equivariant cohomologies.
\end{RMK}

When the coefficient category $\cE$ admits a complex orientation, i.e.~$H^*(BS^1; E) = E^*[u]$ for the monoidal unit $E$, we can prove the following Gysin sequence by \cite[Theorem 1.7]{Li-Zhang26-1}:

\begin{Thm}\label{thm: Gysin}
    Let $U \subset T^*M$ be an open subset and $L\in \bR\cup\infty$.  When $\cE$ is a complex-oriented compactly generated stable category, there is a $C_\ell$-Gysin sequence between the homotopy $S^1$ and $C_\ell$ orbits and fixed points
    \begin{gather*}
    F_L(U; E)_{hC_\ell} \to F_L(U; E)_{hS^1} \xrightarrow{[\ell](u)} F_L(U; E)_{hS^1}[2],\\
    F_L(U; E)^{hS^1}[-2] \xrightarrow{[\ell](u)} F_L(U; E)^{hS^1} \to F_L(U; E)^{hC_\ell} .
    \end{gather*}
    where $u$ is the first Chern class in $\widetilde{H}^2(BS^1; E)$ and $[\ell]$ is the degree $\ell$ map on $S^1$. Similarly, there is a Gysin sequence between the duals
    \begin{gather*}C_L(U; E)_{hC_\ell} \to C_L(U; E)_{hS^1} \xrightarrow{[\ell](u)} C_L(U; E)_{hS^1}[2],\\
    C_L(U; E)^{hS^1}[-2] \xrightarrow{[\ell](u)} C_L(U; E)^{hS^1} \to C_L(U; E)^{hC_\ell} .
    \end{gather*}

    The same statement holds for $C^{G}_{(S,T]}(U,E)$
\end{Thm}

We also have the Künneth formula for the Chiu--Tamarkin invariants following from \cite[Theorem 1.8]{Li-Zhang26-1}:

\begin{Thm}\label{thm:Künneth tamarkin}
Let $U \subset T^*M$ and $W \subset T^*N$ be open subsets and $\cE$ be a complex-oriented compactly generated stable category. Then there is a Künneth isomorphism
\[F(U \times W; E)^{hS^1} \simeq F(U; E)^{hS^1} \widehat\star_{1_\sT^{hS^1}} F(W; E)^{hS^1}.\]
Here $\widehat\star_{1_\sT^{hS^1}}$ means the convolution valued in $D(E^{hS^1})_{comp}$.

\end{Thm}

\section{Geometry of the traces of Tamarkin categories}\label{section: Geometry of the trace of Tamarkin category}

In this section, we realize the (algebraic) circle action on the Hochschild homology of sheaf categories and Tamarkin categories as geometric circle rotations on closed Reeb orbits.

\subsection{Geometric circle action on traces of sheaf categories}\label{section: Geometric circle action}
Now, we will further investigate traces of Tamarkin categories, and show that the circle action on traces in fact comes from a geometrically defined circle action by loop rotations.

First, we recall that Hamiltonian isotopies admit sheaf quantizations. For a Hamiltonian isotopy $\varphi_H: I \times T^*M \to T^*M$ defined by the Hamiltonian function $H: I \times T^*M \to \bR$, we define the Lagrangian movie to be
\[\Lambda_H = \{(t, -H(t, \varphi_H^t(x, \xi)), x, -\xi, \varphi_H^t(x, \xi)) \mid t\in I, (x, \xi) \in T^*M\} \subset T^*(I \times M \times M).\]

\begin{Thm}[{\cite[Theorem 3.7 \& Proposition 4.8]{GKS2012} \cite[Proposition 3.9]{kuo2021wrapped}}]\label{thm:GKS}
    For a homogeneous Hamiltonian isotopy $\varphi_H$ on $T^*M$, there is a unique sheaf $K(\varphi_H) \in \Sh(I \times M \times M; \cE)$ such that $\dot{\ms}(K(\varphi_H)) \subset \Lambda_H$ and $K(\varphi_H)_0 = 1_{\Delta_M}$. Moreover, for Hamiltonian isotopies $\varphi_H$ and $\varphi_{H'}$ with $H' \leq H$, there is a canonical continuation morphism $K(\varphi_{H'})|_{I_{\geq 0}} \to K(\varphi_H)|_{I_{\geq 0}}$ for $I_{\geq 0}=I\cap [0,\infty)$.
\end{Thm}

Then, we can recall the wrapping formula of the microlocal kernels $Q_Z$ and $P_Z$, first established in \cite[Theorem 1.2]{kuo2021wrapped}, in our setting:

\begin{Thm}[{\cite[Theorem 1.2]{kuo2021wrapped}}]\label{thm:wrapping}
Let $U \subset T^*M$ be an open set, and let $H_n$ be any increasing sequence of compactly supported smooth Hamiltonians which are supported on $U$ and such that $H_n(u) \to \infty$ as $n \to \infty$ for any $u \in U$. Then
\begin{equation}
    Q_Z = \varinjlim\nolimits_{n\in \bN}K(\varphi_{H_n})_1,
\end{equation}
where $K(\varphi_{H_n})$ is the sheaf quantization of the Hamiltonian flow $\varphi_{H_n}$, and the morphism $1_{\Delta_M} \to Q_Z$ is induced by the continuation morphisms $1_{\Delta_M} \to K(\varphi_{H_n})_1$.
\end{Thm}
\begin{RMK}\label{rmk: idempotence}
The universal property of $Q_Z$ implies that $Q_Z$ is an idempotent such that $Q_Z \circ Q_Z = Q_Z$, and under the wrapping formula, the continuation maps intertwine the compositions of $K(\varphi_H)_1$ and the compositions of $Q_Z$.
\end{RMK}

We will need the following microsupport estimation, which relates the microsupport of the internal Hom of sheaf quantizations with Reeb orbits. Recall that for a homogeneous Hamiltonian isotopy $\varphi_H$ on $\dot T^*M$, we define $\Spec(\varphi_H) \subset \bR$ to be the spectrum of periods of all periodic orbits of $\varphi_H$.

\begin{Prop}[{\cite[Section 4.3]{Ike_2019}, \cite[Theorem 4.3]{Kuo-Li-Duality}}]\label{prop:microsupport orbit}
Let $H$ be a positive homogeneous Hamiltonian on $T^*M$ and $K(\varphi_H) \in \Sh(I \times M \times M; \cE)$ be the sheaf quantization of the Hamiltonian flow $\varphi_H$. Then
\[\dot{\ms}(\pi_{M\times M*}\sHom(K(\varphi_H), 1_{\Delta \times I})) \subset T^*_+\bR|_{\Spec(\varphi_H)}.\]
Moreover, when the spectrum $\Spec(\varphi_H)$ is discrete, for sufficiently small $\epsilon$, there is a fiber sequence
\[\Hom(K(\varphi_H)_{s+\epsilon}, 1_\Delta) \to \Hom(K(\varphi_H)_{s-\epsilon}, 1_\Delta) \to \dot{\pi}_{\dot T^*M,*}\mhom(K(\varphi_H)_s, 1_\Delta),\]
where microsupport of $\mhom$ is the union of $s$-periodic orbits of $\varphi_H$ in the conormal of the diagonal.
\end{Prop}

The geometric $S^1$-action will be defined by rotating the periodic orbits of the Hamiltonian flow, which arise as singular supports of the sheaf $\sHom(K(\varphi_H), 1_{\Delta \times I})$; such rotation is realized by the simultaneous action of the Hamiltonian flow on both components. In order to realize this geometric rotation, we need to consider the parametric family of sheaf kernels on $\bR \times M \times M$.

We introduce the following notations. Let $z_i: \bR^n \to \bR$ be the $i$-th projection map, and by abuse of notation $z_i: \bR^n \times M \times M \to \bR \times M \times M$ the $i$-th projection map as well.

\begin{Lemma}\label{lem: projection map}
Let $H$ be a time-independent positive homogeneous Hamiltonian on $T^*M$ and $K(\varphi_H) \in \Sh(\bR \times M \times M; \cE)$ be the sheaf quantization of the Hamiltonian flow $\varphi_H$. Then there is a natural isomorphism
\begin{align*}
z_1^*(\pi_{M\times M*}\sHom(K(\varphi_H), 1_{\Delta\times \bR})) \xrightarrow{\sim} \pi_{M\times M*}\sHom(z_1^*K(\varphi_H) \circ_{\bR^2} z_0^*K(\varphi_H), z_1^*K(\varphi_H)).
\end{align*}
\end{Lemma}
\begin{proof}
By \autoref{prop:microsupport orbit}, we can show that the sheaf $\sHom(z_1^*K(\varphi_H) \circ_{\bR^2} z_0^*K(\varphi_H), z_1^*K(\varphi_H))$ is $z_1$-non-characteristic. Hence
\[{\ms}(z_{1*}\pi_{M\times M*}\sHom(z_1^*K(\varphi_H) \circ_{\bR^2} z_0^*K(\varphi_H), z_1^*K(\varphi_H))) \subset 0_{\bR}.\]
Then the pull-back along $z_1: \bR^2 \rightarrow \bR$ is an isomorphism by \cite[Corollary 2.7.7]{KS90}.
\end{proof}

\begin{Lemma}\label{lem: action map}
Let $H$ be a time-independent positive homogeneous Hamiltonian on $T^*M$ and $K(\varphi_H) \in \Sh(\bR \times M \times M; \cE)$ be the sheaf quantization of the Hamiltonian flow $\varphi_H$. Then the natural action by $K(\varphi_H)_\bR$ induces another isomorphism
\begin{gather*}
s^*( \pi_{M\times M*}\sHom(K(\varphi_H), 1_{\Delta\times \bR})) \xrightarrow{\sim} \pi_{M\times M*}\sHom(z_1^*K(\varphi_H) \circ_{\bR^2} z_0^*K(\varphi_H), z_1^*K(\varphi_H)),
\end{gather*}
where $s: \bR \times \bR \to \bR$ is the addition map.
\end{Lemma}
\begin{proof}
By \autoref{prop:microsupport orbit}, we can show that the sheaf $\sHom(z_1^*K(\varphi_H) \circ_{\bR^2} z_0^*K(\varphi_H), z_1^*K(\varphi_H))$ is $s$-non-characteristic. Hence
\[{\ms}(s_{*}\pi_{M\times M*}\sHom(z_1^*K(\varphi_H) \circ_{\bR^2} z_0^*K(\varphi_H), z_1^*K(\varphi_H))) \subset 0_{\bR}.\]
Then the pull-back along $s: \bR^2 \rightarrow \bR$ is an isomorphism by \cite[Corollary 2.7.7]{KS90}.
\end{proof}

The above lemmas can be viewed as the projection map and action map by $\varphi_H$ on the sheaf Hom respectively. They define a topological $\bR$-action:

\begin{Prop}\label{prop: R action}
Let $H$ be a time-independent positive homogeneous Hamiltonian on $T^*M$ and $K(\varphi_H) \in \Sh(\bR \times M \times M; \cE)$ be the sheaf quantization of the Hamiltonian flow $\varphi_H$. Then the simplicial diagram, regarded as a cartesian sheaf on the simplicial space $B_\bullet\bR$
\begin{align*}
\cdots & \substack{\longrightarrow \\[-0.8em] \longrightarrow \\[-0.8em] \longrightarrow \\[-0.8em] \longrightarrow}\, \pi_{M\times M*}\sHom(z_1^*K(\varphi_H) \circ_{\bR^3} z_2^*K(\varphi_H) \circ_{\bR^3} z_0^*K(\varphi_H), z_1^*K(\varphi_H) \circ_{\bR^3} z_2^*K(\varphi_H)) \\
&\,\substack{\longrightarrow \\[-0.8em] \longrightarrow \\[-0.8em] \longrightarrow}\, \pi_{M\times M*}\sHom(z_1^*K(\varphi_H) \circ_{\bR^2} z_0^*K(\varphi_H), z_1^*K(\varphi_H)) \,\substack{\longrightarrow \\[-0.8em] \longrightarrow}\, \pi_{M\times M*}\sHom(K(\varphi_H), 1_{\Delta\times\bR})
\end{align*}
defines a natural $\bR$-action on $\pi_{M\times M*}\sHom(K(\varphi_H), 1_{\Delta\times\bR})$, where the projection map is defined by Lemma \ref{lem: projection map}, and the action map is induced by continuation map by left action of $K(\varphi_H)_\bR$.
\end{Prop}
\begin{proof}
We exhibit the $\bR$-action as a consequence of the parametric sheaf quantization \autoref{thm:GKS} of contact isotopies \cite[Remark 4.9]{GKS2012}. In fact, by writing the Hamiltonian flow $\varphi_H$ as parametric Hamiltonian flow defined on the simplicial diagram $B_\bullet\bR$, the parametric sheaf quantization theorem of contact isotopies implies that we have a coherent diagram of sheaf quantizations
\begin{align*}
\cdots \substack{\longrightarrow \\[-0.8em] \longrightarrow \\[-0.8em] \longrightarrow \\[-0.8em] \longrightarrow }\, z_1^*K(\varphi_H) \circ_{\bR^3} z_2^*K(\varphi_H) \circ_{\bR^3} z_3^*K(\varphi_H) \substack{\longrightarrow \\[-0.8em] \longrightarrow \\[-0.8em] \longrightarrow }\, z_1^*K(\varphi_H) \circ_{\bR^2} z_2^*K(\varphi_H) \,\substack{\longrightarrow \\[-0.8em] \longrightarrow }\, K(\varphi_H) 
\end{align*}
Consider the action of (parametric) integral kernels $K(\varphi_H)$ on the internal Homs by the addition 
\[\pi_{M\times M*}\sHom(K(\varphi_H), 1_{\Delta \times \bR}) \mapsto \pi_{M\times M*}\sHom(z_1^*K(\varphi_H) \circ_{\bR^2} z_2^*K(\varphi_H), z_1^*K(\varphi_H)),\]
we obtain the desired diagram, which defines a topological $\bR$-action.
\end{proof}

\begin{Prop}\label{prop: trivial Z action}
Let $H$ be a time-independent positive homogeneous Hamiltonian on $T^*M$ with discrete spectrum and $K(\varphi_H) \in \Sh(\bR \times M \times M; \cE)$ be the sheaf quantization of the Hamiltonian flow $\varphi_H$. Then the $\bZ$-action induced by the addition map
\begin{align*}
s^*: \pi_{M\times M*}\sHom(K(\varphi_H), 1_{\Delta\times \bR}) \to \pi_{M\times M*}\sHom(K(\varphi_H) \circ K(\varphi_H), K(\varphi_H))
\end{align*}
admits a canonical trivialization.
\end{Prop}
\begin{proof}
Consider the $\bZ$-action defined by the restriction of the $\bR$-action
\begin{align*}
&\substack{\longrightarrow \\[-0.8em] \longrightarrow \\[-0.8em] \longrightarrow \\[-0.8em] \longrightarrow}\,\pi_{M\times M*}\sHom(K(\varphi_H) \circ K(\varphi_H) \circ K(\varphi_H), K(\varphi_H) \circ K(\varphi_H)) \\
&\substack{\longrightarrow \\[-0.8em] \longrightarrow \\[-0.8em] \longrightarrow}\,\pi_{M\times M*}\sHom(K(\varphi_H) \circ K(\varphi_H), K(\varphi_H)) \,\substack{\longrightarrow \\[-0.8em] \longrightarrow}\, \pi_{M\times M*}\sHom(K(\varphi_H), 1_{\Delta\times\bR}).
\end{align*}
Our main observation is that the right action by the kernel $K(\varphi_H)^{- n}$ on each simplex commutes with the left actions by the addition map. Therefore, after applying the right actions $K(\varphi_H)^{- n}$ simultaneously on each simplex, we have a diagram
\begin{align*}
\substack{\longrightarrow \\[-0.8em] \longrightarrow \\[-0.8em] \longrightarrow \\[-0.8em] \longrightarrow}\,\pi_{M\times M*}\sHom(K(\varphi_H), 1_{\Delta\times \bR}) \substack{\longrightarrow \\[-0.8em] \longrightarrow \\[-0.8em] \longrightarrow}\,\pi_{M\times M*}\sHom(K(\varphi_H), 1_{\Delta\times \bR}) \,\substack{\longrightarrow \\[-0.8em] \longrightarrow}\, \pi_{M\times M*}\sHom(K(\varphi_H), 1_{\Delta\times\bR}),
\end{align*}
where the addition map $s: \bR \times \bR \to \bR$ now induces the identity. Hence we obtain a trivialization of the $\bZ$-action.
\end{proof}

\begin{Thm}\label{thm:geometric circle action}
Let $H$ be a time-independent positive homogeneous Hamiltonian on $T^*M$ with discrete spectrum and $K(\varphi_H) \in \Sh(\bR \times M \times M; \cE)$ be the sheaf quantization of the Hamiltonian flow $\varphi_H$. Then the simplicial diagram
\begin{align*}
\cdots & \substack{\longrightarrow \\[-0.8em] \longrightarrow \\[-0.8em] \longrightarrow \\[-0.8em] \longrightarrow}\, \pi_{M\times M*}\sHom(z_1^*K(\varphi_H) \circ_{\bR^3} z_2^*K(\varphi_H) \circ_{\bR^3} z_0^*K(\varphi_H), z_1^*K(\varphi_H) \circ_{\bR^3} z_2^*K(\varphi_H)) \\
&\,\substack{\longrightarrow \\[-0.8em] \longrightarrow \\[-0.8em] \longrightarrow}\, \pi_{M\times M*}\sHom(z_1^*K(\varphi_H) \circ_{\bR^2} z_0^*K(\varphi_H), z_1^*K(\varphi_H)) \,\substack{\longrightarrow \\[-0.8em] \longrightarrow}\, \pi_{M\times M*}\sHom(K(\varphi_H), 1_{\Delta\times \bR}).
\end{align*}
defines a geometric $\bR$-action that descends to a $S^1$-action.
\end{Thm}
\begin{proof}Apply \autoref{prop: R action} and \autoref{prop: trivial Z action}.
\end{proof}

\begin{Coro}\label{cor:geometric circle action}
Let $Z \subset T^*M$ be a closed conic subset. Then there is a natural $S^1$-action on the dual Hochschild homology
\[\begin{tikzcd}
\HH_*(\Sh_Z(M))^\vee \xrightarrow{\sim} \Hom(Q_Z, \Delta_*\omega_M).
\end{tikzcd}\]
\end{Coro}
\begin{proof}
By \autoref{thm:GKS}, we know that for any $H \leq H'$, there is a natural continuation morphism $K(\varphi_{H})_1 \to K(\varphi_{H'})_1$. This induces the continuation morphism
\[\Hom(K(\varphi_{H'})_1, \Delta_*\omega_M) \to \Hom(K(\varphi_H)_1, \Delta_*\omega_M).\]
The continuation morphism is compatible with the $S^1$-action because for parametric families of positive Hamiltonians, the continuation morphism is canonical.
Then using \autoref{thm:wrapping}, considering smooth Hamiltonians supported on $U = T^*M \setminus Z$ such that $H_n(u) \to \infty$ for any $u \in U$, we get a geometric $S^1$-action (the idempotence of $Q_Z$ \autoref{rmk: idempotence} shows the commutativity of the simplicial diagram of sheaves over $BS^1$)
\[\Hom(Q_Z, \Delta_*\omega_M) = \lim_{n \to \infty}\Hom(K(\varphi_{H_n})_1, \Delta_*\omega_M).\]

Finally, it suffices to notice that $\Hom(\HH_*(\Sh_Z(M)), 1) = \Hom(\pi_!\Delta^*Q_Z, 1) = \Hom(Q_Z, \Delta_*\omega_M)$ by standard adjunction.
\end{proof}

\subsection{Geometric and algebraic circle action on traces}
\newcommand{\RelTr}{\operatorname{RelTr}}

We set $\Delta^n=\{z\in [0,1]^{n}\mid \sum_{i=1
}^n z_i\leq 1\}$ and $\Sigma^n=\{z\in \Delta^{n}\mid \sum_{i=1
}^n z_i = 1\}$. Define face maps to be
\[d_i: \Delta^{n-1} \to \Delta^n, \quad (z_1, \dots, z_{n-1}) \mapsto (z_1, \dots, z_{i-1}, 0, z_{i}, \dots, z_{n-1}),\]
and the degeneration maps to be 
\[s_i: \Delta^n \to \Delta^{n-1}, \quad (z_1, \dots, z_n) \mapsto (z_1, \dots, z_{i-1} + z_{i}, \dots, z_n).\]
These maps define a cosimplicial topological space. Let the cyclic rotation be
\[t_n: \Delta^n \to \Delta^n, \quad (z_1, \dots, z_{n-1}, z_n) \mapsto (z_n, z_1, \dots, z_{n-1}).\]
These maps now define a cocyclic topological space. These maps admit natural linear extensions, which we denote using the same notations $d_i: \bR^{n-1} \to \bR^n$, $s_i: \bR^n \to \bR^{n-1}$, and $t_n: \bR^n \to \bR^n$, and they restrict to maps between $\Sigma^n$ as well.

\begin{Prop}\label{prop: circle action Ham}
Let $\varphi_H: \bR \times \dot T^*M \to \dot T^*M$ be a homogeneous Hamiltonian flow on $T^*M$. Then there is a cyclic object $\Tr(K(\varphi_H)_1,-): \boldsymbol{\Lambda} \to \cE$ such that
\[\Tr(K(\varphi_H)_1, n) = \Gamma(\Sigma^n, (\id_{\bR^n}\times \pi_{n})_! (\id_{\bR^n}\times \widetilde{\Delta}_n)^* (K(\varphi_H) \boxtimes \cdots \boxtimes K(\varphi_H)),\]
such that the face maps are induced by the pull-back along $d_i$, the degeneracy maps are induced by the pull-back along $s_i$, and the cyclic rotation is induced by $t_n: (\bR \times M \times M)^n \to (\bR \times M \times M)^n$.
\end{Prop}
\begin{proof}
We consider the following sequence of sheaves on $\bR^n$, which we call the {\emph relative trace}
\[\RelTr(K(\varphi_H), \bR^n)\coloneqq (\id_{\bR^n}\times \pi_{n})_! (\id_{\bR^n}\times \widetilde{\Delta}_n)^* \left(K(\varphi_H) \boxtimes \cdots \boxtimes K(\varphi_H) \right)\in \Sh(\bR^n,\cE).\]
Recall that the parametric sheaf quantization \autoref{thm:GKS} of contact isotopies \cite[Remark 4.9]{GKS2012} shows $K(\varphi_H)$ defines an $\bR$-action. By writing the Hamiltonian flow $\varphi_H$ as parametric Hamiltonian flow defined on the simplicial space $B_\bullet\bR$, the parametric sheaf quantization theorem of contact isotopies implies that we have a coherent diagram of sheaf quantizations
\begin{align*}
\cdots \substack{\longrightarrow \\[-0.8em] \longrightarrow \\[-0.8em] \longrightarrow \\[-0.8em] \longrightarrow }\, z_1^*K(\varphi_H) \circ_{\bR^3} z_2^*K(\varphi_H) \circ_{\bR^3} z_3^*K(\varphi_H) \substack{\longrightarrow \\[-0.8em] \longrightarrow \\[-0.8em] \longrightarrow }\, z_1^*K(\varphi_H) \circ_{\bR^2} z_2^*K(\varphi_H) \,\substack{\longrightarrow \\[-0.8em] \longrightarrow }\, K(\varphi_H) 
\end{align*}
Then, since we have that
\begin{align*}
&\,(\id_{\bR^n}\times \pi_{n})_! (\id_{\bR^n}\times \widetilde{\Delta}_n)^* \left(K(\varphi_H) \boxtimes \cdots \boxtimes K(\varphi_H) \right) \\
=&\,(\id_{\bR^n}\times \pi_{n})_! (\id_{\bR^n}\times \widetilde{\Delta}_n)^* \left(z_1^*K(\varphi_H) \boxtimes_{\bR^n} \cdots \boxtimes_{\bR^n} z_n^*K(\varphi_H) \right) \\
=&\, (\id_{\bR^n} \times \pi_M)_!(\id_{\bR^n} \times \Delta_M)^*(z_1^*K(\varphi_H) \circ_{\bR^n} \dots \circ_{\bR^n} z_n^*K(\varphi_n)),
\end{align*}
we can apply the functor $(\id_{\bR^n} \times \pi_M)_!(\id_{\bR^n} \times \Delta_M)^*$ on all simplices. For any (linear extension of) simplicial structure map $d_i$ or $s_i$, we have
\[d_i^*\RelTr(K(\varphi_H), \bR^{n})\xrightarrow{\simeq } \RelTr(K(\varphi_H), \bR^{n-1}), \quad s_i^*\RelTr(K(\varphi_H), \bR^{n-1})\xrightarrow{\simeq } \RelTr(K(\varphi_H), \bR^n).\]
Similarly, for the cyclic rotation map $t_n$, we have
\[t_n^*\RelTr(K(\varphi_H), \bR^n)  \xrightarrow{\simeq }  \RelTr(K(\varphi_H), \bR^n).\]
Thus, we obtain a cartesian sheaf on the cocyclic space by the coherence of pull-back functors of sheaves. 

Given the relative trace, which defines a cartesian  sheaf on the cocyclic space for $\bR^n$, we can form a cocyclic object by taking global sections on each simplex. Then by proper base change formula and coherence of six-functors, we conclude that
\[\Tr(K(\varphi)_1, n) = \Gamma(\Sigma^n, \RelTr(K(\varphi_H), \bR^n))\]
defines a cocyclic object $\boldsymbol{\Lambda} \to \cE$.
\end{proof}

\begin{Prop}Let $\varphi_H: \bR \times \dot T^*M \to \dot T^*M$ be a homogeneous Hamiltonian flow on $T^*M$. The cocyclic object $\Tr(K(\varphi_H)_1,-)$ is a cartesian cocylic object.
\end{Prop}
\begin{proof}
Let $s_n: \bR^n \to \bR, (z_1, \dots, z_n) \mapsto \sum_{i=1}^n z_i$. Consider the identification
\[\RelTr(K(\varphi_H), \bR^n)\simeq s_n^*\pi_{M!}\Delta_M^*K(\varphi_H).\]
Since $s_{n-1} = s_n \circ d_i$ and $s_n = s_{n-1} \circ s_i$, any simplicial structural map $d_i$ or $s_i$ that induces
\[\Gamma(\Sigma^{n-1}, d_i^*\RelTr(K(\varphi_H),\bR^n)) \to \Gamma(\Sigma^n, \RelTr(K(\varphi_H),\bR^n))\]
is conjugated to the following isomorphism
\[\Gamma(s_{n-1}(\Sigma^{n-1}),\pi_{M!}\Delta_M^*K)\to \Gamma(s_n(\Sigma^n),\pi_{M!}\Delta_M^*K). \]
This shows that $\Tr(K(\varphi_H)_1,-)$ is a cartesian cocyclic object $\widetilde{\boldsymbol{\Lambda}} \to \cE$.
\end{proof}

Our first theorem in this section will show that the circle action in \autoref{thm:geometric circle action} can be realized as the above cyclic trace:

\begin{Thm}
Let $\varphi_H: \bR \times \dot T^*M \to \dot T^*M$ be a homogeneous Hamiltonian flow on $T^*M$. Then under the duality
\[ (\pi_{M!} \Delta_M^*K(\varphi_H)_1)^\vee \simeq \pi_{M\times M*}\sHom(K(\varphi_H)_1, \Delta_{M*}\omega_M),\]
the circle action in \autoref{thm:geometric circle action} and \autoref{prop: circle action Ham} are equivalent.
\end{Thm}
\begin{proof}
Consider the inclusion of simplicial spaces where $\Sigma^n \hookrightarrow \bR^n$ for $n \in \bN$. On each simplex $\Sigma^n$, we have the duality
\begin{align*}
&\,(\id_{\Sigma^n} \times \pi_n)_!(\id_{\Sigma^n}\times \widetilde{\Delta}_n)^*(z_1^*K(\varphi_H) \otimes \dots \otimes z_n^*K(\varphi_H))^\vee \\
\simeq &\,\pi_{M\times M*}\sHom(z_1^*K(\varphi_H) \circ_{\Sigma^n} \dots \circ_{\Sigma^n} z_n^*K(\varphi_H), \pi_{\Sigma^n}^*\Delta_{M*}\omega_M)
\end{align*}
that is compatible with the face and degeneracy maps. Then, we have natural commutative diagrams:
\[\begin{tikzcd}[column sep=16pt]
\Hom(z_1^*K(\varphi_H) \circ \dots \circ z_n^*K(\varphi_H), 1_{\Delta\times \Sigma^n}) \ar[r] \ar[d] & \Hom(z_n^*K(\varphi_H)\circ z_1^*K(\varphi_H) \circ \dots \circ z_{n-1}^*K(\varphi_H), 1_{\Delta\times \Sigma^n}) \ar[d]\\
\Hom(z_1^*K(\varphi_H) \circ \dots \circ z_n^*K(\varphi_H), 1_{\Delta \times \Sigma^n}) \ar[r] & \Hom(z_n^*K(\varphi_H) \circ z_1^*K(\varphi_H) \circ \dots \circ z_n^*K(\varphi_H), z_n^*K(\varphi_H))
\end{tikzcd}\]
where the top horizontal map is induced by cyclic rotation, the bottom horizontal map is induced by the left action of $z_n^*K(\varphi_H)$, and the right vertical map is induced by the right action of $z_n^*K(\varphi_H)$, which commutes with the face and degeneracy maps by the left action.

This shows we have an isomorphism of cartesian sheaves on the cosimplicial spaces $B\boldsymbol{\Sigma}$ formed by the simplices $\Sigma^n$. Moreover, since the inclusion of cosimplicial spaces $B\boldsymbol{\Sigma} \to B\bR$ is a homotopy equivalence, and the cyclic permutation agrees with the geometric circle action, this shows we have an isomorphism of sheaves on $\widetilde{\boldsymbol{\Lambda}} \xrightarrow{\sim} B\boldsymbol{\Sigma}/S^1 \xrightarrow{\sim} BS^1$, where $\widetilde{\boldsymbol{\Lambda}}$ is the groupoid completion of the cyclic category, i.e. the localization of $\boldsymbol{\Lambda}$ with respect to all its morphisms.
\end{proof}

\begin{Coro}\label{cor:algebraic circle action}
For all $Z\subset T^*M$, let $H_n$ be a cofinal sequence on $T^*M\setminus Z$, then we have an equivalence of cocyclic objects (where the right hand side is the cyclic trace in \autoref{thm:S1-action-sheaf} \cite[Theorem 1.1]{Li-Zhang26-1})
\[\varinjlim_{n\to \infty}\Tr(K(\varphi_{H_n})_1,-) \simeq \Tr(Q_Z,-).\]
\end{Coro}
\begin{proof}
For this construction, it is clear by \autoref{thm:GKS} (more precisely \cite[Remark 4.9]{GKS2012} and \cite[Proposition 3.9 \& Remark 3.10]{kuo2021wrapped}) that if $H\leq H'$, and the canonical continuation maps induces a morphism between cocyclic objects
\[\Tr(K(\varphi_{H})_1,-) \to \Tr(K(\varphi_{H'})_1,-).\]

Similarly, for compactly supported autonomous $H\geq 0$ with $H|_Z=0$, by \autoref{thm:wrapping} there exists a canonical morphism
\(c_H:K(\varphi_H)_1 \to Q_Z\)
that is compatible with compositions.
Then using \autoref{thm:wrapping}, considering smooth Hamiltonians supported on $U = T^*M \setminus Z$ such that $H_n(u) \to \infty$ for any $u \in U$, we have a map of cocyclic objects (the idempotence of $Q_Z$ \autoref{rmk: idempotence} shows that the commutativity of the cocyclic diagram)
\[\varinjlim_{n\to\infty}\Tr (K(\varphi_{H_n})_1,-)\to \Tr(Q_Z,-).\]

Finally, notice that the colimit in $\Fun_{\simeq}(\boldsymbol\Lambda,\cE)$ is taken pointwise and the forgetful functor $\Fun_{\simeq}(\boldsymbol\Lambda,\cE)\simeq \Fun(BS^1,\cE)\to \cE$ is conservative. Then we have an equivalence of cocyclic objects.
\end{proof}

Consequently, we can conclude our main theorem, that the geometric $S^1$-action on the trace of the category of sheaves coincides with the algebraic $S^1$-action on the trace of the category of sheaves constructed in \autoref{thm:S1-action-sheaf} \cite[Theorem 1.1]{Li-Zhang26-1}.

\begin{Thm}\label{thm:geometric algebraic action}
Let $Z \subset T^*M$ be a closed conic subset. Then the geometric $S^1$-action in \autoref{cor:geometric circle action} coincides with the algebraic $S^1$-action in \autoref{cor:algebraic circle action} on the dual Hochschild homology $\HH_*(\Sh_Z(M))^\vee$ in \autoref{thm:S1-action-sheaf} \cite[Theorem 1.1]{Li-Zhang26-1}.
\end{Thm}

\subsection{Geometric circle action on traces of Tamarkin categories}\label{sec: geometric circle tamarkin}
We first recall the sheaf quantization results of \cite{GKS2012}.

\begin{Thm}\label{GKS tamarkin version}Let $(I,0)$ be a pointed contractible manifold. Given a complete Hamiltonian isotopy $\varphi: I \times T^*M \to T^*M$, 
    there is a unique $\cK(\varphi) \in \Sh\left(I, \sT\left(T^*(M\times M) \right) \right)$ such that 
    $\cK(\varphi)|_0 \cong 1_{\Delta_{M}}^\sT$, and  
    \begin{equation}\label{equatio: GKS T-linear total graph}
        \underline{ss}(\cK(\varphi)) = \{(z, - H(z,{\varphi}_{z}(q,p))  ,(q,-p),{\varphi}_{z}(q,p) ) : z\in I, (q,p)\in T^*M\}.
    \end{equation}

    Moreover, if $\varphi, \varphi'$ are generated by complete Hamiltonians with $\varphi' \leq  \varphi$, then
    there is a canonical map $\cK(\varphi')|_{1}\rightarrow \cK(\varphi)|_{1}$.
\end{Thm}
\begin{RMK}In the compactly supported case, this is proven in {\cite[Theorem 3.7 \& Proposition 4.8 \& Proposition A.6]{GKS2012}} as a corollary of their quantization for homogeneous Hamiltonian and some geometric observations. See also \cite[Corollary 2.3.2]{guillermou2019sheaves}. The proof for a complete flow will be explained in \autoref{appendix: existence kernel}. Notice that the continuation statement is not explained in the appendix, but follows from the same argument as \cite[Corollary 5.23]{Hochschild-Kuo-Shende-Zhang}.
    
\end{RMK}

\begin{Coro}[{\autoref{appendix: Hamiltonian embedding}}]Let $\varphi: I\times U\to T^*M$ be an isotopy of Hamiltonian embedding (which is not necessarily defined on $T^*M$),  
    there is an equivalence $\sT(U)\xrightarrow{\sim} \sT({\varphi}_1(U))$. And when $\varphi_z$ is globally defined $T^*M$, the functor is induced by $\cK(\varphi)_1\circ_\sT - $.    
\end{Coro}

We can recall the wrapping formula of the microlocal kernels $Q_U$ and $P_U$, first established in \cite[Theorem 1.2]{kuo2021wrapped}, in our setting:

\begin{Thm}[{\cite[Corollary 6.6]{Hochschild-Kuo-Shende-Zhang}}]\label{theorem: wrapping formula Tamarkin Q}
Let $U \subset T^*M$ be an open set, and let $H_n$ be any increasing sequence of non-negative smooth Hamiltonians which are supported on $U$ and such that $H_n(u) \to \infty$ as $n \to \infty$ for any $u \in U$. Then
\begin{equation}
    Q_U = \varinjlim\nolimits_{n\in \bN}\cK(\varphi_{H_n})_1,
\end{equation}
where $\cK(\varphi_{H_n})$ is the sheaf quantization of the Hamiltonian flow $\varphi_{H_n}$, and the morphism $1_{\Delta_M} \to Q_U$ is induced by the continuation morphisms $1_{\Delta_M} \to \cK(\varphi_{H_n})_1$.
\end{Thm}

Similar to the previous subsection, we can construct a geometric circle action on the Hochschild cohomology. The geometric $S^1$-action will be defined by rotating the periodic orbits of the Hamiltonian flow, which arise as singular supports of the sheaf $\Hom^\sT(\cK(\varphi_H), 1_{\Delta \times I}^\sT)$; such rotation is also realized by the simultaneous action of the Hamiltonian flow on both components. 

\begin{Thm}\label{thm:geometric circle action tamarkin}
Let $H$ be a time-independent positive Hamiltonian on $T^*M$ with discrete spectrum and $\cK(\varphi_H) \in \Sh(I \times M \times M; \sT)$ be the sheaf quantization of the Hamiltonian flow $\varphi_H$. Then the simplicial diagram
\begin{align*}
\cdots & \substack{\longrightarrow \\[-0.8em] \longrightarrow \\[-0.8em] \longrightarrow \\[-0.8em] \longrightarrow}\, \pi_{M\times M*}^\sT\sHom(z_1^*\cK(\varphi_H) \star_{\bR^3} z_2^*\cK(\varphi_H) \star_{\bR^3} z_0^*\cK(\varphi_H), z_1^*\cK(\varphi_H) \star_{\bR^3} z_2^*\cK(\varphi_H)) \\
&\,\substack{\longrightarrow \\[-0.8em] \longrightarrow \\[-0.8em] \longrightarrow}\, \pi_{M\times M*}^\sT\sHom(z_1^*\cK(\varphi_H) \star_{\bR^2} z_0^*\cK(\varphi_H), z_1^*\cK(\varphi_H)) \,\substack{\longrightarrow \\[-0.8em] \longrightarrow}\, \pi_{M\times M*}^\sT\sHom(\cK(\varphi_H), 1_{\Delta_M}^\sT).
\end{align*}
defines a geometric $\bR$-action that descends to a $S^1$-action.
\end{Thm}
\begin{proof}
The proof is essentially the same as \autoref{thm:geometric circle action}.
\end{proof}

\begin{Coro}\label{cor:geometric circle action tamarkin}
Let $Z \subset T^*M$ be a closed subset and $U = T^*M \setminus Z$. Then there is a natural $S^1$-equivariant morphism where on both sides the $S^1$-actions are obtained by loop rotations
\[\begin{tikzcd}
\HH_*^\sT(\sT(U))^\vee \xrightarrow{\sim} \Hom^\sT(P_U, \Delta_*\omega_{M}^\sT).
\end{tikzcd}\]
\end{Coro}
\begin{proof}
Using the 6-operations and orientation as in \cite[Corollary 7.3]{Hochschild-Kuo-Shende-Zhang}, we have 
\[\HH_*^\sT(\sT(U))^\vee \simeq \HH^*_\sT(\sT(U), \omega_M^\sT) \simeq \Hom^\sT(P_U, \Delta_*\omega_{M}^\sT).\]
The rest of the proof is then similar to \autoref{cor:geometric circle action}.
\end{proof}

Again, we can show that the geometric circle action defined by loop rotation agrees with the algebraic circle action on the Hochschild homology.

\begin{Prop}\label{prop: circle action Ham Tamarkin}
Let $\varphi_H: \bR \times  T^*M \to  T^*M$ be a complete Hamiltonian flow on $T^*M$. Then there is a cartesian cyclic object $\Tr(\cK(\varphi_H)_1,-): \boldsymbol{\Lambda} \to \sT$ such that
\[\Tr(\cK(\varphi_H)_1, n) = \Gamma^\sT(\Sigma^n, (\id_{\bR^n}\times \pi_{n}^\sT)_! (\id_{\bR^n}\times \widetilde{\Delta}_n^\sT)^* (\cK(\varphi_H) \boxtimes \cdots \boxtimes \cK(\varphi_H)),\]
such that the face maps are induced by the pull-back along $d_i$, the degeneracy maps are induced by the pull-back along $s_i$, and the cyclic rotation is induced by $t_n: (\bR \times M \times M)^n \to (\bR \times M \times M)^n$.
\end{Prop}
\begin{proof}
The proof is essentially the same as \autoref{prop: circle action Ham}.
\end{proof}

\begin{Coro}\label{cor:algebraic circle action Tamarkin}
For all $Z\subset T^*M$, let $H_n$ be a cofinal sequence on $Z$, then we have an equivalence of cocyclic objects (where the right hand side is the cyclic trace in \autoref{thm:S1-action} \cite[Theorem 1.5]{Li-Zhang26-1})
\[\varinjlim_{n\to \infty}\Tr(\cK(\varphi_{H_n})_1,-) \simeq \Tr(Q_U,-).\]
\end{Coro}
\begin{proof}
The proof is the same as \autoref{cor:algebraic circle action}.
\end{proof}

In conclusion, we can identify both circle actions on the dual trace of Tamarkin categories:

\begin{Thm}\label{thm: geometric circle tamarkin}
Let $Z \subset T^*M$ be a closed subset and $U = T^*M \setminus Z$. Then the geometric $S^1$-action in \autoref{cor:geometric circle action tamarkin} coincides with the algebraic $S^1$-action in \autoref{cor:algebraic circle action Tamarkin} and \autoref{thm:S1-action} \cite[Theorem 1.5]{Li-Zhang26-1} on the dual Hochschild homology $\HH_*(\sT(U))^\vee$.
\end{Thm}
\begin{proof}
The proof is essentially the same as \autoref{thm:geometric algebraic action}.
\end{proof}

\begin{RMK}By the generating function formula of \autoref{appendix:generating function}, we can also define an $S^1$-equivariant generating function cohomology using the same cyclic construction.   After evaluating at the $[\ell]$-place of the cyclic category by \autoref{thm:Z/ell-action}, we recover the $C_\ell$-equivariant generating function cohomology (both the domain version and the Hamiltonian version) introduced in \cite{large_non_equeezing_GF} from the cyclic trace construction introduced here.
\end{RMK}

\subsection{Conley--Zehnder index and spectrality}\label{section:CZ index and spectrality}
We have seen that the cyclic trace of the Tamarkin category admits a natural $\bR$-filtration whose associated graded is the cohomology of Morse--Bott Reeb orbits. In this section, we relate the degrees of the object with Conley--Zehnder indices of the orbits, using the work of Ike \cite[Section 4.4]{Ike_2019}.

Let $\omega_{T^*M}$ be the standard symplectic form on $T^*M$. Recall that we say a bounded domain $U \subset T^*M$ has restricted contact boundary (RCT) if there is a vector field $X$ defined on $\overline{U}$ and transversely outward pointing on $\partial U$  such that \[L_X\omega_{T^*M} = \omega_{T^*M}.\]
Under this condition, the restriction $\iota_X\omega|_{\partial U}$ defines a contact form on $\partial U$. The dynamics of the Reeb vector field $R$ on $\partial U$ is called Morse-Bott non-degenerate if for any action value $T>0$, the space of parameterized $T$-periodic orbits
\[N_T=\{x\in \partial U\mid R_T(x)=x\}\]
is a smooth submanifold of $\partial U$ and moreover
\[T_xN_T=\ker(d_xR_T-\id),\quad \forall x\in N_T.\]
We say that a submanifold $N$ (for example an $S^1$-orbit) is Morse--Bott non-degenerate if it is a component of the critical point set and the dynamics is non-degenerate in a neighborhood of $N$.

First, we recall the simple relation between Maslov index of Morse--Bott Lagrangian intersections and the Conley--Zehnder index of closed Hamiltonian orbits. Let $\mu_{CZ}$ be the Conley--Zehnder index of a closed orbit \cite{RobbinSalamon}, $\mu_{RS}$ be the Robbin--Salamon Maslov index of Lagrangian intersections \cite{RobbinSalamon}, and $I$ be Seidel's grading of Lagrangian intersections \cite{Seidel-graded} (which we generalize to Morse--Bott intersections following Asano \cite[Appendix B.3]{Ike_2019}).

\begin{Lemma}[{\cite[Proposition 2.2]{CieliebakFloerHoferWysocki-application} \cite[Lemma 3.4]{Bourgeois-Oancea-automomous}}]
Let $H$ be a positive Hamiltonian on $T^*M$. Suppose the time-$1$ orbit $\gamma$ is Morse--Bott non-degenerate. Then 
\begin{align*}\mu_{CZ}(\gamma) &= -\mu_{RS}(\gamma, \Gamma_{\Id}, \Gamma_{\varphi_H^1}) + 1/2 = I(\gamma,\Gamma_{\Id}, \Gamma_{\varphi_H^1}) - n+1\\
&= \mu_{RS}(\gamma, \Gamma_{\varphi_H^1}, \Gamma_\Id) - 2n + 1/2=n - I(\gamma, \Gamma_{\varphi_H^1}, \Gamma_{\Id}).
\end{align*}
\end{Lemma}
\begin{proof}
We identify the cotangent bundle $T^*(M \times M)$ with the symplectic product $T^*M^- \times T^*M$ by the antipodal map of the first factor. For the first line, the first equality is a combination of \cite[Proposition 2.2]{CieliebakFloerHoferWysocki-application} and \cite[Lemma 3.4]{Bourgeois-Oancea-automomous}. The second equality is a generalization of the definition in \cite[Section 2d]{Seidel-graded} by \cite[Equation (B.9)]{Ike_2019}. For the second line, the equalities follow from \cite[Section 2d]{Seidel-graded}.
\end{proof}

Then the result of Asano and Ike \cite[Section 4 \& Appendix B]{Ike_2019} computes the degree of $\mhom$ in terms of the Maslov grading. In our situation, it says:

\begin{Prop}[{Tomohiro Asano and Yuichi Ike \cite[Proposition 4.15 \& Theorem B.6]{Ike_2019}}]\label{prop: mhom CZ index}
Let $H$ be a positive Hamiltonian on $T^*M$ and $K(\varphi_H) \in \Sh(I \times M \times M; \sT)$ be the sheaf quantization of the flow $\varphi_H$. Suppose the time-$1$ orbits $\gamma_j$ with action $L$ are Morse--Bott non-degenerate for $j \in J$.
Then
\[\mhom(\cK(\varphi_H)_1, T_{L*}1_\Delta) \simeq \bigoplus_{j\in J}or_{\gamma_j}[\mu_{CZ}(\gamma_j)-n],\]
where $or_{\gamma}$ is the orientation sheaf of the maximal bundle of $\gamma^*T(M \times M)$ where the local generating function of $\varphi_{H,1}$ is negative definite on $\gamma$.
\end{Prop}
\begin{proof}
We know $\dot{\ms}(\cK(\varphi_H)_1) \cap \dot{\ms}(T_{L*}1_\Delta) = \bigcup_{j\in J}\gamma_j \times L \subset T^*(M \times M) \times \bR$. By \cite[Proposition 4.15 \& Theorem C.6]{Ike_2019}, we can show that
\[\mhom(\cK(\varphi_H)_1, T_{L*}1_\Delta) \simeq \bigoplus_{j\in J}or_{\gamma_j/T^-\gamma_j}[-I(\gamma_j)].\]
We explain the notations: Consider an open neighborhood of $\gamma_j$ such that $\Gamma_\Id \times L$ is sent to the zero section, and the Legendrian lift $\widetilde\Gamma_{\varphi_H^1}$ is sent to the graph of a Morse--Bott function. Let $T^-\gamma_j$ be the maximal subbundle of the normal bundle $\gamma_j \times L$ inside $\Gamma_\Id \times L$ such that the Hessian of the function is negative definite, and $or_{\gamma_j/T^-\gamma_j}$ is its orientation sheaf.

We can explicitly compute the vector bundle as follows. On the time-1 orbit closed orbit $\gamma$ starting at $(q,p)$, we can compute the primitive of the exact Lagrangians $\Gamma_{\varphi_H^1}$ and $\Gamma_\Id$ is a neighborhood of $\gamma$. One can compute directly that the primitive of $\Gamma_{\varphi_H^1}$ is
\[S_H(q,p)=A_H(\gamma) = \int_{\gamma} \lambda_{T^*M} - H(t, \gamma(t))dt,\]
and the primitive of $\Gamma_\Id$ is $0$ \cite[Proposition 9.19]{mcduff2017introST}. 
Hence $T^-\gamma$  is exactly the the maximal subbundle of the normal bundle $\gamma \times L$ inside $\Gamma_\Id \times L$ such that the Hessian of a local formula of $S_H$, which is a local generating function of $\varphi_{H,1}$ in a neighborhood of $\gamma$, is negative definite.
\end{proof}
\begin{RMK}
Note that the coherent orientation on the sheaf theoretic trace should not agree with the coherent orientation in the Morse--Bott model of symplectic homology \cite[Section 4.4]{Bourgeois-Oancea-automomous}, since the sheaf theoretic trace is not isomorphic to  symplectic homology when the manifold is not Spin \cite[Appendix B]{guillermouviterbo_gammasupport}.
\end{RMK}

We now show that $C_T^{S^1}(U, E)$ is generated by constant orbits in $U$ and Morse--Bott families of (parametrized) closed Reeb orbits of action at most $T$, by approximating the Hochschild homology given by the colimit of wrappings by a finite wrapping at any finite action level. 

Recall that $\Spec(\partial U, \lambda)$ is the action or length spectrum of Reeb orbits on the RCT boundary $\partial U$\footnote{We also include $0$, i.e. the length of constant  Reeb orbit, for convenience.}. Under the RCT boundary condition, the minimal action $\inf \Spec(\partial U, \lambda) \setminus\{0\}$ is finite and positive. 

When $U$ has RCT boundary $\partial U$, there is a Hamiltonian function $H$ such that $\partial U = H^{-1}(1)$ is a regular hypersurface. Consider a neighbourhood of $\partial U$ of the form $(0, 1+\epsilon) \times \partial U$ where $H(r, y) = r$ for $(r, y) \in (0, 1+\epsilon) \times \partial U$. 
Following \cite[Section 3]{Bourgeois-Oancea-automomous} and \cite[Section 2]{Cieliebak-Frauenfelder-Oancea2010}, we consider a class of Hamiltonian functions $H_{\alpha}: T^*M \to \bR$ for any $\alpha > 0$ and $\epsilon \ll \alpha^2$ such that
\begin{enumerate}
    \item $H_{\alpha}|_{\partial U \times [1/2, 1-\epsilon]} = 2\alpha r$;
    \item $H_{\alpha}|_{U \setminus \partial U \times [1/2-\epsilon, 1]} \equiv \alpha$, $H_{\alpha}|_{T^*M \setminus U} \equiv 0$;
    \item $H_\alpha|_{\partial U \times [1/2-\epsilon,1/2]} = \alpha-h(r)$, $H_\alpha|_{\partial U \times [1-\epsilon,1]} = h(r)$, $h'(r) > 0$.
\end{enumerate}
Recall that the action of the periodic orbits are exactly given by the values of $H_\alpha$ on the orbits.
When $\alpha \to \infty$, this is a cofinal sequence of Hamiltonians, and thus
\[Q_U = \varinjlim\nolimits_{\alpha \to \infty} \cK(\varphi_{H_{\alpha}})_1.\]
We now prove the following result based on the dynamics analysis from \cite[Proposition 2.5]{Cieliebak-Frauenfelder-Oancea2010}:
\begin{Prop}\label{prop:compute limit at finite action}
Let $U \subset T^*M$ be a domain with RCT boundary $\partial U$ with Morse--Bott non-degenerate Reeb dynamics. Then we have
\[\Gamma_{(-\infty,-L]}(\bR, \sHom^\star(Q_U, 1_\Delta)) \simeq \Gamma_{(-\infty,-L]}(\bR, \sHom^\star(\cK(\varphi_{H_{L+\varepsilon}})_1, 1_\Delta)).\]

Suppose the Reeb orbits with action $L$ are $\gamma_j$ for $j \in J$. We have
\[\mhom(Q_U, T_{L*}1_\Delta) \simeq \mhom(\cK(\varphi_{H_{L+\varepsilon}})_1, T_{L*}1_\Delta) \simeq \bigoplus_{j\in J}or_{\gamma_j}[\mu_{CZ}(\gamma_j)-n].\]
\end{Prop}
\begin{proof}
First, by \autoref{eq: C_T formula}, we will show that 
\[\Gamma_{(-\infty,-L]}(\bR, \Hom^\star(Q_U, 1_\Delta)) \simeq \Gamma_{(-\infty,-L]}(\bR, \Hom^\star(\cK(\varphi_{H_{L+\varepsilon}})_1, 1_\Delta)).\]
Clearly, by the $\sT$-linear 6-operations and \autoref{lemma: hom computation} \cite[Proposition 7.6]{Hochschild-Kuo-Shende-Zhang}, it is equivalent to the statement that
\[\Hom_{E}(\cK(\varphi_{H_{L+\varepsilon}})_1,\nT_{L*}1^\sT_{\Delta_M}) \simeq \varprojlim_{\alpha>L}\Hom_{E}(\cK(\varphi_{H_{\alpha}})_1,\nT_{L*}1^\sT_{\Delta_M}).\]
Consider the cofinal sequence $H_{\alpha}$ for $\alpha>L$ such that
$Q_U = \varinjlim_{\alpha>L} \cK(\varphi_{H_{\alpha}})_1.$
Consider the parametric family of Hamiltonians $H = H_\alpha$.
Let $c$ and $c_\alpha$ be the infimum and supremum of the Reeb length spectrum $\Spec(\partial U, \lambda) \cap (0, \alpha)$. Then the 1-periodic orbits of the Hamiltonian $H_{\alpha}$ are either constant orbits or closed Reeb orbits. Since the action of periodic orbits are just values of $H_\alpha$ on the orbits, they can be classified into the following types:
\begin{enumerate}
    \item Constant orbits in $T^*M \setminus U$ with action $0$;
    \item Constant orbits in $U \setminus \partial U \times (1/2, 1)$ with action $\alpha$;
    \item Reeb orbits around $\partial U \times \{1/2\}$ with action $[\alpha + c/2, \alpha + c_\alpha/2]$;
    \item Reeb orbits around $\partial U \times \{1\}$ with action $[c, c_\alpha]$.
\end{enumerate}
Denote the parametric sheaf quantizations by $\cK(\varphi_H) \in \Sh(\bR_\alpha \times \bR_z \times M \times M; \sT)$. Then 
\[ \Hom_E(\cK(\varphi_{H_{\alpha}})_1,\nT_{L*}1^\sT_{\Delta_M})= \Hom_E(\cK(\varphi_{H})_{\alpha,1},\nT_{L*}1^\sT_{\Delta_M}).\]
We know that the intersection of microsupports of the above sheaves corresponds to closed orbits of action $L$ in \autoref{prop:microsupport orbit}:
\[\dot{\ms}(\cK(\varphi_{H_\alpha})_1) \cap \dot{\ms}(T_{L*}1_{\Delta_M}^\sT) = \bigcup_{j\in J}\gamma_j \times L,\]
where $\gamma_j$ are Reeb orbits of length $L$. They are excluded by our above construction and orbit classification. Since $U$ is relatively compact, we know that $\cK(\varphi_{H})_{1}$ is the constant sheaf supported on $\Delta_M \times [0, +\infty) \times \bR_\alpha$ outside a compact subset. Therefore, by microlocal Morse lemma \cite[Proposition 5.4.19]{KS90}, for any $\alpha, \alpha' > L$, we have
\[\Hom_{E}(\cK(\varphi_{H})_{\alpha,1},\nT_{L*}1^\sT_{\Delta_M}) \xrightarrow{\sim} \Hom_{E}(\cK(\varphi_{H})_{\alpha',1},\nT_{L*}1^\sT_{\Delta_M}).\]
This completes the proof of the first part. The second part now follows from \autoref{prop: mhom CZ index}.
\end{proof}

Finally, we can justify why the geometric circle action defined in \autoref{thm:geometric circle action} and \autoref{thm: geometric circle tamarkin} are indeed rotations of the $1$-periodic orbits.

\begin{Prop}
Let $H$ be a homogeneous Hamiltonian on $\dot T^*M$ with Morse--Bott non-degenerate periodic orbits. 
Suppose the $1$-periodic orbits are $\gamma_j$ for $j \in J$. Then the circle action on
\[\mhom(K(\varphi_{H})_1, 1_\Delta) \simeq \bigoplus_{j\in J}or_{\gamma_j}[\mu_{CZ}(\gamma_j)-n]\]
is induced by the reparametrization of the $1$-periodic orbits.
\end{Prop}
\begin{proof}
Consider the action by the continuation map for any $z \geq 0$:
\[\Hom(K(\varphi_H)_z, 1_\Delta) \to \Hom(K(\varphi_H)_{s+z}, K(\varphi_H)_{s}).\]
When $z = 0$, the action is induced by the integral kernel $1_\Delta$ and it is clear that the continuation map induces the identity.
For general $z > 0$, by \autoref{prop:microsupport orbit}, we have a diagram of fiber sequences
\[\begin{tikzcd}[column sep=16pt]
    \Hom(K(\varphi_H)_{z+\epsilon}, 1_\Delta) \ar[r] \ar[d, "\rotatebox{90}{$\sim$}"] & \Hom(K(\varphi_H)_{z}, 1_\Delta) \ar[r] \ar[d] & \dot{\pi}_{\dot T^*M*}\mhom(K(\varphi_H)_{z}, 1_\Delta) \ar[d, "\rotatebox{90}{$\sim$}"] \\
    \Hom(K(\varphi_H)_{s+z+\epsilon}, K(\varphi_H)_{s}) \ar[r] & \Hom(K(\varphi_H)_{s+z}, K(\varphi_H)_{s}) \ar[r] & \dot{\pi}_{\dot T^*M*}\mhom(K(\varphi_H)_{s+z}, K(\varphi_H)_{s}).
\end{tikzcd}\]
Taking microlocalization, one can see that $K(\varphi_H)_s$ acts as contact transformation on the support of $\mhom(K(\varphi_H)_z, 1_\Delta)$ \cite[Theorem 7.2.1]{KS90}, and we have that 
\[\mhom(K(\varphi_H)_s \circ K(\varphi_H)_z, K(\varphi_H)_s \circ 1_\Delta) \simeq \varphi_{H*}^s\mhom(K(\varphi_H)_z, 1_\Delta).\]
We know the support of $\mhom(K(\varphi_H)_z, 1_\Delta)$ is given by closed $z$-periodic orbits of $\varphi_H$ by \autoref{prop:microsupport orbit}:
\[\gamma_z = \dot{\ms}(K(\varphi_H)_z) \cap \dot{\ms}(1_\Delta) = \dot{\ms}(K(\varphi_H)_{s+z}) \cap \dot{\ms}(K(\varphi_H)_z).\]
Therefore, we can conclude that the action of $K(\varphi_H)_s$ induces rotations on $z$-periodic orbits.
\end{proof}

\begin{Prop}
Let $H$ define a complete Hamiltonian flow on $T^*M$ with Morse--Bott non-degenerate $1$-periodic orbits. 
Suppose the $1$-periodic orbits with action $L$ are $\gamma_j$ for $j \in J$. Then the circle action on
\[\mhom(\cK(\varphi_{H})_1, T_{L*}1_\Delta) \simeq \bigoplus_{j\in J}or_{\gamma_j}[\mu_{CZ}(\gamma_j)-n]\]
is induced by the reparametrization of the orbits.
\end{Prop}
\begin{proof}
The proof is essentially the same as above.
\end{proof}

In conclusion, the circle action on the (dual of) trace is also induced by rotations of Reeb orbits:

\begin{Thm}\label{cor: mhom circle action trace}
Let $U \subset T^*M$ be a domain with RCT boundary $\partial U$ with Morse--Bott non-degenerate Reeb dynamics. 
Suppose the Reeb orbits with action $L$ are $\gamma_j$ for $j \in J$. Then the circle action on
\[\mhom(Q_U, T_{L*}1_\Delta) \simeq \mhom(\cK(\varphi_{H_{L+\varepsilon}})_1, T_{L*}1_\Delta) \simeq \bigoplus_{j\in J}or_{\gamma_j}[\mu_{CZ}(\gamma_j)-n]\]
is induced by the reparametrization of the Reeb orbits.
\end{Thm}

\subsection{Perfectness of the traces}
In this section, we show a perfectness result of the Hochschild homology of Tamarkin categories on all finite action levels. This will be used later when comparing the $S^1$ and $C_\ell$-equivariant symplectic capacities.

The main result of this section is the following finiteness result:
\begin{Thm}\label{theorem: perfectness of CT}For a relatively compact open set $U\subset T^*M$ with Morse-Bott non-degenerate RCT boundary and $E$ be a commutative ring spectrum, we have, for any $L\in \bR$ and any invertible local system $\mathscr{L}_M^\sT$, that $C_L(U; \mathscr{L}_M^\sT)$ 
is a perfect $E$-module.
\end{Thm}

We prove the finiteness result in \autoref{theorem: perfectness of CT} by approximating the Hochschild homology given by the colimit of wrappings by a finite wrapping at any finite action level as in \autoref{prop:compute limit at finite action}.

\begin{proof}[Proof of \autoref{theorem: perfectness of CT}]
For large enough $\alpha > 0$, by \autoref{prop:compute limit at finite action}, we have
\begin{align*}
\Gamma_{(-\infty,-L]}(\bR, \HH_\sT^*(\sT(U))) 
\simeq \Gamma_{(-\infty, -L]}(\bR, \pi_*^\sT\sHom^\star(\cK(\varphi_{H_\alpha})_1, 1_{\Delta_M}^\sT)).
\end{align*}
Since $\cK(\varphi_{H_\alpha})_1$ is the sheaf quantization of Hamiltonian isotopies, it is constructible and always has perfect microstalks. 

The rest of the proof follows from \cite[Proposition B.12]{guillermouviterbo_gammasupport}; see also \cite[Proposition 7.13]{Hochschild-Kuo-Shende-Zhang}. Consider the fiber sequence
\[\Gamma_{(a,b]}(\bR, \HH_\sT^*(\sT(U))) \to \Gamma_{(a,c]}(\bR, \HH_\sT^*(\sT(U))) \to \Gamma_{(b,c]}(\bR, \HH_\sT^*(\sT(U))). \]
We assume all non-constant time-1 orbits of $H_\alpha$ are non-degenerate. We order action of non-constant closed orbits of action smaller than $L$ by $0<L_1<\cdots<L_n<L$.
Take $\varepsilon >0$ small enough, and $a=-\infty, b=-L_{i}-\varepsilon, c=-L_{i}+\varepsilon$. 
Then 
\[\Gamma_{I}(\bR, \HH_\sT^*(\sT(U)))=\Hom(1_{I}, \Hom^\star(\cK(\varphi_{H_\alpha})_1, 1_{\Delta_M}^\sT).\]
They are perfect modules by \autoref{prop:compute limit at finite action}. For $I=(-\infty,\varepsilon]$, we use \cite[Theorem 7.20]{Hochschild-Kuo-Shende-Zhang} that $\Gamma_{(-\infty,\epsilon]}(\bR, \HH_\sT^*(\sT(U))) \simeq \Gamma(\overline{U}; E)$ to conclude perfectness since $\overline{U}$ is a compact manifold with boundary. 
Finally, we argue by induction, using the action window sequence above, to conclude that for $I=(-\infty,-L]$, also perfect.
\end{proof}

\begin{Coro}Under the condition of \autoref{theorem: perfectness of CT}, we have, for any $L\in \bR$ and any invertible local system $\mathscr{L}_M^\sT$, that $F_L(U; \mathscr{L}_M^\sT)$
is a perfect $E$-module.
\end{Coro}
\begin{proof}Using the 6-operation and orientation on $M$, we have
\begin{align*}
C_L(U, \mathscr{L}_M^\sT \star \omega_M^\sT)&= \Gamma_{(-\infty,-L]}(\bR,\HH_\sT^*(\sT(U), \mathscr{L}_M^\sT \star \omega_M^\sT))\\
&\simeq \Hom(i^*_L\HH^\sT_*(\sT(U), \mathscr{L}_M^\sT), E[-d]) = \Hom(F_L(U, \mathscr{L}_M^\sT), E[-d]).
\end{align*}
On the category of perfect $E$-modules, taking linear dual is an self-equivalence. Then $F_L(U, \mathscr{L}_M^\sT)$ is a perfect $E$-module.
\end{proof}

By Lemma \ref{lem:perfectness}, we can immediately get the following corollaries, on finiteness of fixed points:

\begin{Coro}\label{corollary: equivariant perfectness for costalk}
Let $U\subset T^*M$ be a relatively compact subset with Morse-Bott restricted contact boundary and $E$ be a ring spectrum satisfying the conditions in \autoref{lem:perfectness}. Then for all $L\geq 0$, we have $C_L(U; \mathscr{L}_M^\sT)^{hS^1} \in \Perf(E\fpw{u})$.
\end{Coro}
\begin{Coro}\label{corollary: equivariant perfectness for stalk}
Let $U\subset T^*M$ be a relatively compact subset with restricted contact boundary and $E$ be a ring spectrum satisfying the conditions in \autoref{lem:perfectness}. Then for all $L\geq 0$, we have $F_L(U; \mathscr{L}_M^\sT)^{hS^1} \in \Perf(E\fpw{u})$.
\end{Coro}

\subsection{Fixed point localization}
We discuss the properties of the Tate construction $C_T^{S^1}(U; E)^{tS^1}$, which, under the perfectness assumption, is an $E\fls{u}$-module obtained by localization of the homtopy fixed point $C_T^{S^1}(U, E)^{hS^1}$. The following result is analogous to \cite[Theorem 1.1]{zhaoperiodic}, \cite[Theorem 1.1]{AlbersCieliebakFrauenfelder2016} and \cite[Theorem 4.1]{CoteKartal24}.

\begin{Thm}\label{thm: localization}
Let $E$ be a complex-oriented ring spectrum. For a bounded open set $U\subset T^*M$ with RCT boundary and Morse--Bott non-degenerate Reeb dynamics, we then have $E\fpw{u}$-module equivalence
\[C^{S^1}_0(U ,E) \simeq C_{2n-*}(\overline{U},\partial U,E)\fpw{u}\simeq C^*( {U},E)\fpw{u},\]
and the $E\fpw{u}$-module 
$C^{S^1}_{(0,T]}(U ,E) $
is perfect and acyclic with respect to the $\ell$-series $[\ell](u)$ of the formal group law of $E$ for all $\ell$. 

Consequently,
inverting the $\ell$-series of the formal group law $[\ell](u)$, we have an equivalence 
\[C^{S^1}_T(U ,E)\otimes_{E\fpw{u}} E\fpw{u}[([\ell](u))^{-1}, \ell \in \bN] \simeq C^*(\overline{U} ,E)\otimes E\fpw{u}[([\ell](u)^{-1}, \ell \in \bN].\]
\end{Thm}
\begin{proof}
For the first statement, we know from \cite[Theorem 7.20]{Hochschild-Kuo-Shende-Zhang} that
\[C_0(U ,E)\simeq C_{2n-*}(\overline{U},\partial U,E)\simeq  C^*(U ,E).\]
On the other hand, by \autoref{cor:geometric circle action tamarkin}, we see that the circle action on $C_0(U ,E)$ is the trivial action of $S^1$ on $ {U}$ as we  identify $ {U}$ with constant Hamiltonian orbits in $ {U}$. Therefore, 
\[C^{S^1}_0(U ,E)\simeq (C_0(U ,E))^{hS^1}\simeq C^*(\overline{U} ,E)\fpw{u}.\]

Next, let us consider $C^{S^1}_{(0,T]}(U ,E)$. For the perfectness statement, we use \autoref{prop:compute limit at finite action}. Moreover, suppose the Reeb orbit spectrum is given by $0 < L_1 < \cdots < L_i < \cdots$. Then on each action window, we have 
\[C^{S^1}_{(L_i-\varepsilon,L_{i}+\varepsilon]}(U ,E) \simeq \dot\pi_*\mhom(Q_U, T_{L_i*}1_{\Delta}^\sT) \simeq \bigoplus_{j\in J}\Gamma(\gamma_j, or_{\gamma_j})[\mu_{CZ}(\gamma_j)-n]\]
We notice that the circle action is induced by reparameterization on non-constant Hamiltonian orbits given by $H_\alpha$ for large enough $\alpha$, which correspond to Reeb orbits of $\partial U$, by \autoref{cor: mhom circle action trace}. Then by the $C_\ell$-Gysin sequence \cite[Proposition 2.25]{Li-Zhang26-1}, we have that on a $\ell$-multiple cover of a simple orbit $\gamma$, the action of $[\ell](u)$ on the $E\fpw{u}$-module (which is the Thom spectrum of $or_\gamma$)
\[\Thom_\infty(or_\gamma)[\mu_{CZ}(\gamma)-n] \simeq \Gamma(\gamma, or_{\gamma})[\mu_{CZ}(\gamma)-n]\]
is trivial; see also \cite[Lemma 2.10]{CoteKartal24}. Then, this implies that $C^{S^1}_{(L_i-\varepsilon,L_{i}+\varepsilon]}(U ,E)$ is acyclic with respect to $S=\langle [\ell](u) \mid \ell \geq 1 \rangle$ the multiplicative subset generated by the Euler classes $[\ell](u)$, and hence $C^{S^1}_{(0,T]}$ is a $S$-acyclic module. 
Then the fiber sequence 
\[C^{S^1}_0(U ,E) \rightarrow  C^{S^1}_T(U ,E)\to C^{S^1}_{(0,T]}(U ,E)\]
induces equivalence after localization at all $[\ell](u)$ for $\ell \in \bN$.\end{proof}

\begin{Coro}
Let $\Bbbk$ be a discrete ring containing $\bQ$. For a bounded open sets $U\subset T^*M$ with Morse-Bott restricted contact boundary, we have an equivalence
\[C^{S^1}_T(U ,\Bbbk)\otimes_{\Bbbk\fpw{u}} \Bbbk\fls{u} \simeq C^*( {U} ,\Bbbk)\otimes \Bbbk\fls{u}.\]
\end{Coro}

\begin{Coro}Let $\Bbbk$ be a field containing $\bQ$. If $U\subset \bR^{2n}$ is bounded star-shaped with respect to the origin. Then for all $T\geq 0$, as a finitely generated module over the discrete valuation ring $\Bbbk\fpw{u}$, the free part of $H^*C^{S^1}_T(U ,\Bbbk)$ is of rank $1$.    
\end{Coro}
\begin{proof}The $\Bbbk\fpw{u}$-module $H^*C^{S^1}_T(U ,\Bbbk)$ is perfect since $U\subset \bR^{2n}$ is star-shaped (which has a RCT boundary.) Then it has a finite rank free part, and the rank of the free part equals the $\Bbbk\fls{u}$-dimension of $C^{S^1}_T(U ,\Bbbk)\otimes_{\Bbbk\fpw{u}} \Bbbk\fls{u}$. Then \autoref{thm: localization} shows it is equal to $\Bbbk\fls{u}$-dimension of
\[H^*( {U},\Bbbk)\otimes \Bbbk\fls{u},\]
which equals $1$ since $U$ is star-shaped. 
\end{proof}

\begin{Prop}\label{prop: unfiltered CT=SH}Let $U\subset \bR^{2n}$ be a bounded star-shaped domain. For a complex-oriented commutative ring spectrum $E$, we have
\[H^*C_\infty^{S^1}(U,E)\coloneqq \varinjlim_{T\to \infty}H^*C_T^{S^1}(U,E) \simeq E\fls{u}[([\ell](u))^{-1},\ell\geq 2],
\]
where $[\ell](u)$ stands for the $\ell$-series of the formal group law of $E$.
\end{Prop}
\begin{proof}By \cite[Corollary 2.6]{Capacities2021}, we have the natural identification $C^{S^1}_T(U,E) \simeq C^{S^1}_{r^2T}(rU,E)$. Then for $r>1$ and for every star-shaped domain $U$, we have the following commutative diagram 
\[\begin{tikzcd}
{C^{S^1}_T(U,E)} \arrow[rr, "\simeq"] \arrow[rd] &                         & {C^{S^1}_{r^2T}(rU,E)} \arrow[ld] \\
                                                   & {C^{S^1}_{r^2T}(U,E).} &                                    
\end{tikzcd}\]
We should be careful that star-shaped implies the right restriction morphism is well-defined.

Now, as $U$ is bounded and $0\in U$, we can find a ball $B=B^{2n}_A$ and $r>1$ such that 
\[B \subset U \subset rB \subset rU.\]
Then the above commutative diagram shows that the following compositions
\[\begin{tikzcd}{C^{S^1}_T(U,E)} \arrow[r] & {C^{S^1}_T(B,E)} \arrow[r, "\simeq"] & {C^{S^1}_{r^2T}(rB,E)} \arrow[r] & {C^{S^1}_{r^2T}(U,E)}\end{tikzcd}\]
and
\[\begin{tikzcd}
{C^{S^1}_T(B,E)} \arrow[r, "\simeq"] & {C^{S^1}_{r^2T}(rB,E)} \arrow[r] & {C^{S^1}_{r^2T}(U,E)} \arrow[r] & {C^{S^1}_{r^2T}(B,E)}
\end{tikzcd}\]
are identified with the natural persistence continuation maps.

Because $\bR\to \bR,T\to r^2 T$ is cofinal in $T\to \infty$. Then we have
\[C^{S^1}_\infty(B,E)\xrightarrow{\simeq } C^{S^1}_\infty(U,E),\qquad C^{S^1}_\infty(U,E)\xrightarrow{\simeq } C^{S^1}_\infty(B,E).\]

Therefore, it remains to compute $C^{S^1}_\infty(B,E)$. However, by the ball computation in \cite[Section 3.9]{CyclicZHANG} or \autoref{prop:chiu result on loop sheaf equiv}, we have (take $A=1$ for simplicity)
\[C^{S^1}_T(B,E)\simeq u^{-n\lfloor T \rfloor}E\fpw{u}.\]
Moreover, we have the formula for the continuation map for $t<T$: we set
\[c(t,T)=\prod_{j=\lf t\rf+1}^{\lf T\rf}[j](u)^n\in E^*\fpw{u} \in E^*\fpw{u},\qquad d(t,T)=\prod_{j=\lf t\rf+1}^{\lf T\rf}\frac{[j](u)}{u}^n\in E^*\fls{u}.\]
Then 
\begin{align*}
    C^{S^1}_t(B,E)\simeq u^{-n\lfloor t \rfloor}E\fpw{u} &\to C^{S^1}_T(B,E)\simeq u^{-n\lfloor T \rfloor}E\fpw{u}, \\u^{-n\lfloor t \rfloor}f(u)&\mapsto c(t,T)(u)^{-n\lfloor T \rfloor}f(u)= d(t,T)\; u^{-n\lfloor t \rfloor}f(u).
\end{align*}
We therefore have
\[C^{S^1}_\infty(U,E) \simeq  C^{S^1}_\infty(B,E) =\varinjlim_{\ell\geq 0} \frac{1}{c(0,\ell)} E\fpw{u}  \simeq E\fpw{u}[([\ell](u)^{-1},\ell\in \bN]=E\fls{u}[([\ell](u)^{-1},\ell\geq 2].\qedhere\]
\end{proof}

In the end, we show that the unfiltered version of Chiu-Tamarkin is equivalent to the $S^1$-equivariant symplectic cohomology.

\begin{Coro}
Let $U\subset \bR^{2n}$ be a bounded star-shaped domain. For a discrete commutative ring $\Bbbk$, we have
\[H^*C_\infty^{S^1}(U,\Bbbk) \simeq (\Bbbk \otimes_\bZ \bQ)\fls{u}\simeq SH_{S^1}^*(\overline{U},\Bbbk).
\]
\end{Coro}
\begin{proof}
From the computation above, we see that, $C^{S^1}_\infty(U,\Bbbk)$ coincides with (a version of) $SH^*_{S^1}(\overline{U},\Bbbk)$, see \cite[Section 4.3 and 5.1]{AlbersCieliebakFrauenfelder2016} or \cite[Section 8.1]{zhaoperiodic} (notice that their long exact sequence and computation still works for $\Bbbk$).
\end{proof}
\begin{RMK}\label{rem: versions of SH}
    There are versions of the notion of $S^1$-equivariant symplectic (co)homology which we write $SH_{S^1}^*$ for the unfiltered version. We believe $C_T^{S^1}(M)^{hS^1}$ here should correspond to $HB^*(M) = \underrightarrow{H}\underleftarrow{B}^*(M)$ in \cite[Section 4.3]{AlbersCieliebakFrauenfelder2016} (more precise, $H^*(\varprojlim_{a<0}C^T_a)$ therein), $SH^*_{S^1}(M)$ in \cite[Section 2d]{Seidel-Equivariant}, or $HF_{eq}^*(M;T)$ in \cite{Seidel2018Connections},  since we do not localize by $u$ in the definition. The colimit of localizations $\varinjlim_{T\to \infty}C_T^{S^1}(M)^{tS^1}$ should correspond to $\underrightarrow{H}\underleftarrow{T}^*(M)$ in \cite{AlbersCieliebakFrauenfelder2016} or $HP_{S^1,loc}^*(M)$ in \cite{zhaoperiodic}, which has the feature of Jones-Petrack theory \cite{JonesPetrack1990}. On the other hand, the localization of the colimit $(\varinjlim_{T\to \infty}C_T^{S^1}(M))^{tS^1}$ should correspond to $HF^{*,\infty}_{S^1}(M)$ in \cite{GanatraS1CY}, $\underleftarrow{H}\underrightarrow{T}^*(M)$ in \cite{AlbersCieliebakFrauenfelder2016} or $HP_{S^1}^*(M)$ in \cite{zhaoperiodic}, which has the feature of the Goodwillie theory \cite{Goodwillie85}.
\end{RMK}

\subsection{Fundamental class}\label{section: fundamental class}
We now consider a particular cohomology class which plays the role of the fundamental class in usual cohomology theory.  We assume that $M$ is $E$-oriented now.

For an open inclusion $U \subset T^*M$, we recall the notion in \autoref{sec:Tamarkin trace} or \cite{Capacities2021,CyclicZHANG} that 
\[ F^{S^1}(U) = \HH_*^\sT(\sT(U)) \in \Fun(BS^1, \sT).\]

We can regard it as an object of 
    \begin{gather*}
 \Fun(BS^1,\sT) \hookrightarrow \Fun(BS^1,\Fun(\larbR,D(E))\simeq \Fun(\larbR,\Fun(BS^1,D(E)) .\\
F^{S^1}(U) \mapsto [T\mapsto F^{S^1}(U)_{T}=i_T^*F^{S^1}(U)].       
    \end{gather*}
Precisely, we fully-faithfully embed $\sT$ into $\Sh_{\geq}(\bR)$, and then into $\Fun(\larbR,D(E))$, where $F^{S^1}(U)$ is mapped to $F^{S^1}(U)_{T}=i_T^*F^{S^1}(U)$ for each $T \in \bR$, regarded as objects in $\Fun(BS^1, D(E))$  via the pull-back functor $\Fun(BS^1, \Fun(\larbR,D(E))) \to \Fun(BS^1, D(E))$.

For $G\leq S^1$, we can define by restriction the action an object $F^{G}(U,E)\in\Fun(\larbR,\Fun(BG,D(E))$.

The open inclusions $U \subset V \subset T^*M$ induce strongly continuous functors $\sT(U) \rightarrow \sT(V) \rightarrow \sT(T^*M)$ as explained in \autoref{sec:tamarkin domain}, and therefore a $G$-equivariant morphism 
\[F^G(U,E) \rightarrow F^G(V,E) \rightarrow F^G(T^*M,E).\]

Direct computation in \autoref{sec:Tamarkin trace} or \autoref{sec: geometric circle tamarkin} shows that, for all $T\geq 0$ and any $G\leq S^1$, we have $F^G(T^*M,E)_{T} \simeq \Gamma_c(M,E)$ where the right hand side admits the trivial $G$-action, and when $T<0$, for $G \leq S^1$, we have $F^G(T^*M,E)_{T}\simeq 0$. 

Then by taking orientation of $M$, we have the natural map 
\[F^G(T^*M,E)_{T} \simeq \Gamma_c(M,E) \xrightarrow{or} E[-d],\]
where $E$ is the coefficient ring equipped with the trivial $G$-action and the orientation morphism is $G$-equivariant.

\begin{Def}\label{definition: fundamental class}
Let $U \subset T^*M$ be an open set. We \emph{define} $\eta^G_T(U,E)$ to be the composition map
\[\eta^G_T(U,E) \coloneqq \left[F^G(U,E)_{T} \rightarrow F^G(T^*M,E)_{T}\rightarrow E[-d]\right] .\]
\end{Def}

Now, we regard $\eta^G_T(U,E)$ as a morphism in (recall \autoref{sec:Tamarkin trace})
\[\Ext^*_{G}(F^G(U,E)_{T},E[-d])=H^*\HOM_{G}(F^G(U,E)_{T},E[-d])=H^*C^{G}_T(U,E),\] where $\HOM_{G}$ means the mapping spectrum of $\Fun(BG,D(E))$. 

\begin{RMK}Recall that $H^*C^{G}_T(U,E)$ is a left  $H^*A=H^*A_G=H^*\End_G(E)$-module via $A\cong  \End_G(E[-d])$.

\end{RMK}

The fundamental class $\eta^G_T(U,E)$ is compatible with open inclusion, and invariant under Hamiltonian action on $T^*M$. Moreover for $U\subset T^*\bR^n$, by \cite[Corollary 2.6]{Capacities2021}, we have 
\[\eta^G_T(U,E)=\eta^G_{r^2T}(rU,E)\]
via the natural identification \[F^G(U,E)_{T}\simeq F^G(rU,E)_{r^2T}.\]
All of them follow naturally from the corresponding properties of the trace of Tamarkin categories; see, for example, \cite[Section 2.4]{Capacities2021} (where the argument works for $S^1$ as well).

We may want to define symplectic capacities similar to what we did in \cite{Capacities2021,CyclicZHANG} over $\bQ$. However, the naive construction seems to give us no more than $\bQ$-coefficients. 

Notice that $\bR^n$ is contractible, then it is oriented for $E$, so we can define the fundamental class for $U\subset \bR^{2n}=T^*\bR^n$.
\begin{Prop}\label{prop: ellpsoid embedding from complex-oriented ring spectrum}Let $\mathbf{a},\mathbf{a'}\in \bR^n_{>0}$ be two non-decreasing sequences of positive real numbers. Let $E(\mathbf{a}),E(\mathbf{a'}) \subset \bR^{2n}$ be two ellipsoids. Then there exist morphisms
\[ 
({C_T^{S^1}(E(\mathbf{a}),E)},\eta_T^{S^1}) \to ({C_T^{S^1}(E(\mathbf{a'}),E)},\eta_T^{S^1})\]
of persistence modules over $A=C^*(BS^1,E)$ that preserve the fundamental classes  for all complex-oriented ring spectrum $E$ if and only if $c_k^{EH}(E(\mathbf{a}))\leq c_k^{EH}(E(\mathbf{a'}))$ for all $k\geq 1$, for $c_k^{EH}$ defined in terms of $\bQ$-coefficient cohomologies. 
\end{Prop}\begin{RMK}By \autoref{appendix: Hamiltonian embedding}, if there exists a Hamiltonian embedding $E(\mathbf{a})\lhook\joinrel\xrightarrow{Ham} E(\mathbf{a'})$, then there exist such diagrams over all complex-oriented ring spectra. 

Consequently, naively considering complex-oriented coefficients does not give more embedding obstructions between ellipsoids than Ekeland-Hofer capacities. In fact, one can do further computation to see that the following proof should be compatible with the Tate Frobenius, from the cyclotomic structure, on $F^{S^1}(E(\mathbf{a}),E) $.
\end{RMK}
\begin{proof}Let $I_\mathbf{a}(T)=\sum_i \lf T/a_i\rf$, then it is known that $c_k^{EH}(E(\mathbf{a}))\leq c_k^{EH}(E(\mathbf{a'}))$ if and only if $I_\mathbf{a}(T)\geq I_\mathbf{a'}(T)$ for all $T\geq 0$. 

Similar to \autoref{prop: unfiltered CT=SH}, we have 
\[C^{S^1}_T(E(\mathbf{a}),E)\simeq c_{\mathbf{a}}(T)^{-1}A,\qquad c_{\mathbf{a}}(T)=\prod_{i=1}^{n}\prod_{j=1}^{\lf T/a_i\rf}[j](u)\in H^*A=E^*\fpw{u} ,\]
and $\eta_T^{S^1}$ corresponds to $1\in H^*A=E^*\fpw{u}$. It is clear that the identification is compatible with the persistence structure.

Therefore,  an $\eta_T^{S^1}$-preserving map of persistence modules
\begin{equation}\label{equation: ellipsoid fundamental class comparison for E}
C^{S^1}_T(E(\mathbf{a}),E)\simeq c_{\mathbf{a}}(T)^{-1}A \to C^{S^1}_T(E(\mathbf{a'}),E)\simeq c_{\mathbf{a'}}(T)^{-1}A
\end{equation}
exists if and only if
\[c_{\mathbf{a'}}(T)|c_{\mathbf{a}}(T)~\text{in} ~H^*A=E^*\fpw{u}.\]

Now, if \autoref{equation: ellipsoid fundamental class comparison for E} is true for all $E$, then in particular $E=\bQ$. In this case, $c_{\mathbf{a}}(T)=u^{I_{\mathbf{a}}(T)}$, and the division relation gives $I_\mathbf{a}(T)\geq I_\mathbf{a'}(T)$ for all $T\geq 0$.

Conversely, if for $I_\mathbf{a}(T)\geq I_\mathbf{a'}(T)$ for all $T\geq 0$, we show that \autoref{equation: ellipsoid fundamental class comparison for E} exists for the universal complex-oriented ring $E=MU$, the ring spectrum of complex cobordism. Then morphisms in \autoref{equation: ellipsoid fundamental class comparison for E} exist for all $E$. In fact, it is clear that multiplication by
\[r(T)=c_{\mathbf{a}}(T) c_{\mathbf{a'}}(T)^{-1}~\text{in} ~ MU^*\fls{u}\]
is a natural candidate if $r(T)$ located in $MU^*\fpw{u}$. We then show when $I_\mathbf{a}(T)\geq I_\mathbf{a'}(T)$, we do have $r(T)\in MU^*\fpw{u}$, which is enough to conclude the proof. 

To see this, we use \autoref{prop: decomposition of ell series} to see that 
\[c_{\mathbf{a}}(T)=\prod_{i=1}^{n}\prod_{j=1}^{\lf T/a_i\rf}\prod_{m\mid j}\bar{\psi}_m=\prod_{d\geq 1 }(\bar{\psi}_m)^{I_\mathbf{a}(T/d)},\]
and then one has
\[r(T)=\prod_{d\geq 1 }(\bar{\psi}_m)^{I_\mathbf{a}(T/d)- I_\mathbf{a'}(T/d)}.\]
Therefore, as $I_\mathbf{a}(T)\geq I_\mathbf{a'}(T)$ for all $T\geq 0$, we know $r(T)\in MU^*\fpw{u}$.
\end{proof}

\section{Equivariant symplectic capacities}\label{section: capacities}

In this section, we recall the construction of equivariant symplectic capacities in \cite{Capacities2021,CyclicZHANG} and introduce some new constructions, prove the comparison result and the product inequalities for these capacities in general. 

By \autoref{prop: ellpsoid embedding from complex-oriented ring spectrum}, we may not gain too much knowledge if we only consider discrete rings in the following consideration. Therefore, from now on, we assume that $M$ is $\bZ$-orientable with a fixed orientation, and the coefficient ring $\Bbbk$ is a discrete commutative ring; in particular, it will be $\Bbbk=\bZ,\bQ
$ or finite fields, in particular, the conditions of \autoref{lem:perfectness} are satisfied. 

Let $A=C^*(BG,\Bbbk)$ for $G=S^1,C_\ell$. In this case, we have: when $G=S^1$, $H^*A\simeq \Bbbk\fpw{u}$ for all $\Bbbk$; when $G=C_\ell$ with $\ell\geq 2$ and $\Bbbk=\bF_{p_\ell}$ where $p_\ell$ is the minimal prime factor, $H^*A\simeq \Bbbk\fpw{u}[\theta]$ with a relation $Cu=\theta^2$ for a constant $C$ relies on $\ell$.

\subsection{Symplectic capacities from trace}
We start from various definitions of equivariant symplectic capacities.

\subsubsection{Symplectic capacities from division index}\label{section: capacity from division index}
Now, we recall the following definition of $S^1$-equivariant symplectic capacities introduced in \cite{CyclicZHANG,Capacities2021}.

\begin{Def}\label{def: entire S1 capacities}For an open set $U$ and $k\in \bN$ we define
\begin{align*} c_k(U) \coloneqq \inf
\left\lbrace
  L \geq 0:\begin{aligned} \eta^{S^1}_L (U,\bQ) \in u^kH^*C^{S^1}_L(U,\bQ)
  \end{aligned}
\right\rbrace.
\end{align*}

\end{Def}
We also define in \cite{Capacities2021}  a version of equivariant symplectic capacities by approximating $S^1$ by $C_\ell$ in a similar, but a slightly more involved     way. Here, let us present the details. We set $p_\ell$ to be the smallest prime divisor of $\ell \geq 2$. 
\begin{Def}\label{def: discrete S1 capacities}For an open set $U$ and $k\in \bN$ we define
\begin{align*} 
d_k(U)  \coloneqq \inf
\left\lbrace
  T \geq 0:\begin{aligned} \exists\, p\geq 3\text{ prime such that }\forall\,\ell \in \bO, \, p_\ell \geq p,\,  \eta^{C_\ell}_T(U,\bF_{p_\ell}) \in u^kH^{*}C^{C_\ell}_T(U,\bF_{p_\ell})
  \end{aligned}
\right\rbrace.
\end{align*}
\end{Def}

\begin{RMK}In fact, those capacities are invariant under isotopies of Hamiltonian embedding $\varphi_s: U\to T^*M$ that are not necessarily compactly supported, but have a complete flow. We can prove in the same way as we did in references based on results in \autoref{appendix: existence kernel}, especially \autoref{appendix: Hamiltonian embedding}.    
\end{RMK}
\begin{RMK} In \cite{Capacities2021}, the capacities $d_k$ is denoted by $c_k$. In \cite[Section 3.9]{CyclicZHANG}, the capacities $c_k$ is denoted by $\overline{c}_k$.
\end{RMK}
Both $c_k$ and $d_k$ are symplectic capacities as proved in \cite{Capacities2021} and \cite[Section 3.9]{CyclicZHANG}.

Moreover, based on the Gysin sequence in \autoref{thm: Gysin}, we can also prove that
\begin{equation}\label{equation: non-equivariant by Gyin}
 c_1(U)=\inf\{T\geq 0: \eta^{C_1}_T(U,\bQ)=0\}\eqqcolon c(U),
\end{equation}
where we could understand $\eta^{C_1}_T(U,\bQ)$ as obtained  from the restriction map along the trivial subgroup $C_1 \rightarrow \bR/\bZ$.

In \cite[Corollary 1.3]{Hochschild-Kuo-Shende-Zhang}, we proved that $c(U)=c_{Viterbo}(U)$. Notice that the arguments also applies to arbitrary open sets rather than open sets with RCT boundary. We only need the RCT boundary condition for cohomology comparison therein. 

\subsubsection{Symplectic capacities from image index}\label{section: capacity from image index}
In this subsection, we introduce a different construction, potentially different from the one in the previous section but closer to the Gutt-Hutchings capacities \cite{gutthutchings2018capacities}.

Recall from \autoref{definition: fundamental class}, we define 
\[\eta^{S^1}_T(U,\bQ) \coloneqq \left[F^{S^1}(U,\bQ)_{T} \rightarrow F^{S^1}(T^*M,\bQ)_{T}\rightarrow \bQ_{{S^1}}[-d]\right].\]

After taking homotopy fixed points, we pass to $u$-complete-$\bQ\fpw{u}$ modules to obtain  
\[H^*(\eta^{S^1}_T(U,\bQ)^{hS^1})\coloneqq \left[H^*(F^{S^1}(U,\bQ)_{T} ^{hS^1})\rightarrow \bQ\fpw{u}[-d]\right].\]

\begin{Lemma}\label{lemma: image index property}Let $U\subset T^*M$ be an open set. There exists a non-decreasing function $I_U :[0,\infty)\to \mathbb{N}\cup\{\infty\}$ such that $\im(H^*(\eta^{S^1}_T(U,\bQ)[d])^{hS^1})=(u^{I_U(T)}) \subset \bQ\fpw{u}$, where $(u^{\infty})$ denotes the zero ideal. 

The number $I_U(T)$ satisfies the following properties: 
\begin{enumerate}
\item $I_{T^*M}(T)=0$ for $T\geq 0$; 
\item if $U\subset V\subset T^*M$, then we have $I_U(T)\geq I_V(T)$ for $T\geq 0$; 
\item If $U\subset T^*\bR^n$, then for $r>0$, we have $I_{U}(T)=I_{rU}(r^2T)$ for $T\geq 0$;
\item if $\varphi:\bR_z\times U\to T^*M$ is an isotopy of Hamiltonian embedding, then we have $I_U(T)= I_{\varphi_z(U)}(T)$ for $T\geq 0$. 
\end{enumerate}
\end{Lemma}
\begin{proof}Ideals of the one variable formal power series ring $\bQ\fpw{u}$ is classified by $(u^m)$ for $\mathbb{N}\cup\{\infty\}$ with the convention $0=(u^{\infty})$. Then the existence of $I_U(T)$ follows since $\im((\eta^{S^1}_T(U,\bQ)[d])^{hS^1})=(u^{I_U(T)}) $ is a submodule, equivalently an ideal, of $\bQ\fpw{u}$. The number $I_U(T)$ is clearly non-decreasing by the persistence of $\eta^{S^1}_T(U,\bQ)$: that $\eta^{S^1}_{T'}(U,\bQ) = \eta^{S^1}_T(U,\bQ) \circ j_{T,T'}^*$ for the natural map 
\[j_{T,T'}^*: H^*(F^{S^1}(U;\bQ)_{T'}) \to H^*(F^{S^1}(U;\bQ)_T).\]
For all the remaining properties, they follow from corresponding properties of the fundamental class $\eta^{S^1}_T(U,\bQ)$ as stated in \cite{Capacities2021} and \cite[Section 3.9]{CyclicZHANG}. Plus \autoref{appendix: Hamiltonian embedding} in this article.\end{proof}

\begin{Lemma}\label{lem: symplectic ball}Let $B_{A}^{2n} \subset \bR^{2n}$ be the symplectic ball of dimension $2n$, then we have $I_{B_{A}^{2n}}(T)=n\lfloor T /A\rfloor$.    
\end{Lemma}
\begin{proof}The result follows from the ball computation in \cite{Capacities2021} and \cite[Section 3.9]{CyclicZHANG}.
\end{proof}

\begin{Def}\label{definition of capacities of image index}Let $U\subset T^*M$ be an open set and $k\geq 1$. We then define 
\[\widehat{c}_k(U)\coloneqq \inf\{T\geq 0\mid I_U(T)\geq k\}.\]    
\end{Def}
\begin{Thm}The function $U\mapsto \widehat{c}_k(U)$ defines a sequence of symplectic capacities, which are non-decreasing with respect to $k$, and $\widehat c_k(B_{A})=\lceil k/d\rceil A$. 
\end{Thm}

Here, we similarly present a definition of the discrete version of $\widehat{c}_k$. Again, for $\ell\geq 2$, we set $p_\ell$ the minimal prime factor of $\ell$. We write $H^*A_\ell=\End_{C_\ell}(\bF_{p_\ell})\simeq \bF_{p_\ell}\fpw{u}[\theta]$  (with a relation $Cu=\theta^2$ for a constant $C$ depending on $\ell$) for all $\ell \geq 2 $.

\begin{Def}Let $U\subset T^*M$ be an open set and $k\geq 1$. We then define 
\begin{align*}
    \widehat{d}_k(U)\coloneqq \inf \left\lbrace
  T \geq 0:\begin{aligned} \exists\, p \geq 3\text{ prime so that }\forall\,\ell \in \bN_{\geq 2}, \, p_\ell \geq p,\,  \im(H^*(\eta^{C_\ell}_T(U,\bF_{p_\ell})[d])^{hC_\ell}) \subset  u^k \bF_{p_\ell}\fpw{u}[\theta]
  \end{aligned}
\right\rbrace.
\end{align*}
\end{Def}

By similar discussions, we can see that:
\begin{Thm}The function $U\mapsto \widehat{d}_k(U)$ defines a sequence of symplectic capacities, which are non-decreasing with respect to $k$, and $\widehat{d}_k(B_{A})=\lceil k/d\rceil A$. 
\end{Thm}

\subsection{Comparison of symplectic capacities}\label{section: comparison of capacities} In this section, we collect comparison results for the varieties of capacities we introduced. They involve both comparisons between discrete v.s. continuous circle for both classes of capacities, and comparisons between different classes.

\subsubsection{Discrete v.s.~continuous circle}
In this subsection, we compare the symplectic capacities $c_k$ and $d_k$ (respectively, $\widehat c_k$ and $\widehat d_k$).

\begin{Thm}\label{thm: c_k=cbar_k}For an open set $U$, if the underlying derived abelian group of $F^{S^1}(U,\bZ)_{T} $ is perfect  for all $T\geq 0$, then we have
\[c_k(U)=d_k(U),\qquad \widehat c_k(U)=\widehat d_k(U).\]
\end{Thm}

\begin{proof}Let $T\geq 0$. By \autoref{lem: perfectness over Z}, we can take $N>0$ such that the $F^{S^1}(U,\bZ)_{T}^{hS^1}\otimes^{\mathbf L} \bZ[1/N]$ has
the stated decompositions over $\bZ[1/N]$ and all nonzero coefficients of the fundamental
class are units. Choose a prime $p>\max\{2,N\}$. Then for all $\ell$ with $p_\ell >p$,  the natural Gysin splitting in the proof of \autoref{Coro: finite restriction for large prime determine the free part} identifies
\begin{equation}\label{eq: finite field reduction-fundamental class}
    \eta^{C_\ell}_T(U,\bF_{p_\ell})=j_\ell^*(\eta^{S^1}_T(U,\bZ))\otimes^{\mathbf L} \bF_{p_\ell},\qquad \eta^{S^1}_T(U,\bQ)=\eta^{S^1}_T(U,\bZ)\otimes {\bQ}.
\end{equation} 

Consequently, the first result follows from:
\begin{align*}
    &\eta^{S^1}_T(U,\bQ) \in u^kH^*C_T^{S^1}(U,\bQ)\\
    \iff  &\eta^{S^1}_T(U,\bZ) \otimes^{\mathbb{L}} \bZ[1/N]\in u^kH^*C_T^{S^1}(U,\bZ) \otimes^{\mathbb{L}} \bZ[1/N]\\
     \iff  &\eta^{C_\ell}_T(U,\bF_{p_\ell}) \in u^kH^*C_T^{C_\ell}(U,\bF_{p_\ell}),
\end{align*}
and the second result follows from:
\[\im(H^*(\eta^{C_\ell}_T(U,\bF_{p_\ell})[d])^{hC_\ell}) \simeq \im(H^*(\eta^{S^1}_T(U,\bZ)\otimes^{\mathbb{L}}{\bF_{p_\ell}}[d])^{hS^1}) \otimes_{\bF_{p_\ell}\fpw{u}}\bF_{p_\ell}\fpw{u}[\theta]=u^{I_U(T)}\bF_{p_\ell}\fpw{u}[\theta].  \qedhere\]
\end{proof}

\begin{Coro}When $U$ has RCT boundary and non-degenerate Reeb dynamics, we have $c_k(U)=d_k(U)$ and  $\widehat c_k(U)=\widehat d_k(U)$.    
\end{Coro}
\begin{proof}Combine \autoref{thm: c_k=cbar_k} and \autoref{theorem: perfectness of CT}.    
\end{proof}

\begin{Coro}\label{corollary: capacity comparison, discrete vs continuous circle}When $M=\bR^n$ and $U$ has RCT boundary, we have $c_k(U)=d_k(U)$ and  $\widehat c_k(U)=\widehat d_k(U)$.    \end{Coro}
\begin{proof}The argument is geometric and works for both classes: If $\partial U$ has non-degenerate Reeb dynamics, then it follows from the previous corollary. If $\partial U$ does not have non-degenerate Reeb dynamics, then we can perturb $U$ a little to make sure $\partial U$ has non-degenerate Reeb dynamics. Then the result follows that symplectic capacities are continuous with respect to such perturbation (c.f. \cite[Proposition 3]{EkelandHoferII}).
\end{proof}
\begin{RMK}The only reason we need $M=\bR^n$ here is we need the conformal property to make sure capacities are continuous with respect to those perturbations (since rescaling gives the Liouville vector field of $T^*\bR^n$). The result still holds for those $M$ such that the conformal property of $c_k$ is true.
\end{RMK}

Here we notice that \autoref{corollary: capacity comparison, discrete vs continuous circle} and \autoref{equation: non-equivariant by Gyin} give for all open $U\subset \bR^{2n}$ that
\[d_1(U)=c_1(U)=c(U);\]
On the other hand, by \cite{Hochschild-Kuo-Shende-Zhang}, for connected bounded open $U$ with RCT boundary, we have  the second equality below
\begin{equation}\label{eq: c1=c viterbo}
  c_1(U)=c(U)=c_{SH}(\overline{U})=c_{FH}(U)=c_{FHW}(U)=c_{V}(U), \end{equation}
and the third equality with Floer–Hofer capacity is given by \cite[Proposition 5.7]{Hermann2000} and the last two equalities with Floer-Hofer-Wysocki capacity and Viterbo capacity are given by \cite[Proposition 3.2 and Theorem 1.5]{Hermann2004}.

\subsubsection{\texorpdfstring{Various (cohomological) $S^1$-equivariant capacities}{}}
We now compare the capacities defined using division indices and image indices. For $\widehat{c}_k$ and $ {c}_k$, we have the following comparison result.
\begin{Prop}\label{prop: c=hat c, direct summand criterion}Let $U\subset T^*M$ be an open set and $T\geq 0$. Consider the fiber sequence
\[ F^{S^1}_{T}(U,\bQ)^{hS^1}\xrightarrow{\eta_T^{S^1}}   F^{S^1}_{T}(T^*M,\bQ)^{hS^1}=\Gamma_c^{*}(M,\bQ)\fpw{u}\xrightarrow{\iota_T^{S^1}} F^{S^1,out}_{T}(U,\bQ)^{hS^1},\]
and
set $\beta_T^{S^1}=H^*\iota_T^{S^1}(1)$. We have
\begin{align*}
    \widehat{c}_k(U)&=\inf\{T\mid u^{k-1}\beta_T^{S^1}\neq 0\},\\
    c_k(U)&=\inf\{T\mid u^{k-1}\beta_T^{S^1}\neq 0\hspace{-6pt} \mod u^k\}.
\end{align*}
Consequently, we have that for all $k\geq 1$ 
\[\widehat{c}_k(U)\leq {c}_k(U).\]
Moreover, when $\widehat{c}_k(U)<\infty$, we have equalities if and only if
for all $u^{k-1}\beta_T^{S^1} \neq 0$ we have
\[u^{k-1}\beta_T^{S^1} \notin u^kH^*F^{S^1,out}_{T}(U,\bQ)^{hS^1},\]
if and only if $\im(H^{n}\iota_T^{S^1})$ defines a direct summand $\bQ\fpw{u}/(u^{k'})$ for $k' \geq k$.

\end{Prop}
\begin{proof}We ignore $(U,\bQ)$ from the notation to simplify the discussion. 

First, assume that $U \subset \bR^{2n}$ for $M=\bR^n$, where $\Gamma_c(M,\bQ)\simeq \bQ[-n]$. We use the fiber sequence
\begin{equation}\label{equation: in/out triangle in capacity comparison} F_T^{S^1}(U,\bQ)^{hS^1}\xrightarrow{\eta_T^{S^1}} \bQ\fpw{u}[-n] \xrightarrow{\iota_T^{S^1}} F_T^{S^1,out}(U,\bQ)^{hS^1}.  
\end{equation}
Consider $\widehat{c}_k(U)$ first. The definition gives $\im(H^*\eta_T^{S^1})=(u^{I_U(T)})$. By definition of $\beta_T^{S^1}$, we have 
\[H^*\iota_T^{S^1}(u^k)=u^k\beta_T^{S^1},\]
and the cohomology long exact sequence of \autoref{equation: in/out triangle in capacity comparison} shows that
\[u^{k-1}\beta_T^{S^1} \neq 0 \iff u^{k-1} \notin \im(H^*\eta_T^{S^1}) \iff I_U(T)\geq k.\]
This proves the formula for $\widehat{c}_k(U)$.

Next consider ${c}_k(U)$. By applying $\Hom_{\bQ\fpw{u}}(F_T^{S^1}(U,\bQ)^{hS^1}[n],-)$ to the exact triangle 
\[\bQ\fpw{u}[-2k]\xrightarrow{u^{k}}\bQ\fpw{u} \xrightarrow{q_k} \bQ\fpw{u}/(u^{k}),\]
we see that $\eta_T^{S^1}=u^kx_T^{S^1}$ if and only if $q_k \circ \eta_T^{S^1}=0$. Then, we apply $\Hom_{\bQ\fpw{u}}(-,\bQ\fpw{u}/(u^{k}))$ to \autoref{equation: in/out triangle in capacity comparison}
to see that $q_k\circ \eta_T^{S^1}=0$ if and only if there is some $p\in \Hom_{\bQ\fpw{u}}(F_T^{S^1,out}(U,\bQ)^{hS^1},\bQ\fpw{u}/(u^{k}))$ such that $p\circ \iota_T^{S^1} =q_k$. Suppose such $p$ exists, we have $p(\beta_T^{S^1})=p\circ \iota_T^{S^1}(1)=q_k(1)$, so $p(u^{k-1}\beta_T^{S^1})=q_k(u^{k-1})\neq 0$; and at the same time $p(u^k) = 0$. In particular, this implies that 
\[u^{k-1}\beta_T^{S^1}\notin u^kH^*F^{S^1,out}_{T}(U,\bQ)^{hS^1}.\]
Conversely, suppose the above inequality holds. Then we in particular have 
$\beta_T^{S^1}\notin uH^*F^{S^1,out}_{T}(U,\bQ)^{hS^1}.$
Consider the base change from $\bQ\fpw{u}$ to $\bQ$. We have, for $\bar{\beta}_T^{S^1}$ the image of $\beta_T^{S^1}$ under the base change, that $\bar{\beta}_T^{S^1}\neq 0$. Then we have an injection after base change
\[H^*\bar\iota_{T}^{S^1}: \bQ\hookrightarrow H^*F^{S^1,out}_{T}(U,\bQ)^{hS^1}/uH^*F^{S^1,out}_T(U,\bQ)^{hS^1}.\]
In particular, we can show that $\im(H^*\iota_T^{S^1})$ is a direct summand of $H^*F^{S^1,out}_{T}(U,\bQ)$. Now, since we also know that $u^{k-1}\beta_T^{S^1} \neq 0 \hspace{-6pt} \mod u^k$, we can lift $\bar p$ to $p$ such that $p\circ \iota_T^{S^1} =q_k$.

Now, the inequality follows: if $u^{k-1}\beta_T^{S^1}\notin  u^kH^*F^{S^1,out}_{T}(U,\bQ)^{hS^1}$, then clearly $u^{k-1}\beta_T^{S^1}\neq 0$. It is also clear that the reverse implication gives the equality.

When $U \subset T^*M$ is a general open set, we note that $i: \bQ\fpw{u}[-n] \rightarrow H_c^*(M,\bQ)\fpw{u}$ is a direct summand. Therefore, we can consider the diagram
\[\begin{tikzcd}
{F_T^{S^1}(U,\bQ)}\ar[r,"{\eta_T^{S^1}}"] \ar[d,"\simeq"] & {\bQ\fpw{u}[-n]} \ar[r] \ar[d,"i", shift left=1] & F_T^{S^1,nu}(U,\bQ) \ar[d,"i", dashed, shift left=1] \\ 
{F_T^{S^1}(U,\bQ)}\ar[r] & {H_c^*(M,\bQ)\fpw{u}} \ar[r,"{\iota_T^{S^1}}"] \ar[u,"p", shift left=1] & F_T^{S^1,out}(U,\bQ) \ar[u,"p", dashed, shift left=1]
\end{tikzcd}\]
and conclude that $F_T^{S^1,nu}(U,\bQ)$ (which by definition is the cofiber of $\eta_T^{S^1}$) is a direct summand of $F_T^{S^1,out}(U,\bQ)$. The above argument still works for the top fiber sequence, which shows that
\[\eta_T^{S^1} \in u^kH^*C_T^{S^1}(U, \bQ) \iff u^{k-1}\beta_T^{S^1} \notin u^kH^*F_T^{S^1,nu}(U, \bQ).\]
Then we also have $u^{k-1}\beta_T^{S^1} \notin u^kH^*F_T^{S^1,out}(U, \bQ).$
\end{proof}

\begin{Coro}\label{coro: c=hat c, dynamically convex case}Let $U\subset \bR^{2n}$ be a bounded non-degenerate star-shaped domain. If $H^{n-1}F^{S^1}_{T}(U,\bQ)^{hS^1}$ has no $u$-torsion, then $\widehat{c}_k(U)= {c}_k(U)$ for all $k$.

In particular, for $n\geq 2$, if $\mu_{CZ}(\gamma) \geq n$ for all Reeb orbits $\gamma$ in $\partial U$ (for example when $\partial U$ is dynamically convex in the sense of \cite{HWZ1998}), then $\widehat{c}_k(U)= {c}_k(U)$ for all $k$.
\end{Coro}
\begin{proof}When $U\subset \bR^{2n}$ is bounded, we have all $\widehat{c}_k(U)\leq c_k(U)<\infty$. We write $I=I_U(T)$ here. As $U$ is non-degenerate and star-shaped, by \autoref{corollary: equivariant perfectness for costalk}, we have $C_T^{S^1}(U; \bQ)^{hS^1}$ is perfect and rank $1$ as a $\bQ\fpw{u}$-module. So is $H^{*}F^{S^1}_{T}(U,\bQ)^{hS^1}$ as in \autoref{corollary: equivariant perfectness for stalk}. Therefore, we have
\[H^{n+*}F^{S^1}_{T}(U,\bQ)^{hS^1}\simeq u^{I}\bQ\fpw{u}\oplus \operatorname{Tor}_u,\]
where $H^{n+*}(\eta_T^{S^1}(U,\bQ)^{hS^1})(u^{I})=u^{I} \in \bQ\fpw{u}$.

Moreover, because $H^{n-1}F^{S^1}_{T}(U,\bQ)^{hS^1}$ has no $u$-torsion, then there is no non-trivial homotopy between $\eta_T^{S^1}(U,\bQ)|_{\operatorname{Tor}_u}$ and $0$. Therefore, $\eta_T^{S^1}(U,\bQ)|_{\operatorname{Tor}_u}=0$.
Consequently, we have
\[H^{n+*}F^{S^1,out}_{T}(U,\bQ)^{hS^1}= H^{n+*}\cof(\eta_T^{S^1}(U,\bQ))\simeq \bQ\fpw{u}/(u^{I})\oplus \operatorname{Tor}_u[1],\]
and $\im({H^*\iota_T^{S^1}}:\bQ\fpw{u}\rightarrow H^{n+*}F^{S^1,out}_{T}(U,\bQ)^{hS^1})=\bQ\fpw{u}/(u^{I})$ is a direct summand, and $u^{k-1}\beta_T\notin u^kH^{n+*}F^{S^1,out}_{T}(U,\bQ)^{hS^1}$ for all $k$. Then we conclude the equality by \autoref{prop: c=hat c, direct summand criterion}.

Under the condition $\mu_{CZ}(\gamma) \geq n$, we have $H^{q}F^{S^1}_{T}(U,\bQ)^{hS^1}=0$ for $q<n$: We can deduce it either from \autoref{section:CZ index and spectrality}, or reduce to non-equivariant $C_T$ using the Gysin sequence \autoref{thm: Gysin} and then deduce it from the comparison isomorphism with symplectic cohomology \cite{Hochschild-Kuo-Shende-Zhang}. This particularly implies that $H^{n-1}F^{S^1}_{T}(U,\bQ)^{hS^1}=0$ has no torsion.
\end{proof}

\begin{RMK}
In general, we should not expect $\widehat{c}_k= c_k$ in general.
    Nevertheless, we will show in \autoref{prop: toric c = hat_c in general} that for a toric domain $X_\Omega$, we have $\widehat{c}_k(X_\Omega)= c_k(X_\Omega)$ for all $k$.
\end{RMK}

Recall that for any $S^1$-module $M$, let $J=\bQ\fls{u}/u\bQ\fpw{u}$, we have
\[M_{hS^1}= M^{hS^1} \otimes_{\bQ\fpw{u}} J.\]
Here, we remark that $J=\bQ\fls{u}/u\bQ\fpw{u}$ is isomorphic to $H_{*}(BS^1,\bQ)$ such that $[\bC P^i]=u^{-i}$.

Under perfectness of $F^{S^1}(U,\bQ)^{hS^1}_T$, we shall compute $\widehat{c}_k(U)$ using homotopy orbits rather than homotopy fixed points, which are more directly connected to other known symplectic capacities defined using equivariant symplectic homology defined using homotopy orbits.

\begin{Lemma}\label{lemma: hat c definition non-deg case}Let $U\subset T^*M$ be open and we assume that $F^{S^1}(U,\bQ)^{hS^1}_T$ is perfect for all $T\in \bR$.
Then we have 
\[\widehat{c}_k(U)=\inf\{T\geq 0\mid \eta_T^{S^1}(U,\bQ)\otimes u^{1-k}=0\in H^*(C_T^{S^1} (U,\bQ)\otimes_{\bQ\fpw{u}} J)\}.\]
\end{Lemma}
\begin{proof}Consider the $\bQ\fpw{u}$-linear map 
\[\theta_k:\bQ\fpw{u}\to J[2-2k],\quad f(u)\mapsto u^{1-k}f(u),\]
where $\ker(\theta_k)=(u^k)$.

Notice that $J$ is injective as a $\bQ\fpw{u}$-module; we know that its derived hom equals the usual hom. Then for the fundamental class 
\[\eta_T^{S^1}(U,\bQ): F^{S^1}(U,\bQ)^{hS^1}_T\to  \bQ\fpw{u}[-n],\]
we have 
\[I_U(T)\geq k\iff \im (H^*\eta_T^{S^1}(U,\bQ)) \subset(u^k)\iff \theta_k[-n]\circ \eta_T^{S^1}(U,\bQ)=0.\]

That means we have
\[\widehat{c}_k(U)=\inf \{T\geq 0\mid\theta_k[-n]\circ \eta_T^{S^1}(U,\bQ)=0\} .\]

Now, as we assume perfectness of $F^{S^1}(U,\bQ)^{hS^1}_T$ for all $T\in \bR$, then it is dualizable. Consequently:
\[C_T^{S^1}(U,\bQ)\otimes_{\bQ\fpw{u}} J=\Hom_{\bQ\fpw{u}}(F^{S^1}(U,\bQ)^{hS^1}_T,\bQ\fpw{u}[-n])\otimes_{\bQ\fpw{u}} J=\Hom_{\bQ\fpw{u}}(F^{S^1}(U,\bQ)^{hS^1}_T,J[-n]),\]
and under the identification, we have
$\eta^{S^1}_T(U,\bQ)\otimes u^{1-k}$ corresponds to $\theta_k[-n]\circ \eta^{S^1}_T(U,\bQ).$
Therefore
\[\widehat{c}_k(U)=\inf \{T\mid\eta^{S^1}_T(U,\bQ)\otimes u^{1-k} =0 \in H^*(C_T^{S^1}(U,\bQ)\otimes_{\bQ\fpw{u}} J)\}. \qedhere\]
\end{proof}

Now, we state the following:
\begin{Prop}\label{prop: GH=CT capacity level}Let $U\subset T^*M$ be a relatively compact open set with restricted contact boundary. Under the \autoref{conjecture:SH=HH, S1}: $(C_T^{S^1}(U,\bQ),\eta_T)\simeq (SC^{S^1}_{(-\infty,T)}(\overline{U}),\mathbf{1}_{SH})$ for the symplectic cohomology unit $\mathbf{1}_{SH}$ and $T\notin \Spec(\partial U)$, we have that for all $k\geq 1$ that
\[\widehat{c}_k(U)=c_k^{GH}(\overline{U}).\]
\end{Prop}
\begin{proof}We can assume that $\partial U$ is non-degenerate. Otherwise, we can approximate the domain via non-degenerate open domains.

The Gutt-Hutchings capacities are defined using the homotopy orbits of the $S^1$-action on the positive symplectic cohomology cochain. Under the conjectural filtered isomorphism, the homotopy orbits are isomorphic to
\[SC^{S^1}_{(\epsilon,T)}(\overline{U})_{hS^1}=C^{S^1}_{(\epsilon,T)}(U,\bQ) \otimes_{\bQ\fpw{u}} J =C_{(\epsilon,T)}(U,\bQ)_{hS^1}.\]
Denoting by $\delta: H^*(C_{(\epsilon,T)}(U,\bQ)_{hS^1}) \to H_{1+2n-*}(\overline{U},\partial U; \bQ)\otimes J$ the connecting homomorphism, then we can rewrite the definition of Gutt-Hutchings capacity \cite[Definition 4.1]{gutthutchings2018capacities} (the commutativity of $u$ and $\delta$ follows from \cite[Proposition 3.1]{gutthutchings2018capacities}) using Chiu-Tamarkin complex
\[c_k^{GH}(\overline{U})= \inf\{T > 0 \mid \exists\,\alpha \in H^*C_{(\epsilon,T)}(U; \bQ)_{hS^1}, \,\text{s.t.}\,\delta(u^{k-1}\alpha)=u^{k-1}\delta(\alpha) = [\overline{U},\partial U]\otimes 1\}.\]

Then we know that
\[u^{k-1}\delta(\alpha) = [\overline{U},\partial U]\otimes 1 \iff \delta(\alpha)=[\overline{U},\partial U]\otimes u^{1-k} \text{ in }H_{2n}(\overline{U},\partial U; \bQ)\otimes J,\]
and we notice that $\mathbf{1}_{SH}$ is mapped to $[\overline{U},\partial U]$ under $\delta$.

Therefore, the action fiber sequence $C^{*}(U; \bQ)_{hS^1} \to C_T(U; \bQ)_{hS^1} \to C_{(\epsilon,T)}(U; \bQ)_{hS^1}$ shows that
\[c_k^{GH}(\overline{U})= \inf\{T > 0 \mid \mathbf{1}_{SH} \otimes u^{1-k} = 0 \in C_T(U,\bQ)_{hS^1}\}.\]
This is exactly the formula of $\widehat{c}_k(U)$ in \autoref{lemma: hat c definition non-deg case} (notice that $F^{S^1}(U,\bQ)^{hS^1}_T$ is perfect for all $T$ by non-degeneracy).
\end{proof}

Finally, we can also provide a characterization of the capacities $c_k$. They are also related to existing works in symplectic cohomology \cite{Zhou2022kdilations}, but we will omit the comparison.

\begin{Lemma}\label{lemma: c definition non-deg case }Let $U\subset T^*M$ be an open set.
Let $J_k = \bQ\left<1,u^{-1},\dots,u^{1-k}\right>$. Then we have 
\[{c}_k(U)=\inf\{T\geq 0\mid \eta_T^{S^1}(U,\bQ)\otimes u^{1-k}=0\in H^*(C_{T}^{S^1} (U,\bQ)\otimes_{\bQ\fpw{u}} J_k)\}.\]
\end{Lemma}
\begin{proof}The proof is similar to \autoref{lemma: hat c definition non-deg case}. But notice that  we do not need perfectness here since $J_k\simeq \bQ\fpw{u}/(u^k)[2k-2]$ is perfect.
\end{proof}

\begin{RMK}
One can show that (and compare it with the definition in \cite[Definition 3.4]{Zhou2022kdilations})
\[c_k(U) = \inf\{T > 0 \mid \exists\,\alpha \in H^*C^+_T(U,\bQ)_{hS^1}, u^k\alpha = 0, \delta(u^{k-1}\alpha) = [\overline{U},\partial U] \otimes 1\}.\]
\end{RMK}

\subsection{Product inequality for symplectic capacities}\label{section: product inequality}
We start from some discussion of the fundamental class of products, and then obtain the product inequality.

Let $U\subset T^*M$ and $W\subset T^*N$ be open sets, and assume that $M$ and $N$ are both $\bZ$-orientable. By the Künneth formula, we have the following commutative diagram
\[
\begin{tikzcd}
{F^{S^1}(U,\bQ)\star F^{S^1}(W,\bQ)} \arrow[d, "\simeq"] \arrow[r] & {F^{S^1}(T^*M,\bQ)\star F^{S^1}(T^*N,\bQ)} \arrow[d, "\simeq"] \arrow[r] & {\bQ_{[0,\infty)}[-m]\star\bQ_{[0,\infty)}[-n]} \arrow[d, "\simeq"] \\
{F^{S^1}(U\times W,\bQ)} \arrow[r]                           & {F^{S^1}(T^*(M\times N),\bQ)} \arrow[r]                            & {\bQ_{[0,\infty)}[-m-n]}                                           
\end{tikzcd}
\]
where the last two horzantal arrows are induced by taking orientation.

In the above diagram, the sheaf convolution and morphism is a priori taken over $\bQ$. Then we equip $S^1$-action on the first line by the diagonal action of $S^1$-action on trace of each factor, and equip the second line by the $S^1$-action of the trace of $U\times W$. For the last $\bQ$, we equip the trivial $S^1$-action. Under the given action, one can check using the equivariant compatibility of the Künneth formula that the diagram is commutative in $S^1$-equivariant category.

As a direct application of the diagram, we have the stabilization property.

\begin{Thm}[Stabilization invariance]\label{prop: stabilization invariance}For all $U\subset T^*M$ and for all $k\in \bN$, we have
\[c_k(U\times \bR^{2n})=c_k(U),\quad \widehat{c}_k(U\times \bR^{2n})=\widehat{c}_k(U).\]
\end{Thm}
\begin{proof}In fact, since $F^{S^1}(\bR^{2n},\bQ)\simeq \bQ_{[0,\infty)}[-n]$ with the trivial $S^1$-action by the description of \autoref{thm:S1-action}, we have $F^{S^1}(U\times \bR^{2n},\bQ)\simeq F^{S^1}(U,\bQ)[-n]$ and ${\eta^{S^1}_T(U\times  \bR^{2n},\bQ)}={\eta^{S^1}_T(U,\bQ)}[-n]$ by the Künneth formula \autoref{thm:Künneth tamarkin}. Then we can conclude by the definition of capacities.    
\end{proof}

Now we turn to the general product inequality.

\begin{Lemma}\label{lemma: very simple restriction of stalk of convolution}Let $F \to F',G\to G'\in \Sh_{\geq } (\bR,\mathcal{E})$ for a stable category $\cE$. For any $T,a,b\in \bR$ with $a+b\leq T$, there exists a morphism
\[i_{a,b}^*: (F\star G)_T \to  F_a\otimes G_b.\]

The constructed morphism for fixed $a,b,T$ is functorial with respect  $F\to F',G\to G'\in \Sh_{\geq } (\bR,\mathcal{E})$. In particular, there is a natural morphism
\[i_{a+b\leq T}^*: (F\star G)_T \to \varprojlim_{a+b\leq T} F_a\otimes G_b.\]
\end{Lemma}
\begin{proof}We first prove the diagram for $b=T-a$. By definition, we have
\[(F\star G)_T\simeq \Gamma_c(\{(x,y)\mid x+y=T\},F\boxtimes G).\]  
Therefore, as $(a,T-a)$ is in the set $\{(x,y)\mid x+y=T\}$, we have the (closed) restriction map
\[(F\star G)_T\simeq \Gamma_c(\{(x,y)\mid x+y=T\},F\boxtimes G\}) \to \Gamma_c(\{(a,T-a)\},F\boxtimes G\})=F_a\otimes G_{T-a}. \]

Next, notice that $a+b\leq T$ implies $b\leq T-a$, and $G\in \Sh_{\geq}(\bR)$, there exists a continuation map $G_{T-a}\to G_b$. Therefore, the desired morphism is the composition
\[(F\star G)_T\to F_a\otimes G_{T-a}  \to  F_a\otimes G_{b}.\]

It is evident that the construction is functorial with respect to both factors.
\end{proof}

As a simple corollary of \autoref{lemma: very simple restriction of stalk of convolution}, we have a commutative diagram:
\begin{Lemma}\label{lemma: fundamental class product decomposition naive form}Let $U \subset T^*M,W \subset T^*N$ and $a,b,T \geq 0$ with $a+b\leq T$. We have a 
commutative diagram of $S^1$-modules
\[\begin{tikzcd}
{(F^{S^1}(U,\bQ)\star F^{S^1}(W,\bQ))_{T}} \arrow[d, "{\eta^{S^1}_T(U\times W,\bQ)}" left] \arrow[r, "i_{a,b}^*"] & {F^{S^1}(U,\bQ)_{a}\otimes F^{S^1}(W,\bQ)_{b}} \arrow[d, "{\eta^{S^1}_a(U,\bQ)\otimes \eta^{S^1}_b(W,\bQ)}"] \\
{\bQ[-m-n]} \arrow[r,"{\simeq}"]  & {\bQ[-m] \otimes \bQ[-n].}
\end{tikzcd}\]
In particular, there is a commutative diagram of $S^1$-modules
\[\begin{tikzcd}
{(F^{S^1}(U,\bQ)\star F^{S^1}(W,\bQ))_{T}} \arrow[d, "{\eta^{S^1}_T(U\times W,\bQ)}" left] \arrow[r,"i_{a+b\leq T}^*"] & {\displaystyle\varprojlim_{a+b\leq T}F^{S^1}(U,\bQ)_{a}\otimes F^{S^1}(W,\bQ)_{b}} \arrow[d, "{\varprojlim \eta^{S^1}_a(U,\bQ)\otimes \eta^{S^1}_b(W,\bQ)}"] \\
{\bQ[-m-n]} \arrow[r,"{\simeq}"]  & {\bQ[-m-n].}
\end{tikzcd}\]
\end{Lemma}
\begin{proof}
This follows from the commutativity of the diagram
\[\begin{tikzcd}
{(F^{S^1}(U,\bQ)\star F^{S^1}(W,\bQ))_{T}} \arrow[d] \arrow[r, "i_{a,b}^*"] & {F^{S^1}(U,\bQ)_{a}\otimes F^{S^1}(W,\bQ)_{b}} \arrow[d] \\
{(F^{S^1}(T^*M,\bQ) \star F^*(T^*N, \bQ))_T} \arrow[r,"{\simeq}"]  & {F^{S^1}(T^*M, \bQ)_a \otimes F^{S^1}(T^*N, \bQ)_b}
\end{tikzcd}\]
and the product formula of fundamental classes on both sides.
\end{proof}

\begin{Coro}\label{coro:index inequality}Let $U \subset T^*M, W \subset T^*N$. Recall the definition of $I_U(T), I_W(T)$ in \autoref{lemma: image index property}. If $U$ has non-degenerate RCT boundary, we then have
\begin{equation}\label{eq: index inequality}I_{U\times W}(T)\geq \sup_{a,b\geq 0,a+b\leq T}\{I_U(a)+I_W(b)\}\end{equation}
\end{Coro}
\begin{proof}After taking the homotopy fixed point to the diagram of \autoref{lemma: fundamental class product decomposition naive form}, we obtain that 
\begin{align*}
   ~\im(H^*(\eta_T^{S^1}(U\times W,\bQ)^{hS^1}))
   \subset\im(H^*(\eta^{S^1}_a(U,\bQ)^{hS^1}~\widehat\otimes_{A} ~ \eta^{S^1}_b(W,\bQ)^{hS^1})).
\end{align*}

Since $U$ has a non-degenerate RCT boundary, we have $F^{S^1}_a(U,\bQ)^{hS^1}$ is a perfect $A$-module by \autoref{theorem: perfectness of CT}. Then
we can use \autoref{prop:eilenberg-moore} and \autoref{cor: complete eilenberg-moore} to see that
\begin{align*}
   ~&\im(H^*(\eta^{S^1}_a(U,\bQ)^{hS^1} ~\widehat{\otimes}_{A} ~ \eta^{S^1}_b(W,\bQ)^{hS^1}))\\
   = 
   ~&\im(H^*(\eta^{S^1}_a(U,\bQ)^{hS^1}\otimes_{A} \eta^{S^1}_b(W,\bQ)^{hS^1}))\\
   =~&(u^{I_U(a)}) \otimes_{H^*A} (u^{I_W(b)})=(u^{I_U(a)+I_W(b)}).
\end{align*}
On the other hand, by definition, we have 
\[\im(H^*(\eta_T^{S^1}(U\times W,\bQ)^{hS^1}))= (u^{I_{U\times W}(T)}).\]
Then the inclusion relation implies that
\[I_{U\times W}(T)\geq I_U(a)+I_W(b).\]

The inequality is true for all $a,b\geq 0$ and $a+b\leq T$. Then we conclude the proof.
\end{proof}

\begin{Thm}[Product inequalities]\label{thm: product inequality}For $U \subset T^*M,W \subset T^*N$, we have
\begin{enumerate}[fullwidth]
    \item If $U$ has RCT boundary, then \[    \widehat{c}_k(U \times W) \leq \min_{i+j=k, i,j\geq 0}\{\widehat{c}_i(U) + \widehat{c}_j(W)\}.\]
    \item For any open $U,W$, we have
\[ {c}_k(U \times W) \leq \min_{i+j=k, i,j\geq 0}\{ {c}_i(U) +  {c}_j(W)\}.\]
\end{enumerate}
\end{Thm}
\begin{proof}For the inequality of $\widehat{c}_k$, it follows directly from the index inequality \autoref{coro:index inequality}.
For the inequality of ${c}_k$, we simply observe that if $\eta^{S^1}_a(U,\bQ)=u^i x^{S^1}_a$, and $ \eta^{S^1}_b(W,\bQ)=u^j y^{S^1}_b$, then by \autoref{lemma: fundamental class product decomposition naive form}, $\eta^{S^1}_T(U\times W,\bQ) = i_{a,b}^*(u^{i+j} x^{S^1}_a \otimes y^{S^1}_b) = u^{i+j}i_{a,b}^*(x^{S^1}_a\otimes y^{S^1}_b)$.
\end{proof}

\begin{RMK}\label{remark: index equality iff capacity formula}Due to the definition and the proof, we see that the product inequality~(1) becomes equality if the index inequality \autoref{eq: index inequality} becomes equality for all $T$.  
\end{RMK}

\begin{RMK}As we mentioned in the introduction, we will construct an example that makes the first inequality strict. In fact, we will see later that the same example also makes the second inequality strict.
\end{RMK}

Based on the above discussion, it is crucial to consider morphisms in \autoref{lemma: fundamental class product decomposition naive form} more carefully and for all $a+b\leq T$ in detail. This is harder to control without any geometric condition. Therefore, in below, we will assume $U$ and $W$ have non-degenerate RCT boundaries.

\begin{Lemma}\label{lemma: stalk formula of sheaf convolution}For $F,G\in \Sh_{\geq } (\bR)$. Assume that there exist $A,B\in \bR$ such that $F|_{(-\infty,A)}=0$ and $G|_{(-\infty,B)}=0$ and they are weakly constructible with respect to a locally finite stratification, then we have for all $T\in \bR$ that
\[(F\star G)_T\simeq   \varprojlim_{a+b\leq T} F_a\otimes G_b,\]
and this is actually a finite limit.

Such identification is functorial with respect to $F\to F'$ and $G\to G'$ for $F',G'$ satisfying the same support and constructibility condition.
\end{Lemma}
\begin{proof}We first assume that $F\simeq M_{[c,\infty)}$ and $G\simeq N_{[d,\infty)}$ for $M, N\in \cC$. In this case, a direct computation shows that 
\[(F\star G)_{T}\quad=\quad \begin{cases}M\otimes N, &T\geq c+d\\ 0, &T<c+d    
\end{cases}\quad =\varprojlim_{a+b\leq T} F_a\otimes G_b,\]
and the right-hand side is essentially a finite limit.

This verifies this particular case. In general, we will reduce to this case.

For $H\in \Sh_{\geq}(\bR)$, the definition of microsupport gives for $t<T$ that $H_T\simeq \Gamma({T},H)\simeq \Gamma((t,T],H)$. Notice that $F \star G \in \Sh_{\geq}(\bR)$ by \autoref{lem: convolution estimates}. 

Under the vanishing condition, we then have that $s:\bR^2\to \bR$ is proper on $\supp(F\boxtimes G)$, and 
\[(F \star G)_T=s_!(F\boxtimes G)_T=\Gamma(\{(a,b)\mid a\geq A,b\geq B,a+b\leq T\},F\boxtimes G).\]

Since $\{(a,b)\mid a\geq A,b\geq B,a+b\leq T\}$ is compact, we may assume both $F$ and $G$ are constructible with respect to a finite stratification by considering $F\otimes 1_{(A,T-B+1)}$ here. 

However, for $H\in \Sh_{\geq}(\bR)$ that is constructible with respect to a finite stratification, $H$ can be written as a finite limit of sheaves $M_{[c,d)}$ for $M\in \cC$. Then $M_{[c,d)}$ can be written as $\fib(M_{[c,\infty)}\to M_{[d,\infty)})$. Therefore, $H$ can be written as a finite limit of sheaves $M_{[c,\infty)}$.

Since our final result commutes with finite limits/colimits, the result follows from the simple computations we explained at the beginning. In particular, $\varprojlim_{a+b\leq T} F_a\otimes G_b$ is equivalent to a finite limit since $F_a$ and $G_b$, and the continuation maps $F_a\to F_{a'}$ (also for $G_b \to G_{b'}$) only change at finitely many points.

From the proof, it is clear that the limit formula is functorial in both $F$ and $G$.\end{proof}

\begin{RMK}
The main reason we need the constructibility assumption is that our capacities are defined in terms of stalks of $F^{S^1}(U, \bQ)$ at $T \geq 0$, and taking stalks do not commute with arbitrary limits. This could be avoided by taking costalks as in the APT correspondence \cite{APT-Tatsuki-Zhang}, in which case we have a colimit formula of the Tamarkin convolution, which is exactly the Day convolution. However, we prefer not to reorganize all the capacity package from the beginning using costalks.
\end{RMK}

\begin{Coro}\label{coro:  stalk formula of fundamental class convolution}For $U\subset T^*M,W \subset T^*N$, assume that $F^{S^1}(U,\bQ)$ and $F^{S^1}(W,\bQ)$ satisfy the condition of \autoref{lemma: stalk formula of sheaf convolution}. Then we have a finite limit 
\[\eta_T^{S^1}(U\times W,\bQ) \simeq \varprojlim_{a+b\leq T} \eta^{S^1}_a(U,\bQ)\otimes \eta^{S^1}_b(W,\bQ),\]
and
\[(\eta_T^{S^1}(U\times W,\bQ) )^{hS^1}\simeq   \varprojlim_{a+b\leq T}\eta^{S^1}_a(U,\bQ)^{hS^1}~\widehat\otimes_A ~\eta^{S^1}_b(W,\bQ)^{hS^1}.\]
\end{Coro}
\begin{proof}By \autoref{lemma: stalk formula of sheaf convolution}, the result is true non-equivariantly. However, we can track the proof to see that the $S^1$-action is compatible with this limit. Then we have the first equality.

To see the second equality, we notice that the homotopy fixed point is a limit, so it commutes with the given limit formula. Moreover, the homotopy fixed point is a monoidal functor.
\end{proof}
The first application here is the product formula for $c_1$ and $\widehat{c}_1$.
\begin{Thm}[Product equality for the first capacity]\label{Product formula for c1 RCT}For $U \subset T^*M,W\subset T^*N$, assume that $F^{S^1}(U,\bQ)$ and $F^{S^1}(W,\bQ)$ satisfy the condition of \autoref{lemma: stalk formula of sheaf convolution} and $F^{S^1}(U,\bQ)_T$ is perfect for all $T\geq 0$. 
Then we have 
\[\widehat{c}_1(U \times W) = \min\{\widehat{c}_1(U), \widehat{c}_1(W)\},\qquad {c}_1(U \times W) = \min\{{c}_1(U), {c}_1(W)\}.\]
In particular, the equality holds if $U, W$ have non-degenerate RCT boundaries.
\end{Thm}
\begin{proof}Let $T\geq 0$. As $F^{S^1}(U,\bQ)$ and $F^{S^1}(W,\bQ)$ are vanish on $(-\infty,0)$, then the limit over $a+b\leq T$ is equivalent to the limit over $a+b\leq T$ where $a\leq T$ and $b\leq T$. Then by the limit formula for the convolution \autoref{lemma: stalk formula of sheaf convolution}, we have a commutative diagram of $S^1$-modules
\[\begin{tikzcd}
                                      & {F^{S^1}(U,\bQ)_T\otimes F^{S^1}(W,\bQ)_T} \arrow[d, "r_t"] \arrow[rd, "\eta\otimes\id"] \arrow[ld, "\id \otimes \eta" above] &                                       \\
{F^{S^1}(U,\bQ)_T[-n]} \arrow[rd, "\eta" below] & {F^{S^1}(U\times W,\bQ)_T} \arrow[l] \arrow[r] \arrow[d, "\eta"]                                                   & {F^{S^1}(W,\bQ)_T[-m]} ,\arrow[ld, "\eta"] \\
                                      & \bQ[-m-n]                                                                                                                  &                                      
\end{tikzcd}\]
where $\eta$ denotes the corresponding fundamental classes.
Now, we prove the formulas. 

For $\widehat{c}_1(X)$, we can take the homotopy fixed points and base change to $\bQ = \bQ\fpw{u}/u\bQ\fpw{u}$. Therefore, the above diagram implies that
\[I_{U\times W}(T)=0 \iff I_{U}(T)=I_{W}(T)=0,\]
which implies that
$\widehat{c}_1(U \times W) = \min\{\widehat{c}_1(U), \widehat{c}_1(W)\}.$

For ${c}_1(X)$, we can forget the $S^1$-structure, and use the Gysin sequence to see that $\eta^{S^1}_T(X)\in uH^*C_T^{S^1}(X,\bQ)$ if and only the non-equivariant fundamental class $\eta_T(X)=0$ as in \autoref{equation: non-equivariant by Gyin}. Notice that the forgetful functor is monoidal. Therefore, the diagram again shows that 
\[\eta_T(U\times W)\neq 0 \iff \eta_T(U )\neq 0 \text{ and }\eta_T( W)\neq 0,\] 
which implies that
$ {c}_1(U \times W) = \min\{ {c}_1(U),  {c}_1(W)\}.$
\end{proof}

\begin{RMK}
We give an alternative proof to explain more clearly the reason of using the limit formula. By the limit formula for convolutions \autoref{lemma: stalk formula of sheaf convolution}, we have
\[F^{S^1}(U, \bQ)_T \simeq \varprojlim_{a+b\leq T}F^{S^1}(U,\bQ)_a \otimes F^{S^1}(V, \bQ)_b.\]
The homotopy fixed point of the diagram and perfectness of $F^{S^1}(U,\bQ)$ implies that
\[F^{S^1}(U, \bQ)_T^{hS^1} \simeq \varprojlim_{a+b\leq T}F^{S^1}(U,\bQ)_a^{hS^1} \otimes F^{S^1}(V, \bQ)_b^{hS^1}.\]
We can now take the tensor product on both sides with $\bQ\fpw{u}/u\bQ\fpw{u}$, and obtain an isomorphism of $\bQ$-modules
\[F(U, \bQ)_T \simeq \varprojlim_{a+b\leq T}F(U,\bQ)_a \otimes F(V, \bQ)_b.\]

For $\widehat c_1(U\times W)$, we know that $\sup_{a+b\leq T}\{I_U(a) + I_W(b)\}$ if and only if 
\[\im\big(\varprojlim\nolimits_{a+b\leq T} H^*\eta_a(U) \otimes H^*\eta_b(W)\big) = \bQ[-m-n].\]
Now \autoref{lem: finite lim index} implies that there is a canonical map
\[\begin{tikzcd}
\mathbf{R}^i\varprojlim\nolimits_{a+b\leq T}H^*F(U,\bQ)_a \otimes H^{j-*}F(W,\bQ)_b \ar[r, Rightarrow] \ar[d,"\varprojlim\nolimits_{a+b\leq T} H^*\eta_a(U) \otimes H^{j-*}\eta_b(W)" left] & H^{i+j}F(U \times W; \bQ)_T \ar[d,"\eta_{T}(U\times W)"] \\
\mathbf{R}^i\varprojlim\nolimits_{a,b\geq 0,a+b\leq T}\bQ[-m-n] \ar[r, Rightarrow] & \bQ[-m-n].
\end{tikzcd}\]
Since $\varprojlim\nolimits_{a+b\leq T}H^*F(U,\bQ)_a \otimes H^{j-*}F(W,\bQ)_b$ is always a direct summand of $H^jF(U \times W; \bQ)_T$ as $\bQ$-modules, we can conclude that
\[\im(H^*F(U \times W; \bQ)_T) = \im\big(\varprojlim\nolimits_{a+b\leq T} H^*\eta_a(U) \otimes H^*\eta_b(W)\big) = \bQ[-m-n].\]
\end{RMK}

The second application is an equality for a disk factor in large radius.
\begin{Thm}[Stabilization invariance]\label{coro: disk stablization}For $U \subset T^*M$, assume that $F^{S^1}(U,\bQ)$  satisfies the condition of \autoref{lemma: stalk formula of sheaf convolution}. Let $W=D(r)\subset T^*\bR$. Then if $\widehat{c}_k(U),c_k(U)\leq r$, we have
\[\widehat{c}_k(U\times D(r))=\widehat{c}_k(U), \quad c_k(U\times D(r))=c_k(U).\]
In particular, the equality holds if $U$ has RCT boundary.
\end{Thm}
\begin{proof}
The computation for general balls \autoref{lem: symplectic ball} shows that $F^{S^1}(D(r),\bQ)$ is known to satisfy the condition of \autoref{lemma: stalk formula of sheaf convolution}, and moreover,
\[\eta_b^{S^1}(D(r),\bQ): F^{S^1}(D(r),\bQ)_b\to F^{S^1}(T^*\bR,\bQ)_b\]
is an isomorphism for $0\leq b< r$. Then \autoref{coro:  stalk formula of fundamental class convolution} shows that when $0\leq T< r$ (in which case $0 \leq a,b < r$ for $0\leq a+b< T$), we have
\[F_T^{S^1}(U\times D(r),\bQ) \simeq \varprojlim_{a+b\leq T} F^{S^1}(U,\bQ)_a\otimes F^{S^1}(D(r),\bQ)_b \simeq F^{S^1}(U,\bQ)_T\]
is an equivalence. Then the result follows by the definition of capacities.    
\end{proof}

Next, we discuss the algebraic obstruction to the product inequality being an equality in general. First, we have a direct criterion as a direct corollary of \autoref{prop: c=hat c, direct summand criterion}. Recall the exact triangle
\[H^*F^{S^1}(U \times W, \bQ)_T^{hS^1} \to H_c^*(M \times N,  \bQ)\fpw{u} \to H^*F^{S^1,out}(U \times W, \bQ)_T^{hS^1},\]
and $\beta_T^{S^1}(U\times W, \bQ)$ is the image of $\mathbf{1}_{[M\times N]}$ in the positive Chiu--Tamarkin complex.
\begin{Prop}\label{prop: product equality 1}
For $U, W$ open domains such that $F^{S^1}(U, \bQ)_T^{hS^1}$ and $F^{S^1}(W,\bQ)_T^{hS^1}$ are perfect $\bQ\fpw{u}$-modules for $T \in \bR$. Then, for $T_0 = \min_{i+j=k}\{\widehat c_i(U) + \widehat c_j(W)\}$ and all $0\leq T < T_0$,
\[u^{k-1}\beta_{T}^{S^1}(U \times W, \bQ) = 0,\]
then the product inequality in \autoref{thm: product inequality}~(1) becomes an equality.
\end{Prop}

We can also provide a more specific criterion using the formula of finite homotopy limits:

\begin{Prop}\label{prop: product equality summand}
For $U, W$ open domains such that $F^{S^1}(U, \bQ)_T^{hS^1}$ and $F^{S^1}(W,\bQ)_T^{hS^1}$ are perfect $\bQ\fpw{u}$-modules for $T \in \bR$. Then, if 
\[\varprojlim_{a+b\leq T}H^*F^{S^1}(U, \bQ)_a^{hS^1} \otimes H^*F^{S^1}(W, \bQ)_b^{hS^1} \to H^*F^{S^1}(U \times W, \bQ)^{hS^1}_T\]
is a direct summand, then the product inequality in \autoref{thm: product inequality}~(1) becomes an equality.
\end{Prop}
\begin{proof}
This follows from \autoref{cor: complete eilenberg-moore} and \autoref{lem: finite lim index}. In fact, we have
\[\begin{tikzcd}
{\mathbf{R}^i\varprojlim_{a+b\leq T}H^*F^{S^1}(U, \bQ)_a^{hS^1} \otimes H^{j-*}F^{S^1}(W, \bQ)_b^{hS^1}} \ar[r, Rightarrow] \ar[d,"{\varprojlim H^*\eta_a^{S^1}(U,\bQ) \otimes H^*\eta_b^{S^1}(W,\bQ)}" left] & {H^{i+j}F^{S^1}(U \times W, \bQ)^{hS^1}_T} \ar[d,"{H^*\eta_T^{S^1}(U\times W, \bQ)}"] \\
{\mathbf{R}^i\varprojlim_{a,b\geq 0,a+b\leq T}\bQ\fpw{u}} \ar[r, Rightarrow] & {\bQ\fpw{u}.}
\end{tikzcd}\]
Hence by \autoref{lem: finite lim index} we have
\[\im(H^*\eta_T^{S^1}(U\times W, \bQ)) \subset \im\big(\varprojlim\nolimits_{a+b\leq T} H^*\eta_a^{S^1}(U,\bQ) \otimes H^*\eta_b^{S^1}(W,\bQ)\big).\]
Suppose $\varprojlim_{a+b\leq T} H^*F^{S^1}(U,\bQ)_a^{hS^1} \otimes H^*F^{S^1}(W,\bQ)_b^{hS^1}$ is a direct summand of $F^{S^1}(U \times W, \bQ)_T^{hS^1}$. Then the images are the same and thus
\[I_{U\times W}(T) = \sup_{a,b\geq 0;a+b\leq T}\{I_U(a) + I_W(b)\}.\]
The index equality now implies the capacity equality.
\end{proof}

\section{Chain model for toric domains}\label{section: chain model overall}In this section, the coefficient ring $\Bbbk$ shall be a field containing $\bQ$, or $\bZ$, depending on the requirement.

Let $U\subset T^*M$ be open and $Z=T^*M\setminus U$.  Recall that the cyclic trace of $\Sh_{Z}(M,\sT)$ is denoted by $F^{S^1,out}(U,\Bbbk)$. We have the tautological fiber sequence
\[F^{S^1}(U,\Bbbk)\to F^{S^1}(T^*M,\Bbbk)\to F^{S^1,out}(U,\Bbbk)  \text{  in  } \Fun(BS^1,\sT),\]
and its homotopy fixed point at  $T\geq 0$, 
\[F^{S^1}(U,\Bbbk)_T^{hS^1}\to \Gamma_c(M,\Bbbk)\fpw{u} \to F^{S^1,out}(U,\Bbbk)_T^{hS^1} \text{  in  } D(\Bbbk\fpw{u})_{comp}.\]

We similarly define the action window version for $-\infty \leq \epsilon<T<\infty$: 
\[F_{(\epsilon,T]}^{S^1,out}(U,\Bbbk)\coloneqq \fib(F^{S^1,out}(U,\Bbbk)_T\to F^{S^1,out}(U,\Bbbk)_\epsilon),\]
which gives
\[ F_{(\epsilon,T]}^{S^1,out}(U,\Bbbk)^{hS^1}\simeq \fib(F^{S^1,out}(U,\Bbbk)_T^{hS^1}\to F^{S^1,out}(U,\Bbbk)_\epsilon^{hS^1}).\]
Here, we use the notation $F_{(-\infty,T]}^{S^1,out}(U,\Bbbk)\coloneqq F^{S^1,out}(U,\Bbbk)_T$.

The main goal of this section is to formulate a chain model for $F^{S^1}(X_\Omega,\Bbbk)_T^{hS^1} $ for a toric domain $X_\Omega \subset \bR^{2n}$.

\subsection{Geometric setup}

Let $V=\bR^n$ and identify $T^*V$ with $\mathbb C^n$ in the standard way.  Let 
\[\Omega\overset{open}{\subset}\bR_{\geq0}^n,\qquad
  U_\Omega=\bR_{\geq0}^n\setminus\Omega.
\]

Then for $\mu(\zeta)=(\pi|\zeta_1|^2,\dots,\pi|\zeta_n|^2):\bC^n\to \bR_{\geq0}^n$, we shall take 
\[X_\Omega=\mu^{-1}(\Omega),\qquad X_\Omega^{out}=\mu^{-1}(U_\Omega),\]
where we call $X_\Omega$ the open toric domain associated with $\Omega$.

\begin{Def}
For $\Omega$ and $U_\Omega$, we impose the following \textbf{geometric assumptions and notations}:
\begin{enumerate}[fullwidth]
    \item We say an open set $\Omega$ is a domain if  $\Omega$ is connected.
    
    \item We say an open $\Omega$ is continuously/Lipschitz/H\"older/smooth star-shaped domain if there exists a continuous/Lipschitz/H\"older/smooth function
    \[\rho: S^{n-1}_+=S^{n-1}\cap\bR_{\geq0}^n\to [0,\infty) \]
    such that 
    \[\Omega=\{\zeta=ru\mid u\in S^{n-1}_+, ~r\geq 0,~ r\rho(u)<1\}.\]
    In this case, $ U_\Omega\cong \{\rho>0\}\times [1,\infty) , ru/\rho\mapsfrom (u,r)$ is locally contractible since $\{\rho>0\}$ is an open set in $S^{n-1}_+$.

    We note that for a continuously star-shaped domain $\Omega$, $\Omega$ is bounded if and only if the corresponding $\rho>0$ on $S^{n-1}_+$.

    \item Let $\Omega$ be a domain. For $m\in\bZ^n$ and $T\in\bR$, set
\begin{gather*}
    \Omega(T;m)=\{\zeta\in \Omega\mid T+\langle m,\zeta\rangle\geq 0\},\quad U_\Omega(T;m)=\{\zeta\in U_\Omega\mid T+\langle m,\zeta\rangle\geq 0\},\\
    I(m)=m_1+\cdots+m_n.
\end{gather*}
  We say $\Omega$ is homologically star-shaped domain if $U_\Omega$ is locally contractible and for every $T\in \bR$, there exists a number $C(T)$ such that if $I(m)\geq C(T)$, then $H^*_c(U_\Omega(T;m),\Bbbk)=0$.
\end{enumerate}
\end{Def}

\begin{eg}Here, we give some examples of Lipschitz star-shaped domains.
\begin{enumerate}
    \item Let $\Omega=\bR^{n}_{\geq0}$, then we can take $\rho=0$ to see that $\Omega$ is smooth star-shaped.
    \item Let $\widetilde{\Omega}$ be a generic union of finitely many open half-spaces containing the origin, then $\Omega=\widetilde{\Omega}\cap \bR^{n}_{\geq0}$ is generally a non-smooth Lipschitz star-shaped domain.\end{enumerate}
\end{eg}

\begin{Lemma}If $\Omega$ is bounded smooth star-shaped domain, then $X_\Omega$ is a smooth star-shaped domain.     
\end{Lemma}
\begin{proof}Let $j(\zeta)=|\zeta|\rho(\zeta/|\zeta|)>0$ for $\zeta\neq 0$. Then $j$ is smooth, homogeneous, and $\{j=1\}=\partial \Omega$ (for the boundary in $\bR^n_{\geq0}$). Therefore, for $h(z)=j(\mu(z))$ for $z\neq 0$, we have $h$ is smooth at $z\neq 0$ and homogeneous of degree $2$. Moreover, we have $\partial X_\Omega=\{h=1\}$ is smooth.

We then use the Euler identity for quadratic property of $h$ to see that for $R_z$ the radial vector field of $\bC^n$ we have
\[dh_z(R_z)=2h(z)=2,\quad z\in \partial X_\Omega=\{h=1\},\]
which implies the transversality of $R_z$ and $\partial X_\Omega$.    
\end{proof}

\begin{Prop}\label{proposition: regularity for lipschitz star-shaped domain}For a Lipschitz star-shaped $\Omega$ domain, we have
\begin{enumerate}
    \item $H_c^*(\Omega,\Bbbk)=0$ and $H_c^*(U_\Omega,\Bbbk)=0$.
    \item $\Omega$ is homologically star-shaped.
\end{enumerate}    
\end{Prop}
\begin{proof}We write $\rho$ for the defining function of $\Omega$. 
We notice that the following excision triangle
\begin{gather*}
    \Gamma_c(\Omega,\Bbbk)\to \Gamma_c(\bR^n_{\geq0},\Bbbk)\simeq 0\to \Gamma_c(U_\Omega,\Bbbk),\\ \Gamma_c(\Omega(T;m),\Bbbk)\to \Gamma_c(\{\zeta\in\bR^n_{\geq0} \mid T+\langle m,\zeta\rangle\geq 0\},\Bbbk)\to \Gamma_c(U_\Omega(T;m),\Bbbk)
\end{gather*}

The statement (1) is clear since $U_\Omega\cong \{\rho>0\}\times [1,\infty)$ and then by Künneth we have \[\Gamma_c(U_\Omega,\Bbbk)\simeq \Gamma_c(\{\rho>0\},\Bbbk)\otimes \Gamma_c([1,\infty),\Bbbk)=0,\]
and then $\Gamma_c(\Omega,\Bbbk)\simeq 0$ by the excision triangle.

For (2), we only need to show that for every $T$, there exists a number $C(T)>n$ such that if $I(m)\geq C(T)$, then $H^*_c(U_\Omega(T;m),\Bbbk)=0$. In fact, we have that $\Gamma_c(\{\zeta\in\bR^n_{\geq0} \mid T+\langle m,\zeta\rangle\geq 0\},\Bbbk)\simeq \Bbbk$ only when $T\geq 0$ and $m\in \bZ^n_{<0}$, and otherwise vanishes. When $I(m)>n$, by definition of $I$ we must have $m\notin \bZ^n_{<0}$. Consequently, the excision sequence shows that $H^*_c(\Omega(T;m),\Bbbk)=0$ if and only if $H^*_c(U_\Omega(T;m),\Bbbk)=0$.

Let $j(\zeta)=|\zeta|\rho(\zeta/|\zeta|)\geq 0$ for $\zeta\neq 0$ be the gauge function of $\Omega$. Since $j$ is homogeneous, we can continuously extend $j$ to $0$, which is still denoted by $j$. Moreover, $\rho$ is Lipschitz; one can check that $j$ is also Lipschitz on $\bR^n_{\geq0} $. We set $L\geq 1$ to be the Lipschitz constant of $j$. 

Let $T^+=\max\{T,0\}$, then we take $C(T)=n(2LT^+ +1)$. If $I(m)>C(T)$, we show that $U_\Omega(T;m)\cong  E\times [0,\infty)$ for a relative open set $E$. If so, the vanishing result we need follows from Künneth.

Now, under the condition $I(m)>C(T)$, there exists $i$ such that $m_i\geq 2LT^++1$, which implies that $m_i/2L>T$. Without loss of generality, we assume $i=1$, and we write $a=1/(2L)$. First, we have a change of coordinates
\[\bR^n_{\geq0} \cong [a,\infty)\times\bR^{n-1}_{\geq0} , \zeta\mapsto (t,v)=\Big(\zeta_1+a,\frac{(\zeta_2,\dots,\zeta_n)}{\zeta_1+a}\Big).\]
Now, we consider the projection
\[U_\Omega(T;m)\eqqcolon A\subset \bR^n_{\geq0} \cong [a,\infty)_t\times\bR^{n-1}_{\geq0} \to\bR^{n-1}_{\geq0}.  \]
Then its fiber over $v\in \bR^{n-1}_{\geq0}$ under the $t$-coordinate is
\[A_v=\{t\geq a\mid j(t-a,tv)\geq 1,\quad T+\langle m,(t-a,tv)\rangle \geq 0\}.\]
We then claim that $A_v$ is either empty or of the form $[x,\infty)$. We prove it by showing that if $t\in A_v$, then for all $s\geq t$, we have $s\in A_v$. Equivalently, we need to show that
\[j(s-a,sv)\geq 1,\quad T+\langle m,(s-a,sv)\rangle \geq 0\]
when $t$ satisfies the similar inequalities. 

Here $t>0$, then $\lambda=s/t\geq 1$. We then have
\[(s-a,sv)=\lambda(t-a,tv)+a(\lambda-1)e_1.\]
Therefore, by Lipschitz condition, we have
\[j(s-a,sv)\geq \lambda j(t-a,tv)-La(\lambda-1)\geq 1,\]
and
\[T+\langle m,(s-a,sv)\rangle = \lambda(T+\langle m,(t-a,tv))+(\lambda-1)(am_1-T)\geq 0.\]

Therefore, we then show that the projection
\[A\subset \bR^n_{\geq0} \cong [a,\infty)_t\times\bR^{n-1}_{\geq0} \to\bR^{n-1}_{\geq0}\]
is a fiber bundle onto its image whose fiber are $[x,\infty)$. However, such a bundle is trivial. And then $U_\Omega(T;m)=A\cong E \times [0,\infty)$ and we conclude the proof.
\end{proof}
\begin{RMK}The same proof actually gives that: for a continuously star-shaped domain $\Omega$, if there exists $L\geq 0$ such that the gauge function $j$ satisfies $j(\zeta+he_i)\geq j(\zeta)-Lh$ for all $h>0$ and $i=1,\dots,n$, then $\Omega$ is homologically star-shaped. Such examples include monotone open sets, which satisfy that for $\eta$ with $0\leq \eta\leq \zeta$ and $\zeta\in \Omega$, then $\eta\in \Omega$.    
\end{RMK}
\begin{RMK}
    There exist some non-Lipschitz homologically star-shaped domains, and some of them can have a Hölder gauge function. But not all Hölder star-shaped domains are homologically star-shaped.
\end{RMK}

\subsection{Statement of the chain model}\label{section: chain model statement}

We start from preparing some algebraic constructions:

For $-\infty\leq \epsilon<T<\infty$, we set $\Omega(-\infty;m)=\varnothing$ and
\[G_{m,(\epsilon,T]}=C^{*+n+2I(m)}_c({\Omega(T;m)},\Omega(\epsilon;m) ;\Bbbk  ).\]

We consider the following $\Bbbk $-linear map associated with natural inclusions
\[t_{i,m}:G_{m,(\epsilon,T]} \to G_{m-e_i,(\epsilon,T]}.\]

Let $\bigwedge^\bullet\Bbbk ^n$ be the exterior algebra of $\Bbbk ^n$ with the basis $e_{I}=e_{i_1}\wedge \cdots\wedge e_{i_p}$. We set $|e_i|=1$ and $\iota_{e_i}:\bigwedge^\bullet\Bbbk ^n\to \bigwedge^\bullet\Bbbk ^n$ the inner product $\iota_{e_{i_r}}(e_{I})=(-1)^{r-1}e_{I\setminus \{i_r\}}$. We also write the parity function $(-1)^N$ on $\bigwedge^\bullet\Bbbk ^n$, i.e. $(-1)^N|_{\wedge^q\Bbbk ^n}=(-1)^q$. 

Under the notation, we then consider the graded $\Bbbk\fpw{u}$-module
\[I_{(\epsilon,T]}^{S^1}(X_\Omega)\coloneqq \bigoplus_{m\in \bZ^n}\left( G_{m,(\epsilon,T]}  \fpw{u} \otimes \bigwedge\nolimits^\bullet \Bbbk ^n\right),\]
with the differential
\begin{equation}\label{equation: differential chain model}d|_m=d_{\mathrm{sing}}\otimes (-1)^N+\sum_{i=1}^n  (m_iu+t_i)\otimes \iota_{e_i},\end{equation}
where $d_\mathrm{sing}$ should be understood as the usual singular differential tensored with the identity map on $\Bbbk\fpw{u}$. We write $\widehat{I}_{(\epsilon,T]}^{S^1}(X_\Omega)$ the $u$-adic completion of $I_{(\epsilon,T]}^{S^1}(X_\Omega)$. 

If $\epsilon =-\infty$, we write $I_{T}^{S^1}(X_\Omega)=I_{(\epsilon,T]}^{S^1}(X_\Omega)$.

\begin{RMK}\label{remark: toric chain model positive lattice trivial}Notice that when $T>\epsilon>0$ and $m\in \bZ_{\geq 0}^n$, we have 
\[\Omega=\Omega(T;m)=\Omega(\epsilon;m),\quad U_\Omega=U_\Omega(T;m)=U_\Omega(\epsilon;m),\]
and then the cochain corresponding to $m\in \bZ_{\geq 0}^n$ contributes acyclic terms in this case.
\end{RMK}
\begin{RMK}\label{remark: abuse notation for the chain model for locally closed sets}Here, compare to \autoref{remark: coend formula for general locally closed V}, for any locally closed subset $V\subset \bR^n_{\geq 0}$, we will use the notation $X_V$ and $I_{(\epsilon,T]}^{S^1}(X_V)$, where we simply replace $\Omega$ by $V$ everywhere in the formula of $I_{(\epsilon,T]}^{S^1}(X_\Omega)$, in a slightly general way even though it may happen that $X_V$ is not a toric domain.
\end{RMK}

The main result in this section is the following theorem.
\begin{Thm}\label{thm: chain model toric domain}For $-\infty \leq \epsilon<T<\infty$. Let $\Omega\subset \bR^n$ be an open set, and $X_\Omega$ be its toric domain, we have

\begin{enumerate}[fullwidth]
\item The chain complex $\widehat{I}^{S^1}_{(\epsilon,T]}(X_\Omega)$ over $\Bbbk\fpw{u}$ is a chain model of $F_{(\epsilon,T]}^{S^1}(X_\Omega)^{hS^1}$, i.e. $F_{(\epsilon,T]}^{S^1}(X_\Omega)^{hS^1}\simeq  \widehat{I}^{S^1}_{(\epsilon,T]}(X_\Omega)$ in $D(\Bbbk\fpw{u})_{comp}$. 

\item If $(\epsilon_2,T_2]\subset (\epsilon_1,T_1]$, then the natural morphism $F_{(\epsilon_1,T_1]}^{S^1}(X_\Omega)^{hS^1}\to F_{(\epsilon_2,T_2]}^{S^1}(X_\Omega)^{hS^1}$ is represented by $\widehat{I}_{(\epsilon_1,T_1]}^{S^1}(X_\Omega)\to\widehat{I}_{(\epsilon_2,T_2]}^{S^1}(X_\Omega)$. 

\item If $\Omega_1\subset \Omega_2 \subset \bR^n_{\geq 0}$ are two open sets, then the natural map $F_{(\epsilon,T]}^{S^1}(X_{\Omega_1})^{hS^1}\to F_{(\epsilon,T]}^{S^1}(X_{\Omega_2})^{hS^1}$ is represented by $\widehat{I}^{S^1}_{(\epsilon,T]}(X_{\Omega_1})\to \widehat{I}^{S^1}_{(\epsilon,T]}(X_{\Omega_2})$, and we have the following Mayer-Vietoris short exact sequence
\[0\to \widehat{I}^{S^1}_{(\epsilon,T]}(X_{\Omega_1\cap \Omega_2})\to \widehat{I}^{S^1}_{(\epsilon,T]}(X_{\Omega_1})\oplus \widehat{I}^{S^1}_{(\epsilon,T]}(X_{\Omega_2})\to \widehat{I}^{S^1}_{(\epsilon,T]}(X_{\Omega_1\cup \Omega_2})\to 0.\]

\item If $\Omega$ is a homologically star-shaped domain (for example, a Lipschitz star-shaped domain by \autoref{proposition: regularity for lipschitz star-shaped domain}), then $I^{S^1}_{(\epsilon,T]}(X_\Omega)$ is $u$-adic complete and $\widehat{I}^{S^1}_{(\epsilon,T]}(X_\Omega)\simeq I^{S^1}_{(\epsilon,T]}(X_\Omega)$. In particular, we only need to consider this direct sum chain complex in this case.

\item  If we assume in addition that $\Bbbk$ is a field containing $\bQ$ and if $\Omega$ is a Lipschitz star-shaped domain, there exists a filtration on $\widehat{I}^{S^1}_{(\epsilon,T]}(X_\Omega)\simeq F_{(\epsilon,T]}^{S^1}(X_\Omega)^{hS^1}$, whose associated spectral sequence has $E_1$-page as follows 
\[E_1^{p,q}=\bigoplus_{m\notin \bZ^n_{\geq 0}, I(m)= p,j} H^{q+n+3p-j}_c({\Omega(T;m)},\Omega(\epsilon;m) ;\Bbbk )   \otimes \wedge^j \ker(m) \Rightarrow H^{p+q} F_{(\epsilon,T]}^{S^1}(X_\Omega,\Bbbk)^{hS^1},\]
and the spectral sequence is $E_n$-degenerate. 
\end{enumerate} 
\end{Thm}

\begin{RMK}\label{remark: non-equivariant chain model}We set \[I_{(\epsilon,T]}(X_\Omega)=I_{(\epsilon,T]}^{S^1}(X_\Omega)\otimes_{\Bbbk\fpw{u}}\Bbbk, \textnormal{ equipped with the differential } d_0=d\otimes_{\Bbbk\fpw{u}}\Bbbk,\]
or equivalently require $u=0$ in $(I_{(\epsilon,T]}^{S^1}(X_\Omega),d)$. Therefore, we have $I_{(\epsilon,T]}(X_\Omega)\simeq F_{(\epsilon,T]}(X_\Omega,\Bbbk)$ in $D(\Bbbk)$. This is clear since the functor forgetting the $S^1$-action is given by $\Fun(BS^1,D(\Bbbk))\simeq D(\Bbbk\fpw{u})_{comp}\to D(\Bbbk)$ where $M \to M\otimes_{\Bbbk\fpw{u}}\Bbbk $.
\end{RMK}

\begin{RMK}\label{remark: in and out for chain model}We also have versions for $F_{(\epsilon,T]}^{S^1,out}(X_\Omega)^{hS^1}\simeq \widehat{O}_{(\epsilon,T]}^{S^1}(X_\Omega)$, where $\widehat{O}_{(\epsilon,T]}^{S^1}(X_\Omega)\coloneqq \widehat{I}_{(\epsilon,T]}^{S^1}(X_{U_\Omega})$ (recall \autoref{remark: abuse notation for the chain model for locally closed sets}), and all the other results are true by the same argument. But in this case, we need to assume that $U_\Omega$ is locally contractible.

One can also argue in the following way: by \autoref{thm: chain model toric domain}~(3), we take $\Omega_1=\Omega$ and $\Omega_2=\bR^n_{\geq 0}$, and we have the commutative diagram of fiber sequences
\[\begin{tikzcd}
{F_{(\epsilon,T]}^{S^1}(X_{\Omega})^{hS^1}} \arrow[r] \arrow[d, "\simeq"] & {F_{(\epsilon,T]}^{S^1}(\bC^n)^{hS^1}} \arrow[d, "\simeq"] \arrow[r] & {F_{(\epsilon,T]}^{S^1,out}(X_{\Omega})^{hS^1}} \arrow[d, "\simeq"] \\
{\widehat{I}^{S^1}_{(\epsilon,T]}(X_{\Omega})} \arrow[r]                  & {\widehat{I}^{S^1}_{(\epsilon,T]}(\bC^n)} \arrow[r]                  & C,                                                               
\end{tikzcd}\]
where the first row is the in/out tautological sequence. In particular, the cofiber $C$ is equivalent to $F_{(\epsilon,T]}^{S^1,out}(X_\Omega)^{hS^1}$. However, one chain model of $C$ is clearly given by $\widehat{O}_{(\epsilon,T]}^{S^1}(X_\Omega)$ due to the open/closed excision of the compactly supported cohomology.
\end{RMK}

\begin{RMK}For $G= C_\ell$, the similar method for $S^1$ will give similar chain model. In fact, we will see in the proof that compare to the
$S^1$-case, the major difference is simply taking $hG$ for the trivial $G$-module $\Bbbk$ (restricted from $S^1$-action) rather than taking $hS^1$. So, we will not explain too much on the finite cyclic group case.
\end{RMK}

\begin{RMK}As we explain in the introduction, the construction is mainly motivated by understanding Irie's chain model of $S^1$-equivariant symplectic homology for toric domain in $\bR^4$ \cite{Irie_ChainModel}. We will explain in \autoref{section: Irie's model} how our chain model (when $n=2$ and $\bQ=\Bbbk$) recovers Irie's chain model.
\end{RMK}

Here, we give some quick applications of the theorem.
\begin{Coro}For $T\in \bR$. Let $\Omega\subset \bR^n$ be an open set, and $X_\Omega$ be its toric domain, we have
\[C^{S^1}_T(X_\Omega,\Bbbk)\simeq \Big( \prod_{m\in\bZ^n} \left[C^{BM}_{-*+2n+2I(m)}(\Omega(T;m),\Bbbk)\fpw{u}\otimes (\bigwedge\nolimits^\bullet \Bbbk ^n)^\vee\right],D\Big),\]
where $|e_S^\vee| =-|S|$ in $(\bigwedge\nolimits^\bullet \Bbbk ^n)^\vee$ and
\[D|_m=d_{BM}\otimes\id
+
\sum_{i=1}^n
(m_i u +t_i^\vee )\otimes(e_i^\vee\wedge -).
\]
\end{Coro}
\begin{proof}
    Using \autoref{lemma: linear dual of S1-action} and the equivalence of \autoref{thm: fixed point as a monoidal functor}, we have 
\[C_T(X_\Omega,\Bbbk)\simeq \Hom_{\Bbbk\fpw{u}}(\widehat{I}^{S^1}_{T}(X_\Omega),\Bbbk\fpw{u}[-n]).\]  
However, the completion is a left-adjoint functor and $\Bbbk\fpw{u}$ is complete already, so we can remove the completion $\wedge$ from the hom-computation. Then we conclude using the Poincar\'e duality $\Hom_{\Bbbk}(C_c^*(\Omega,\Bbbk),\Bbbk[-n])\simeq C^{BM}_{-*+n}(\Omega,\Bbbk)$.
\end{proof}

\begin{Coro}\label{coro: chain model compute SH}Let $\Omega\subset \bR^n_{\geq0}$ be a bounded smooth star-shaped domain. For $T \notin \Spec(\partial X_\Omega)$, we have
\[SH^*_{(-\infty,T)}(\overline{X}_\Omega,\bQ)\simeq H^* \Big( \prod_{m\in\bZ^n} \left[C^{BM}_{-*+2n+2I(m)}(\Omega(T;m),\bQ)\otimes (\bigwedge\nolimits^\bullet \bQ ^n)^\vee\right],D\Big),\]
where $|e_S^\vee| =-|S|$ in $(\bigwedge\nolimits^\bullet \bQ ^n)^\vee$ and 
\[D|_m=d_{BM} +\sum_{i}^n t_i^\vee \otimes(e_i^\vee\wedge -).\]
\end{Coro}
\begin{proof}Recall \autoref{remark: non-equivariant chain model} and $I_{(-\infty,T]}(X_\Omega)$ the non-equivariant chain model. By \cite{Hochschild-Kuo-Shende-Zhang}, we have
\[SH^*_{(-\infty,T)}(\overline{X_\Omega},\bQ)\simeq HH^*_\sT(\sT(U),(-\infty,T])= H^*\Hom_\bQ(F (X_\Omega,\bQ)_T,\bQ[-n]),\]
where the second equality is a formal consequence of the definition of $F$.

Therefore, we have
\[I_{(-\infty,T]}(X_\Omega)\simeq F^{S^1}(X_\Omega,\bQ)_T^{hS^1}\otimes_{\bQ\fpw{u}}\bQ \simeq F (X_\Omega,\bQ)_T\]
by \autoref{thm: chain model toric domain}, and we conclude by a similar dual computation.
\end{proof}

\subsection{Sheaf quantization of Hamiltonian rotation}
Let $z \in \bR$. The $2$-dimensional rotation $\varphi_z(u)=\exp{(-2i\pi z)}u$ on $\bC_u$ is the Hamiltonian flow of the Hamiltonian function $H(u)=\pi|u|^2 $. Here, we identify $\bC_u$ with $T^*\bR_q$ by $u=q+ip$.

In \cite[Proposition 3.10]{chiu2017}, Chiu constructs a sheaf quantization of Hamiltonian rotation in all dimensions, particularly for the $2$-dimensional  $\varphi_z$, say 
$\mathcal{S} \in \Sh(\bR_z\times \bR_{q}\times \bR_{q'}; \sT)$.

The construction of $\mathcal{S}$ uses generating functions, which is summarized in Appendix \ref{appendix:generating function}.
Restricting to $z \in (-1/4,0)$, one can see that 
\[ 
\mathcal{S}|_{\{z\}\times\bR_q \times \bR_{q'}}\cong {\Bbbk }_{\{(z,q,q',t):t+ S_{H}(z,q,q')\geq 0\}},   
\]
where $S_H$ is a generating function of the Hamiltonian rotation: 
\begin{equation}\label{equation: formula of S_H}
S_H(z,q,q')= \frac{q^2+{q'}^2}{2\tan(2\pi z)} -\frac{qq'}{\sin(2\pi z)}.    
\end{equation}

One can extend it to the entire $\bR_z$ using convolution relative to $z$, where a detailed discussion on a more general formula will be explained in Appendix \ref{appendix:generating function}, especially \autoref{eg: rotation quantization using V-type GF}. 

When restricted to $ (-N/4,N/4)$ for $N\in \bN$, we have the following formula \cite[Equation (3.2)]{Capacities2021}:
\begin{equation}\label{equation: formula of GF sheaf quantization}
\mathcal{S}|_{\{z\}\times\bR_{q} \times \bR_{q'}}\cong  \begin{cases}
     \pi_{(q_2,\dots,q_N)!}{\Bbbk }_{\{(z,q_1,\dots,q_{N+1},t):t+\sum_{j=1}^N S_{H}(z/N,q_j,q_{j+1})\geq 0\}}, &z\in (-N/4,0)\\
      {\Bbbk }_{\Delta_{\bR_{q}}\times [0,\infty)}, &z=0\\
     \pi_{(q_2,\dots,q_N)!}{\Bbbk }_{\{(z,q_1,\dots,q_{N+1},t):t+\sum_{j=1}^N S_{H}(z/N,q_j,q_{j+1})\geq 0\}}[N], &z\in (0,N/4),
 \end{cases}  
\end{equation}
where we set $(q_1,q_{N+1})=(q,q')$.
\begin{RMK}In \cite[Equation (3.2)]{Capacities2021}, the formula at $z>0$ does not involve the shifting $N$, which is a mistake. However, we do not need the $z>0$ part in loc. cit., which does not invalidate main results therein.
\end{RMK}

In \textnormal{\cite[Formula (38)]{chiu2017}}, the first part of the following proposition is stated; in \cite[Proposition 3.9]{Capacities2021} and \cite[Proposition 4.10]{zhangthesis}, we described  the sheaf $\cE$ on $\bR_{z\leq 0}$ in terms of the sheaf quantization $\mathcal{S}$. Here, we provide an explicit formula on the sheaf $\cE$ on entire $\bR_z$; see \cite[Proposition 3.11]{Capacities2021} and \cite[Proposition 4.12]{zhangthesis}.
\begin{Prop}[{\cite[Propositions 3.9 \& 3.11]{Capacities2021} \cite[Proposition 4.10 \& 4.12]{zhangthesis}}]\label{prop:chiu result on loop sheaf} There exists a (unique) sheaf ${\mathcal{E}}\in \Sh(\bR_z)$ such that we have an isomorphism in $\Sh(\bR_z\times \bR_t)$:
\begin{equation}\label{loop sheaf of GF model of rotation non-equi}
 \pi_{q!}\Delta_q^{*} \mathcal{S} \cong \cE \boxtimes {\Bbbk }_{[0,\infty)},   
\end{equation}
where the sheaf in $\Sh(\bR_z)$ is
\begin{equation}\label{E sheaf of GF model of rotation}
{\mathcal{E} }\cong  \bigoplus_{k\in \bZ}\Bbbk _{(k,k+1]}[1+2k].
\end{equation}
\end{Prop}
\begin{proof}By \autoref{equation: formula of GF sheaf quantization}, we find the problem is reduced to compute (where we write $q_1 = q_{N+1}$)
\[\Gamma_c\left(\left\{(q_1,\dots,q_{N}):t+\sum\nolimits_{j=1}^N S_{H}(z/N,q_j,q_{j+1})\geq 0\right\},\Bbbk \right).\]
By Morse lemma, we only need to compute at the unique critical value of the quadratic form $Q_z(q_1,\dots,q_{N}) = \sum_{j=1}^NS_H(z/N, q_j, q_{j+1})$ at the point $t=0$. We can diagonalize the quadratic form under the coordinate change $x_j = \sum_{i \in C_N}q_{i+1}\omega^{ij}_N$. Under the basis $(x_0, x_1, \dots, x_{N/2}) \in \bR \times \bC \times \dots \times \bC \times \bR$ for $N$ even and $(x_0, x_1, \dots, x_{(N-1)/2}) \in \bR \times \bC \times \dots \times \bC$ for $N$ odd, we have \cite[Equation (3.17)]{Capacities2021}
\[Q_z(x_0,x_1, \dots, x_{\lfloor N/2\rfloor}) = \lambda_0(z)x_0^2 + \sum_{j=1}^{\lfloor N/2\rfloor}\lambda_j(z)|x_j|^2.\]

The computation is essentially carried out in \cite[Section 3.3]{Capacities2021} in $z<0$. For the convenience of the readers and later applications, we review it here: Consider the superlevel set of the quadratic form $Q_z(q_1,\dots,q_{N}) = Q_z(x_0,x_1, \dots, x_{\lfloor N/2\rfloor})$.
By Morse lemma, it is isomorphic to $\Bbbk [-\operatorname{ind}^{\geq0}(Q_z)]$. Therefore, we only need to count the number of negative eigenvalues of $Q$. It turns out that the diagonal entries $\lambda_k(z) \geq 0$ iff $k \leq \lfloor -z \rfloor$. The final stalk is $\cE_{z}\simeq \Bbbk [-1-2\lfloor -z\rfloor]$.

For $z>0$, the computation is the same as in loc. cit. where $\lambda_k(z) \geq 0$ iff $k \geq -\lfloor -z \rfloor$, so the number of negative eigenvalues will be $N+1+2\lfloor -z\rfloor$ for the subdivision parameter $N$. However, the shifting $N$ is canceled with the shifting $-N$ that appears in the definition of $\cS$ in \autoref{equation: formula of GF sheaf quantization}. So the final stalk is $\cE_{z}\simeq \Bbbk [-1-2\lfloor -z\rfloor]$ and does not depend on $N$.
\end{proof}

When considering the $S^1$-action, one can do a similar weight computation as the $C_\ell$-case in \cite[Section 3.3, after Proposition 3.11]{Capacities2021} \cite[after Proposition 4.12]{zhangthesis}.

\begin{Prop}\label{prop:chiu result on loop sheaf equiv} Recall that by \autoref{thm:geometric circle action tamarkin}, the sheaf $\pi_!\Delta^*\mathcal{S}$ has a $S^1$-action described by a cartesian cyclic object, say $\pi_!\Delta^*\mathcal{S}:\boldsymbol{\Lambda} \to \Sh(\bR_z\times \bR_t)$ as well. There exists a (unique) sheaf with circle action ${\mathcal{E}}: \boldsymbol{\Lambda} \to \Sh(\bR_z)$ such that we have an isomorphism of cartesian cyclic objects $\boldsymbol{\Lambda} \to \Sh(\bR_z\times \bR_t)$:
\begin{equation}\label{loop sheaf of GF model of rotation equi}
\pi_{q!}\Delta^{*} \mathcal{S} \cong \cE \boxtimes {\Bbbk }_{[0,\infty)},   
\end{equation}
where the underlying non-equivariant sheaf in $\Sh(\bR_z)$ is the sheaf $\cE\cong  \bigoplus_{k\in \bZ}1_{(k,k+1]}[1+2k]$ in \autoref{prop:chiu result on loop sheaf}, and extension classes between blocks are described by the $\Bbbk$-coefficient Chern class of the line bundle $\mathcal{O} (-k) \to BS^1=\bC P^\infty$. In particular, over $\bZ$, the extension class between $\bZ_{(k-1,k]}[-1+2k]$ and $\bZ_{(k,k+1]}[1+2k]$ is described by $-ku\in \bZ\fpw{u}$ over $\bZ$.
\end{Prop}
\begin{figure}[htbp]
    \centering
\begin{tikzpicture}[x=1.8cm,y=1cm,>=stealth]
  \draw[->] (-3.6,0) -- (3.6,0) node[right] {$z$};

  \foreach \x in {-3,-2,-1,0,1,2,3}
  {
    \draw (\x,0.08) -- (\x,-0.08);
    \node[below] at (\x,-0.08) {$\x$};
  }

  \node[above=8pt] at (-2.5,0) {$\bZ[-5]$};
  \node[above=8pt] at (-1.5,0) {$\bZ[-3]$};
  \node[above=8pt] at (-0.5,0) {$\bZ[-1]$};
  \node[above=8pt] at (0.5,0) {$\bZ[1]$};
  \node[above=8pt] at (1.5,0) {$\bZ[3]$};
  \node[above=8pt] at (2.5,0) {$\bZ[5]$};

  \node[above=20pt] at (-2,0) {$2u$};
  \node[above=20pt] at (-1,0) {$u$};
  \node[above=20pt] at (0,0) {$0$};
  \node[above=20pt] at (1,0) {$-u$};
  \node[above=20pt] at (2,0) {$-2u$};

  \draw[->] (-2.2,0.6) -- (-1.8,0.6);
  \draw[->] (-1.2,0.6) -- (-0.8,0.6);
  \draw[->] (-0.2,0.6) -- (0.2,0.6);
 \draw[->] (0.8,0.6) -- (1.2,0.6);
  \draw[->] (1.8,0.6) -- (2.2,0.6);

  \node at (-3,0.6) {$\cdots$};
  \node at (3,0.6) {$\cdots$};
\end{tikzpicture}
    \caption{The constructible sheaf $\cE$: the information above the axis describes the stalk on corresponding strata and their extension class.}
    \label{fig: construction of E}
\end{figure}
\begin{proof}
From \autoref{thm:geometric circle action tamarkin}, we know that the circle action on $\pi_!\Delta^*\mathcal{S}$ is described by the following cyclic object, where the simplex $n$ is sent to
\[\pi_{n!}\widetilde{\Delta}_n^*s_{t!}^n(\mathcal{S}^{\boxtimes n}) \simeq \cE_n \boxtimes \Bbbk_{[0,+\infty)},\]
where the sheaf $\cE_n \in \Sh(\bR_z \times \bR_t)$ is given by
\[(\cE_n)_{z} = \Gamma_c\Big(\Big\{(q_1, \dots, q_{nN}):  \sum\nolimits_{j=1}^{nN}S_H(z/nN, q_j, q_{j+1}) \geq 0\Big\}, \Bbbk \Big),\]
at $z<0$ and at $z\geq 0$ by the same formula but shifted by $N$.

We first consider the $z<0$ case. Since the face maps, degeneracy maps and permutation maps are induced by the projection formula for $\bR_z^n \to \bR^{n-1}_z$ and respectively the permutation for $\bR^n_z$, the circle action is induced by a circle action on $\cE: \boldsymbol{\Lambda} \to \Sh(\bR_z \times \bR_t)$. We will define a cyclic anima (homotopy type) $\mathcal{W}^N(z): \boldsymbol{\Lambda} \to \mathrm{An}$ such that for any $\ell \in \boldsymbol{\Lambda}$,
\[\mathcal{W}^{N\ell}_\ell(z) = \Big\{(q_1, \dots, q_{N\ell}): \sum\nolimits_{j=1}^{N\ell}S_H(z/N\ell, q_j, q_{j+1}) \geq 0\Big\}.\]
Diagonalize the form $Q_z(q_1, \dots, q_{N\ell}) = \sum_{j=1}^{N\ell}S_H(z/N\ell, q_j, q_{j+1})$ by the coordinate change $x_j = \sum_{i \in C_{N\ell}}q_{i+1}\omega^{ij}_{N\ell}$ (where $\omega_{N\ell}$ is the primitive $N\ell$-th root of unity). We can write
\[Q_z(x_0,x_1, \dots, x_{\lfloor N\ell/2\rfloor}) = \lambda_0(z)x_0^2 + \sum_{j=1}^{\lfloor N\ell/2\rfloor}\lambda_j(z)|x_j|^2.\]
First, the $C_\ell$-action on $\mathcal{W}_\ell^{N\ell}(z)$ is generated by $(q_j)_{1\leq j\leq N\ell} \mapsto (q_{j-N})_{1\leq j\leq N\ell}$. After the coordinate change, the $C_\ell$-action is by $\omega_\ell \cdot (x_k)_{0\leq k \leq \lfloor N\ell/2\rfloor} = (\omega_\ell^kx_k)_{0\leq k \leq \lfloor N\ell/2\rfloor}$ for any $\ell \in \bN$ \cite[Equation (3.20)]{Capacities2021} \cite[Equation (4.20)]{zhangthesis}. 
Therefore, we can consider the cyclic space $\mathcal W_\bullet(z)$ such that the $C_\ell$-action is induced by $\omega_\ell \cdot (x_k)_{0\leq k\leq \lfloor -z\rfloor} = (\omega_\ell^k x_k)_{0\leq k\leq \lfloor -z\rfloor}$. Consider now the space $W(z)$ with the circle action by 
\[t \cdot (x_k)_{0\leq k \leq \lfloor N/2\rfloor} = (t^kx_k)_{0\leq k \leq \lfloor N/2\rfloor}. \]
Note that we have an equivalence of anima induced by an $S^1$-map (and thus an equivalence after taking the fiber product $- \times_G EG$):
\[\mathcal{W}_{\ell}^{N\ell}(z) \simeq W(z) \coloneqq \{(x_0, x_1,\dots, x_{\lfloor -z\rfloor}, 0, \dots, 0): (x_0, x_1,\dots, x_{\lfloor -z\rfloor}) \in \bR \times \bC \times \dots \times \bC\}.\]
We show that the induced circle agrees with the cyclic anima $\mathcal{W}^N_\bullet(z) \simeq W_\bullet (z)$ above.
The compactly supported cohomology is isomorphic to $(\cE_n)_{z}$ as a cyclic object. Now, it suffices to compute the equivariant cohomology of the circle action on $\mathcal{W}^N(z) \simeq W(z)$.

We now review the argument in \cite[after Remark 4.13]{zhangthesis} that determines the circle action: For $z<0$, $S^1$ acts on $\bR^{1+2\lfloor |z| \rfloor}$ such that it acts on $\bR$ trivially, and acts on $\bC^{\lfloor |z| \rfloor} \cong \bR^{2\lfloor |z|\rfloor}$ with weight $(1, 2, \dots, \lfloor |z| \rfloor)$. The morphism of $S^1$-equivariant complexes for $z'<z<0$
\[\Gamma_c(\bR^{1+2\lfloor |z'|\rfloor}; \Bbbk ) \to \Gamma_c(\bR^{1+2\lfloor |z|\rfloor}; \Bbbk )\]
is induced by the top Chern class on the vector bundle $\bR^{1+2\lfloor |z'|\rfloor} \times_{S^1} S^\infty \to \bR^{1+2\lfloor |z|\rfloor} \times_{S^1} S^\infty$, by $(\lfloor |z'|\rfloor !/\lfloor |z|\rfloor!)u^{\lfloor |z'|\rfloor -\lfloor |z|\rfloor}$, as in the figure over $\bZ$.

In the $z>0$ case, the computation is basically the same, but one needs to be careful about the dimension at $z$, which should be $N+1+2\lfloor -z\rfloor$, as explained in \autoref{prop:chiu result on loop sheaf}, and the weight of the corresponding representations, which should be $(-\lfloor -z\rfloor,   \dots, (N-1)/2)$ by a similar computation.
\end{proof}

\subsection{Microlocal projector for toric domains}
Recall that the Chiu--Tamarkin complex is defined as $F(X; \Bbbk ) = \pi_!\Delta^*(s_{t!}P_X)$. In this subsection, we recall the computation of the microlocal projector $P_X$.
Let
\begin{equation}\label{sheafquantization of Hamiltonian torus action}
 \mathcal{T}\coloneqq \mathcal{S}^{\boxstar n}=s_{t!}^n(\mathcal{S}^{\boxtimes n}) \in \Sh(\bR_{{z}}^n\times V_{1}\times V_{2};\sT ),   
\end{equation}
where $ V_{i}=\bR^n_{\bq_i}$.
The microsupport estimates show that $\mathcal{T}$ is a sheaf quantization of the standard torus action on $T^*V = T^*\bR_{\bq}$. 

Recall that the Fourier--Sato--Tamarkin transform is defined by \cite[Definition 1.13]{Capacities2021}
\[\widehat{\mathcal{T}} := \mathcal{T} \star_{z} \Bbbk _{\{(z,\zeta,t):t+\left<z,\zeta\right>\geq 0\}} [n]= (\pi_z \times s_{t})_!(\mathcal{T} \boxtimes_z \Bbbk _{\{(z,\zeta,t):t+\left<z,\zeta\right>\geq 0\}})[n] \in \Sh(\bR^n_\zeta \times V_1 \times V_2; \sT).\]
Then we have the following computation; see \cite[Propositions 2.8, 3.3 \& Equation (3.7)]{Capacities2021}. 
\begin{Prop}Let $\widetilde{\Omega}$ is an open set of $\bR^n$ such that $\Omega=\widetilde{\Omega}\cap \bR^n_{\geq 0}$. Then we have
\begin{equation}
    P_{X_{\Omega}}\coloneqq \widehat{\mathcal{T}}\circ {\Bbbk }_{\widetilde{\Omega}} ,\qquad Q_{X_{\Omega}}\coloneqq \widehat{\mathcal{T}}\circ {\Bbbk }_{\bR^n_{{\zeta}}\setminus\widetilde{\Omega}},
\end{equation}
\end{Prop}
To make sense of the notation, we shall show that $(P_{X_{\Omega}}, Q_{X_{\Omega}}) $  does not depend on $\widetilde{\Omega}$.

\begin{Prop}\label{prop: extension independence in toric kernel}
We have $\supp(\widehat{\mathcal{T}}) \subset \bR^n_{\zeta\geq 0}$, and then, for $\Omega=\widetilde{\Omega}\cap \bR^n_{\geq 0}$, 
\[P_{X_{\Omega}}\cong \widehat{\mathcal{T}}\circ {\Bbbk }_{{\Omega}} ,\qquad Q_{X_{\Omega}}\coloneqq \widehat{\mathcal{T}}\circ {\Bbbk }_{\bR^n_{\geq 0}\setminus{\Omega}}.\]  
\end{Prop}
\begin{proof}We will prove that $\supp(\widehat{\mathcal{T}}) \subset \bR^n_{\zeta\geq 0}$ in 2 methods. The first uses the microsupport estimation of the Fourier transform, it appears in \cite[After Proposition 3.3]{Capacities2021} already, which we review here; and the second one uses a direct cohomology computation based on microlocal cut-off, which is new and helpful for our later discussion.

(1) Recall the microsupport estimate of $\widehat{\mathcal{T}}$ \cite[Equation (2.4)]{Capacities2021}, that if $(\zeta,z,\bq,\bp,\bq',\bp',t,\tau)\in \msnz(\widehat{\mathcal{T}})$, then we have $\zeta=\mu(\bq,\bp)\in \bR^n_{\zeta\geq 0}$. So, if $\zeta\notin \bR^n_{\zeta\geq 0}$ and $(\zeta,z,\bq,\bp,\bq',\bp',t,\tau)\in SS(\widehat{\mathcal{T}})$, we have $(\bp,\bp',\tau)=0$. Accordingly, for any $\zeta\notin \bR^n_{\zeta\geq 0}$, we have $\ms(\widehat{\mathcal{T}}|_{(\zeta,\bq,\bq')})\subset \{\tau= 0\}$. So $\widehat{\mathcal{T}}|_{(\zeta,\bq,\bq')}\cong M_\bR$ is a constant sheaf over $\bR$ for some $M\in D(\Bbbk )$.
As $\widehat{\mathcal{T}}|_{(\zeta,\bq,\bq')}\in \sT(\pt)$, we have $\widehat{\mathcal{T}}|_{(\zeta,\bq,\bq')}\cong M_\bR\cong M_\bR\star \Bbbk _{[0,\infty)} \cong 0 $ for any $\zeta\notin \bR^n_{\zeta\geq 0}$.

(2) Let $U_i=\{\zeta\in \bR_\zeta^n\mid\zeta_i<0\}$ be the open half-space of $\zeta_i$ coordinate. Notice that $ \bR_\zeta^n\setminus \bR_{\geq 0}^n \times V_1 \times V_2\times \bR_t=\cup_{i=1}^n U_i \times V_1 \times V_2\times \bR_t$, we only need to show that $\widehat{\mathcal{T}}|_{U_i\times V_1\times V_2\times  \bR_t}=0$ for all $i=1,\dots,n$.

By $\ms(\mathcal{T})\subset \{\zeta_i\leq 0, \forall i=1,\dots, n\}$ and microlocal cut-off lemma \cite[Proposition 5.2.1]{KS90}, for $K_i=\Bbbk_{\{z_i\geq 0\}}\in \Sh(\bR_z^n)$, we have $\mathcal{T}\star_z K_i=0$ for all $i=1,\dots, n$. Therefore, we have
\[0=\widehat{\mathcal{T}\star_z K_i}|_{U_i\times V_1\times V_2\times  \bR_t}\simeq \widehat{\mathcal{T}}\star (H_i|_{U_i\times \bR_t}),\qquad H_i=\pi_{r!}\Bbbk_{\{(r,\zeta,t)\in \bR_r\times \bR_\zeta\times \bR_t\mid r\geq 0,\, t+r\zeta_i\geq 0 \}}.\]
By the formula for $H_i$, we have 
\[H_i|_{U_i\times \bR_t}=\Bbbk_{U_i\times [0,\infty)},\]
and then
\[0=\widehat{\mathcal{T}\star_z K_i}|_{U_i\times V_1\times V_2\times  \bR_t}\simeq \widehat{\mathcal{T}}\star (H_i|_{U_i\times \bR_t})\simeq \widehat{\mathcal{T}}|_{U_i\times V_1 \times V_2 \times \bR_t},\]
which concludes.
\end{proof}

\begin{Prop}\label{prop: cyclic trace for toric domain-first formula}
For $\Omega\subset \bR^n_{\geq 0}$ open, we have
\[F(X_\Omega) \simeq \pi_{z,\zeta!}(\pi_{\zeta,t}^*\cE^{\boxtimes n} \otimes \Bbbk _{\{(z,\zeta,t): \zeta \in \Omega, t+\left<z,\zeta \right> \geq 0\}})[n].\]
Moreover, this can be enhanced to an $S^1$-equivariant isomorphism of cartesian cyclic objects in $\sT$ such that
\[(F^{S^1}(X_\Omega))_\ell\simeq F_\ell(X_\Omega) \simeq \pi_{z,\zeta!}(\pi_{\zeta,t}^*\widetilde{\Delta}_\ell^*\cE^{\boxtimes n\ell} \otimes \Bbbk _{\{(z,\zeta,t): \zeta \in \Omega, t+\left<z,\zeta \right> \geq 0\}})[n].\]
\end{Prop}
\begin{proof}
Consider the following diagram of maps:
\[\begin{tikzcd}[sep=scriptsize]
V \times \bR_t \ar[r, "\Delta_q"]& V_1 \times V_2 \times \bR_t & \bR^n_z \times V_1 \times V_2 \times \bR_t^2 \ar[l,"\pi_z\times s_t"] \\
& \ar[u,"\pi_\zeta"] \bR^n_\zeta \times V_1 \times V_2 \times \bR_t \ar[d,"\pi_{q,q'}"] & \bR^n_z \times \bR^n_\zeta \times V_1 \times V_2 \times \bR^2_t \ar[l,"\pi_z \times s_t"] \ar[d,"\pi_{q,q'}"] \ar[u,"\pi_\zeta"]\\
& \bR^n_\zeta & \bR^n_z \times \bR^n_\zeta \times \bR_t \ar[l,"\pi_z \times s_t"].
\end{tikzcd}\]
By the projection formula, we can compute that
\begin{align*}
F(X_\Omega) &\simeq \pi_{q!}\Delta_q^*\pi_{\zeta!}((\pi_z \times s_t)_!(\pi_\zeta^*\mathcal{T} \otimes \pi_{q,q'}^*\Bbbk _{\{(z,\zeta,t): t+\left<z,\zeta \right> \geq 0\}}) \otimes \pi_{q,q'}^*\Bbbk _{\Omega})[n] \\
& \simeq \pi_{q!}\Delta_q^*\pi_{\zeta!}((\pi_z \times s_t)_!(\pi_\zeta^*s_{t!}^n(\mathcal{S}^{\boxtimes n}) \otimes \pi_{q,q'}^*\Bbbk _{\{(z,\zeta,t): t+\left<z,\zeta \right> \geq 0\}}) \otimes \pi_{q,q'}^*\Bbbk _{\Omega})[n] \\
&\simeq \pi_{q!}\Delta_q^*\pi_{\zeta!}(\pi_{z} \times s_{t})_!(\pi_\zeta^*s_{t!}^n(\mathcal{S}^{\boxtimes n}) \otimes \pi_{q,q'}^*\Bbbk _{\{(z,\zeta,t): \zeta \in \Omega, t+\left<z,\zeta \right> \geq 0\}})[n].
\end{align*}
Consider now the diagrams of maps
\[\begin{tikzcd}[sep=scriptsize]
V \times \bR_t \ar[d,"\Delta"] & \bR^n_z \times \bR^n_\zeta \times V \times \bR_t \ar[d,"\Delta"] \ar[l,"\pi_{z,\zeta}"] & \\
V_1 \times V_2 \times \bR^2_t & \bR^n_z \times \bR^n_\zeta \times V_1 \times V_2 \times \bR_t \ar[d,"\pi_z"] \ar[l,"\pi_{z,\zeta}"] & \bR^n_z \times \bR^n_\zeta \times V \times \bR_t^2 \ar[ul,"s_t"] \ar[d,"\Delta"] \\
& \ar[ul,"\pi_\zeta"] \bR^n_\zeta \times V_1 \times V_2 \times \bR_t  & \bR^n_z \times \bR^n_\zeta \times V_1 \times V_2 \times \bR^2_t \ar[l,"\pi_z \times s_t"] \ar[ul,"s_t"]
\end{tikzcd}
\]
and also consider the following diagram of maps:
\[
\begin{tikzcd}[sep=scriptsize]
\bR^n_z \times \bR^n_\zeta \times \bR_t & \ar[l,"\pi_q"]\bR^n_z \times \bR^n_\zeta \times V \times \bR_t^2 \ar[d,"\Delta"] \ar[r,"\pi_{\zeta,t}"] & \bR^n_z \times V \times \bR_t \ar[d,"\Delta"] &  (\bR_z \times \bR_q \times \bR_t)^n \ar[l,"s_t^n"] \ar[d,"\Delta^n"]\\
& \ar[ul,"\pi_{q,q'}"]\bR^n_z \times \bR^n_\zeta \times V_1 \times V_2 \times \bR^2_t  \ar[r,"\pi_{\zeta,t}"] & \bR^n_z \times V_1 \times V_2 \times \bR_t & (\bR_z \times \bR_q \times \bR_{q'} \times \bR_t)^n \ar[l,"s_t^n"] 
\end{tikzcd}
\]
Then, by proper base change formula, and the observation that pull-backs commute with tensor products, we can use proper base change formula again and find that
\begin{align*}
F(X_\Omega) &\simeq \pi_{q!}\pi_{z,\zeta!}s_{t!}\Delta_q^*(\pi_\zeta^*s^n_{t!}(\mathcal{S}^{\boxtimes n}) \otimes \pi_{q,q'}^*\Bbbk _{\{(z,\zeta,t): \zeta \in \Omega, t+\left<z,\zeta \right> \geq 0\}})[n]\\
&\simeq \pi_{q!}\pi_{z,\zeta!}s_{t!}(s_{t!}^n(\Delta_q^*\mathcal{S})^{\boxtimes n} \otimes \pi_q^*\Bbbk _{\{(z,\zeta,t): \zeta \in \Omega, t+\left<z,\zeta \right> \geq 0\}})[n]
\end{align*}
Finally, by the projection formula again and \autoref{prop:chiu result on loop sheaf}, we can conclude that
\begin{align*}
F(X_\Omega) &\simeq \pi_{z,\zeta!}s_{t!}\pi_{q!}(s_{t!}^n(\Delta_q^*\mathcal{S})^{\boxtimes n} \otimes \pi_q^*\Bbbk _{\{(z,\zeta,t): \zeta \in \Omega, t+\left<z,\zeta \right> \geq 0\}})[n]\\
&\simeq \pi_{z,\zeta!}s_{t!}(s_{t!}^n(\pi_{q!}\Delta_q^*\cS)^{\boxtimes n} \otimes \Bbbk _{\{(z,\zeta,t): \zeta \in \Omega, t+\left<z,\zeta \right> \geq 0\}})[n]\\
&\simeq \pi_{z,\zeta!}s_{t!}(s_{t!}^n(\cE \boxtimes 1_{[0,\infty)})^{\boxtimes n} \otimes \Bbbk _{\{(z,\zeta,t): \zeta \in \Omega, t+\left<z,\zeta \right> \geq 0\}})[n].
\end{align*}
Since $1_{[0,\infty)}$ is the monoidal unit under the Tamarkin convolution, $s_{t!}^n(1_{[0,\infty)}^{\boxtimes n}) = 1_{[0,\infty)}$, and we conclude that
\[F(X_\Omega)\simeq \pi_{z,\zeta!}(\pi_{\zeta,t}^*\cE^{\boxtimes n} \otimes \Bbbk _{\{(z,\zeta,t): \zeta \in \Omega, t+\left<z,\zeta \right> \geq 0\}})[n].\]
\vspace{1pt}

Finally, for the statement on $S^1$-equivariant isomorphism, all the proofs remain the same. In fact, we only need to notice that in the second step above, the proper base change formula now gives
\begin{align*}
F^{S^1}(X_\Omega)_\ell &\simeq \pi_{\ell!}\pi_{z,\zeta!}s_{t!}\widetilde\Delta_\ell^*(\pi_\zeta^*s^n_{t!}(\mathcal{S}^{\boxtimes n})^{\boxtimes \ell} \otimes \pi_{q,q'}^*\Bbbk _{\{(z,\zeta,t): \zeta \in \Omega, t+\left<z,\zeta \right> \geq 0\}}^{\boxtimes \ell})[n]\\
&\simeq \pi_{\ell!}\pi_{z,\zeta!}s_{t!}(s_{t!}^n(\widetilde\Delta_\ell^*\mathcal{S})^{\boxtimes n\ell} \otimes \pi_q^*\Bbbk _{\{(z,\zeta,t): \zeta \in \Omega, t+\left<z,\zeta \right> \geq 0\}})[n]
\end{align*}
The rest of the computation is the same. The coherence between the isomorphisms with the morphisms in the cyclic category follows from six functor formalism.
\end{proof}

\subsection{An exodromy equivalence}\label{section: lattice sheaf general view}
Based on \autoref{prop: cyclic trace for toric domain-first formula}, we shall study $\cE^{\boxtimes n} $, which is a sheaf that is constructible with respect to the lattice stratification of $\bR^n$ and has (fully) negative microsupport. In this section, we discuss a general setting for those kinds of sheaves and apply it to $\cE^{\boxtimes n} $ later. In fact, all those results should be well known, though we did not find nice enough references. 

To strengthen the full power of the construction, we shall use the language of derived coend. The language is a fancy way to express the (slight generalization of) Koszul complex, which we shall explain a bit later. 

\subsubsection{Exodromy equivalence}
To begin with, we set some notation for high-dimensional intervals. Let $x,y\in (\bR\cup\{\infty\})^n$ with $x\leq y$ under the product order of $\bR\cup\{\infty\}$. We set 
\[ (x,y)\coloneqq \prod_{i=1}^n(x_i,y_i)\]
as well as its closed/open-closed versions. 

We consider the following Lagrangian sets.
\[L_1\coloneqq  0_\bR\cup \bigcup_{k\in\bZ} \{k\}\times (-\infty,0] \subset T^*\bR,\qquad L_n\coloneqq (L_1)^n\subset T^*\bR^n\]

\begin{figure}[htbp]
    \centering
\begin{tikzpicture}[x=1.15cm,y=1cm,>=stealth]
  \draw[->,gray] (-3.9,0) -- (4.05,0) node[right,black] {$q$};
  \draw[->,gray] (0,-3.1) -- (0,0.85) node[above,black] {$p$};
  \draw[blue!70!black,very thick] (-3.65,0) -- (3.65,0);
  \foreach \k in {-3,...,3} {
    \draw[blue!70!black,very thick] (\k,0) -- (\k,-2.7);
    \fill[blue!70!black] (\k,0) circle (1.6pt);
    \node[above=5pt] at (\k,0) {$\k$};
  }
\end{tikzpicture}
    \caption{The picture of $L_1\subset T^*\bR$}
    \label{fig: picture of lattice strata in dim=1}
\end{figure}
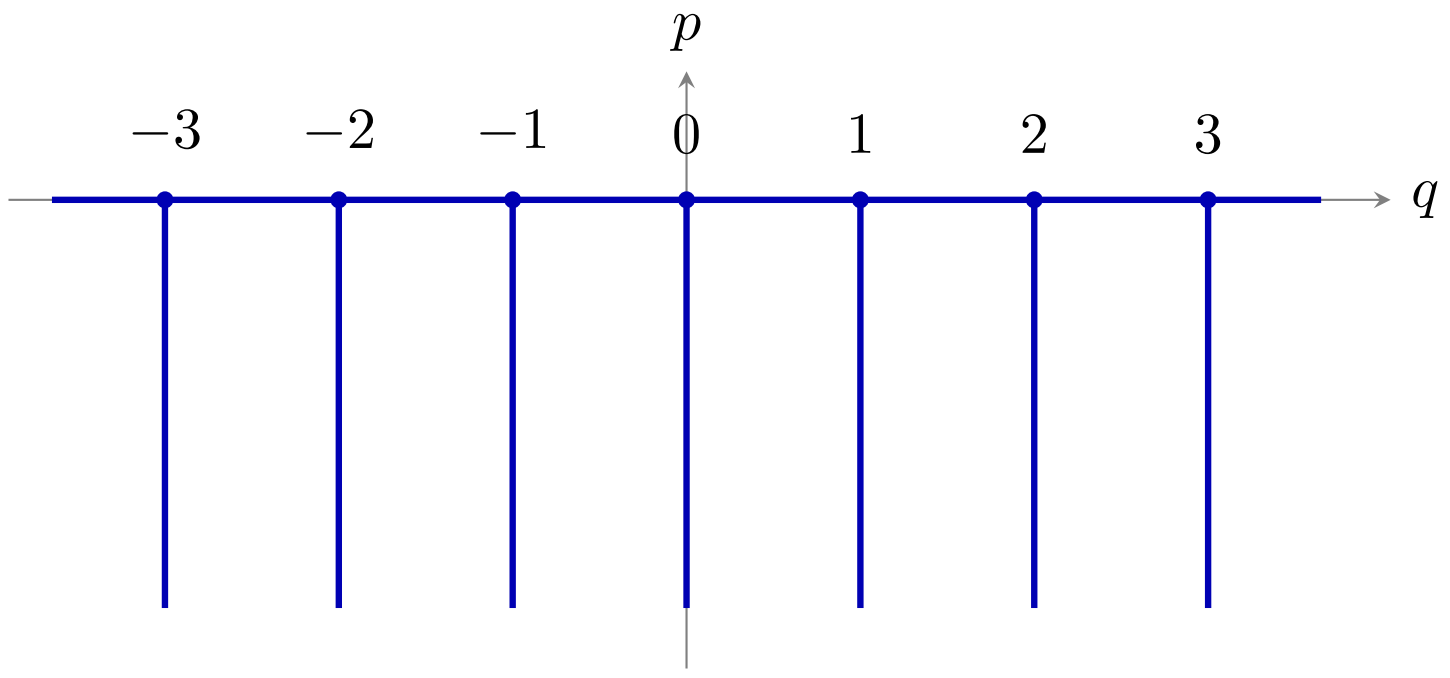

Let $\cC$ be a presentable stable category. The category $\Sh_{L_n}(\bR^n,\cC)$ is a full subcategory of $\mathcal{S}$-constructible sheaves on $\bR^n$ for the integer-lattice stratification of $\bR^n$:
\[
\mathcal{S}
=
\left\{
\operatorname{RelInt}([k,k+\varepsilon])
\;\middle|\;
k\in\bZ^n,\quad \varepsilon\in\{0,1\}^n
\right\},
\]
where when $\varepsilon_i=0$ (resp. $\varepsilon_i=1$), it means that the $i$-th factor is $\{k_i\}$ (resp. $(k_i,k_i+1)$).

Moreover, negative microsupport assigns the propagating direction positively. We have the exodromy equivalence 
\[ \Sh_{L_n}(\bR^n,\cC)\simeq \Fun(\rarbZ^n,\cC),\]
where $\rarbZ^n$ is the category associated with the poset $\bZ^n$ with the product order of $\bZ$. Notice that $\Fun(\rarbZ^n,D(\Bbbk))\simeq \operatorname{IndCoh}([\mathbb{A}^n/\mathbb{G}_m^n])$ (follows from the discussion in \cite[Tag 0EKL]{StacksProject}); then it follows from toric mirror symmetry, for example \cite{ToricMSrevisted-Shende,SunFunctorialityToricCCC} for ring coefficients.) On the other hand, the equivalence also follows from the general exodromy equivalence applied to the Bondal--Ruan stratification \cite[Example 5.3.10]{HainePortaTeyssier-conicality}. Here, we shall describe it over a general coefficient $\cC$.

On one side, the exodromy functor can be given by
\[\Ex_\cC: \Sh_{L_n}(\bR^n,\cC) \to \Fun(\rarbZ^n,\cC)\]
is given by
\[\Ex_\cC(F)(m)=\Gamma((m,\infty),F),\qquad  \Ex_\cC(F)(l\leq m)=[\Gamma((l,\infty),F)\to \Gamma((m,\infty),F)].\]

\medskip

In the following, we write a quasi-inverse of $\Ex_\cC$. 

For any $M\in \Fun(\rarbZ^n,\cC)$, we define a bifunctor
\[B_M: (\rarbZ^n)^\op\times \rarbZ^n \to \Sh_{L_n}(\bR^n,\cC),\quad (l,m)\mapsto  j_{l!}(M(m))\]
where $j_l$ is the open inclusion $(l,\infty)\to \bR^n $ and we simply take $j_{l!}X$ the constant sheaf supported on $(l,\infty)$ with value $X$. Here, we notice that $B_M$ takes values in $\Sh_{L_n}(\bR^n,\cC)$ by the usual microsupport estimation \cite[Proposition 5.4.1 \& 5.4.8]{KS90}.
Then, we have an object defined through the coend of the above bifunctor \cite[Section 5]{GrothPontoShulman}
\[\Real_\cC(M)\coloneqq \int^{\rarbZ^n}B_M = \int^{\rarbZ^n} j_{m!}(M(m))\in \Sh_{L_n}(\bR^n,\cC).\]
Let $\mathrm{TwAr}(\rarbZ^n)$ be the twisted arrow category where objects are $m \leq m'$ (and morphisms are $m\leq m' \mapsto n \leq n'$ if $n \leq m \leq m'\leq n'$). Then the coend can be computed by the formula
\[\int^{\rarbZ^n}B_M = \colim_{m\leq l \in \mathrm{TwAr}(\rarbZ^n)}j_{l!}(M(m)).\]
The claim that the colimit is contained in $\Sh_{L_n}(\bR^n)$ comes from the cocompleteness of $\Sh_{L_n}(\bR^n)$.

It is clear that $\Real$ defines a functor
\[\Real_\cC: \Fun(\rarbZ^n,\cC) \to \Sh_{L_n}(\bR^n,\cC).\]
In fact, by the universal property of coend and the Yoneda lemma, we have
\begin{align*}
\Map(M,\Ex_\cC(F))&=\int_{\rarbZ^n}\Map(M(m),\Gamma((m,\infty),F))\\
   &=\int_{\rarbZ^n}\Map(j_{m!}M(m),F) =\Map(\Real_\cC(M),F),
\end{align*}
which supplies the functoriality of $\Real_\cC$ and shows that $\Real_\cC$ is the left adjoint of $\Ex_\cC$. 

\begin{Prop}The functor $\Real_\cC$ is a/the quasi-inverse of $\Ex_\cC$.
\end{Prop}
\begin{proof}We first assume that $\cC=\Sp$. In this case, we have $j_{m!}(M(m)))=\bS_{(m,\infty)}\otimes M(m)$ for the sphere spectrum $\bS$.

Let $h_m$ be the Yoneda functor on $ \rarbZ^n$ of $m\in \rarbZ^n$. On the one hand, we have 
\begin{align*}
    \Ex_\Sp(\Real_\Sp(M)) (m)&=\Gamma\Big((m,\infty),\int^{\rarbZ^n} \bS_{(m,\infty)}\otimes M(m)\Big)    =\Hom_{\Sh_{L_n}}\Big(\bS_{(m,\infty)},\int^{\rarbZ^n} \bS_{(m,\infty)}\otimes M(m)\Big)\\
   &=\int^{\rarbZ^n} \Hom_{\Sh_{L_n}}(\bS_{(m,\infty)}, \bS_{(m,\infty)}\otimes M(m)) =\int^{\rarbZ^n} h_m\otimes M(m) 
     =M(m).
\end{align*}

The third equality is true because $\bS_{(m,\infty)}$ is a compact object in $\Sh_{L_n}(\bR^n,\Sp)$, and the last equality is the co-Yoneda lemma.

On the other hand, the evaluation functor induces
\[j_{m!}(\Gamma((m,\infty),F))\to F,\]
and then we have a natural map
\[\Real_\Sp(\Ex_\Sp(F)) = \int^{\rarbZ^n} j_{m!}(\Gamma((m,\infty),F)) \to F.
\]
Thus, we conclude equivalence by the stalk-wise test by hypercompleteness of $\bR^n$ and $\cC$: For any $x\in \bR^n$, we have
\[  \left(\int^{\rarbZ^n} j_{m!}\Gamma((m,\infty),F))\right)_x\simeq  \int^{\rarbZ^n}\left( j_{m!}\Gamma((m,\infty),F)))\right)_x \to F_x.\]
But the left-hand side is computed by
\[\varinjlim_{m< x}\Gamma((m,\infty),F)\simeq F_{\lceil x\rceil}.\]
However, by $\ms(F)\subset L_n$, the object $F_{\lceil x\rceil}$ is isomorphic to $F_x$ by the microlocal Morse lemma \cite[Proposition 5.4.19]{KS90}. We conclude the spectra-valued case.

In general, we have
\[\Fun(\rarbZ^n,\cC)\simeq\Fun(\rarbZ^n,\Sp)\otimes\cC\simeq\Sh_{L_n}(\bR^n,\Sp)\otimes\cC\simeq \Sh_{L_n}(\bR^n,\cC),\]
where the first equivalence follows from \cite[Lemma 4.8.1.16, Proposition 4.8.1.17]{HA}, and the third equivalence follows from \cite[Proposition 2.7]{Zhang_K_Remark}.

Here, the middle equivalence is $ \Real_\Sp\otimes \,\cC$ whose quasi-inverse is $\Ex_\Sp\otimes \,\cC$, which is then the right adjoint of $ \Real_\Sp\otimes \,\cC$. But one can check that $\Real_\cC\simeq \Real_\Sp\otimes \,\cC$. Then, combining that $\Real_\cC$ is the left adjoint of $\Ex_\cC$, one has that $\Ex_\cC\simeq \Ex_\Sp\otimes\, \cC$ by uniqueness of adjunction.\end{proof}

\subsubsection{Coend of cubes} Here, we explain how to compute the $\rarbZ^n$-coend explicitly, and we will see its relation to the Koszul complex. We need some notation and results about $n$-cubes in a category, which we refer to in \cite{BeckertGroth_Cubical_homotopy}.

\begin{Prop}\label{prop: coend of cube as total cofiber}Let $\cC$ be a cocomplete stable category, and $B:(\rarbZ^n)^\op\times \rarbZ^n\to \cC$ be a bifunctor. Let $e_i$ be the unit vector in $\bZ^n$ and for $I\subset\{1,\ldots,n\} $, write $e_I\coloneqq\sum_{i\in I}e_i$.

For $I=\{i_1<\cdots<i_p\}$ and $i=i_r\in I$, define
\[
\partial_i^0=B(m+e_{I\setminus\{i\}}\le m+e_I,m):
B(m+e_I,m)
\longrightarrow
B(m+e_{I\setminus\{i\}},m)
\]
and 
\[
\partial_i^1=B(m+e_I,m \le m+e_{i}):
B(m+e_I,m)
\longrightarrow
B(m+e_I,m+e_i).
\]

We then define 
\[C_B(I)\coloneqq  
\bigoplus_{v\in\bZ^n}B(v+e_{[n]\setminus I},v).\]
When $i\notin I$, we denote
\[\delta_i=\partial^0_i-\partial^1_i: C_B(I) \to C_{B}(I\cup \{i\}).\]

Then 
\begin{enumerate}
    \item The assignment $I\subset \{1,\ldots,n\} \mapsto C_I(B)$ and $\delta_i$ induces a functor
    \[C_{B}:\mathcal{P}([n]) \to \cC,\]
i.e. an $n$-cube of \cite{BeckertGroth_Cubical_homotopy}.
\item 
Then there is a natural equivalence
\[
\int^{m\in\rarbZ^n}B(m,m)
\simeq
\operatorname{tcof}(C_B)=\cof\left(
\varinjlim_{I\subsetneq[n]} C_B(I)
\to C_B([n]) \right),
\]
where the right-hand side is the total cofiber (its definition is presented in the following proof).
\end{enumerate}
\end{Prop}

\begin{proof}
The first statement is clear since $B$ is a bifunctor. We prove the second by induction on $n$.

We start from $n=1$. Since $\rarbZ$ is the free category generated by the elementary arrows 
$
m\leq m+1$, $
 m\in\bZ, 
$
then the coend is computed by the following cofiber
\[
\int^{m\in\rarbZ}B(m,m)
\simeq
\cof\left(
\bigoplus_{m\in\bZ}B(m+1,m)
\xrightarrow{\delta_1}
\bigoplus_{m\in\bZ}B(m,m)
\right).
\]
By \cite[Remark 8.5(i)]{BeckertGroth_Cubical_homotopy}, the cofiber of a $1$-cube is
canonically its total cofiber. 

For $n\geq 2$, we use the Fubini theorem for coends:
\[
\int^{m\in\rarbZ^n}B(m,m)
\simeq
\int^{m_1\in\rarbZ}
\cdots
\int^{m_n\in\rarbZ}
B(m,m).
\]

Then we can apply the one-dimensional formula successively in the $n$
coordinate directions. Thus, taking the coend in the $i$-th coordinate is the same as applying the cone construction of $\delta_i$ in the $i$-th direction of the cube $C_B$. And then the coend over the entire $\rarbZ^n$ is an iterated cone. Therefore, \cite[Theorem 8.25]{BeckertGroth_Cubical_homotopy} concludes the proof.
\end{proof}
 
We are particularly interested in a linear algebra situation, which gives a reformulation of the coend formula in terms of the Koszul complex.
\begin{Coro}\label{coro: coend computing in chain}
    In the situation of \autoref{prop: coend of cube as total cofiber}, we assume further that $\cC=D(A)$ is the derived category of a dg-$\Bbbk$-algebra $A$. 
    
    We set
\[
C_B^p\coloneqq \bigoplus_{|I|=p} C_B(I)=
\bigoplus_{{m\in\bZ^n, |I|=p}}
B(m+e_{[n]\setminus I},m).
\]

For $I\subset[n]$ with $|I|=p$, we write $ [n]\setminus I=\{i_1<\cdots<i_{n-p}\}, $ and then define 
\[ \partial_I \coloneqq \sum_{r=1}^{n-p} (-1)^{r-1}\delta_{i_r}: C_B(I)\longrightarrow C_B^{p+1}, \]
and then define $\partial: C_B^{p}\to C_B^{p+1}$ with component $\partial_I$.

Then we have
\begin{enumerate}[fullwidth] \item We have $ \partial^2=0$. Hence \[ C_B^0\xrightarrow{\partial} C_B^1\xrightarrow{\partial}\cdots \xrightarrow{\partial}C_B^{n-1} \xrightarrow{\partial}C_B^n \] is a complex in $D(A)$.  \item There is a natural equivalence\footnote{The degree shifting comes since we put $C^0_B$ as the $0$-th term.} \[ \int^{m\in\rarbZ^n}B(m,m) \simeq \operatorname{Tot} \left( C_B^0\xrightarrow{\partial} C_B^1\xrightarrow{\partial}\cdots \xrightarrow{\partial}C_B^{n-1} \xrightarrow{\partial}C_B^n  \right)[n]=\left(\bigoplus_p C_B^p[n-p],d=(-1)^{n-p}d_{C^p_B}+\partial\right). \] 

\end{enumerate}
\end{Coro}
\begin{proof}
The total cofiber is computed by iterated cones in $n$ directions, which can be computed by the given totalization.
\end{proof}

\subsection{Cyclic trace of toric domains}\label{section: trace for toric domain}
In this section, we describe a coend formula  for cyclic trace based on the lattice constructible feature of $\cE^{\boxtimes n}$ for the sheaf $\cE$ in \autoref{prop:chiu result on loop sheaf equiv} using results of \autoref{section: lattice sheaf general view}.

\begin{figure}[htbp]
    \centering
    \begin{tikzpicture}[
        x=1.65cm,
        y=1.4cm,
        >=stealth,
        stalk/.style={font=\small},
        maplabel/.style={
            font=\scriptsize,
            fill=white,
            inner sep=1pt
        }
    ]
        \draw[step=1, gray!35, thin] (-3,-3) grid (3,3);
        \draw[->] (-3.4,0) -- (3.5,0)
            node[right] {$z_1$};
        \draw[->] (0,-3.4) -- (0,3.5)
            node[above] {$z_2$};

        \fill (0,0) circle (1.2pt)
            node[below left, fill=white, inner sep=1pt] {$O$};

        \foreach \i in {-3,...,2} {
            \foreach \j in {-3,...,2} {
                \pgfmathtruncatemacro{\stalkshift}{2*(\i+\j+1)}
                \node[stalk] at ({\i+0.5},{\j+0.5})
                    {$\mathbb Z[\stalkshift]$};
            }
        }

        \foreach \k/\c in {-2/2,-1/{},0/0,1/{-},2/{-2}} {
            \foreach \j in {-3,...,2} {
                \draw[->]
                    ({\k-0.18},{\j+0.5})
                    -- ({\k+0.18},{\j+0.5})
                    node[midway, above=3pt, maplabel]
                    {$\ifnum\k=0 0\else \c u_1\fi$};
            }

            \foreach \i in {-3,...,2} {
                \draw[->]
                    ({\i+0.5},{\k-0.18})
                    -- ({\i+0.5},{\k+0.18})
                    node[midway, right=3pt, maplabel]
                    {$\ifnum\k=0 0\else \c u_2\fi$};
            }
        }
    \end{tikzpicture}
    \caption{The sheaf $\cE^{\boxtimes 2}$.
    Labels inside the squares indicate the stalks;
    arrows indicate the extension classes across the integer walls.}
    \label{fig:E-external-product}
\end{figure}
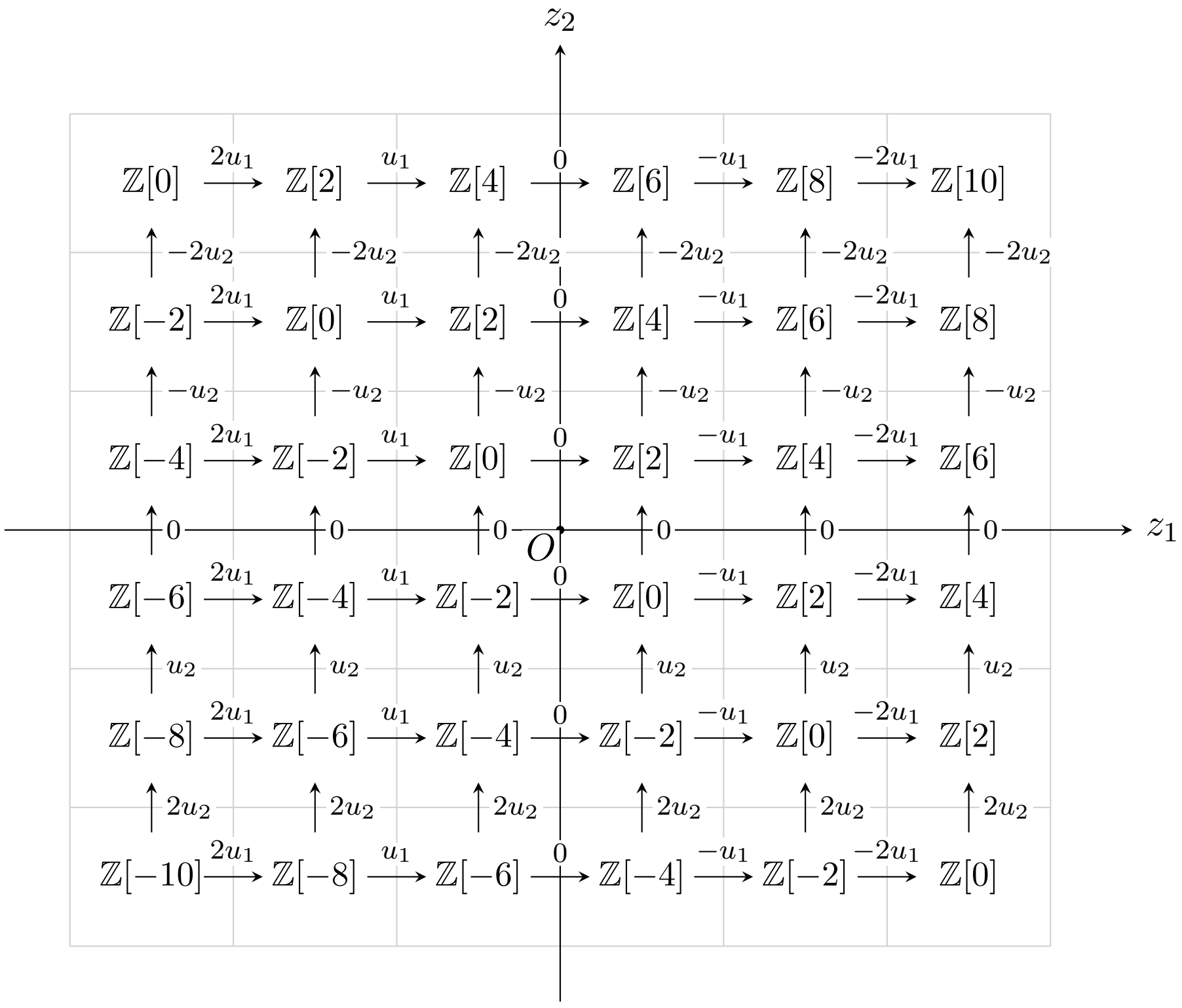

\begin{Lemma}\label{lemma: coend formula for E^n}We have that $\cE^{\boxtimes n} \in \Fun(BS^1,\Sh_{L_n}(\bR^n_z))\simeq \Fun(BS^1,\Fun(\rarbZ^n,D(\Bbbk)))$. Here, for $m\in \rarbZ$, we write $\cE^n(m)\coloneqq \Ex(\cE^{\boxtimes n})(m)\simeq \cE^{\boxtimes n}_{m+(1,\cdots,1)}$.

Consequently, we have
\[\cE^{\boxtimes n}\simeq \int^{m\in\rarbZ^n} 1_{(m,\infty)}\otimes \cE(m) .\]    
\end{Lemma}
\begin{proof}By the constructibility and the microsupport estimation of \autoref{prop:chiu result on loop sheaf equiv}, we see that $\cE$  gives $\ms(\cE)\subset L_1$. So, $\ms(\cE^{\boxtimes n})\subset L_n$.  The sheaf $\cE$ admits the diagonal $S^1$-action, and then it comes with a $S^1$-action $\cE^{\boxtimes n}$ by diagonal action.
\end{proof}

\begin{Lemma}\label{lemma: cyclic trace coend second}We have the following equivalence of $S^1$-modules over $\Bbbk$
\[
F^{S^1}(X_\Omega)_{T}\simeq \int^{m\in\rarbZ^n}\left(\Gamma_c\left({\{(z,\zeta)\mid \zeta\in \Omega, T+\langle z,\zeta\rangle \geq 0\}};\Bbbk_{(m,\infty)}\right)\otimes \cE(m)[n]\right)\in \Fun(BS^1,D(\Bbbk)).
\]
\end{Lemma}
\begin{proof}By \autoref{prop: cyclic trace for toric domain-first formula}, we have 
\[F^{S^1}(X_\Omega)_T \simeq \Gamma_c(\{(z, \zeta) \mid \zeta \in \Omega, T + \left<z, \zeta \right> \geq 0\}; \cE^{\boxtimes n})[n].\]
However, the coend is a colimit, and $\Gamma_c$ commutes with colimits. Then we conclude by \autoref{lemma: coend formula for E^n} and a proper base change.
\end{proof}

\begin{Lemma}\label{lemma: cyclic trace coend integration}For all $m\in \bZ^n$, recall that we denote
\[
  \qquad \Omega(T;m)=\{\zeta\in \Omega\mid T+\langle m,\zeta\rangle\geq 0\}.
\]
Then we have
\begin{align*}&\,\Gamma_c\left({\{(z,\zeta)\mid \zeta\in \Omega, T+\langle z,\zeta\rangle \geq 0\}};\Bbbk_{(m,\infty)}\right)\otimes \cE(m) \\
\simeq &  \,\Gamma_c\left(\Omega(T;m) ;\Bbbk  \right)[-n]\otimes \cE(m) \in \Fun(BS^1,D(\Bbbk)).
\end{align*}
\end{Lemma}
\begin{proof}

We only need to show that, in $D(\Bbbk)$,
\begin{align*}
        &\,\Gamma_c({\{(z,\zeta)\mid \zeta\in \Omega, T+\langle z,\zeta\rangle \geq 0\}};\Bbbk_{(m,\infty)})\simeq\Gamma_c( \Omega(T;m) ;\Bbbk  )[-n].
\end{align*}

We denote 
\[Y_\Omega={\{(z,\zeta) \in (m,\infty)\times  \Omega\mid T+\langle z,\zeta\rangle \geq 0\}},\qquad \pi_{z}: Y_\Omega\to  \Omega.\]

For each $\zeta\in  \Omega \subset \bR^n_{\geq 0}$, we will study the fiber
\[\pi_z^{-1}(\zeta)\cong\{u\in (0,\infty)\mid T+\langle m,\zeta\rangle + \langle u,\zeta\rangle  \geq 0\}.  \]
When $T+\langle m,\zeta\rangle  \geq 0$, the fiber is homeomorphic to $(0,\infty)^n$ with compactly supported cohomology $\Bbbk [-n]$; and when $T+\langle m,\zeta \rangle  < 0$, we have $ \langle u,\zeta\rangle \geq D\coloneqq -T-\langle m,\zeta\rangle >0 $. Using the following change of variables
\[s=\langle u,\zeta\rangle,\quad \theta=u/\|u\|_1,\]
we see that 
\[\pi_z^{-1}(\zeta) \to [D,\infty)\times \Delta^\circ_{n-1}, z\mapsto (s,\theta)\]
defines a homeomorphism, and the fiber has vanishing compactly supported cohomology. 

Therefore, we have 
\[\pi_{z!}\Bbbk_{Y_\Omega}\simeq \Bbbk_{\{ \zeta\in  \Omega\mid T+\langle m,\zeta\rangle \geq 0\}}[-n]=\Bbbk_{ \Omega(T;m)}[-n]. \]
The result now follows from the cohomology Fubini theorem.
\end{proof} 

Now, we recall by definition of $\cE^{\boxtimes n}$, for $m\in \bZ^n$, we have
\[\cE(m)\simeq \Bbbk[n+2I(m)],\]
where $\Bbbk$ is equipped with the trivial $S^1$-action.

In summary, by \autoref{lemma: cyclic trace coend second} and \autoref{lemma: cyclic trace coend integration}, we have
\begin{Prop}\label{prop: trace for toric domain derived final}Let $\Omega\subset \bR^n_{\geq 0}$ be an open set, for $T\in \bR$, we have
\begin{align*}
    F^{S^1}(X_\Omega)_{T}\simeq &\int^{m\in\rarbZ^n} \Gamma_c({ \Omega(T;m)} ;\Bbbk  )\otimes \cE(m)\\
    \simeq &\int^{m\in\rarbZ^n} \Gamma_c({ \Omega(T;m)} ;\Bbbk  )[n+2I(m)]\in \Fun(BS^1,D(\Bbbk)),
\end{align*}
where $S^1$ acts on $\Gamma_c({ \Omega(T;m)} ;\Bbbk  )$ trivially.
\end{Prop}
\begin{RMK}Warning: the non-triviality of the $S^1$-action is recorded by the coend during its colimit process.    
\end{RMK}

\begin{RMK}\label{remark: coend formula for general locally closed V}
    So far, we did not pose any condition on the open set $\Omega$. Results in \autoref{section: trace for toric domain} are true for all toric domains. 
    
    We can define $F^{S^1}(X_V) \in \Fun(BS^1,\sT)$ for a locally closed set $V\subset \bR^n_{\geq0}$ using formula in \autoref{prop: cyclic trace for toric domain-first formula} simply by replacing $\Omega$ by $V$. Then the formula in \autoref{prop: trace for toric domain derived final}  is still true by replacing $\Omega(T;m)$ by $V(T;m)$, and the only thing we do not know is whether it is the equivariant trace of some category.
\end{RMK}

In the end, we prepare a coordinate boundary removal lemma for later applications. 
\begin{Lemma}\label{lemma: remove coordinate boundary in toric domain}Let $V\subset \bR^n_{\geq 0}$ locally closed and $T\in \bR$, we write $F^{S^1}(X_V)$ defined in \autoref{remark: coend formula for general locally closed V}. Now, if we assume that $V\subset \partial\bR^n_{\geq 0}=\cup_{i=1}^n \bR^n_{\geq 0}\cap H_i$ for $H_i=\{\zeta\mid \zeta_i=0\}$, we have that $F^{S^1}(X_V)_T\simeq 0$ for all $T\in \bR$.
\end{Lemma}
\begin{proof}We use the following coend formula as explained in \autoref{prop: trace for toric domain derived final}
\begin{align*}
    F^{S^1}(X_V)_{T}\simeq \int^{m\in\rarbZ^n} \Gamma_c({ V(T;m)} ;\Bbbk  )\otimes \cE(m).
\end{align*}

We set $V_1=V\cap H_1$, and $V^\circ=V\setminus V_1$. Then $V_1$ is closed in $V$ and $V^\circ$ is open in $V$. Then using this open-closed decomposition, we have a fiber sequence
\[F^{S^1}(X_{V^\circ})_{T}\to F^{S^1}(X_{V})_{T}\to F^{S^1}(X_{V_1})_{T}.\]
Therefore, we only need to show that $F^{S^1}(X_{V_1})_{T}\simeq 0$. If so, we see that $V^\circ$ is still a locally closed set satisfying a similar condition, but with no $H_1$ part. Then we simply repeat this process for $H_2$.

Now, we show $F^{S^1}(X_{V_1})_{T}\simeq 0$. Now, we notice that since $V_1\subset H_1$, then we have
\[V_1(T;m)=V_1(T;m+e_1)\]
for all $T$. Then the induced map
\[t_1:\Gamma_c({ V_1(T;m+e_1)} ;\Bbbk  )\to \Gamma_c({ V_1(T;m)} ;\Bbbk  )  \]
is the identity.

Recall the idea we use in \autoref{prop: coend of cube as total cofiber}, for $m=(m_1,m')$, we use the Fubini theorem for coends to see that
\[\int^{m\in\rarbZ^n} \Gamma_c({ V_1(T;m)} ;\Bbbk  )\otimes \cE(m)\simeq \int^{m'\in\rarbZ^{n-1}} \int^{m_1\in \rarbZ}\Gamma_c({ V_1(T;m)} ;\Bbbk  )\otimes \cE(m).\]
However, by $t_1=\id$, we have
\[
\int^{m_1\in\rarbZ}
\Gamma_c(V_1(T;(m_1,m'));\Bbbk)\otimes\cE(m_1,m')
\simeq
\Gamma_c(V_1(T;(0,m'));\Bbbk)
\otimes
\operatorname*{colim}_{m_1\in\rarbZ}\cE(m_1,m').
\]
However, we have
\[\operatorname*{colim}_{m_1\in\rarbZ}\cE(m_1,m')\simeq \operatorname*{colim}_{m_1\in\rarbZ}\Bbbk[n+2m_1+2I(0,m')] \simeq 0.\]
In fact, forgetting the circle action, it is clear the colimit is $0$. Then we conclude because forgetting the circle action functor is conservative. We then conclude the proof here.
\end{proof}

\subsection{Chain model and proof of \autoref{thm: chain model toric domain}~(1) - (4)}\label{section: chain model of homotopy fixed point}

In this section, we construct the chain model of $ F^{S^1}(X_\Omega)_{T}^{hS^1}$. In particular, we finish the proof of (1) - (4) of \autoref{thm: chain model toric domain}.

\begin{proof}[Proof of \autoref{thm: chain model toric domain}~(1)] We only consider the case $ F^{S^1}(X_\Omega)_{T}^{hS^1} $. The relative case follows directly.

By standard sheaf theory \cite[Chapter III.1]{BredonSheaf} or \cite[Section 3.5]{HainePortaTeyssier2023}, we know that for a locally contractible locally compact Hausdorff topological space $X$, we have
\[\Gamma_c(X;\Bbbk )\simeq C^*_c(X;\Bbbk ),\]
where $C^*_c$ stands for the compactly supported singular cochain. 

Therefore, by \autoref{prop: trace for toric domain derived final}, we have
\begin{align*}
    F^{S^1}(X_\Omega)_{T}\simeq  \int^{m\in\rarbZ^n} C^{*+n+2I(m)}_c({\Omega(T;m)} ;\Bbbk  )\in \Fun(BS^1,D(\Bbbk)),
\end{align*}

Now, we shall use \autoref{thm: fixed point as a monoidal functor}. It follows that we have
\[F^{S^1}(X_\Omega)_{T}^{hS^1}\simeq  \widehat{\int} ^{m\in\rarbZ^n} C^{*+n+2I(m)}_c({ \Omega(T;m)} ;\Bbbk  )\fpw{u}\in D(\Bbbk\fpw{u})_{comp}.\]

Importantly, the notation $\widehat{\int} $ indicates that this coend is the complete colimit taking in $D(\Bbbk\fpw{u})_{comp}$, which is, by definition, the completion of the incomplete coend in $D(\Bbbk\fpw{u})$
\begin{equation}\label{equation: underived coend in chain}
  \int^{m\in\rarbZ^n} C^{*+n+2I(m)}_c({ \Omega(T;m)} ;\Bbbk  )\fpw{u}\in D(\Bbbk\fpw{u}).  
\end{equation}

Now, for this incomplete coend, we shall use \autoref{coro: coend computing in chain} to compute it. To do so, we consider the following bifunctor that defines the coend
\begin{align*}
    B:&\,(\rarbZ^n)^\op \times \rarbZ^n\to D(\Bbbk\fpw{u})\\ &\,(l,m) \mapsto  C^{*+n+2I(m)}_c({ \Omega(T;l)} ;\Bbbk  )\fpw{u} 
\end{align*}
with the morphisms
\[(\partial_i^0)_m=t_{i,m},\qquad (\partial_i^1)_m=-(m_i+1)u.\]

Recall the exterior algebra collection: Let $\bigwedge^\bullet\Bbbk ^n$ be the exterior algebra of $\Bbbk ^n$ with the basis $e_{I}=e_{i_1}\wedge \cdots\wedge e_{i_p}$, and we set $|e_i|=1$ and $\iota_{e_{i_r}}:\bigwedge^\bullet\Bbbk ^n\to \bigwedge^\bullet\Bbbk ^n$ the inner product $\iota_{e_{i_r}}(e_{I})=(-1)^{r-1}e_{I\setminus \{i_r\}}$. We also write the parity function $(-1)^N$ on $\bigwedge^\bullet\Bbbk ^n$, i.e. $(-1)^N|_{\wedge^q\Bbbk ^n}=(-1)^q$. 

Then by \autoref{coro: coend computing in chain}, we have that the incomplete coend \autoref{equation: underived coend in chain} is computed by the graded $\Bbbk\fpw{u}$-module
\[I_{T}^{S^1}(X_\Omega)\coloneqq \bigoplus_{m\in \bZ^n}\left( C^{*+n+2I(m)}_c({ \Omega(T;m)} ;\Bbbk  )\fpw{u} \otimes \bigwedge\nolimits^\bullet \Bbbk ^n\right),\]
with the differential
\[d=d_{\mathrm{sing}}\otimes (-1)^N+\sum_i  (m_iu+t_i)\otimes \iota_{e_i},\]
where $d_{\mathrm{sing}}$ should be understood as the usual singular differential tensor with $\Bbbk\fpw{u}$. 

We set $\widehat I_T^{S^1}(X_\Omega)$ the $u$-adic completion of $I_T^{S^1}(X_\Omega)$. Therefore, by definition of $\widehat{\int} $, we have
\[F_{T}^{S^1}(X_\Omega)^{hS^1}\simeq(\widehat I_T^{S^1}(X_\Omega),d).\qedhere\]

In the end, the relative version follows from: (i) by definition a fiber of the absolute case, and (ii) the completion is an exact functor.
\end{proof}

\begin{proof}[Proof of \autoref{thm: chain model toric domain}~(2) and (3)] The proof for (2) and (3) is clear. Both of them are induced by corresponding open inclusions between $\Omega(T;m)$. Notice that since we are working on compact support cochains, the induced maps are given by extension by zero rather than restriction.
\end{proof}

Next, we shall discuss the completeness of $I_{(\epsilon,T]}^{S^1}(X_\Omega)$. In particular, we prove (4) of \autoref{thm: chain model toric domain}.
Here, we consider the following algebraic lemma.
\begin{Lemma}\label{lemma: sufficient condition for completeness}Let $\{C_m\}_{m\in M}$ be a set of derived $\Bbbk$-modules equip with the trivial $S^1$-action. We assume that there exists $N$ such that $H^q(C_m)=0$ for all $q<N$ and all $m\in M$, then the $\Bbbk\fpw{u}$-module $\oplus_m C_m^{hS^1}=\oplus_m C_m\fpw{u}$ is $u$-adic complete, i.e. $\oplus_m C_m^{hS^1} \to \widehat\oplus_m C_m^{hS^1}$ is an equivalence in $D(\Bbbk\fpw{u})$.
\end{Lemma}
\begin{proof}As the $S^1$-action is trivial, we have $C_m^{hS^1}=C_m\fpw{u}=\prod_j C_m[-2j]$.
We set $P=\bigoplus_m C_m^{hS^1}=\bigoplus_m\prod_j C_m[-2j] \in D(\Bbbk)$. Because direct sums/products are exact in $D(\Bbbk)$, we have that 
\[H^q(P)=\bigoplus_m\prod_j H^{q-2j}(C_m)=0,\quad \forall q<N.\]

However, the Milnor sequence \cite[Tag 0CQE]{StacksProject} and the previous vanishing implies that
\[H^q\big(\varprojlim_r P[-2r]\big)=H^q\big(\varprojlim (\cdots\xrightarrow{u} P[-4] \xrightarrow{u} P[-2]\xrightarrow{u} P )\big)=0,\quad \forall q<N.\]
Then we conclude that $P$ is $u$-adic complete by (a graded version of) \cite[Tag 091P-(3)]{StacksProject}.
\end{proof}

\begin{proof}[Proof of \autoref{thm: chain model toric domain}~(4)] 
Now, for $I_{(\epsilon,T]}^{S^1}(X_\Omega)$, again we assume that $\epsilon=-\infty$ and we fix $T$, we set
\[C_m=C^{*+n+2I(m)}_c({\Omega(T;m)} ;\Bbbk  )\otimes \bigwedge\nolimits^\bullet \Bbbk ^n.\]
Then we have 
\[H^q(C_m)=\bigoplus_{a=0}^nH^{q-a+n+2I(m)}_c({\Omega(T;m)} ;\Bbbk)\otimes \bigwedge\nolimits^a \Bbbk ^n.\]

Since $\Omega$ is homologically star-shaped, by definition, we have for every $T$, there exists $C(T)$ such that if $I(m) \geq C(T)$ then $H^{*}_c({\Omega(T;m)} ;\Bbbk)=0$ by excision. Therefore, we have
\[\forall q<N=-n-2C(T), ~\forall m\in \bZ^n,\qquad H^q(C_m)=0.\]
Then $I_{(\epsilon,T]}^{S^1}(X_\Omega)$ is complete by the previous Lemma.
\end{proof}
\begin{RMK}In fact, one can argue that $\Omega$ being homologically star-shaped is also the necessary condition for $I_{(\epsilon,T]}^{S^1}(X_\Omega)$ being $u$-adic complete.
\end{RMK}

\subsection{Comparison with Irie's chain model}\label{section: Irie's model}
In this section, we compare our chain model and Irie's conjectural chain model in dimension $4$ \cite{Irie_ChainModel}. Here, we will explain Irie's construction in the proof of \autoref{prop: compare with Irie Model}. To compare with Irie's model, we shall use $O_{(\epsilon,T]}^{S^1}(X_\Omega)$ in \autoref{remark: in and out for chain model} rather than $I_{(\epsilon,T]}^{S^1}(X_\Omega)$, and we will see \autoref{lemma: remove coordinate boundary in toric domain} that we can compute $O_{(\epsilon,T]}^{S^1}(X_\Omega)$ using $U_\Omega(T;m)\cap \bR^n_{>0}$ as Irie did.

\begin{Prop}\label{prop: compare with Irie Model}Let $0<\epsilon<T<\infty$ and let $\Omega$ be a bounded smooth star-shaped domain. When $n=2$ and $\Bbbk=\bQ$ (where Irie's model is defined), the chain complex $O_{(\epsilon,T]}^{S^1}(X_\Omega)$ is chain homotopy equivalent to (a cohomology version of) Irie's conjectural model for positive $S^1$-equivariant symplectic homology for a toric domain in $\bR^4$.
\end{Prop}

\begin{proof}For $m\neq 0$, we set $\iota_m=\sum_{i=1}^2 m_i\iota_{e_i}$. Then the Koszul complex $
\left(
\bQ\fpw{u}\otimes\bigwedge^\bullet \bQ^2,
u\iota_m
\right)
$
admits an explicit deformation retract onto $\bigwedge^\bullet \ker(m)= H^*(S^1)$ where we view $m: \bQ ^2 \to \bQ $ as a covector. 

Precisely, pick $v=(v_{1,m},v_{2,m} ) $ such that $m_1v_{1,m}+m_2v_{2,m}=1$, and set $w:=m_2e_1-m_1e_2$. Then $\ker(m)=\bQ w$ and $
\bQ^2=\bQ v\oplus \bQ w$.
Let $i_m:(\bigwedge^\bullet \ker(m),0)
\longrightarrow
(\bQ\fpw{u}\otimes\bigwedge^\bullet \bQ^2, u\iota_m)$ be the natural inclusion. 
Define
\[
p_m:
\Big(\bQ\fpw{u}\otimes\bigwedge\nolimits^\bullet \bQ^2, u\iota_m\Big)
\longrightarrow
\Big(\bigwedge\nolimits^\bullet \ker(m), 0\Big)
\]
by
\[
p_m(\eta)=\eta,
\qquad
p_m(u^k\eta)=0\quad (k>0),
\qquad
p_m(u^k v\wedge\eta)=0,
\]
for $\eta\in\bigwedge^\bullet\ker(m)$. 
Finally, define the homotopy by
\[
h_m(u^k\eta)
=
\begin{cases}
-u^{k-1}v \wedge\eta,&k\geq1,\\
0,&k=0,
\end{cases}
\qquad
h_m(u^k v \wedge\eta)=0.
\]
For the differential $\delta = u \iota_m$, a direct computation gives
\[
\delta h_m+h_m\delta
=
i_mp_m-\id.
\]
Thus $(i_m,p_m,h_m)$ defines a deformation retract.

Take the direct sum over $m \notin\bZ_{\geq0}^2$, and denote the resulting complexes by
\begin{gather*}
\mathcal{O}_{(\epsilon,T]}(X_\Omega):=\bigoplus_{m\notin\bZ_{\geq0}^2}
C^{*+2+2I(m)}_c({U_\Omega(T;m)},U_\Omega(\epsilon;m) ;\bQ  )\fpw{u}
\otimes\bigwedge\nolimits^\bullet \bQ^2,\\
\mathcal{H}_{(\epsilon,T]}(X_\Omega)
:=
\bigoplus_{m\notin\bZ_{\geq0}^2}
C^{*+2+2I(m)}_c({U_\Omega(T;m)},U_\Omega(\epsilon;m) ;\bQ  )
\otimes\bigwedge\nolimits^\bullet \ker(m),
\end{gather*}
and denote the resulting maps by $p,i,h$. Then we have homotopy equivalences
\[p: (\mathcal{O}_{(\epsilon,T]}(X_\Omega), d_{\mathrm{sing}}\otimes (-1)^N + u\iota_m\Big)
\longrightarrow 
(\mathcal{H}_{(\epsilon,T]}(X_\Omega), d_{\mathrm{sing}} \otimes (-1)^N).\]
Recall that the pull-back of the natural inclusion is denoted by
\[t_i: C^*_c({U_\Omega(T;m)},U_\Omega(\epsilon;m) ;\bQ  ) \to C^*_c({U_\Omega(T;m-e_i)},U_\Omega(\epsilon;m-e_i) ;\bQ  ), \quad i=1,2.\]
We now regard $\tau=t_1\otimes \iota_{e_1}+t_2\otimes \iota_{e_2}$ as a perturbation on the differential $d=d_{\mathrm{sing}} \otimes (-1)^N + u \iota_m$.

Every homogeneous element of $\mathcal O_{(\epsilon,T]}(X_\Omega)$
has finite lattice support and bounded $u$-degree. Since $h$ lowers the $u$-degree, $\tau h$ is locally nilpotent.
Thus the homological perturbation lemma applies \cite{Homological-perturbation}, and it 
gives a chain homotopy equivalence
\[
\Phi
=
p(1-\tau h)^{-1}
:
(\mathcal{O}_{(\epsilon,T]}(X_\Omega), d_{\mathcal{O}})
\longrightarrow
(\mathcal{H}_{(\epsilon,T]}(X_\Omega), d_{\mathcal{H}})
\]
for the perturbed differential $d_{\mathcal{O}} = d_{\mathrm{sing}}\otimes (-1)^N + u\iota_m + \tau$ on the complex $\mathcal{O}_{(\epsilon,T]}(X_\Omega)$ and some perturbed differential $d_{\mathcal{H}} = d_{\mathrm{sing}}\otimes (-1)^N + p(1-\tau h)^{-1}\tau i$ on the complex $\mathcal{H}_{(\epsilon,T]}(X_\Omega)$.

The chain complex $\mathcal{O}_{(\epsilon,T]}(X_\Omega)$ is homotopy equivalent to $ {O}^{S^1}_{(\epsilon,T]}(X_\Omega)$ (where the sum is taken for all $m\in \bZ^n$) since the singular cochain part contribute acyclic terms when $m \in\bZ_{\geq0}^2$.

On the other hand, we consider $\mathcal{H}_{(\epsilon,T]}(X_\Omega)$. Since $\tau i$ has $u$-degree zero, we have $h\tau i=0$. Therefore
\[
d_{\mathcal{H}}
=
d_{\mathrm{sing}}\otimes (-1)^N+p\tau i.
\]
Moreover, $
p\tau i(w_m)
=
m_2t_1-m_1t_2
$. Consequently,
\[
\mathcal{H}_{(\epsilon,T]}(X_\Omega)
=
\bigoplus_{m\notin\bZ_{\geq0}^2}
C^{*+2+2I(m)}_c({U_\Omega(T;m)},U_\Omega(\epsilon;m) ;\bQ  )
\otimes
\bigl(
\bQ\cdot1_m\oplus\bQ\cdot w_m
\bigr)
\]
has transferred differential
\[
d_{\mathcal{H}}(\alpha\otimes1_m)
=
d_{\mathrm{sing}}\alpha\otimes1_m
,\quad 
d_{\mathcal{H}}(\alpha\otimes w_m)
=
-d_{\mathrm{sing}}\alpha\otimes w_m
+
m_2t_1\alpha \otimes 1_{m-e_1}-m_1 t_2 \alpha
\otimes 1_{m-e_2}.
\]

After identifying
\[
1\longleftrightarrow e_0,
\qquad
w_m\longleftrightarrow e_1,
\]

this is precisely the cohomological version of Irie's chain model
\cite{Irie_ChainModel} (notice that our $m$ should be $-m$ therein).
\end{proof}
\begin{RMK}In fact, the strategy here for simplifying the chain model works for all $n\geq 2$. It may help for further computation in machine.
\end{RMK}

We now compare multiplication by the equivariant variable $u$ with the
chain map $U$ introduced by Irie \cite[Proposition 3.1]{Irie_ChainModel}.

\begin{Prop}The chain homotopy $\Phi$ constructed above ``interwines" $u$-action on $\mathcal{O}_{(\epsilon,T]}^{S^1}(X_\Omega)$ and the $u$-action on $\mathcal{H}_{(\epsilon,T]}(X_\Omega)$ constructed by Irie on a direct summand. 
\end{Prop}
\begin{RMK}Here, we say ``interwine" since on certain sub-complexes, we need to adjust $\Phi$ a bit. We will see this more precisely in the proof.    
\end{RMK}
\begin{proof}It is clear that the chain homotopy equivalence $O_{(\epsilon,T]}^{S^1}(X_\Omega)\simeq \mathcal{O}_{(\epsilon,T]}(X_\Omega)$ commutes with $u$-action. We only consider the latter.

We follow the same homological perturbation strategy. Let
\[
\mu_u:
\mathcal{O}_{(\epsilon,T]}(X_\Omega)
\longrightarrow
\mathcal{O}_{(\epsilon,T]}(X_\Omega)[2]
\]
be multiplication by $u\in \bQ\fpw{u}$. Transfer this map through the above deformation
retract to a map 
\[U_{\mathrm{tr}}: \mathcal{H}_{(\epsilon,T]}(X_\Omega)
\longrightarrow
\mathcal{H}_{(\epsilon,T]}(X_\Omega)[2]
\]
where $
U_{\mathrm{tr}}
\coloneqq 
\Phi\mu_u\Psi,
$
where 
$\Phi = p(1-\tau h)^{-1}$
and
$
\Psi=(1-h\tau )^{-1}i
$
\cite{Homological-perturbation}. Since $h\tau i=0$, we simply have $\Psi=i$ and $U_{\mathrm{tr}} = p(1-\tau h)^{-1}\mu_ui = p\tau h \mu_u i$.

The transferred $u$-action can be computed explicitly. Recall $v = (v_{1,m} ,v_{2,m})$ and $h_m(u^k\eta) = -u^{k-1}v \wedge \eta$, and $\bigwedge^\bullet \ker(m)=\bQ 1_m\oplus \bQ \omega_m$.

We have
\begin{align*}
    U_{\mathrm{tr}}(\alpha\otimes 1_m)
&=-
v_{1,m}t_1(\alpha)\otimes 1_{m-e_1}
-v_{2,m}t_2(\alpha)\otimes 1_{m-e_2},\\
U_{\mathrm{tr}}(\alpha\otimes w_m)
&=
-v_{1,m-e_1}t_1(\alpha)\otimes w_{m-e_1} 
-v_{2,m-e_2}t_2(\alpha)\otimes w_{m-e_2} .
\end{align*}

On the other hand, Irie's $U$-map is \cite[Proof of Proposition 3.1]{Irie_ChainModel}
\[
U_{\mathrm{Irie}}
=
t_1+t_2.
\]
We consider the weight $
w=w(m,\eta)
\coloneqq 
I(m)-|\eta|,
\eta\in\bigwedge^{|\eta|}\ker(m)$ on the complex $\mathcal{H}_{(\epsilon,T]}(X_\Omega)$.
Let\[
L|_{C^*(U_\Omega(T;m),U_\Omega(\epsilon;m))\otimes\wedge^j\ker(m)}
=
(I(m)-j)\id.
\]
The differential $d_{\mathcal{H}}$ preserves this weight, so the complex $\mathcal{H}$ splits by weights.
Define a degree-one map \[
K:\mathcal{H}_{(\epsilon,T]}(X_\Omega)
\longrightarrow
\mathcal{H}_{(\epsilon,T]}(X_\Omega)
\]
by
\[
K(1_m)
=
\bigl(
v_{2,m}-v_{1,m}
\bigr)w_m,
\qquad
K(w_m)=0.
\]
Then we have
\[
U_{\mathrm{Irie}}
+
U_{\mathrm{tr}}L
=
d_{\mathcal{H}}K+Kd_{\mathcal{H}}.
\]

Therefore, on the weight-$w\geq 0$ summand of the chain complex, we have
\[
U_{\mathrm{Irie}}
\simeq
-w U_{\mathrm{tr}}.
\]
On the $w< 0$ summand, we should be careful that $U_{\mathrm{Irie}}$ only differs from $U_{\mathrm{tr}}$ by a suitable rescaling (depending on $w$), we leave the computation to readers.
\end{proof}

\subsection{Spectral sequence and proof of \autoref{thm: chain model toric domain}~(5)}\label{subsection: spectral sequence}
In this section, we build a filtration on the chain model $I_{(\epsilon,T]}^{S^1}(X_\Omega,\Bbbk)$ of $F_{(\epsilon,T]}^{S^1}(X_\Omega)$ and study its convergence under the assumptions that $\Omega$ is \textit{homologically star-shaped} and $ \Bbbk$ is a field containing $\bQ$, which proves \autoref{thm: chain model toric domain}~(5).

Here, we notice that $I_{(\epsilon,T]}^{S^1}(X_\Omega,\Bbbk)$ is complete as proved in \autoref{thm: chain model toric domain}-(4), so we can work on this direct sum model. Many of the computations here are straightforward generalizations of \cite{Irie_ChainModel}.  

\begin{Lemma}For any $p$, we denote
\[F^pI_{(\epsilon,T]}^{S^1}(X_\Omega)= \bigoplus_{m\in \bZ^n,I(m)\leq p}\left( C^{*+n+2I(m)}_c({\Omega(T;m)},\Omega(\epsilon;m) ;\Bbbk )\fpw{u} \right)\otimes \bigwedge\nolimits^\bullet \Bbbk ^n.\]

Then $F^\bullet I_{(\epsilon,T]}^{S^1}(X_\Omega)$ is an increasing filtration of $I_{(\epsilon,T]}^{S^1}(X_\Omega)$, which is exhaustive and Hausdorff.    
\end{Lemma}
\begin{proof}It is clear that $d_{\mathrm{sing}}$ preserves $m$, while $t_i\otimes \iota_i$ maps $m$ to $m-e_i$, which decreases $I$. Then under the total differential $d$, we have $d (F^p I_{(\epsilon,T]}^{S^1}(X_\Omega))\subset  F^p I_{(\epsilon,T]}^{S^1}(X_\Omega)$. Therefore, $F^p I_{(\epsilon,T]}^{S^1}(X_\Omega)$ is a subcomplex. 
The filtration is clearly increasing, exhaustive, and Hausdorff.    
\end{proof}
We now have the spectral sequence of the filtration $(E_r,d_r)$.
\begin{Lemma}The spectral sequence is (strongly) convergent:
\[E_r\Rightarrow H^nI_{(\epsilon,T]}^{S^1}(X_\Omega)\simeq H^nF^{S^1}_{(\epsilon,T]}(X_\Omega).\]
\end{Lemma}
\begin{proof}We already see the filtration is exhaustive and Hausdorff. To prove the convergence, we show the filtration is bounded from below in the sense of \cite[Definition~5.4.2(2)]{Weibel1994}: For a homogeneous element $\alpha u^k\otimes \eta\in F^pI^q_{(\epsilon,T]}(X_\Omega)$, we have
\[q=|\alpha|+2k+|\eta|-n-2I(m)\geq  -n-2p.\]
Therefore, if $2p<-n-q$, we have $F^pI^q_{(\epsilon,T]}=0$. Then we conclude convergence by \cite[Theorem~5.5.1(2)]{Weibel1994}.
\end{proof}

We then compute the $E_0$ and $E_1$ pages of the spectral sequence.

\begin{Lemma}\label{lem: E_0 page}We have
\[E_0^{p,\bullet-p}=\operatorname{gr}^pI_{(\epsilon,T]}^{S^1}(X_\Omega)=\bigoplus_{I(m)= p} \left(C^{*+n+2p}_c({\Omega(T;m)},{\Omega(\epsilon;m)} ;\Bbbk )\fpw{u}\otimes \bigwedge\nolimits^\bullet \Bbbk ^n \right),\]
with the differential
\[d_0(\alpha u^k\otimes \eta)=(-1)^{|\eta|} d_{\mathrm{sing}}(\alpha u^k)\otimes \eta+\alpha u^{k+1}\otimes \left(\sum_{i=1}^n m_i \iota_{e_i}\eta\right).\]    
\end{Lemma}
\begin{proof}Recall that the differential on $I_{(\epsilon,T]}^{S^1}(X_\Omega)$ is $d=d_{\mathrm{sing}}\otimes (-1)^N+\sum_{i=1}^n (m_iu+t_i)\otimes \iota_{e_i}$. On the associated graded, we have that $t_{i,m}$ vanishes since $t_{i,m}F^p\subset F^{p-1}$, which gives the required differential formula.    
\end{proof}

\begin{Lemma}\label{lem: E_1 page}For $m\in \bZ^n$, set $\iota_m\coloneqq \sum_{i=1}^n m_i \iota_{e_i} $. We have 
\begin{align*}
    E_1^{p,q-p}=&\bigoplus_{m\neq 0, I(m)= p} \bigoplus_{j=0}^{n-1}H^{q+n+2p-j}_c({\Omega(T;m)},{\Omega(\epsilon;m)} ;\Bbbk )   \otimes \wedge^j \ker(m).
\end{align*}
Moreover, we have that $u E_1=0$ and the differential formula, 
\[d_1([\alpha] \otimes \eta)= 
    \sum_{i=1}^n(t_i)_*[\alpha]\otimes \iota_{e_i}\eta
.\]
\end{Lemma}
\begin{proof}
We compute the $E_1$-page using the differential in \autoref{lem: E_0 page} on each summand:
\begin{align*}
d_0: & \,C_c^{*+n+2p}(\Omega(T; m), \Omega(\epsilon; m); \Bbbk)\fpw{u} \otimes \bigwedge\nolimits^\bullet \Bbbk^n \\
&\, \to C_c^{*+n+2p}(\Omega(T; m), \Omega(\epsilon; m); \Bbbk)\fpw{u} \otimes \bigwedge\nolimits^{\bullet} \Bbbk^n.\end{align*}

When $m=0$, we have 
\(d_0(\alpha u^k \otimes \eta) = (-1)^{|\eta|}d_{\mathrm{sing}}(\alpha)u^k \otimes \eta.\) 
For $T\geq 0$, we have $\Omega(T;m)=\Omega$ satisfying $H^*_c(\Omega(T;m),\Bbbk)= 0$ under the homologically star-shaped condition; for $\epsilon\in \bR$, we either have $\Omega(\epsilon;m)=\Omega$ for $\epsilon\geq 0$ or $\Omega(\epsilon;m)=\varnothing$ for $\epsilon<0$, and both of them implies $H^*_c(\Omega(\epsilon;m),\Bbbk)=0$. Consequently, it does not contribute to the $E_1$-page of the spectral sequence.

When $m\neq 0$, we pick $v\in \Bbbk^n$ such that $m(v)=1$, and we split $\Bbbk^n=\Bbbk v\oplus \ker (m)$. Without the singular cochains, the cohomology of the Koszul complex
\[\Big(\Bbbk\fpw{u} \otimes \bigwedge\nolimits^\bullet \Bbbk ^n,m_1 u, \cdots, m_n u\Big)\]
has differential $d_{\mathrm{Kos}}(u^k \otimes \eta) = u^{k+1} \sum_{i=1}^n m_i \iota_{e_i}\eta$ and is isomorphic to 
\[\bigwedge \nolimits^\bullet \ker(m).\]
By the Künneth formula, for the differential $d_0 = (-1)^{N}d_{\mathrm{sing}} \otimes 1 + 1 \otimes d_{\mathrm{Kos}}$, we then have that its cohomology is
\[H^{s(q)}_c({\Omega(T;m)} ;\Bbbk )   \otimes \wedge^j \ker(m),\]
for some $s(q)$. 
Here, at degree $q$, the degree for singular cochains satisfies the equation
\[q=s(q)+j-n-2p,\]
which gives $s(q)=q+n+2p-j$.

To see that $uE_1=0$: we pick $\alpha$ closed and $\eta\in\Lambda^\bullet\ker (m)$, then 
\[d_0(\alpha \otimes v\wedge \eta)=u \alpha\otimes \eta,\]
which tells us that for any representative $\alpha\otimes \eta$ on the $E_1$ page, multiplying by $u$ gives a $d_0$ differential. It means that the action of $ u$ on $E_1$ is $0$. The differential formula is then clear.
\end{proof}

Lastly, we state the following degeneration result.
\begin{Prop}\label{prop: degeneration of SS}Let $-\infty\leq \epsilon<T$, then we have the spectral sequence of $I_{(\epsilon,T]}^{S^1}(X_\Omega)$ degenerates at $E_n$.
\end{Prop}
\begin{proof}Recall that by \autoref{lem: E_1 page} we have
\[
    E_1^{p,q}=\bigoplus_{m\neq 0, I(m)= p} \bigoplus_{j=0}^{n-1}H^{q+n+3p-j}_c({\Omega(T;m)},{\Omega(\epsilon;m)} ;\Bbbk )   \otimes \wedge^j \ker(m).
\]

For $\alpha u^k \otimes \eta$ in the $m$-summand, we consider the weight
\[w=I(m)-k-|\eta|.\]
By the differential formula, we see that $d$ does not change $w$, and the filtration is also preserved with respect to $w$. Let
\[I(w)=\bigoplus_{\substack{m\in \bZ^n,k\geq 0,n\geq j\geq 0,\\I(m)-k-j=w}}C_c^*({\Omega(T;m)},{\Omega(\epsilon;m)} ;\Bbbk )  u^k \otimes\wedge^j\Bbbk^n,\quad F^pI(w)=F^pI_{(\epsilon,T]}^{S^1}(X_\Omega)\cap I(w).\]

Then we have
\[I_{(\epsilon,T]}^{S^1}(X_\Omega)=\bigoplus_w I(w),\quad F^pI_{(\epsilon,T]}^{S^1}(X_\Omega) =\bigoplus_w F^pI(w)\]
We write the corresponding fixed-weight spectral sequence as $E_r^{p,q}(w)$ for $r\geq 1$, see \autoref{fig:spectral sequence degeneration}. 

Now, we go to the $E_1$ page. We have
\[E_1^{p,q}(w) = \bigoplus_{\substack{m\neq 0,n\geq j\geq 0,\\I(m)=p,p-j=w}} \bigoplus_{j=0}^{n-1}H_c^{q+n+3p-j}(\Omega(T;m),\Omega(\epsilon; m); \Bbbk) \otimes \wedge^j \ker(m).\]
We have that for fixed $w$, we have $p=w+j$ where $j$ is the grading of the exterior algebra, and then
\[ E_1^{p,q}(w)=0,\qquad p\notin[w,w+n-1].
\]
On the other hand, the weight is preserved by $d$. So, we have for all $r\geq 1$ that
\[ E_r^{p,q}(w)=0,\qquad p\notin[w,w+n-1].
\]

Then we conclude $d_r=0$ for $r\geq n$ since $d_r$ has bigrading $(-r,r+1)$.    
\end{proof}
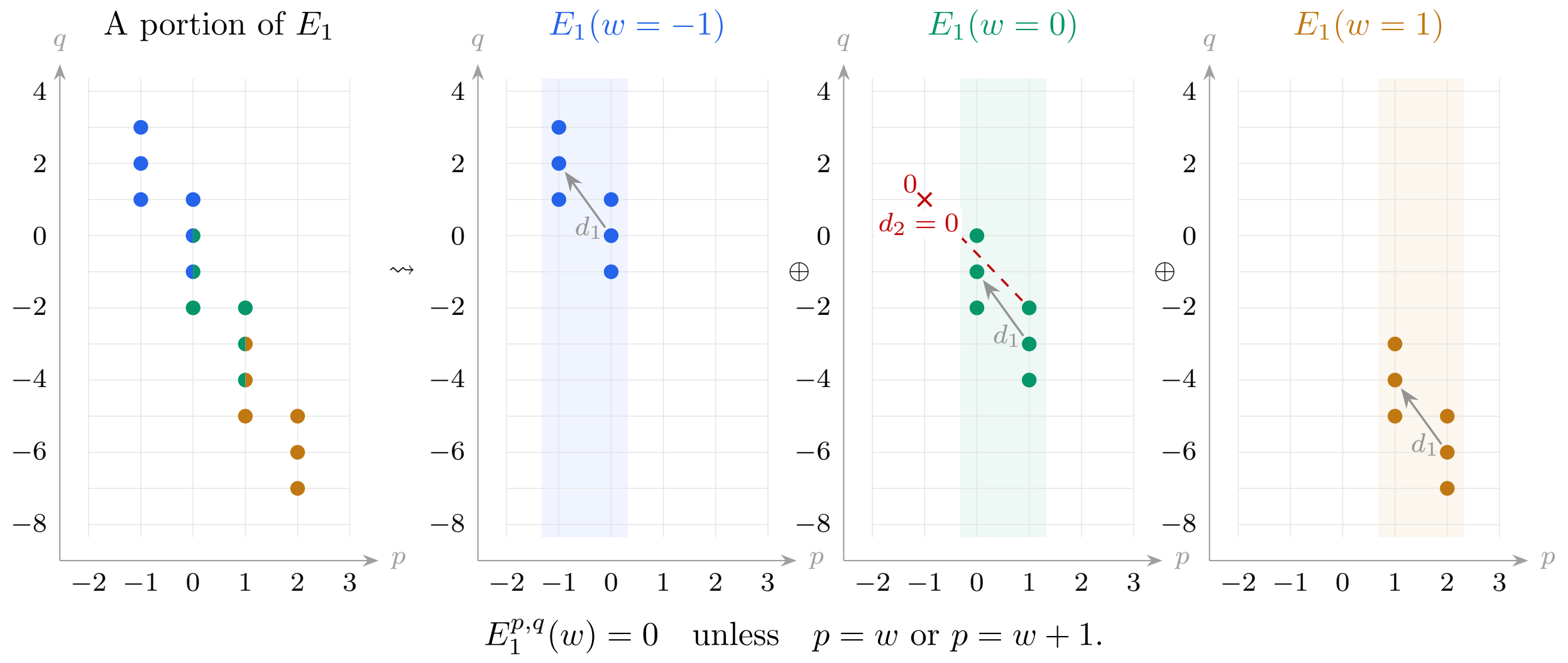
\begin{figure}[htbp]
    \centering
\begin{tikzpicture}[
  x=0.55cm,y=0.38cm,
  >=Stealth,
  every node/.style={font=\scriptsize},
  lattice/.style={gray!22,line width=0.25pt},
  axis/.style={gray!75,line width=0.4pt,->},
  dOne/.style={->,gray!85,line width=0.65pt,shorten <=3pt,shorten >=3pt},
  dTwo/.style={->,red!75!black,dashed,line width=0.65pt,
    shorten <=3pt,shorten >=4pt}
]
\definecolor{weightblue}{HTML}{2563EB}
\definecolor{weightgreen}{HTML}{059669}
\definecolor{weightorange}{HTML}{C07812}

\foreach \p in {-2,...,3}
  \draw[lattice] (\p,-8.35)--(\p,4.35);
\foreach \q in {-8,...,4}
  \draw[lattice] (-2,\q)--(3,\q);
\draw[axis] (-2.55,-9)--(3.55,-9) node[right] {$p$};
\draw[axis] (-2.55,-9)--(-2.55,4.75) node[above] {$q$};
\foreach \p in {-2,...,3}
  \node[below] at (\p,-9) {$\p$};
\foreach \q in {-8,-6,...,4}
  \node[left] at (-2.55,\q) {$\q$};
\node[font=\small] at (0.5,5.8) {A portion of $E_1$};

\foreach \q in {1,2,3}
  \fill[weightblue] (-1,\q) circle[radius=2.2pt];
\fill[weightblue] (0,1) circle[radius=2.2pt];
\fill[weightgreen] (0,-2) circle[radius=2.2pt];
\fill[weightgreen] (1,-2) circle[radius=2.2pt];
\fill[weightorange] (1,-5) circle[radius=2.2pt];
\foreach \q in {-7,-6,-5}
  \fill[weightorange] (2,\q) circle[radius=2.2pt];
\foreach \q in {-1,0}{
  \begin{scope}[shift={(0,\q)}]
    \fill[weightblue] (0,2.2pt)
      arc[start angle=90,end angle=270,radius=2.2pt]--cycle;
    \fill[weightgreen] (0,2.2pt)
      arc[start angle=90,end angle=-90,radius=2.2pt]--cycle;
  \end{scope}
}
\foreach \q in {-4,-3}{
  \begin{scope}[shift={(1,\q)}]
    \fill[weightgreen] (0,2.2pt)
      arc[start angle=90,end angle=270,radius=2.2pt]--cycle;
    \fill[weightorange] (0,2.2pt)
      arc[start angle=90,end angle=-90,radius=2.2pt]--cycle;
  \end{scope}
}
\node at(4.0
,-1) {$\rightsquigarrow$};

\foreach \offset/\ww/\col in
  {8/-1/weightblue,15/0/weightgreen,22/1/weightorange}{
  \begin{scope}[shift={(\offset,0)}]
    \fill[\col!7] (\ww-0.32,-8.35)
      rectangle (\ww+1.32,4.35);
    \foreach \p in {-2,...,3}
      \draw[lattice] (\p,-8.35)--(\p,4.35);
    \foreach \q in {-8,...,4}
      \draw[lattice] (-2,\q)--(3,\q);
    \draw[axis] (-2.55,-9)--(3.55,-9) node[right] {$p$};
    \draw[axis] (-2.55,-9)--(-2.55,4.75) node[above] {$q$};
    \foreach \p in {-2,...,3}
      \node[below] at (\p,-9) {$\p$};
    \foreach \q in {-8,-6,...,4}
      \node[left] at (-2.55,\q) {$\q$};
    \node[font=\small,text=\col] at (0.5,5.8) {$E_1(w=\ww)$};

    \pgfmathtruncatemacro{\sourceq}{-3-3*\ww}
    \pgfmathtruncatemacro{\targetq}{-1-3*\ww}
    \draw[dOne] (\ww+1,\sourceq)
      --node[pos=0.42,below,yshift=-2pt,inner sep=1pt] {$d_1$}
      (\ww,\targetq);

    \ifnum\ww=0
      \draw[dTwo] (1,-2)
        --node[pos=0.63,above left,fill=white,inner sep=1pt] {$d_2=0$}
        (-1,1);
      \draw[red!75!black,line width=0.7pt]
        (-1.13,0.81)--(-0.87,1.19)
        (-1.13,1.19)--(-0.87,0.81);
      \node[red!75!black,above left,inner sep=2pt] at (-1,1) {$0$};
    \fi

    \foreach \j in {0,1}{
      \pgfmathtruncatemacro{\pp}{\ww+\j}
      \foreach \a in {0,1,2}{
        \pgfmathtruncatemacro{\qq}{\a-2-3*\pp+\j}
        \fill[\col] (\pp,\qq) circle[radius=2.2pt];
      }
    }
  \end{scope}
}
\node at (11.6,-1) {$\oplus$};
\node at (18.6,-1) {$\oplus$};

\node[font=\small] at (11.5,-11.1)
  {$E_1^{p,q}(w)=0\quad\text{unless}\quad p=w\ \text{or}\ p=w+1.$};

\end{tikzpicture}
      \caption{Let $n=2$. All $d_r$ preserve $w$; hence $d_r=0$ for $r\ge2$, so $E_2=E_\infty$. Here $w=-1,0,1$.}
    \label{fig:spectral sequence degeneration}
\end{figure}
\begin{RMK}In fact, we can similarly consider the spectral sequence for $\widehat{I}_{(\epsilon,T]}(X_\Omega,\Bbbk)$ using the $I(m)\leq p $ filtration on the complete direct sum. Then one can show that if $H^*_c(\Omega,\Bbbk)=0$, we still have the $E_n$-degeneration as stated here. However, without the homologically star-shaped assumption, which implies completeness, we are not able to show that the spectral sequence converges to $H^*F^{S^1}_{(\epsilon,T]}(X_\Omega,\Bbbk)^{hS^1}$.
\end{RMK}
\begin{RMK}There exist some examples such that the spectral sequences are not degenerate at $ E_ {n-1} $ even in the case that $X_\Omega$ is (strongly) convex. Nevertheless, on $E_1$-page, the reason that $I(m)=w,\dots,w+n-1$ could potentially be geometrically interpreted as follows: the generic torus fibers of the moment map, as a family of parameterized closed orbits (of some positive Hamiltonian, which is not necessarily $H=1$), contribute to the cohomology $\wedge^\bullet\Bbbk^{n-1}\simeq H^*_{S^1}(T^n,\Bbbk)$.
\end{RMK}

\section{Applications of the chain model for toric domains}\label{section: application of the chain model}
We discuss the applications of the chain model for toric domains on symplectic capacities. In this section, we only consider the rational case $\Bbbk=\bQ$.

\subsection{\texorpdfstring{Comparison between $c_k$ and $\widehat{c}_k$ for toric domains}{}}
The goal of this section is the following:
\begin{Prop}\label{prop: toric c = hat_c in general}
  For any toric domain $X_\Omega$, we have  
  \[c_k(X_\Omega)=\widehat{c}_k(X_\Omega).\]
\end{Prop}
We prove the statement using the criterion \autoref{prop: c=hat c, direct summand criterion}. In this subsection, we always take $T\geq 0$.

We first prepare some results. First, we shall identify the class $\beta_T^{S^1}$ using the chain model.

We recall the chain model first. Let $V\subset \bR^n_{\geq0}$ be a locally closed subset. Recall that \autoref{remark: abuse notation for the chain model for locally closed sets}, we consider 
\[I_{T}^{S^1}(X_V) = \bigoplus_{m\in \bZ^n} C_c^{*+n+2I(m)} (V (T;m),\bQ) \fpw{u} \otimes \wedge ^\bullet \bQ ^n.\]
Using this notation, we write $I^{S^1}_{T}(\bC^n)$ for $V=\bR^n_{\geq0}$ and $O^{S^1}_{T}(X_\Omega)$ for $V=U_\Omega$. We say a class $\alpha u^r\otimes e_J$ is at $(m,u^r,e_J)$ if $\alpha$ is in the $C_c^{*+n+2I(m)} (V(T;m),\bQ)$ summand.

As explained in \autoref{remark: in and out for chain model}, the closed inclusion $U_\Omega\subset \bR^n_{\geq0}$ induces the chain map (in fact, before completion)
\begin{equation}\label{equation: chain map in toric capacity comparison}
    I^{S^1}_{T}(\bC^n) \to O^{S^1}_{T}(X_\Omega).
\end{equation}

\begin{Lemma}\label{lemma: chain reprensent beta class for toric}Let $K_T=U_\Omega\cap \Delta_T$, the chain map \autoref{equation: chain map in toric capacity comparison} induced by the closed inclusion identifies the cocycle 
\[z_T=1_{K_T} ~\text{ at }~(m,u^r,e_J)=((-1,\cdots,-1),1,e_{\varnothing})\]
as a representative of $\beta_T^{S^1}$ up to a non-zero rational (which we may simply regard as $1$) after completion. Here, we take $z_T=0$ if $K_T=\varnothing$.    
\end{Lemma}
\begin{proof}It is known that $I^{S^1}_{T}(\bC^n)\simeq \bQ\fpw{u}[-n]$ abstractly. Using the chain model $I^{S^1}_{T}(\bC^n)$, we precisely identify (up to a nonzero rational scalar, which we assume to be $1$) the generator of $1\in \bQ\fpw{u}[-n]$ with the constant cochain $1_{\Delta_T}$ at $(m,u^r,e_J)=((-1,\cdots,-1),1,e_{\varnothing})$, where $\Delta_T=\{\zeta\in \bR^n_{\geq0}\mid \sum \zeta_i\leq T\}$. (We may also see \autoref{prop: concave computation} for a simpler discussion since $\bC^n$ is concave.)

The chain map \autoref{equation: chain map in toric capacity comparison} induced by the closed inclusion therefore identifies the constant cochain $1_{\Delta_T}$ to the constant cochain
\[z_T=1_{K_T} ~\text{ at }~(m,u^r,e_J)=((-1,\cdots,-1),1,e_{\varnothing}),\]
which by definition is a representative of $\beta_T^{S^1}$.
\end{proof}

Next, we shall fix a subcomplex of $O^{S^1}_{T}(X_\Omega)$, and locate the class $\beta_T^{S^1}$ more precisely. We state the following lemma in a slightly more general situation.
\begin{Lemma}\label{lemma: negative speed decomposition for chain model of toric}For $V\subset \bR^n_{\geq0}$ a locally closed subset. The chain complex $I^{S^1}_{T}(X_V)$ has a direct sum decomposition  $I^{S^1}_{T}(X_V)=I_T^-\oplus I_T'$ over $\bQ\fpw{u}$ such that $I_T^-$ is the subcomplex of $I^{S^1}_{T}(X_V)$ consisting of chains $\alpha u^r\otimes e_J$ that are at $(m,u^r,e_J)$ with $m_i<1_{J}(i)$ for all $i=1,\dots, n$ (here, $1_J$ is the characteristic function of the set $J\subset \{1,\dots,n\}$). In particular, $I_T^-$ is $u$-adic complete. 
\end{Lemma}
\begin{proof}Recall the rotation sheaf $\cE$ in \autoref{prop:chiu result on loop sheaf equiv}. Notice that as sheaves of $S^1$-modules, the extension class there between $z\leq 0$ and $z> 0$ part is $0$, so we have the following decomposition of sheaves of $S^1$-modules
\[\cE_-\coloneqq \cE\otimes 1_{(-\infty,0]},~\cE_+\coloneqq \cE\otimes 1_{(0,\infty)}, \cE=\cE_-\oplus \cE_+.\]
Using this decomposition in $\Fun(bS^1,\Sh(\bR))$, we write
\[\cE^{\boxtimes n}=\cE_-^{\boxtimes n}\oplus \mathcal{N},\quad \mathcal{N}\coloneqq \bigoplus_{d\in \{-,+\}^n, d\neq (-,\cdots,-)} \cE^\boxtimes_{d_i}. \] 

Notice that both $\cE_-^{\boxtimes n}, \mathcal{N}$ are constructible with respect to the lattice stratification. Then we can run the whole machine of \autoref{section: chain model overall}. In particular, we first use \autoref{prop: trace for toric domain derived final} to see that 
\[F^{S^1}(X_\Omega,\bQ)_T= F^-_T\oplus F'_T\]
for
\begin{align*}
    F^-_T=&\int^{m\in\rarbZ^n} \Gamma_c({ V(T;m)} ;\Bbbk  )\otimes \cE^n(m) \otimes \chi_{\cE_-^{\boxtimes n}}(m)  \\
    F'_T=&\int^{m\in\rarbZ^n} \Gamma_c({ V(T;m)} ;\Bbbk  )\otimes\cE^n(m)\otimes \chi_{\mathcal{N}}(m),
\end{align*}
where $\chi_{\cE^n_-}:\rarbZ^n\to \Fun(BS^1,D(\bQ))$ means an indicator functor that takes the value $\bQ$ for $m \in \rarbZ^n_{<0}$ and takes the value $0$ otherwise and $\chi_{\mathcal{N}}$ for a similar formula.

Then we may use \autoref{coro: coend computing in chain} to represent $F^-_T$ and $F'_T$. Here, we first denote $I_T^-$ and $I_T'$ the corresponding usual direct sum chain complex with $I^{S^1}_{T}(X_V)=I^-_T\oplus I'_T$. Then we have $F^-_T\simeq \widehat I^-_T$ and $F'_T\simeq \widehat I'_T$. This is the complex $I^-_T$ we need.

To see the description of $I_T^-$, we keep track of \autoref{coro: coend computing in chain}. In fact, the bifunctor $B(m,l)=\Gamma_c({ \Omega(T;m)} ;\Bbbk  )\otimes \cE^n(l) \otimes \chi_{\cE_-^{\boxtimes n}}(l)$, and the component at wedge degree $|J|$ in the chain model is given by $B(m,m-\sum_i 1_{J(i)}e_{i})$. Then it follows that 
\[C_B^p = \bigoplus_{m\in \bZ^n,|J|=p}\Gamma_c(V{T;m}; \Bbbk) \otimes \cE^n(m- \Sigma_i1_{J(i)}e_i) \otimes \chi_{\cE_-^{\boxtimes n}}(m- \Sigma_i1_{J(i)}e_i)\]
and hence $I_T^-$ consists of chains at  $(m,u^r,e_J)$ with $m_i-1_{J}(i)<0$ for all $i$. 

To see the completeness of $I_T^-$, we notice that the condition $m_i-1_{J}(i)<0$ for all $i$ implies that
\[I(m)-|J|\leq -n.\]
Therefore, for $\alpha u^r\otimes e_J\in I_T^-$, we have
\[|\alpha u^r\otimes e_J|=|\alpha|_{sing}+2r+|J|-n-2I(m)\geq 0. \]
So we conclude by \autoref{lemma: sufficient condition for completeness}.
\end{proof}

\begin{proof}[Proof of \autoref{prop: toric c = hat_c in general}] 
We work on the incomplete complex $O^{S^1}_{T}(X_\Omega)$ first. Then we write the decomposition of \autoref{lemma: negative speed decomposition for chain model of toric} (for $V=U_\Omega$) as $O^{S^1}_{T}(X_\Omega)=O^-_T\oplus O'_T$

We consider the chain-level weight $w(\alpha u^r\otimes e_J)=I(m)-r-|J|$ (where $\alpha\in C_c^{*} (U_\Omega(T;m),\bQ)$) we used in \autoref{prop: degeneration of SS}.  Using this weight, we similarly decompose $O^{S^1}_{T}(X_\Omega)$ (be careful that the weight changes for the $u$-action, so it is only a $\bQ$-decomposition):
\[O(w)=\bigoplus_{\substack{m\in \bZ^n,r\geq 0,n\geq j\geq 0,\\I(m)-r-j=w}}C_c^*({U_\Omega(T;m)} ;\bQ )  u^r \otimes\wedge^j\bQ^n,\qquad O^{S^1}_{T}(X_\Omega)=\bigoplus_w O(w).\]
In particular, by the formula of $z_T$, we see that $w(z_T)=-n$.

Given that $u^{k-1}[z_T] \neq 0$, we want to use the weight decomposition to exclude the case $u^{k-1}[z_T]\in u^{k}H^*O^{S^1}_{T}(X_\Omega)$. We shall notice that $w(z_T u^{k-1})=1-k-n$ for $k\geq 1$. Therefore, one may try to show that $H^*O(w)=0$ for $w>-n$. However, this is generally not true. A trickier thing is that we eventually need to work in the completion $\widehat O^{S^1}_{T}(X_\Omega)$, but $O(w)$ is not a $\bQ\fpw{u}$-submodule, so there is no way to discuss the weight decomposition of the entire $\widehat O^{S^1}_{T}(X_\Omega)$ after completion.

However, notice that $z_T$ is at the place $((-1,\cdots,-1),1,e_{\varnothing})$. Therefore, $z_T\in O_T^-$. It motivates us to consider $O_T^-$ in the decomposition (after completion): let $O_T^-(w)\coloneqq O_T^-\cap O(w)$, we have
\[F^{S^1,out}(X_\Omega,\bQ)^{hS^1}\simeq \widehat{O}^{S^1}_{T}(X_\Omega)=O_T^-\oplus \widehat O'_T, \quad O_T^-=\bigoplus_w O_T^-(w),\]
where the first is a decomposition of complete $\bQ\fpw{u}$-modules, and the second is a decomposition of $\bQ$-modules. Here, we notice that $u:O_T^-(w) \to O_T^-(w-1)$. 

On the other hand, by definition of $O_T^-$, for a vector $\alpha u^r\otimes e_J$ inside, we have
\[w(\alpha u^r\otimes e_J)=I(m)-r-|J|=\sum_i (m_i-1_J(i))-r\leq -n-r\leq -n.\]
Therefore, we have $O_T^-(w)=0$ for $w>-n$. 

Now, suppose $u^{k-1}\beta_T^{S^1}=[z_T u^{k-1}]^\wedge ~\in u^k H^*F^{S^1,out}(X_\Omega,\bQ)^{hS^1}$. We must have 
\[u^{k-1}\beta_T^{S^1}\in u^k \left(H^*(F^{S^1,out}(X_\Omega,\bQ)^{hS^1})\cap  H^*( O^-_T(1-k-n))\right)\subset H^*O^-_T(1-n)=0.\] 

Consequently, we have shown that for a toric domain, 
\[u^{k-1}\beta_T^{S^1}\neq 0\iff u^{k-1}\beta_T^{S^1} \notin u^k H^*F^{S^1,out}(X_\Omega,\bQ)^{hS^1},\]
and then we conclude the equality of two capacities by \autoref{prop: c=hat c, direct summand criterion}.
\end{proof}

\subsection{Weakly convex and concave toric domains, and toric annulus}\label{section:convex and concave toric domains}

In this section, we give a computer-friendly formula for $F^{S^1}(X_\Omega,\bQ)_T^{hS^1}$ for weakly convex and concave toric domains using \autoref{thm: chain model toric domain}. Consequently, we prove the Weyl law for all toric domains. We also explain a result concerning toric annuli.

We prepare some notation in advance. \textbf{Warning}: in this section, we will carefully distinguish derived and underived tensor products!

In this section, we consider the $\bZ$-graded ring
\[R=\bQ[x_1,\cdots, x_n],\quad |x_i|=2.\]
Write $x^m = x_1^{m_1} \dots x_n^{m_n}$. For the monomial $x^m\in R$, the degree is $2I(m)\in \bZ$. We always regard $\bQ\fpw{u}$ as an $\bZ$-graded $R$-module via $x_i\mapsto u$. In the following, all modules are $\bZ$-graded unless specified.

For a finitely generated (underived) $R$-module $M$, we consider its Koszul resolution
\[K(M)=M\otimes R\otimes\wedge^\bullet\bQ^n\]
with the differential
\begin{equation}\label{equation: differential Koszul}d=\sum_{i=1}^n  (x_i - t_i)\otimes \iota_{e_i}\end{equation}
where $x_i:M\otimes R \to M\otimes R, ~x_i(m\otimes r)=m\otimes x_ir$ and  $t_i:M\otimes R \to M\otimes R, ~t_i(m\otimes r)=x_i\cdot m\otimes r$. 
By homological algebra, this is a free resolution of $M$. Therefore, it computes the derived tensor product: for any other $R$-module, we have
\begin{equation}\label{equatiuon: definition of derived tensor over polynomial using Koszul}
    M\otimes_{R}^{\mathbf L}N \coloneqq (K(M) \otimes_R N,d\otimes_R \id_N ).
\end{equation}

In this article, we simply rewrite it as a definition of the derived tensor product.
\begin{RMK}
    The formula is similar to the one in \autoref{thm: chain model toric domain}. The reason is deep. In fact, to prove that formula, we use that $\cE^{\boxtimes n}$ can be regarded as an object in $\Fun(\rarbZ^n,D(\bQ\fpw{u}))$. However, the category is equivalent to the derived category $D_{gr}(\bQ\fpw{u}[x_1,\cdots, x_n])$ by the Rees construction. In this sense, \autoref{coro: coend computing in chain} is simply computing the coend using the Koszul complex under the above-mentioned equivalence. In some sense, our computation here is tautological from this point of view; however, the advantage is that a computer algebra system is mature for computing the Koszul resolution of polynomial rings.
\end{RMK}
\subsubsection{Concave toric domains}
We can see that our capacities for concave toric domains agree with the Gutt--Hutchings capacities \cite[Theorem 1.14]{gutthutchings2018capacities}. Note that our notion of concave toric domains is slightly more general than theirs (where they always assume $\Omega$ to be bounded).

\begin{Prop}\label{prop: concave computation}Let $\Omega \subset \bR^n_{\geq0}$ be a concave open set, namely $U_{\Omega}=\mathbb R_{\geq0}^n\setminus\Omega$ is convex. We assume further that $U_\Omega+\bR^n_{\geq 0} \subset U_\Omega$ (for example, when $\Omega$ is bounded).
We set $\nu=-m$ and $e=(1,1,\cdots,1)$. Let 
\[\underline{h}_\Omega(\nu)=\inf_{\zeta\in {U_{\Omega}}}\langle \nu,\zeta\rangle \] and 
\[J_{\Omega,T}=(x^{\nu-e} \mid \nu\in \bZ^n_{>0},h_\Omega(\nu)> T )\subset R\]
be the ideal of $R$. 

Then we have an equivalence of cofiber sequences
\begin{align*}
    &\left[F^{S^1}(X_\Omega,\bQ)_T^{hS^1}\to \bQ\fpw{u}[-n] \to F^{S^1,out}(X_\Omega,\bQ)_T^{hS^1}\right]   \\
    \simeq &\left[J_{\Omega,T}\otimes_R^{\mathbb{L}} \bQ\fpw{u} \xrightarrow{\eta}  \bQ\fpw{u} \to R/J_{\Omega,T}\otimes_R^{\mathbb{L}} \bQ\fpw{u}\right][-n],
\end{align*}
where the fundamental class is identified with $\eta_T^{S^1}(X_\Omega)(x^{\nu-e} u^r)=u^{I(\nu)+r-n}$.

In particular, we have $F^{S^1}(X_\Omega,\bQ)_T^{hS^1}$ is perfect over $\bQ\fpw{u}$ for all $T$.
\end{Prop}
\begin{proof}Here, it is more convenient to work with $F^{S^1,out}(X_\Omega,\bQ)_T^{hS^1}$ and $O_{T}(X_\Omega)$ first. See \autoref{remark: in and out for chain model} for the chain model in this case.

The closed set $U_\Omega$ is convex and $U_\Omega+\bR^n_{\geq 0} \subset U_\Omega$. Consequently, $U_\Omega$ is locally contractible, and moreover we have the following: when $\nu\notin \bZ^n_{>0}$, we then have $C^*_c(U_\Omega(T;-\nu),\bQ)\simeq 0$; when $\nu\in \bZ^n_{>0}$, $U_\Omega(T;-\nu)$ is compact and convex, and non-empty only when $\inf_{\zeta\in {U_{\Omega}}}\langle \nu,\zeta\rangle \leq T$, which implies that $C^*_c(U_\Omega(T;-\nu),\bQ)\simeq \bQ$ only when $\inf_{\zeta\in {U_{\Omega}}} \langle \nu,\zeta\rangle\leq T$. 

We set $h_\Omega(\nu)=\inf_{\zeta\in {U_{\Omega}}} \langle \nu,\zeta\rangle$. Therefore, the chain model is quasi-isomorphic to 
\[ O_T^{S^1}(X_\Omega)\simeq\bigoplus_{\nu\in \bZ^n_{>0},h_\Omega(\nu)\leq T}  \bQ\fpw{u}[n-2I(\nu)] \otimes \wedge^\bullet \bQ ^n.\]

Let $b_\nu$ be the generator in the direct summand $\bQ$ for $\nu \in \bZ^n_{>0}$. 
Then by comparing \autoref{equation: differential chain model} and \autoref{equation: differential Koszul}, the following change of variables, for $I\subset [n]$,
\[b_\nu e_I/\prod_{l}(\nu_l-1)!\prod_{i\in I} \nu_i \leftrightarrow x^{\nu-e}e_I,\]
defines a chain isomorphism
\[ \bigoplus_{\nu\in \bZ^n_{>0},h_\Omega(\nu)\leq T}\left( \bQ\fpw{u}[n-2I(\nu)] \right)\otimes \wedge^\bullet \bQ ^n\cong R/J_{\Omega,T}\otimes_R^{\mathbb{L}} \bQ\fpw{u}\]
by the definition in \autoref{equatiuon: definition of derived tensor over polynomial using Koszul}.

Here, $R/J_{\Omega,T}$ is a finitely generated $R$-module. But $R$ is regular, so finitely generated modules are perfect. Then $R/J_{\Omega,T}\otimes_R^{\mathbb{L}} \bQ\fpw{u}$ is already perfect as $\bQ\fpw{u}$ module. So is $J_{\Omega,T}\otimes_R^{\mathbb{L}} \bQ\fpw{u}$. 
Notice that perfect modules over $\bQ\fpw{u}$ are complete (we use that $|u|=2$ here). We regard $\bQ\fpw{u}$ as a $R=\bQ[x_1,\cdots, x_n]$ algebra via $x_i\mapsto u$, we have 
\[F_T^{S^1,out}(X_\Omega)[n] \simeq O_T^{S^1}(X_\Omega) [n]\simeq R/J_{\Omega,T}\otimes_R^{\mathbb{L}} \bQ\fpw{u} .\]

Therefore, via the cofiber sequence
\[F_T^{S^1}(X_\Omega)\to \bQ\fpw{u}[-n] \to F_T^{S^1,out}(X_\Omega),\]
we see that
\[F_T^{S^1}(X_\Omega)[n] \simeq J_{\Omega,T}\otimes_R^{\mathbb{L}} \bQ\fpw{u},\]
and the fundamental class is represented by
\begin{gather} \label{equation: eta_T concave}
   \eta_T^{S^1}(X_\Omega)=\left[F_T^{S^1}(X_\Omega)\to \bQ\fpw{u}[-n]\right]\simeq \left[J_{\Omega,T}\otimes_R^{\mathbb{L}} \bQ\fpw{u} \to R\otimes_R^{\mathbb{L}} \bQ\fpw{u}\simeq \bQ\fpw{u}\right][-n]\\
   \eta_T^{S^1}(X_\Omega)(x^\nu u^r)=u^{I(\nu)+r}, \notag
\end{gather}
induced by the ideal inclusion $J_{\Omega,T}\subset R$ and the natural action of $R$ on $\bQ\fpw{u}$.
\end{proof}
\begin{Coro}\label{prop: concave capacity}Let $\Omega \subset \bR^n_{\geq0}$ be a concave open set, namely $U_{\Omega}=\mathbb R_{\geq0}^n\setminus\Omega$ is convex. We assume further that $U_\Omega+\bR^n_{\geq 0} \subset U_\Omega$.  We have
\[{c}_k(X_{\Omega})
=
\widehat{c}_k(X_{\Omega})
=
\max_{\substack{\nu\in\bZ_{>0}^n\\ I(\nu)=k+n-1}}    
\inf_{\zeta\in {U_{\Omega}}}
\langle \nu,\zeta\rangle   .
\]

\end{Coro}\begin{proof}Here, we only need to compute $c_k$. By the formula of the fundamental class, we have
\[\eta_T^{S^1}(X_\Omega) \in u^kH^*C_T^{S^1}(X_\Omega) \iff k\leq d_T\coloneqq \min\{I(\nu)-n\mid \nu\in \bZ^n_{>0},\quad h_\Omega(\nu)>T\}.\]
Then we conclude by the following observation:
\[{c}_k(X_\Omega)=\inf\{T\geq 0\mid d_T\geq k\}.\qedhere\]  
\end{proof}
\subsubsection{Weakly convex toric domains}
Our computation of capacities of {\it weakly} convex toric domains extends the previous well known computations for convex domains (for various capacities, for instance \cite[Theorem 1.6]{gutthutchings2018capacities} and \cite[Theorem 3.7]{Capacities2021}). This will be a key new input for the Weyl law.

\begin{Prop}\label{prop: weakly convex cohomology}Let $T\geq 0$. Suppose that $X_{\Omega}$ is a weakly convex toric domain, namely $\Omega$ is a convex open subset of $\mathbb R_{\geq 0}^n$. Let 
\[\overline{h}_\Omega(m)\coloneqq \sup_{\zeta\in \Omega}\langle \zeta,m\rangle.\]
We consider the following $\bQ\fpw{u}\otimes R$-module: its underlying $\bQ\fpw{u}$-module is (where $b_{-m}$ is the generator of the corresponding $\bQ\fpw{u}$-summand)
\[\mathcal{I}_{\Omega,T}=\bigoplus_{m\in \bZ^n, \overline{h}_\Omega(m)\leq T} \bQ\fpw{u}, \quad |b_{-m}|=2I(m),\]
and the $R$-module structure is given by (where $b_{-(m+e_i)}=0$ when $\overline{h}_\Omega(m+e_i)> T$)
\[x_i\cdot b_{-m}=b_{-(m+e_i)}-m_i u b_{-m}.\]
We write $\widehat{\mathcal{I}}_{\Omega,T}$ for its $u$-adic completion.

Then we have
\[F^{S^1}(X_\Omega,\bQ)_T^{hS^1}\simeq \widehat{\mathcal{I}}_{\Omega,T}\otimes^\mathbb{L}_R\bQ ,\]
and the fundamental class is given by (the completion of) the formula 
\[\eta_{T}^{S^1}(b_{-m}\otimes e_I)=\begin{cases}
    1,&m=0,I=\{1,\dots, n\},\\
    0,&\text{else.}
\end{cases}\]
\end{Prop}
\begin{proof}
Here, recall the chain complex $I^{S^1}_T(X_\Omega)$ before completion in \autoref{section: chain model statement}. We write $\Omega^{int}=\Omega\cap \bR^n_{>0}$. We can consider $I^{S^1}_T(X_{\Omega^{int}})$ and then \autoref{lemma: remove coordinate boundary in toric domain} (we may pick a chain model using a similar idea for the boundary part) shows that the natural morphism induced by open inclusion
\[\widehat I^{S^1}_T(X_{\Omega^{int}})\to \widehat I^{S^1}_T(X_{\Omega})\]
is a quasi-isomorphism of $\bQ\fpw{u}$-chain complexes. Therefore, we only need to consider the $I^{S^1}_T(X_{\Omega^{int}})$. 

Now, because $\Omega$ is convex open in $\bR^n_{\geq 0}$, ${\Omega^{int}}$ is convex open in $\bR^n_{>0}$. In particular, $\Omega^{int}\cong \bR^n$ and the chain map
\[C^*_c(\Omega^{int},\bQ) \to \bQ[-n]\]
induced by pairing with the (locally compact) fundamental class of $\bR^n$ is a quasi-isomorphism. 

Now, for $m\in \bZ^n$ and $T\geq 0$, we have two chain maps
\[ C^*_c(\Omega^{int}(T;-m),\bQ)\leftarrow  C^*_c(\Omega^{int},\bQ) \to \bQ[-n],\]
such that, if $\overline{h}_\Omega(m)\leq  T$, then the left-arrow is the identity since $\Omega^{int}(T;-m)=\Omega^{int}$, and, if  $\overline{h}_\Omega(m) > T$, then $C^*_c(\Omega^{int}(T;-m),\bQ)$ is acyclic.

We write $b_{-m}$ with the $\bQ$-generator of the chain complex $\bQ[-n+n+2I(-m)]$, then we have 
\[\bigoplus_{\overline{h}_\Omega(m)\leq  T}
\bQ[-n+n+2I(-m)]\otimes \bQ\fpw{u} \otimes \wedge^\bullet \bQ ^n = {\mathcal{I}}_{\Omega,T}\otimes^\mathbb{L}_R\bQ\]
by the definition in \autoref{equatiuon: definition of derived tensor over polynomial using Koszul}.
Taking the sum over all $m\in \mathbb{Z}^n_{>0}$, we  obtain a chain map
\[I^{S^1}_T(X_{\Omega^{int}})   \to {\mathcal{I}}_{\Omega,T}\otimes^\mathbb{L}_R\bQ ,\]
which is induced by identity for $m$-components with $\overline{h}_\Omega(m)\leq  T$, and $0$ otherwise. And it induces an equivalence in $D(\bQ\fpw{u})$:
\[I^{S^1}_T(X_{\Omega^{int}})  \simeq {\mathcal{I}}_{\Omega,T}\otimes^\mathbb{L}_R\bQ.\]

Here, we notice that $\bQ=R/(x_1,\dots,x_n)$ admits a finite free resolution, so the $u$-adic completion commutes with the tensor $\otimes^\mathbb{L}_R\bQ$. Then we have
\[F^{S^1}(X_{\Omega},\bQ)_T^{hS^1}\simeq \widehat{I}^{S^1}_T(X_{\Omega})\simeq \widehat{I}^{S^1}_T(X_{\Omega^{int}})  \simeq \widehat{\mathcal{I}}_{\Omega,T}\otimes^\mathbb{L}_R\bQ.\]

The formula for the fundamental class can be verified in the following way: we notice that the computation already gives a chain morphism $I^{S^1}_T(X_{\bR^n_{>0}})  \to {\mathcal{I}}_{\bR^n,T}\otimes^\mathbb{L}_R\bQ = \bQ\fpw{u}[-n]$ where the latter complex is generated by $b_0\otimes e_{\{1,2,\cdots,n\}}$, so the fundamental class is computed by completion of
\[\begin{tikzcd}
{I}^{S^1}_T(X_{\Omega}) \arrow[d] & {I}^{S^1}_T(X_{\Omega^{int}}) \arrow[r,"\simeq"] \arrow[d] \arrow[l,"\simeq" above] & {{\mathcal{I}}_{\Omega,T}\otimes^\mathbb{L}_R\bQ} \arrow[d, dashed] &                                                  \\
I^{S^1}_T(\bC^n)                  & I^{S^1}_T(X_{\bR^n_{>0}})  \arrow[r] \arrow[l]                 & {{\mathcal{I}}_{\bR^n,T}\otimes^\mathbb{L}_R\bQ } \arrow[r, equal]    & {\bQ\fpw{u}[-n]}
\end{tikzcd}\]
where the two rows are quasi-isomorphisms, and the columns all represent the fundamental class. In particular, we can compute the fundamental class using the dashed column after completion.
\end{proof}

\begin{Coro}\label{coro: capacity weakly convex toric}Suppose that $X_{\Omega}$ is a weakly convex toric domain, namely $\Omega$ is a convex open subset of $\mathbb R_{\geq 0}^n$. We have
\[
 {c}_k(X_{\Omega})=\widehat{c}_k(X_{\Omega})
=
\min_{\substack{m\in\bZ_{\geq 0}^n\\ I(m)=k}} \sup_{\zeta\in \Omega}\langle m,\zeta\rangle,
\] \end{Coro}
\begin{proof}Here, we choose to compute $\widehat{c}_k$ instead. Let $d_T=\max\{I(m)\mid m\in \bZ^n_{\geq 0},\overline{h}_\Omega(m)\leq T\}$. We set $m!\coloneqq \prod_i m_i!$. 
We show that $I_{X_{\Omega}}(T)=d_T$, i.e.
\[\im (H^*\eta_T^{S^1}(X_\Omega, \bQ))=(u^{d_T}).\]

On the one hand, we pick a generic $z=\sum_{\overline{h}_\Omega(m)\leq T}a_mb_{-m} \otimes e_{\{1,\cdots,n\}}+z_{<n}$, where $z_{<n} \in \widehat{\mathcal{I}}_{\Omega,T} \otimes R \otimes \wedge^{<n}\bQ^n$. If $z$ is closed in $\widehat{\mathcal{I}}_{\Omega,T}\otimes^\mathbb{L}_R\bQ$, then the Koszul differential computes that
\[0=dz=\sum_{1\leq i\leq n,\overline{h}_\Omega(m)\leq T} (-1)^i(a_{m-e_i}-m_iua_m)b_{-m} \otimes \iota_{e_i} e_{\{1,\cdots,n\}}+\text{terms of wedge degree}\neq ~n-1.\]

Then we have
\[a_{m-e_i}=m_iua_m,\quad \overline{h}_\Omega(m)\leq T,\quad m_i>0,\]
and we can compute the evaluation of the fundamental class on $z$
\[\eta_T^{S^1}(X_\Omega, \bQ)(z)=a_0=m!u^{I(m)}a_m.\]

On the other hand, we set
\[z_T=\sum_{m\in \bZ^n_{\geq 0},~\overline{h}_\Omega(m)\leq T}\frac{u^{d_T-I(m)}}{m!}b_{-m}e_{\{1,\cdots,n\}},\]
which one can check using the previous formula to see that $z_T$ is closed in $\widehat{\mathcal{I}}_{\Omega,T}\otimes^\mathbb{L}_R\bQ$, and we can compute the evaluation of the fundamental class on $z_T$
\[\eta_T^{S^1}(X_\Omega, \bQ)(z_T)=u^{d_T}.\]

Therefore, one can verify directly 
\[  I_{X_\Omega}(T)=d_T\geq k \iff \exists\, m\in \bZ^n_{\geq 0}, I(m)=k,\overline{h}_\Omega(m)\leq T,\]
and then the formula follows. \end{proof}

\begin{Coro}Under the setting of \autoref{prop: weakly convex cohomology}, we assume further that $0\in \Omega$. Then we have that $F^{S^1}(X_\Omega,\bQ)_T^{hS^1}$ is a perfect $\bQ\fpw{u}$-module for all $T$.
\end{Coro}
\begin{proof}By \autoref{lem:perfectness}, we only need to show that, forgetting the $S^1$-action, $F (X_\Omega,\bQ)_T$ is perfect as a $\bQ$-module. 

Let $\mathcal{J}_{\Omega,T}=\mathcal{I}_{\Omega,T}/u\mathcal{I}_{\Omega,T}$ with the $R$-module structure $x_i\cdot b_{-m}=b_{-(m+e_i)}$. The corresponding chain model of $F(X_\Omega,\bQ)_T$ is 
\[\mathcal{J}_{\Omega,T}\otimes^\mathbb{L}_R\bQ=\bigoplus_{m\in \bZ^n,\overline{h}_\Omega(m)\leq T} \bQ \quad ,\quad d(b_{-m}\otimes e_J)=\sum_{j\in J}b_{-m-e_j}\otimes \iota_{e_j}e_J.\] 

Let $S = \{m \in \bZ^n \mid \overline{h}_\Omega(m)\leq T\}$. It is clear that $S-\bZ^n_{\geq 0} \subset S$. On the other hand, because $0\in \Omega$ and $\Omega$ is open, we can find a ball such that $B_\epsilon(0)\cap \bR^n_{\geq 0}\subset \Omega$, and then for all $m\in S$,
\[\epsilon m_i =\langle \epsilon e_i,m\rangle \leq \overline{h}_\Omega(m)\leq T\Rightarrow m_i\leq N\coloneqq \lf T/\epsilon \rf +1.\]
Therefore, let $a=(N,\dots,N)$, we have $m\leq a $ in $\bZ^n$.
Consequently, if we set
\[\mathfrak{a}=\bigoplus_{m\in S} \bQ x^{a-m}\subset R,\]
then $\mathfrak{a}$ is an ideal of $R$. Therefore, by the Hilbert basis theorem, we have that $\mathfrak{a}$ is finitely generated. The key observation here is that
\[\mathcal{J}_{\Omega,T}\simeq \mathfrak{a}^\vee[-2Nn]\]
as $R$-modules, where $^\vee$ here means the $\bZ^n$-graded linear dual.

Consequently, since $\bQ=R/(x_1,\dots,x_n)$ is perfect, we have
\[\mathcal{J}_{\Omega,T}\otimes^\mathbb{L}_R\bQ\simeq \mathfrak{a}^\vee[-2Nn]\otimes^\mathbb{L}_R\bQ\simeq \mathbf{R}\Hom_{R}(\bQ,\mathfrak{a}[-2Nn])^\vee.\] 
Now, both $\bQ$ and $\mathfrak{a}$ are finitely generated $R$-modules, then so is $\mathbf{R}\Hom_{R}(\bQ,\mathfrak{a})$. On the other hand, by $\bQ=R/(x_1,\dots,x_n)$, we see that $x_i$-action on $\mathbf{R}\Hom_{R}(\bQ,\mathfrak{a})$ are trivial. So, we $\mathbf{R}\Hom_{R}(\bQ,\mathfrak{a})$ is in fact a $\bQ$-module, which is also finitely generated. 
Therefore, $\mathcal{J}_{\Omega,T}\otimes^\mathbb{L}_R\bQ$ is perfect $\bQ$-module.
\end{proof}

\subsubsection{Toric annulus}
We give some examples of Liouville cobordisms rather than domains. While in Floer theory and symplectic field theory, it is known that under the dynamical convex condition, positive symplectic cohomology only depends on the contact boundary \cite{HoferWysockiZehnder,Cieliebak-Oancea-ES,Bourgeois-Oancea-17,HutchingsNelson}. Here, we provide the first sheaf theoretic explanation on the phenomenon.

\begin{Def}Consider $\Omega_-\subset \Omega  \subset \bR^n_{\geq 0}$, where $\Omega_-$ is bounded and continuously star-shaped. Write $ K=\overline{\Omega}_-\subset \Omega$. Consider $U_{\Omega \setminus K}= \bR^n_{\geq 0}\setminus (\Omega\setminus K)=K\sqcup U_{\Omega}$. In this case, we say $X_{\Omega\setminus K}$ is a toric annulus.
\end{Def}

\begin{Thm}\label{thm: toric annuli homology and capacity}Let $T\geq 0$. For a toric annulus $X_{\Omega\setminus K}$, we have the natural $S^1$-morphism 
\[F^{S^1}(X_{\Omega\setminus K},\bQ)_{T}^{hS^1} \xrightarrow{\simeq } F^{S^1} (X_{\Omega},\bQ)_{T}^{hS^1}   \]
induced by the inclusion $\Omega\setminus K \subset   \Omega$ is an equivalence.

In particular, we have for all $k\geq 1$ that 
\[c_k(X_{\Omega\setminus K})=c_k(X_{\Omega}),\]
i.e., the capacity of the toric annuli only depends on the outer toric domain.
\end{Thm}
\begin{RMK}We thank Yuan Yao for explaining to us the triviality of the result from the holomorphic curve perspective. From the perspective of symplectic field theory, obstructions to symplectic embeddings are usually defined as an obstruction to the existence of Liouville cobordisms between the positive contact boundaries. (One subtle point is whether positive symplectic cohomology only depends on the contact boundary. However, under the dynamical convex condition \cite{HoferWysockiZehnder,Cieliebak-Oancea-ES}, this is clear from symplectic field theory using neck stretching \cite{Bourgeois-Oancea-17,HutchingsNelson}.) Thus it is expected that for the type of capacities we study, having negative boundaries also does not change the capacities.
\end{RMK}
\begin{proof} The natural morphism is given in $\Fun(BS^1,\sT)$. Since the forgetting $S^1$-action functor is conservative, we only need to verify the equivalence stalk-wise and forget the $S^1$-structure.

Then, in terms of the notation in \autoref{remark: abuse notation for the chain model for locally closed sets} and \autoref{remark: non-equivariant chain model}, we can write down the following excision $\bQ$-chain complex short exact sequence by \autoref{thm: chain model toric domain}~(3)
\[0\to I_{T}(X_{\Omega\setminus K})\to I_{T}(X_{\Omega }) \to I_{T}(X_{K}) \to 0,\]
where the chain map $I_{T}(X_{\Omega\setminus K})\to I_{T}(X_{\Omega })$ represents the stalk of the natural map at action $T$.

We only show that $I_{T}(X_{K}) $ is acyclic. The strategy is similar to \autoref{lemma: remove coordinate boundary in toric domain}, but we do it in a slightly different and explicit way.

Recall that 
\[K(T;m)=K\cap \{\zeta\mid T+\langle m,\zeta\rangle\geq 0\}.\]
By definition of $K$, we have $tK\subset K$ for $t\in [0,1]$. Therefore, for $\zeta\in K(T;m)$, we have
\[T+\langle m,t\zeta\rangle=(1-t)T+t(T+\langle m,\zeta\rangle) \geq 0\]
for all $m$. In particular, we have $K(T;m)$ is star-shaped. Therefore, the chain map induced by the constant cochain 
\[\bQ \to C^*_c(K(T;m),\bQ)\]
is a quasi-isomorphism. 
Consequently, we have the following quasi-isomorphism 
\[\Big(\bigoplus_{m\in\bZ^n} \bQ v_m[n+2I(m)]\otimes \wedge^\bullet\bQ^n,d=\sum_{i=1}^n \id\otimes \iota_i \Big)\to I_{T}(X_{K}),\]
where $v_m$ are formal generators with $|v_m|=-n-2I(m)$.

Here, we recall the polynomial algebra and Koszul complex notation introduced in \autoref{section:convex and concave toric domains}.
Then we may use a change of variable
\[  \bQ v_m[n+2I(m)]  \ni v_m\leftrightarrow x^{-m} \in \bQ[x_1^\pm,\dots,x_n^\pm],\]
which induces a chain isomorphism
\[\Big(\bigoplus_{m\in\bZ^n} v_m\bQ[n+2I(m)]\otimes \wedge^\bullet\bQ^n,d=\sum_{i=1}^n \id\otimes \iota_i \Big)\cong \bQ[x_1^\pm,\dots,x_n^\pm]\otimes_{R}^{\mathbb{L}}\bQ[n].\]
In particular, we have that 
\[\bQ[x_1^\pm,\dots,x_n^\pm]\otimes_{R}^{\mathbb{L}}\bQ [n]\to I_{T}(X_{K})\]
is a quasi-isomorphism.

Now, we notice that $\bQ[x_1^\pm,\dots,x_n^\pm]$ is a flat module over $R$. Then the derived tensor product is acyclic (as the underived tensor product is zero by computation). So is $I_{T}(X_{K})$. We thereby conclude the proof.
\end{proof}
\begin{RMK}One may ask if $K(T;m)$ needs to be locally contractible. However, this is unnecessary. In fact, one can finish this argument using the construction explained in \autoref{remark: coend formula for general locally closed V} without picking a chain representative of the derived section $\Gamma^*_c(K(T;m),\bQ)$; however, we only use local contractibility to pick the singular cochain representative.
\end{RMK}
\begin{RMK}If $T<0$, and $\Omega$ is continuously star-shaped, then the same approach shows that $I_T(X_\Omega)=0$, and then we have $I_T(X_{\Omega\setminus K})\simeq I_T(X_K)[-1]$ (as well as its $\widehat I_T^{S^1}$ version). This is somehow compatible with our understanding that we shall regard the closure of $X_{\Omega\setminus K}$ as a Liouville cobordism with negative boundary $\partial K$.
\end{RMK}

For $\mathbf{a} = (a_1, \dots, a_n)$ with $a_i>0$, let $E(\mathbf{a})$ be the corresponding ellipsoid; for $t\in (0,1)$, we set $A_t(\mathbf{a}) = E(\mathbf{a}) \setminus E(t\mathbf{a})$. 

\begin{Prop}\label{prop: ellipsoid stabilization}Let $n\geq 2$ and $N\geq 0 $. With the notation above, we have that if there exists an isotopy of Hamiltonian embedding $\varphi_s:A_t(\mathbf{a})\times \bC^N \to \bC^{n+N}$ such that
\[\varphi_1(A_t(\mathbf{a})\times \bC^N)\subset B^{2n}(A)\times \bC^N,\]
then we have $\lceil k/n \rceil A\geq \widehat{c}_k(E(\mathbf{a}))$ for all $\ k$.
In particular $A\geq \widehat{c}_n(E(\mathbf{a}))$.
\end{Prop}
\begin{proof}By \autoref{thm: toric annuli homology and capacity}, we have that $c_k(E(\mathbf{a}))=c_n(A_t(\mathbf{a}))$ for all $k$. Then stabilization invariance of capacities together with the Hamiltonian invariance and monotonicity of the capacity show  that the existence of such a Hamiltonian embedding implies that 
\[\widehat{c}_k(E(\mathbf{a}))=\widehat{c}_k(A_t(\mathbf{a}))\leq \widehat{c}_k(B^{2n}(A)).\qedhere\]
\end{proof}

\subsubsection{Weyl law}
We show the Weyl law of toric domains using \autoref{prop: concave capacity} and \autoref{coro: capacity weakly convex toric}.

\begin{Thm}\label{thm: weyl law for toric}Let $\Omega\subset \mathbb R_{\geq 0}^n$ be any non-empty open set. Set $a(\Omega)\coloneqq \sup_{x\in {\Omega}}\min_i x_i $. Then we have
\[ ka(\Omega)\leq c_k(X_\Omega) \leq (k+n-1)a(\Omega). \]
In particular, 
\[\lim_{k\to \infty}\frac{\widehat{c}_k(X_\Omega)}{k}=\lim_{k\to \infty}\frac{{c}_k(X_\Omega)}{k} = a(\Omega).\]
\end{Thm}
\begin{proof}We consider the limit formula of $c_k(X_\Omega)/k$. The proof for the formula of $\widehat{c}_k(X_\Omega)/k$ is identical. The proof is the same as \cite{Irie_ChainModel} in dimension $4$. 

Suppose first that $a=a(\Omega)<\infty$.
Fix $\epsilon>0$ and choose $x\in\Omega$ such that
$\min_i x_i>a-\epsilon$.
Choose a convex relatively open neighborhood
$x\in V\subset\Omega$.
By monotonicity and the weakly convex computation,
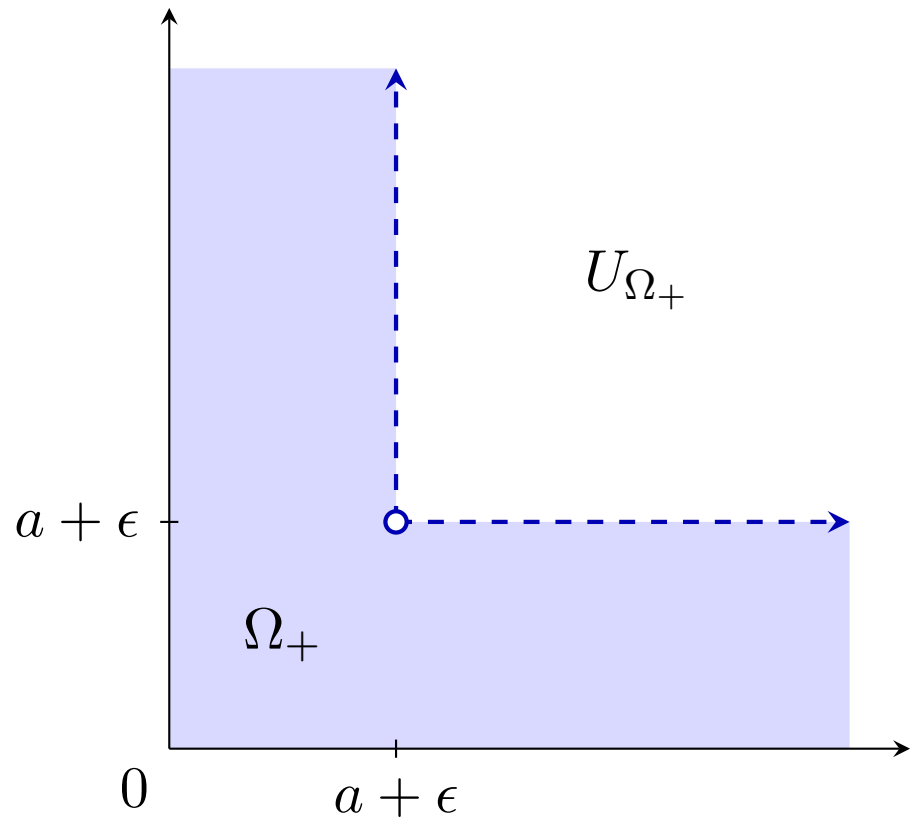
\[
  c_k(X_\Omega)
\geq   c_k(X_V)
\geq k\min_i x_i
>k(a-\epsilon).
\]
\[\begin{tikzpicture}[scale=1, >=stealth]
  \def\b{1.5}
  \def\L{4.5}

  \fill[blue!15]
    (0,0) -- (\L,0) -- (\L,\b)
    -- (\b,\b) -- (\b,\L) -- (0,\L) -- cycle;

  \draw[blue!70!black, thick, dashed, ->]
    (\b,\b) -- (\L,\b);
  \draw[blue!70!black, thick, dashed, ->]
    (\b,\b) -- (\b,\L);

  \draw[->] (0,0) -- (4.9,0) ;
  \draw[->] (0,0) -- (0,4.9) ;
  \node[below left] at (0,0) {$0$};

  \draw (\b,0.06) -- (\b,-0.06)
    node[below] {$a+\epsilon$};
  \draw (0.06,\b) -- (-0.06,\b)
    node[left] {$a+\epsilon$};

  \draw[blue!70!black, thick, fill=white]
    (\b,\b) circle (2pt);

  \node at (0.75,0.75) {$\Omega_+$};
  \node at (3.1,3.1) {$U_{\Omega_+}$};
\end{tikzpicture}\]
On the other hand, set
\[
\Omega_+
=\{\zeta\in\mathbb R_{\ge0}^n:
\min_i\zeta_i<a+\epsilon\}.
\]
Then $\Omega\subset\Omega_+$, and
$U_{\Omega_+}=[a+\epsilon,\infty)^n$ is convex
and upward closed.
Since
\[
\underline h_{\Omega_+}(\nu)
=(a+\epsilon)I(\nu),
\]
the concave computation gives
\[
  c_k(X_\Omega)
\leq c_k(X_{\Omega_+})
=(k+n-1)(a+\epsilon).
\]
Letting $\epsilon\to0$ yields
\[ 
ka\leq  c_k(X_\Omega)\leq(k+n-1)a.
\]
Dividing by $k$ proves the claimed limit.

If $a(\Omega)=\infty$, the same lower bound gives
$  c_k(X_\Omega)=\infty$ for every $k\geq 1$.
\end{proof}

Motivated by the Weyl law, we consider the following two limits of sub-leading terms when $a(\Omega)<\infty$
\begin{equation}\label{equation: refined weyl limits}
      0\leq a_-(\Omega)\coloneqq \liminf_{k\to \infty} \left[ c_k(X_\Omega)-ka(\Omega) \right] \leq \limsup_{k\to \infty} \left[ c_k(X_\Omega)-ka(\Omega) \right] \eqqcolon a_+(\Omega) \leq (n-1)a(\Omega) .
\end{equation}
are Hamiltonian invariants of a toric domain $X_\Omega$, which simply follows from the Hamiltonian invariance of the capacities.

\begin{Prop}\label{prop: weyl limit for convex/concave refinement}Let $\Omega\subset \bR^n_{\geq 0}$ be a bounded open set. If $X_\Omega$ is bounded concave toric domain, then $ a_+(\Omega)=(n-1)a(\Omega)$; If $X_\Omega$ is bounded weakly convex toric domain , then $ a_-(\Omega)=0$.
\end{Prop}
\begin{proof}We prove the concave case; the weakly convex case is similar. By \autoref{prop: concave capacity}, for $\underline{h}_\Omega(\nu)=\inf_{\zeta\in {U_{\Omega}}}\langle \nu,\zeta\rangle $ and $\Delta=\{v\in \bR^n_{\geq 0} \mid I(v)=1\}$ (where we define $I(v) = \sum_{i=1}^nv_i$), we have
\[a(\Omega)=\max_{v\in \Delta} \underline{h}_\Omega(v). \]
We notice that for $v\in  \Delta \cap \partial \bR^n_{\geq 0} $, since $\Omega$ is bounded, $\underline{h}_\Omega(v)=0$. So, there exists $v_*\in  \Delta\cap \bR^n_{>0}$ such that $\underline{h}_\Omega(v_*)=a(\Omega).$ We then pick a sequence $t_j\to \infty$ of real numbers and $m_j\in \bZ^n_{\geq 0}$, by compactness of $\bR^n/\bZ^n$, such that
\[m_j-t_jv_* \to 0,\quad \underline{h}_\Omega(m_j)-I(m_j)a(\Omega)=   \underline{h}_\Omega(m_j)-I(m_j)\underline{h}_\Omega(v_*)\to 0.\]

Let $k_j=I(m_j)-n+1$, we have
\[c_{k_j}-(k_j+n-1)a(\Omega) \to 0.\qedhere\]    
\end{proof}

\begin{Coro}\label{coro: neither convex nor concave}Let $n\geq 2$ and $2\leq a_i<\infty$ for $i=1,\cdots,n$. Let
\[\Omega=\Big\{\zeta\in \prod_{i=1}^n[0,a_i)\mid \#\{i\mid \zeta_i\geq 1\}\leq 1\Big\},\]
whose toric domain is the union of polydisks
\[X_\Omega=\bigcup_{i=1}^n P(1,\dots,a_i,\dots,1).\]
Then, $X_\Omega$ is monotone, such that $X_\Omega$ is not Hamiltonian diffeomorphic to a bounded weakly convex toric domain for all $n\geq 2$, and is not Hamiltonian diffeomorphic to a bounded concave toric domain  when $n\geq 3$.\end{Coro}
\begin{proof}It is clear that $X_\Omega$ is monotone with $a(\Omega)=1$. Notice that $X_\Omega$ is a union of convex toric domains; we can use the Mayer-Vietoris property \autoref{thm: chain model toric domain}~(3) and \autoref{coro: capacity weakly convex toric} (similar to our computation in \autoref{prop: couterexample capacity estimation}) to show that
\[c_k(X_\Omega)=k+1.\]
The computation is straightforward, but the strategy is similar to what we did before, which we omit here. Consequently, by \autoref{equation: refined weyl limits}, we have 
\[a_-(\Omega) = a_+(\Omega)=1.\]
Therefore, by \autoref{prop: weyl limit for convex/concave refinement}, this cannot happen if $X_\Omega$ is Hamiltonian isotopic to a weakly convex one for all $n\geq 2$, or to a concave one for all $n\geq 3$.    
\end{proof}

We can also deduce the asymptotic behavior on Non-Disjoint Union of Cylinders (NDUC) capacities \cite[Definition 4]{GPR22}, as mentioned in the introduction.

\begin{Coro}\label{coro: NDUC}Let $C_r=\{\zeta \in \bR^n_{>=0} : \min_i \zeta_i < r\}$, and $N_n(r)=X_{C_r}$. For $X\subset \bC^{n}$, we define the \textit{Hamiltonian NDUC capacity}
\[c^{N,h}(X)=\inf\{r>0\mid \exists\,\varphi_s:X \lhook\joinrel\xrightarrow{Ham} \bR^{2n}, \quad \varphi_1(X)\subset N_n(r)\}.\]
Then for all open $\Omega\subset \bR^n_{\geq 0}$, we have
\[\lim_{k\to \infty}\frac{\widehat{c}_k(X_\Omega)}{k}=\lim_{k\to \infty}\frac{{c}_k(X_\Omega)}{k}=c^{N,h}(X_\Omega).\]    
\end{Coro}
\begin{proof}On the one hand, the proof of \autoref{thm: weyl law for toric} in \autoref{section:convex and concave toric domains} in fact shows that $X_\Omega\subset N_n(a(\Omega)+\epsilon)$ for all $\epsilon>0$. On the other hand, the Weyl law shows that for $X_\Omega\subset N_n(r)$, we have $a(\Omega)\leq a(C_r) = r$.   
\end{proof}

\subsection{A counterexample to the product formula}\label{section: product counterexample}

In this section, we explain how to construct a counterexample for the product inequality. 

Our method works well for non-smooth domains. Let us start from a non-smooth example, and then explain that smoothing gives a smooth example.
\subsubsection{Algebraic verification}

Let 
\begin{equation}\label{equation: defining function for the product example, concave factor}
q_1(\zeta_1,\zeta_2,\zeta_3)=\frac{3\zeta_1+19\zeta_2+19\zeta_3}{60},\quad q_2(\zeta_1,\zeta_2,\zeta_3)=\frac{17\zeta_1+\zeta_2+\zeta_3}{20}.    
\end{equation}

We denote
\[\Omega=\{\zeta\in \bR^3_{\geq 0}\mid \min\{q_1,q_2\}<1\}=\{q_1<1\}\cup \{q_2<1\},\]
whose supporting function is
\[\underline{h}_\Omega(a,b,c)=\min\{20a,20b,20c,a+3\min\{b,c\}\}.\]
Then we have that $U_\Omega=\bR^3_{\geq 0}\setminus \Omega=\{q_1\geq 1\}\cap \{q_2\geq 1\}$ convex. 

The associated toric domain is concave, which is in fact the union of two $ 6$-dimensional ellipsoids in $\bC^3$:
\[X_\Omega=E(20,60/19,60/19)\cup E(20/17,20,20).\]

Then we can compute its capacities using \autoref{prop: concave computation} and \autoref{prop: concave capacity}, where we have
\begin{align*}c_i(X_\Omega)&=\max_{a+b+c=i+2}\min_{0\leq t\leq 3}\left<(a, b, c), (1, t, 3-t)\right>\\
&=\max_{a+b+c=i+2}(a+3\min(b,c)), \quad i=1,2,3,4,5,\end{align*}
or precisely, as $a, b, c\in \bZ_{>0}$,
\[c_1=4,\,c_2=5,c_3=7,c_4=8,c_5=10.\]

On the other hand, we consider the disk $D(2)=\{z_4\in \bC\mid \pi|z_4|^2<2\}$, which is also a concave toric domain whose capacities are $c_j(D (2))=2j$. Consequently, we have
\[\min_{i+j=5}(c_i(X_\Omega)+c_j(D(2)))=10.\]

\begin{figure}[htbp]
    \centering
\begin{tikzpicture}[
  x={(0.30cm,-0.12cm)},
  y={(-0.22cm,-0.16cm)},
  z={(0cm,0.95cm)},
  font=\small,
  axis/.style={-{Stealth[length=4pt]},line width=0.5pt},
  edge/.style={line width=0.75pt},
  hidden edge/.style={gray!65,densely dashed,line width=0.45pt},
  top edge/.style={gray!85,dash pattern=on 3pt off 2pt,line width=0.65pt},
  every node/.style={inner sep=2pt}
]
  \coordinate (O)  at (0,0,0);
  \coordinate (A)  at (20,0,0);
  \coordinate (C)  at (3,1,0);
  \coordinate (B)  at (0,20,0);
  \coordinate (Ot) at (0,0,2);
  \coordinate (At) at (20,0,2);
  \coordinate (Ct) at (3,1,2);
  \coordinate (Bt) at (0,20,2);
  \coordinate (S)  at (2,0,0);
  \coordinate (P)  at (2,{22/3},0);
  \coordinate (St) at (2,0,2);
  \coordinate (Pt) at (2,{22/3},2);

  \fill[black!11] (B) -- (C) -- (Ct) -- (Bt) -- cycle;
  \fill[black!18] (C) -- (A) -- (At) -- (Ct) -- cycle;
  \fill[black!5] (Ot) -- (At) -- (Ct) -- (Bt) -- cycle;

  \fill[blue!55!cyan,fill opacity=0.22]
    (S) -- (P) -- (Pt) -- (St) -- cycle;

  \draw[hidden edge] (A) -- (O) -- (B);
  \draw[hidden edge] (O) -- (Ot);

  \draw[edge] (B) -- (C) -- (A);
  \draw[edge] (A) -- (At);
  \draw[edge] (B) -- (Bt);
  \draw[edge,line width=0.95pt] (C) -- (Ct);

  \draw[top edge] (Ot) -- (At) -- (Ct) -- (Bt) -- cycle;

  \draw[blue!65!black,line width=0.85pt] (P) -- (Pt);
  \draw[blue!65!black,densely dashed,line width=0.65pt]
    (P) -- (S) -- (St) -- (Pt);
  \node[blue!65!black,above left=2pt] at (2,3.7,2) {$s=2$};

  \draw[axis] (A) -- (23,0,0) node[right] {$s$};
  \draw[axis] (B) -- (0,23,0) node[below left] {$\zeta_1$};
  \draw[axis] (Ot) -- (0,0,3.2) node[above] {$\zeta_4$};
  \node[below=3pt] at (A) {$20$};
  \node[below right=2pt] at (B) {$20$};
  \draw (Ot) -- ++(-0.12cm,0cm) node[left] {$2$};
  \node[gray!70,above left=1pt,font=\scriptsize] at (O) {$0$};

  \node at (1.5,10.5,1) {$q_1=1$};
  \node at (11.5,0.5,1) {$q_2=1$};
  \fill (C) circle[radius=1.1pt];
  \draw[fill=white,line width=0.65pt] (Ct) circle[radius=1.4pt];
  \node[below=11pt,font=\footnotesize] at (C) {$(3,1,0)$};
  \node[above right=3pt,font=\footnotesize] at (Ct) {$(3,1,2)$};
  \node[anchor=west] at (3.25cm,1.7cm) {$0\le\zeta_4<2$};
\end{tikzpicture}

    \caption{The moment domain $W=\Omega \times [0,2)$ and its $s$-slice.}
    \label{fig:moment domain of product example}
\end{figure}
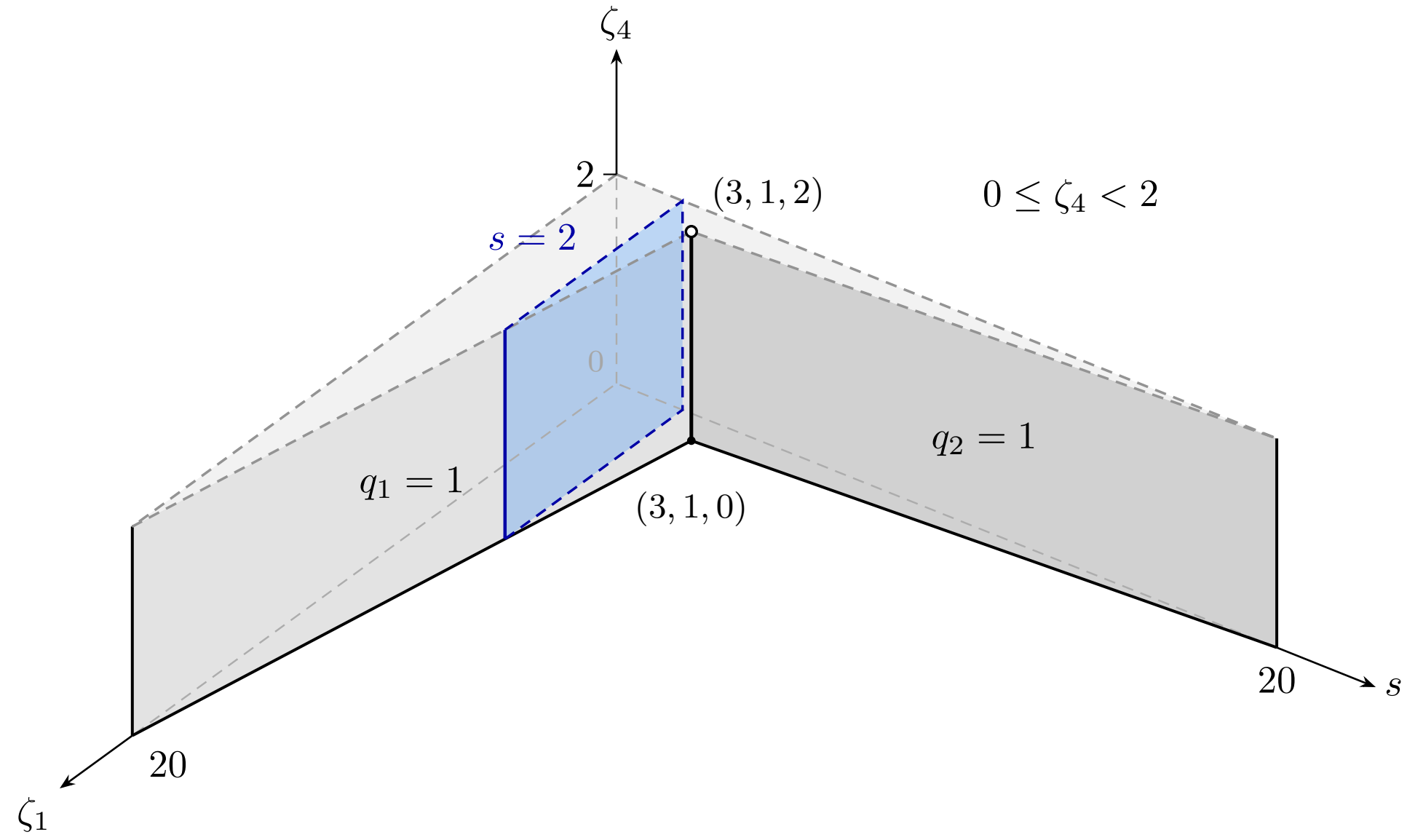

However, we will show that 
\begin{Thm}\label{prop: couterexample capacity estimation}We have
\[c_5(X_\Omega \times D(2))\leq 9.\]
  
\end{Thm}
\begin{proof}

We shall compute $F^{S^1}(X_\Omega \times D(2))_T^{hS^1}$ and show $\eta_T^{S^1}=u^5x_T$ when $9<T<10$.

We use the Mayer-Vietoris property \autoref{thm: chain model toric domain}~(3). Let (see \autoref{fig:moment domain of product example})
\[W_1=\Omega\times \bR_{\geq0},\quad W_2=\bR^3_{\geq0} \times [0,2),\quad W_{12}=W_1\cup W_2, \quad W=W_1\cap W_2.\]
Then we have $X_\Omega \times D(2)=X_W$ and one can take
\[F^{S^1}(X_W)_T^{hS^1}=\fib\left(F^{S^1}(X_{W_1}  )_T^{hS^1}\oplus  F^{S^1}(X_{W_2})_T^{hS^1}   \to F^{S^1}(X_{W_{12}})_T^{hS^1}\right).\]

All three of $W_\bullet$, $\bullet\in \{1,2,12\}$ have upward closed convex complements. Therefore, we use chain model in \autoref{prop: concave computation}. Precisely, there exist ideals 
\[J_{\bullet,T}=\big(x^{\nu-e}\mid \nu\in \bZ_{>0}^4, \, \underline{h}_{W_\bullet}(\nu)=\inf_{\zeta\in {U_{W_\bullet}}}\langle \nu,\zeta\rangle>T\big)\subset R=\bQ[x_1,x_2,x_3,x_4],\] 
such that 
\[F^{S^1}(X_{W_{\bullet}})_T^{hS^1}[4]\simeq J_{\bullet,T} \otimes_R^\mathbb{L} \bQ\fpw{u}.\]

By definition of the Koszul resolution, we have 
\[J \otimes_R^\mathbb{L} \bQ\fpw{u}= K(J)=J\otimes_\bQ\bQ\fpw{u}\otimes \Lambda^\bullet \bQ^4\]
where the basis $\{e_i\}$ of $\bQ^4$ have degree $1$, and we always take $|x_i|=|u|=2$, with the differential as in \autoref{equation: differential Koszul} 
\[d(x^m e_{i_1}\wedge\cdots \wedge e_{i_r})=\sum_{l=1}^{r }(-1)^{l-1}(x_{i_l}-u)x^m e_{i_1}\wedge\cdots \widehat{e_{i_l}}  \wedge \cdots\wedge e_{i_r} .\]

Here, we only list those ideals for $9<T<10$:
\begin{align*}
    J_{1,T}=(x_1^6,x_1^3x_2x_3,x_2^2x_3^2),\quad
    J_{2,T}=(x_4^4),\quad
    J_{12,T}=(x_1^4,x_1x_2x_3,x_2^2x_3^2,x_1^2x_4,x_2x_3x_4,x_4^2).
\end{align*}

Therefore, we have the following specific chain model, which we fix below,
\[F^{S^1}(X_W)_T^{hS^1}[4]=(K(J_{1,T})\oplus K(J_{2,T})\oplus K(J_{12,T})[-1], D)\]
where $K$ means the associated Koszul complex, and the differential is 
\[D(a,b;c)=(da,db;a-b-dc)\]
where $d$ means the internal differential of the Koszul complex.

Similarly, to compute the fundamental class, we can make an identification
\[\bQ\fpw{u} \simeq (\bQ\fpw{u} \oplus \bQ\fpw{u} \oplus \bQ\fpw{u}[-1], D),\]
where the differential on the right hand side is $D(a, b; c) = (0, 0; a-b)$. For the chain map sending $a$ to $(a, a, 0)$, its quasi-inverse sends $(a, b; c)$ to $a$ (with a chain homotopy $H(a, b; c) = (0, c; 0)$). 
Under the identification, we also see that by \autoref{prop: concave computation},
\[\eta_T^{S^1}\coloneqq \eta_T^{S^1}(X_{W})[4]: F^{S^1}(X_W)_T^{hS^1}[4] \to \bQ\fpw{u}\]
is given by
\[\eta_T^{S^1}(a,b;c)=\epsilon(a),\quad \epsilon(x^m)=u^{I(m)},\quad \epsilon(x^m e_{I})=0, |I|>0.\]

Claim: there exists $x_T$ of degree $-10$ such that $\eta_T=u^5x_T$ for $9<T<10$. If so, we conclude $c_5(X_W)=c_5(X_\Omega\times D(2)) \leq 9$. Let $q:\bQ\fpw{u}  \to \bQ\fpw{u}/(u^5)$ be the quotient map. As the proof of \autoref{prop: c=hat c, direct summand criterion} explained, we only need to show that 
$q\eta_T=0$ in  $H^0\Hom_{D(\bQ\fpw{u})}(F^{S^1}(X_W)_T^{hS^1}[4],\bQ\fpw{u}/(u^5))$.

Therefore, we prove that for
\[q\eta_T: F^{S^1}(X_W)_T^{hS^1}[4] \to \bQ\fpw{u} \xrightarrow{q} \bQ\fpw{u}/(u^5),\]
there exists a degree $-1$ map
\[H: F^{S^1}(X_W)_T^{hS^1}[4] \to \bQ\fpw{u}/(u^5)\]
such that
\[q\eta_T=d_{\bQ\fpw{u}/(u^5)}H+HD=HD.\]
I.e., $q\eta_T$ is a null-homotopy via $H$. This computation is eventually finished in \autoref{lemma:counterexample nullhomotopy} below. Then we conclude that $c_5(X_\Omega \times D(2))\leq 9$.
\end{proof}

\begin{Lemma}\label{lemma:counterexample nullhomotopy}Under the notation in the previous proposition, $q\eta_T$ is null-homotopic through a homotopy $H: F^{S^1}(X_W)_T^{hS^1}[4] \to \bQ\fpw{u}/(u^5)[1]$.
\end{Lemma}
\begin{proof}
By the formula of $\eta_T$ and the set of generators in $J_{1,T}$, we see that 
\[ q\eta_T  (x_1^6)=q\eta_T(x_1^3x_2x_3)=0,\quad q\eta_T(x_2^2x_3^2)=[u^4].\]

Therefore, we only need to worry about $x_2^2x_3^2$ in the set of free $\bQ\fpw{u}$ bases of $F_T^{S^1}(X_{W})[4] $. The key observation here is that $D(x_2^2x_3^2,0;0)=(0,0;x_2^2x_3^2)$. 

Consequently, notice that the chain model we are working with is free over $\bQ\fpw{u}$, we set
\[H(0,0;x_2^2x_3^2)=[u^4]\]
and $H$ takes value $0 \in \bQ\fpw{u}/(u^5)$ on ALL other free $\bQ\fpw{u}$ bases of $F_T^{S^1}(X_W)[4] $, and then extend $H$ linearly. As a sanity check, we notice that $|(0,0;x_2^2x_3^2)|=9$ (notice there is a shift $[-1]$ and $|[u^4]|=8$), we therefore have $|H|=-1$. 

Next, we check $q\eta_T=HD$ by evaluating on all monomial generators $x^m\in J_{\bullet,T}$: 
\begin{enumerate}[fullwidth]
\item For monomial generators $x^m\in J_{1,T}$, we have
\[HD(x^m,0;0)=H(0,0;x^m)= 0 \text{ or } [u^4] = q\eta_T(x^m),\]
and for $x^m e_I$ with $|I|\geq 1$, 
\[HD(x^me_I,0;0)=H(d(x^me_I),0;x^me_I)= 0  = q\eta_T(x^me_I).\]
\item For the only monomial generators $x_4^4\in J_{2,T}$, we have 
\[HD(0,x_4^4;0)=H(0,0;-x_4^4)= 0 = q(0)=q\eta_T(0,x_4^4;0).\]

\item  For monomial generators $x\in J_{12,T}$, all free generators are $x^m e_I$.

When $I=\varnothing$, we have $HD(0,0;x^m)=H(0,0;dx^m)=H(0,0;0)$. When $|I|=1$, we have
\[HD(0,0;x^me_i)=H(0,0;-x_ix^m+ux^m)=-H(0,0;x_ix^m)+H(0,0;ux^m).\]
But for all $x^m\in J_{12,T}$, it would not happen that $x_ix^m= x_2^2x_3^2$, so the term $H(0,0;x_ix^m)$ is always zero. For the second term, we may take $x^m=x_2^2x_3^2$, and then $H(0,0;ux^m)=uH(0,0;x^m)=[u^5]=0$. When $|I|\geq 2$, we notice that $D$ decreases the degree of the wedge power, then the $|I|=1$ computation shows that $HD(0,0;x^me_I)=0$. In all cases, we have $HD (0,0; x^me_I)=0=q\eta_T(x^me_I)$.
\end{enumerate}
Then we conclude by linear extension that $q\eta_T=HD$. 
\end{proof}
\begin{RMK}
    A similar computation can be performed for $8<T<9$ to show that $c_5(X_\Omega\times D(2))\geq 9$, which we omit here. In particular, $c_5(X_\Omega \times D(2)) = 9$.
\end{RMK}
Now, we smooth $X_\Omega$ to obtain a smooth concave toric domain, and verify it is still a counterexample.
\begin{Coro}\label{coro: smooth model of the product counterexample}Let $\delta=21/20$ and
    \[j_{sm}=\frac{\delta}{\|(1/q_1,1/q_2)\|_{20}},\quad \Omega_{sm}=\{j_{sm}<1\}.\]
Then $X_{\Omega_{sm}}$ defines a bounded star-shaped concave toric domain with smooth boundary and
    \[c_5(X_{\Omega_{sm}}\times D(2))\leq 9<19/2 \leq  10\delta^{-1}\leq \min_{i+j=5} \left(c_i(X_{\Omega_{sm}})+c_j(D(2)\right). \]
\end{Coro}
\begin{proof}By definition of $\Omega_{sm}$, we have
    \[{\delta}^{-1}\Omega\subset \Omega_{sm} \subset \Omega,\]
    and $X_{\Omega_{sm}}$ defines a bounded star-shaped concave toric domain with smooth boundary. 

By monotonicity, we have
\[c_5(X_{\Omega_{sm}}\times D(2))\leq c_5(X_{\Omega}\times D(2))\leq 9.\]

On the other hand, by monotonicity and the conformal property of the capacity, we have
\[\delta^{-1} c_i(X_{\Omega}) \leq c_i(X_{\Omega_{sm}}),\]
and then
\[19/2 < 10\delta^{-1}\leq \min_{i+j=5} \left(c_i(X_{\Omega_{sm}})+c_j(D(2)\right).\qedhere\]    
\end{proof}

\subsubsection{A potential dynamics explanation}\label{subsection: dynamical explanation of the example}
Here, we present a discussion on the potential dynamical mechanism behind this example. For convenience, we explain the non-smooth version first.

Recall that 
\begin{equation}\label{equation: toric moment map in dim 4}    (\zeta_1,\zeta_2,\zeta_3,\zeta_4)=\mu(z_1,z_2,z_3,z_4)=(\pi|z_1|^2,\pi|z_2|^2,\pi|z_3|^2,\pi|z_4|^2),
\end{equation}
and by \autoref{equation: defining function for the product example, concave factor}, we can set $\Omega=\{q_1<1\}\cup \{q_2<1\}$ and  
\[W_1=\Omega\times \bR_{\geq0},\quad W_2=\bR^3_{\geq0} \times [0,2),\quad W_{12}=W_1\cup W_2, \quad W=W_1\cap W_2,\]
such that 
\[X_{W}=X_\Omega\times D(2)\subset \bC^4.\]
We also set the following continuous Hamiltonian
\begin{equation}\label{equation: continuous Hamiltonian that defines the product counterexample}
    H(z)=\max\Big\{\min(Q_1,Q_2),\frac{\zeta_4}{2}\Big\},
\end{equation}
where $Q_i=q_i\circ \mu$. Then we have $X_W=\{H<1\}$.

Here, $\overline{X_{W}}=\{H\leq 1\}$ is not a smooth manifold with boundary, but a manifold with corners. We consider its topological boundary $\partial X_{W}$. The main observation is that: inside $\partial X_{W}$, there is a particular family of ``Reeb''-orbits, which provides a non-trivial circle bundle whose Euler class then contributes to the capacity index $I_W(9)=5$.

Now, let us explain the observation more precisely.

We set 
\[Z_s=\partial X_{W}\cap \{z \in \bC^4\mid \pi|z_2|^2+\pi|z_3|^2=s\}.\]
We pick variables $(z_2,z_3)$ because the last two capacities of ellipsoids we take are equal. By the definition of $\partial X_{W}$, one can see directly that $ Z_s$ is non-empty when $ 0\leq s\leq 20$.

We first consider \textit{the topology of  }$Z_s$.

We write $S^3(s)=\{z \in \bC^2_{(z_2,z_3)} \mid \pi|z_2|^2+\pi|z_3|^2=s\}$, then
 \[Z_{s}\cong S^3(s)\times \partial (\overline{D(A_s)}\times \overline{D(2)} ) \subset \bC^2_{(z_2,z_3)} \times  \bC^2_{(z_1,z_4)},\]
    for 
\begin{equation}\label{equation: formula of a(s) in geometry of product example}
    A_s=\begin{cases}
        20-19s/3,& 0\leq s \leq 3,\\
        (20-s)/17,& 3\leq s \leq 20.
    \end{cases}
\end{equation}
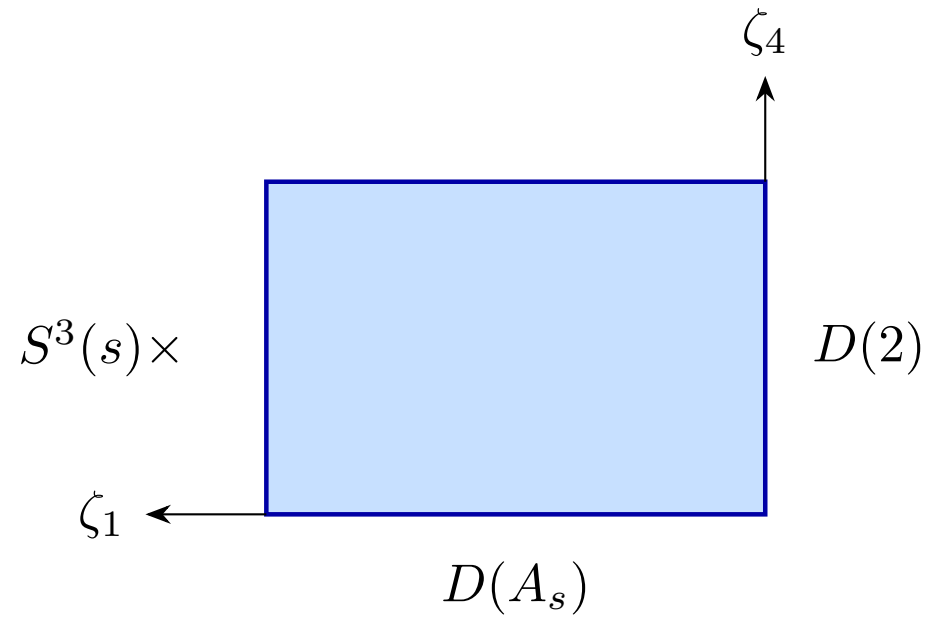
\begin{figure}[htbp]
    \centering
\begin{tikzpicture}[font=\small,>=Stealth]
    \filldraw[fill=blue!55!cyan,fill opacity=0.22,
        draw=blue!65!black,line width=0.85pt]
        (0,0) rectangle (3.3,2.2);

    \draw[->] (0,0) -- (-0.8,0) node[left] {$\zeta_1$};
    \draw[->] (3.3,2.2) -- (3.3,2.9) node[above] {$\zeta_4$};

    \node[anchor=east] at (-0.4,1.1) {$S^3(s)\times$};
    \node[below=5pt] at (1.65,0) {$D(A_s)$};
    \node[right=5pt] at (3.3,1.1) {$D(2)$};
\end{tikzpicture}
    \caption{The $s$-slice of the moment region. Compare to \autoref{fig:moment domain of product example}.}
    \label{fig:s-slice moment domain of product example}
\end{figure}

In particular, for $0<s<20$ (we always assume it in the discussion below), we have
\[Z_s\cong S^3\times S^3.\]
The circle $S^1$ freely acts on $Z_s$ for all $0<s<20$, via weight $(1,2,2,1)$
\[e^{2\pi i t} \cdot (z_1,z_2,z_3,z_4)= (e^{2\pi i t} z_1,e^{4\pi i t} z_2,e^{4\pi i t} z_3,e^{2\pi i t} z_4).\]

After a proper change of coordinates, the quotient $p_s:Z_s\to B_s=Z_s/S^1$ is isomorphic to a non-trivial $S^1$-bundle
\[Z_s\cong S^3(s)\times S^3\xrightarrow{\id\times h} B_s:=S^3 \times S^2,\]
where $h$ is the Hopf bundle, which is non-trivial.

Here, we in particular notice that the Heegaard torus (of the second $S^3$ factor) inside $Z_s$ is given by
\[N_s=S^3(s)\times \partial{D(A_s)}\times \partial{D(2)} \subset \partial X_{W_{12}}\cap \{z \mid \pi|z_2|^2+\pi|z_3|^2=s\}.\]
It means that the non-trivial topology of $Z_s$ comes from the gluing of $\partial X_{W}$ from $\partial X_{W_{1}}$ and $\partial X_{W_{2}}$ along $\partial X_{W_{12}}$. This is compatible with our Mayer-Vietoris computation in the chain model.

Next, we consider the \textit{dynamics of $Z_s$} and explain why the number $9$ comes in.

The $S^1$-action above indicates that $Z_s$ is a union of parameterized $S^1$-orbits 
\[\gamma_z(t)= (e^{2\pi i t} z_1,e^{4\pi i t} z_2,e^{4\pi i t} z_3,e^{2\pi i t} z_4)\] for $z=(z_1,z_2,z_3,z_4)\in Z_s$, parameterized by the moduli $[z]\in B_s=Z_s/S^1$.

For each of the orbits $\gamma_z$ with $z\in Z_s$, its action functional is 
\[\mathcal{A}(\gamma_z)=\int_{\gamma_z} \lambda = \pi|z_1|^2+2\pi|z_2|^2+2\pi|z_3|^2+\pi|z_4|^2=\pi|z_1|^2+\pi|z_4|^2+2s,\]
where $\lambda$ is the radial Liouville form on $\bR^{2n}\cong\bC^n$.

The maximum of the action is
\[\max_{z\in Z_s} \mathcal{A}(\gamma_z)= A_s+2+2s,\]
which is attained at
\[N_s=S^3(s)\times \partial{D(A_s)}\times \partial{D(2)} .\]
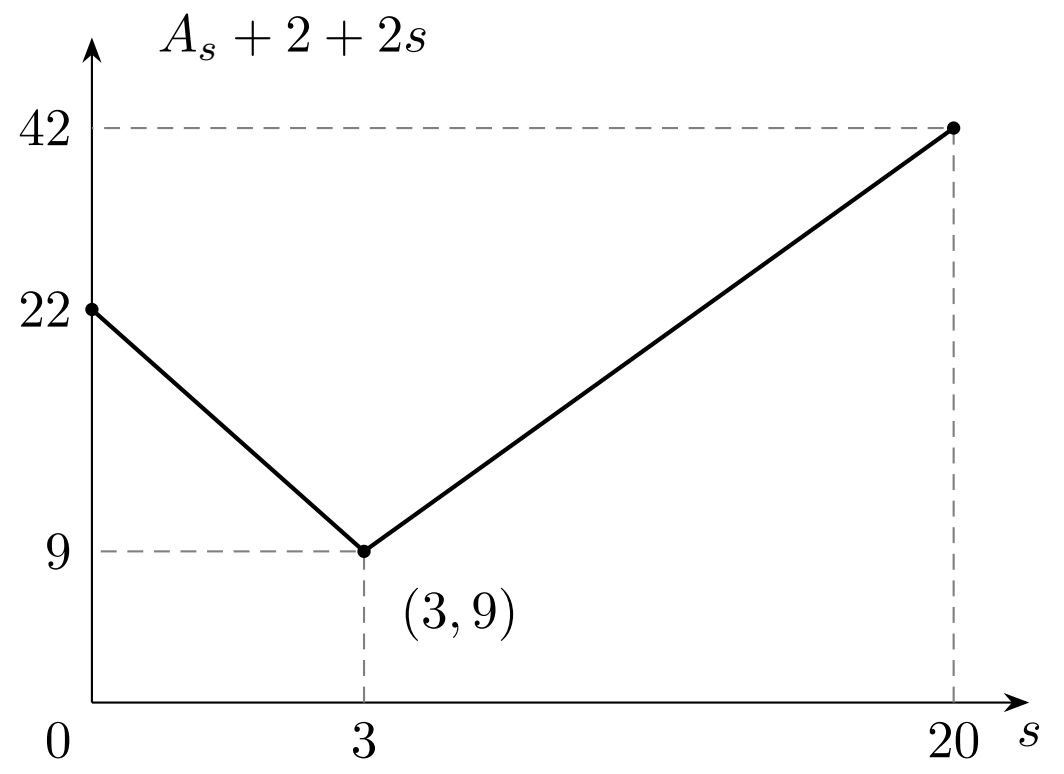
\[\begin{tikzpicture}[font=\small, >=Stealth]
    \coordinate (A) at (0,2.6);
    \coordinate (C) at (1.8,1);
    \coordinate (B) at (5.7,3.8);

    \draw[->] (0,0) -- (6.2,0) node[below] {$s$};
    \draw[->] (0,0) -- (0,4.4)
        node[right=3mm] {$A_s+2+2s$};

    \draw[gray,densely dashed]
        (1.8,0) -- (C) -- (0,1);
    \draw[gray,densely dashed]
        (5.7,0) -- (B) -- (0,3.8);

    \draw[thick] (A) -- (C) -- (B);
    \foreach \p in {A,C,B}
        \fill (\p) circle[radius=1.25pt];

    \node[below left] at (0,0) {$0$};
    \node[left] at (0,1) {$9$};
    \node[left] at (A) {$22$};
    \node[left] at (0,3.8) {$42$};
    \node[below] at (1.8,0) {$3$};
    \node[below] at (5.7,0) {$20$};
    \node[below right=3pt] at (C) {$(3,9)$};
\end{tikzpicture}\]

Therefore,  by the formula of $A_s$, we have the min-max formula
\[\min_{s\in [0,20]}\max_{z\in Z_s} \mathcal{A}(\gamma_z)= 9,\]
attained at $s=3$. Numerically, the number $9$ comes already as the minimal.

Moreover, the family of orbits on the Heegaard torus $N_3\subset Z_3\to B_3$ is special in the sense that $N_3$ is the only family of ``Reeb''/characteristic-orbits in a particular non-smooth sense, which we will not expand here but give an ad hoc explanation: Recall the continuous Hamiltonian $H$ in \autoref{equation: continuous Hamiltonian that defines the product counterexample},
whose Clarke differential \cite[Theorem 10.31]{RockafellarWets1998} at $z\in N_s$ is
\[\partial_C H(z)=\begin{cases}\operatorname{conv}\{dQ_1(z), \pi\Re(\overline{z_4}dz_4)\}, & 0<s<3,\\
\operatorname{conv}\{dQ_1(z),dQ_2(z),\pi\Re(\overline{z_4}dz_4)\},& s=3,\\
\operatorname{conv}\{ dQ_2(z),\pi\Re(\overline{z_4}dz_4)\},& 3<s<20.\end{cases}\]
The normal cone of $\{H=1\}$ at $z\in N_s$ is $\bR_>\partial_C H(z)$ by \cite[Theorem 6.14]{RockafellarWets1998}.

Let $F(z)=\pi|z_1|^2+2\pi|z_2|^2+2\pi|z_3|^2+\pi|z_4|^2$, which is also the Hamiltonian of the circle action of weight $(1,2,2,1)$. 

Then one can check that $F$ is positive/characteristic at $N_s$  in the sense that\footnote{It is clear that the condition is compatible with the smooth case since $\partial_C H(z)=\{dH(z)\}$ then.}
\[dF(z)\in \bR_>\partial_C H(z),\quad z\in N_s\]
if and only if $s=3$.

Now, let us try to summarize what we have so far: We foliate $\partial X_W$ into a family of $S^1$-orbits $Z_s\cong S^3\times S^3$ whose non-trivial topology is obtained along the Heegaard torus $N_s\cong  S^3\times T^2\subset Z_s$. When $s=3$, the family $N_3$ is ``Reeb''/characteristic and has the action $9$. The family $N_3$ potentially contributes more $u$-multiplication at action level $T=9$. (and below $T<9$, the algebraic computation shows that a $u^4$ factor already exists, contributed by an $S^1$-action of weight $(2,2,2,1)$ with the initial position $(z_1,z_2,z_3,0)$).
In particular, by the topological study of $Z_3$, we see that $N_3$ essentially comes from the product structure of $X_W=X_\Omega\times D(2)$.

However, the product formula predicts that the capacity of the product $X_W=X_\Omega\times D(2)$ only records the ``Reeb dynamics'' of both factors; for example, here, $10=c_5(X_\Omega)=\min_{i+j=5}(c_i(X_\Omega)+c_j(D(2))$ indicates that a ``Reeb" orbit in $\partial X_\Omega$ is the only potential $u^5$ contributor at $T=10+\epsilon$. Our example hereby shows that the expectation is not correct.

In the end, we remark that in the smooth model \autoref{coro: smooth model of the product counterexample}, the explanation should still work after suitable perturbation of the family $N_3$ (but still a family since this perturbation is not generic enough).

\begin{RMK}\label{remark: higher dimension example formula in the end}Finally, we mention the high dimensional example presented in \autoref{remark: higher dimension example formula}. We recall it here: Fix $n\geq 2$.  Let $E(a_1,\cdots,a_{n+1})$ be an ellipsoid in $\bR^{2n+2}$ for $n\geq 2$. Let $\sigma>0$ be real and $m\geq 2$ be integer. 

We take $X_\Omega^0\subset \bR^{2n+2}$ to be 
\[E(R,a,\cdots,a)\cup E(b,R,\cdots,R) \subset \bR^{2n+2}.\]

The mechanism behind the example is the same. We match some number for the convenience of the readers. The parameter $m$ here is related to the weight of the $S^1$-action we considered in $\partial (X_\Omega^0\times D(r))$. Precisely, we consider the circle action of weight $(1,m,\dots,m,1)$ on $\bC^{n+2}=\bR^{2n+4}$, which is restricted to $\partial (X_\Omega^0\times D(r))$. The number $s=3$ corresponds to $s=\sigma$ here, and we have a pullback diagram of the orbit space
\[\begin{tikzcd}
Z_s=S^{2n-1}(s)\times \partial(D(A_s)\times D(r)) \arrow[d, "p_s"] \arrow[r, "\operatorname{proj}"] & \partial(D(A_s)\times D(r))\cong S^3 \arrow[d, "h"] \\
B_s=Z_s/S^1 \arrow[r, "\pi_s"]                                               & S^2,                                           
\end{tikzcd}\]
for $h$ the Hopf bundle and
\begin{gather*}
     Z_s=\partial (X_\Omega^0\times D(r)) ~\cap~\{\pi\sum_{j=1}^{n} |z_j|^2=s\},\\
    p_s: Z_s\to B_s=Z_s/S^1,\quad (z_0,\cdots,z_{n+1})\mapsto[z_0,\cdots,z_{n+1}],\\
    \pi_s: B_s\to S^2,\quad[z_0,\cdots,z_{n+1}]\mapsto [z_0,z_{n+1}].       
\end{gather*}
Notice that here, $\pi_s$ is a trivial $S^{2n-1}$ bundle if and only if $mn$ is even.

Therefore, before giving a formula of $a,b,r,R$ using $m,n,\sigma$, we can first estimate capacities using similar method as \autoref{prop: couterexample capacity estimation} in terms of $a,b,r,R$ where $m,n,\sigma$ are given parameters (then some constraints on $a,b,r,R$ are also needed). Therefore, to obtain a counterexample of the product formula in the current form, we only need to solve the inequalities system about $a,b,r,R$:
\begin{equation}\label{equation: inequality system for counterexample}
    \begin{cases}
&A_\sigma =  R\left(1-\frac{\sigma}{a}\right)
=b\left(1-\frac{\sigma}{R}\right)>0,\\
&nA_\sigma\leq r<\sigma,\\
& m\sigma+A_\sigma+r
<\min\{R,(mn+1)r\},\\
    \end{cases}
\end{equation}
Here, we simply pick one particular solution of the system, and then we have
\[
\widehat{c}_{mn+1}(X_\Omega^0\times D(r))
\leq m\sigma+A_\sigma+r
< 
\min_{i+j=mn+1}
\{
\widehat{c}_i(X_\Omega^0)+\widehat{c}_j(D(r))
\}.
\]

\end{RMK}

We end the section by a discussion on the implications on the capacities of symplectic $p$-products \cite[Conjecture 1.6]{HaimKislevOstrover2023}.

\begin{RMK}\label{rmk: HaimKislevOstrover}
We claim that our counterexample also disproves $\widehat{c}_k$-version of this conjecture for sufficiently large $p\gg 2$: On one hand, by definition of $p$-product \cite[Conjecture 1.6]{HaimKislevOstrover2023}, we have
\[X_\Omega\times_p D(2) \subset  X_\Omega \times D(2)\]
for $p\geq2$, then we have $\widehat{c}_5(X_\Omega \times_p D(2)) \leq 9$ for $p\geq 2$; on the other hand, by the H\"older inequality, we have
\[\min_{ i+j=5} 
\left[
\widehat{c}^i(U)^{\frac{p}{p-2}}
+
\widehat{c}^j(W)^{\frac{p}{p-2}}
\right]^{\frac{p-2}{p}}\geq 2^{-2/p} \min_{ i+j=5} 
\{
\widehat{c}^i(U) 
+
\widehat{c}^j(W) 
\}= 2^{-2/p}(9+e),\]
for $e>0$. Then for sufficiently large $p\gg 2$, we have $2^{-2/p}(9+e)>9$. Here, we will not give a precise estimation of $p$ since it depends on the particular $\Omega$. It is clear that our higher-dimensional examples \autoref{remark: higher dimension example formula in the end} also provide more counterexamples for the equality.
\end{RMK}

\appendix
\section{Sheaf quantization for complete Hamiltonian flow}

\subsection{Existence and uniqueness of sheaf quantization}\label{appendix: existence kernel}

We say a Hamiltonian $H\in C^\infty(\bR\times T^*M )$ is complete if the flow of the Hamiltonian vector field $X_H$ is defined for all time $z\in \bR$.

In \cite{Chiu-kernel2021}, Chiu shows that a sheaf quantization for a complete Hamiltonian flow exists under a technical condition. A similar statement for non-negative complete Hamiltonians was explained in \cite[Section 2.3]{Asano-Ike-peg}. Here, we give a proof for all complete Hamiltonians. The result of this appendix is true for all commutative ring spectra.

For a complete Hamiltonian $H\in C^\infty( \bR_z\times T^*M)$, we denote its Hamiltonian flow to be 
\[\varphi=\varphi_H:\bR_z\times T^*M\to T^*M.\]
We denote
\[\Lambda_H=\{(z,-H,q,-p,\varphi_z(q,p))\mid (q,p,z)\in \bR_z\times T^*M\} \subset T^*(\bR_z\times M\times M )\]
be the suspended graph of the Hamiltonian isotopy $\varphi_H^z$. We use the notion of reduced microsupports $\underline{ss}$ for sheaves in the Tamarkin category; see \cite[Remark 5.8]{Hochschild-Kuo-Shende-Zhang}.

\begin{Thm}[Tomohiro Asano\footnote{
Tomohiro Asano obtained a proof of the result in March 2024 and explained it to the first-named author; and the second-named author, without knowing it, independently obtained its proof during the preparation of the paper based on an early version of the result on positive Hamiltonians \cite{Asano-Ike-peg}. We thank Asano for sharing his draft with us.
}]\label{thm: existence of sheafquantization for all complete hamiltonian}If $H\in C^\infty(\bR_z\times T^*M)$ is a complete Hamiltonian, then there exists a unique sheaf $\cK = \cK(H) \in \sT(T^*(\bR_z\times M\times M))$ such that 
\[\underline{ss}(\cK(H))\subset \Lambda_H,\quad i_0^*\cK(H)\simeq 1_{\Delta_{M}}^\sT.\]   

The sheaf $\cK=\cK(H)\in \sT(T^*(\bR_z\times M\times M))$ is called the \textit{sheaf quantization} of $H$.
\end{Thm}
\begin{RMK}\begin{enumerate}[fullwidth]
\item By the results of \cite{GKS2012}, we know the theorem is true if $H$ is compactly supported. However, our proof for a general Hamiltonian depends on the compactly supported result. (Note that for a homogeneous Hamiltonian isotopy, the result \cite{GKS2012} has no compact support conditions. However, to work with inhomogeneous Hamiltonian isotopies in $T^*M$, one needs to lift it to a homogeneous isotopy in $T^*_{>0}(M \times \bR)$, and then extend it to $T^*(M \times \bR)$ in order to apply \cite{GKS2012}. The extension step is where one uses compact supports.)
\item Unlike \cite{GKS2012}, we prove the existence first. The idea is following: we approximate $H$ by a sequence of suitably chosen compactly supported $H_n$. Then we can approximate using the (non-conic) GKS quantization $\cK({H_n})$ (so the proof is not independent of the GKS construction, but a generalization). We shall use the Tamarkin projector to perform cutoff on $T^*M$ directly (rather than cutoff on the base direction by standard shriek extensions), which makes the proof slightly different from GKS's proof. 
Then we finish the uniqueness by a standard microlocal Morse argument.

\item Similar to \cite[Remark 3.9]{GKS2012}, we can extend the result to contractible parameter spaces other than $\bR$.
\end{enumerate}
    
\end{RMK}

We first state a version of the microlocal Morse lemma.
\begin{Lemma}\label{lemma: microlocal morse-Tamarkin}Let $J$ be a contractible manifold and take $0\in J$, we have \[  \sT_{0_{T^*J}\times T^*M }(T^*(J\times M) )\to \sT(T^*M), F\to i^*_0F\] is a categorical equivalence.    
\end{Lemma}
\begin{proof}By the Künneth formula \cite{zhang-cutoff}, the functor can be represented by 
 \[\sT_{0_{T^*J}\times T^*M }(T^*(J\times M) )\simeq\Loc(J) \otimes  \sT(T^*M) \simeq \sT(T^*M) ,\quad F\to i^*_0F,\] 
 which finishes the proof.

Alternatively, one may prove it by the usual microlocal Morse lemma, but be careful on the Tamarkin direction. In fact, carefully applying microlocal Morse lemma shows that for $F, G \in \sT_{0_{T^*J}\times T^*M }(T^*(J\times M) )$,
\[\Gamma^\sT(J, \sHom^\star(F, G)) \simeq \Gamma^\sT(0, \sHom^\star(F, G)),\]
and the essential surjectivity is easy to see.
\end{proof}

Let $J\subset \bR_z$ be open. We naturally identify $\sT(T^*(J\times M\times M))$ with $\Sh(J\times M\times M,\sT)$. The $\sT$-linear 6-functor on $\Sh(J\times M\times M,\sT)$ induces the following relative composition: For $A, B\in \Sh(J\times M\times M,\sT)$, we set for $i,j=1,2,3$ the map $q_{ij}:J\times M^3\to J\times M^2$ to be the projection to corresponding $M$-factors, we define
\[A\circ_{\sT,J} B\coloneqq (q_{13})_!
\bigl(
(q_{12})^*A
\otimes_{\sT}
(q_{23})^*B
\bigr)\in \Sh(J\times M\times M,\sT).\]

The relative composition satisfies that, for each $z\in J$ and $i_z:M^2\to \bR\times M^2 ,(x, x') \mapsto (z,x,x' )$, we have that 
\[i_z^*(A\circ_{\sT,J} B)=i_z^* A\circ_{\sT}i_z^*B.\]

\begin{Lemma}\label{lemma: relative compisition estimation}Let $H,G$ be two complete Hamiltonians, and their composition is denoted by
\[H \# G (z,q,p) \coloneqq H(z,q,p)+G(z,(\varphi^H_z)^{-1}(q,p)).\]

Let $J\subset \bR_z$ an open set and $0 \in J$. Suppose there exists sheaf quantization $\cK(H)$ and $\cK(G)$ in the sense of \autoref{thm: existence of sheafquantization for all complete hamiltonian}.
Then we have 
\[\underline{ss}(\cK(H)|_J\circ_{\sT,J} \cK(K))|_J)\subset \Lambda_{H \# G},\qquad i_0^*(\cK(H)|_J\circ_{\sT,J} \cK(G))|_J)=1_{\Delta_{M}}^\sT\]   
and hence there exists a sheaf quantization $\cK(H \# G)$ in the sense of \autoref{thm: existence of sheafquantization for all complete hamiltonian}.
\end{Lemma}
\begin{proof}By the definition of sheaf quantizations, we have $i_0^*(\cK(H))\simeq i_0^*(\cK(G))\simeq 1_{\Delta_{M}}^\sT$, then by the property of relative composition, we have
  \[i_0^*(\cK(H)\circ_{\sT,J} \cK(G))\simeq 1_{\Delta_{M}}^\sT \circ_{\sT } 1_{\Delta_{M}}^\sT=1_{\Delta_{M}}^\sT.\] 
  This verifies the second condition.

  For the first condition, this is a standard microsupport estimation result similar to \cite[Equation (1.12)]{GKS2012}. By the graph microsupport constraint, the required non-characteristic condition for the singular support estimation \cite[Equation (1.14)]{GKS2012} is true. 
\end{proof}

We now prepare the existence part of \autoref{thm: existence of sheafquantization for all complete hamiltonian}.

We take $J=(-R,R)\subset \bR_z$ for a large enough $R > 0$. We take an exhaustion $U_n$ of $T^*M$ such that $\overline{U_n}\subset U_{n+1}$ are compact. Then by compactness of $\overline{U_n}$, we have
\[C_n=\{\varphi_z(q,p)\mid (q,p,z)\in \overline{U_n\times J}\} \subset T^*M\]
is compact. 

Then we take a cutoff function $\rho_n$ equal to $1$ on a relatively compact neighborhood of $C_n$. We then take $H_n=\rho_n H \in C_c^\infty(\bR_z\times T^*M)$, which is complete, and moreover, 
\[H_n(z,\varphi_{z}^{H_n}(q,p))=H(z,\varphi_{z}^{H}(q,p)),\quad \varphi_{z}^{H_n}(q,p)=\varphi_{z}^{H}(q,p)\]
for $(q,p,z)\in \overline{U_n\times J}$; in fact, those equations are true on $ W\times   J$ for  a neighborhood $W$ of $\overline{U_n}$ since $\rho_n$ equals $1$ on a neighborhood of $C_n$.

For the compactly supported $H_n$, we can take its sheaf quantization in the Tamarkin category using the standard GKS quantization, say $\cK(H_n)$. We set their restriction on $J$ to be $\cK_n\coloneqq \cK(H_n)|_{J}$.

On the other hand, we take the cutoff kernel $P_n\coloneqq P_{U_n}\in \sT(T^*(M\times M))$ for the inclusion functor $\sT_{-U_n \times U_n}(T^*(M \times M))$. By \cite[Proposition 4.5]{zhang-cutoff}, we have an estimation on the reduced microsupport
\begin{equation}\label{equation: estimation of PU}
\underline{ss}(P_{U_n})\subset (-\overline{U_n})\times \overline{U_n} \subset T^*(M \times M). 
\end{equation}

\begin{Lemma}\label{Lemma: continuation of microlocal kernels}For all $n\in \bN$ and $m\geq n$, we have
\[\cK_n\circ_{\sT} P_n\simeq \cK_{m}\circ_{\sT} P_n \in \sT(T^*(J\times M\times M)).\]    
\end{Lemma}
\begin{proof}Consider the isotopy $\psi_z=(\varphi^{H_n}_z)^{-1}\circ \varphi^{H_m}_z$ for $z\in J$, whose Hamiltonian is $D(z,q,p)=(H_m-H_n)(\varphi^{H_n}_z(q,p))$. 
Since $D,H_n,H_m$ are compactly supported, their sheaf quantizations exist, and satisfy 
\[\cK(D)\simeq \cK_n^{-1}\circ_{\sT,J} \cK_m\]
for an inverse $\cK_n^{-1}$ of $\cK_n$ (which is known to be unique) and the $J$-relative composition $\circ_{\sT,J}$.

On the other hand, for a neighborhood $W_n$ of $\overline{U_n}$, we have $D=0$ on the product $W_n \times J$. Then a direct computation shows that
\[\underline{ss}(\cK(D)\circ_{\sT} P_{n})\subset 0_{T^*J}\times T^*(M\times M),\quad i_0^*(\cK(D)\circ_{\sT} P_{U_n})\simeq P_{n}.\]

Therefore, by \autoref{lemma: microlocal morse-Tamarkin}, we have that $\cK(D)\circ_{\sT} P_{n}\simeq 1_{J}\boxtimes P_{n}$, and 
\[
\begin{aligned}
\cK_m\circ_\sT P_n{\simeq}&\,
(\cK_n\circ_{\sT,J}\cK(D))\circ_{\sT}P_n
\simeq
\cK_n\circ_{\sT,J}\bigl(\cK(D)\circ_{\sT}P_n\bigr) \\
\simeq &\,
\cK_n\circ_{\sT,J}(1_J\boxtimes P_n) \simeq \cK_n\circ_{\sT} P_n.
\end{aligned}
\]

Then we finish the proof.
\end{proof}

\begin{proof}[Proof of \autoref{thm: existence of sheafquantization for all complete hamiltonian}]

We prove the existence of sheaf quantizations on $J=(-R,R)$ first. 
We define
\[\cS_n\coloneqq \cK_n\circ_{\sT} P_n.\]
The inclusion $U_n\to U_{n+1}$ induces $P_{n}\to P_{n+1}$ by functoriality of microlocal kernels. Therefore, for all $m\geq n$, we have a morphism
\[\cS_n=\cK_n|_J\circ_{\sT} P_n\simeq \cK_{m}|_J\circ_{\sT} P_n \to \cK_{m}|_J\circ_{\sT} P_{m}=\cS_{m},\]
and we can define
\[\cS_{H,J} \coloneqq \colim_n \cS_n.\]

It is clear that \[i^*_0\cS_{H,J}|_{z=0}=\varinjlim_{n\to\infty} i^*_0\cS_n =\varinjlim_{n\to\infty} P_n=P_{T^*M}=1^\sT_{\Delta_M},\]
since $\bigcup_{n=1}^\infty U_n=T^*M$ (see for example \cite[Proposition 4.8]{zhang-cutoff}).

On the other hand, by the microsupport estimation for colimits \cite[Proposition 3.4]{guillermouviterbo_gammasupport}, we have
\[\underline{ss}(\cS_{H,J})\subset \bigcap_{N\in\bN} \overline{ \bigcup_{n\geq N}\underline{ss}(\cS_{n})}.\]

We define the following (Hamiltonian) diffeomorphism
\begin{align*}
&\Theta_H: T^*(J\times M\times M)\to T^*(J\times M\times M),\\
&\Theta_H(z,\zeta,q,p,q',p')=(z,\zeta-H(z,\varphi^H_z(q',p')),q,p,\varphi^H_z(q',p')).
\end{align*}
By definition of $\cS_n$ and cutoff estimation of $P_U$ \cite{zhang-cutoff}, one can verify that
\[\underline{ss}(\cS_{n})\subset \Theta_H(0_{T^*J}\times \underline{ss}( P_{n}))\subset \Theta_H(0_{T^*J}\times (\Delta^a_{T^*M}\cup (T^*M\setminus U_n)^2)).\]
Then notice that $\cap T^*M\setminus U_n)=\varnothing$, therefore, we have
\[\underline{ss}(\cS_{H,J})\subset  \Theta_H\Big( \bigcap _{N\in\bN} 0_{T^*J}\times (\Delta^a_{T^*M}\cup (T^*M\setminus U_N)^2)\Big)\subset \Theta_H\Big(  0_{T^*J}\times \Delta^a_{T^*M}\Big )\subset \Lambda_H |_{T^*J}.\]
Therefore, $\cS_{H,J}$ is a sheaf quantization of $H$ on $J$.

Next, we prove the uniqueness on $J\subset \bR_z$ open. The proof is standard; we present it here for the convenience of readers.

We start with the ``inverse kernel". For a complete Hamiltonian $H$, we take a sheaf quantization $\cK$ of $H$ and a sheaf quantization $\overline{\cK}$ of the inverse Hamiltonian $\overline{H}(z,q,p)=-H(z,\varphi_z(q,p))$.
By \autoref{lemma: relative compisition estimation}, we have 
\[\underline{ss}(\cK\circ_{\sT,J} \overline{\cK}),\quad \underline{ss}(\overline{\cK}\circ_{\sT,J} \cK)\subset \Lambda_{0}=T_{\Delta_M}(M\times M)\times 0_{T^*J},\]
and
\[i^*_0(\cK\circ_{\sT,J} \overline{\cK})\simeq i^*_0(\overline{\cK}\circ_{\sT,J} \cK )\simeq 1^\sT_{\Delta_M}.\]
Then by \autoref{lemma: microlocal morse-Tamarkin}, we have
\begin{equation}\label{equatiuon:inverse kernel identity}
  \cK\circ_{\sT,J} \overline{\cK}\simeq \cK\circ_{\sT,J} \overline{\cK}\simeq 1^\sT_{\Delta_M\times J}.  
\end{equation}

Now, for the complete Hamiltonian $H$, suppose there are two sheaf quantizations $\cK_1,\cK_2$ of $H$, and we pick a sheaf quantization $\overline{\cK}$ of the inverse Hamiltonian $\overline{H}$. \autoref{equatiuon:inverse kernel identity} is clearly true for both $\cK_1$ and $\cK_2$. Then we have  
\[\cK_1\simeq \cK_1\circ_{\sT,J} \overline{\cK} \circ_{\sT,J} \cK_2  \simeq \cK_2.\]

Lastly, for an increasing sequence $J_n\subset \bR_z$ with $\bigcup_{n\in\bN} J_n=\bR_z$, we can take the sheaf quantizations $S_{H,J_n}$. By the uniqueness, we can glue all $J\mapsto S_{H,J_n}$ to define
\[\cK(H)=\varinjlim_{n\to\infty} \cS_{H,J_n}.\]
The sheaf $\cK(H)$ satisfies the required properties.
\end{proof}

\begin{eg}If $H$ is quadratic at infinity, then we can construct $\cK(H)$ in the following way: By definition, we can write $H=Q\#K$, where $K$ is compactly supported and $Q(q,p)=a(|q|^2+|p|^2)$; for example, we can take $K(z,q,p)=(H-Q)\circ \varphi_Q(z,q,p)$. Then we know the existence of $\cK(K)$ and $\cK(Q)$ (which is a high-dimensional (rescaled) version of \autoref{equation: formula of GF sheaf quantization}). Then by the microsupport estimation lemma, and the uniqueness, we have
$\cK(H)\simeq \cK(Q)\circ_{\sT,\bR_z}\cK(K)$.    
\end{eg}

The trick we use can be extended a bit to isotopy of Hamiltonian embeddings that are not necessarily defined globally. 
\begin{Prop}\label{appendix: Hamiltonian embedding}Let $U\subset T^*M$ be open and $\varphi_s:I\times U\to T^*M$ be a Hamiltonian isotopy of embedding, which is not necessarily extended to $T^*M$. Then there exists kernel $\Phi,\Psi\in \sT(T^*(M\times M))$ such that $\Phi{\circ_\sT} \Psi\simeq P_U$ and $\Psi{\circ_\sT} \Phi\simeq P_{\varphi_1(U)}$. In particular, $\Phi$ induces an equivalence $\sT(U)\simeq \sT(\varphi_1(U))$. 
\end{Prop}
\begin{proof}Let $H_n\in C^\infty(I\times U)$ be the Hamiltonian function of $\varphi_s$. We pick an exhaustion $U_n$ of $U$ and cutoff functions $\rho_n$ as before. Then for $H_n=\rho_n$, we can extend them by $0$ to $T^*M$, which are globally defined on $T^*M$. Then the same construction as before gives $\cK_n$, $P_n$, and we define $\Phi=\varinjlim_n \cK_n \circ_\sT P_n$, and $\Psi= \cK_n^{-1}\circ_\sT P_{\varphi_1(U_n)}$ via \autoref{Lemma: continuation of microlocal kernels}. Then one can check easily that $\Psi{\circ_\sT} \Phi\simeq P_U$ and $\Phi{\circ_\sT} \Psi\simeq P_{\varphi_1(U)}$. 

In particular, convolution functors of $\Phi$ and $\Psi$ descend to an equivalence $\sT(U)\simeq \sT(\varphi_1(U))$.    
\end{proof}

\subsection{Sheaf quantization of Hamiltonians from generating functions}\label{appendix:generating function}

In this section, we give an explicit formula for $\cK(H)$ for $H\in C^\infty(\bR^{2n})$ with $\|D^2H\|_\infty<\infty$ using a generating function, following Chaperon's idea on broken trajectories. The idea is essentially to go back to the original construction of \cite{tamarkin2013}, and also see \cite[Proposition 10.4]{Viterbo_Princeton_Lecture}. In loc. cit., the quantization $\cK(H)$ of $\varphi_z^H$ has been well explained for a single $z$, but not formulated in the entire $\bR_z$-family. To complete the literature, we supply some details on this aspect. 

We will use the so-called generating function type-V following Chaperon \cite{Chaperon1984}; see also \cite[Section 9.2]{mcduff2017introST} for relevant discussion. However, the idea of the proof works if the other type of generating function exists; for example, we explain a construction for Hamiltonian rotation based on a generating function of type-S in \autoref{equation: formula of GF sheaf quantization} and how they are related in \autoref{eg: rotation quantization using V-type GF}. 

To be precise, we give the generating function results we need in detail.

\begin{Prop}\label{prop: existence of function F}Let
$
H\in C^\infty(\mathbb{R}^{2n})
$
be an autonomous Hamiltonian satisfying $
\|D^2H\|_\infty\leq L.
$ We write its Hamiltonian flow as
\[
\varphi_H^z(q,p)
=
\bigl(Q_z(q,p),P_z(q,p)\bigr).
\]
Then there exists $\varepsilon>0$ such that the following statements hold for every
$|z|<\varepsilon$.

\begin{enumerate}[fullwidth]
\item
The map
\[
{\Psi}:
(-\varepsilon,\varepsilon)\times \mathbb{R}^{n}_q\times \bR^n_p
\longrightarrow
(-\varepsilon,\varepsilon)\times \mathbb{R}^{n}_X\times \bR^n_p,
\]
defined by
\[
{\Psi}(z,q,p)
=
\bigl(z,Q_z(q,p),p\bigr),
\]
is a diffeomorphism. 

In particular, denote ${\Psi}^{-1}$ to be
\[{\Psi}^{-1}(z,X,p)=(z,x_z(X,p),p),\]
we have $Q_z(x_z(X,p),p)=X.$

\item Let 
\[
F:
(-\varepsilon,\varepsilon)
\times \mathbb{R}^n_X
\times \mathbb{R}^n_p
\longrightarrow \mathbb{R}
\]
be the function
\[
F(z,X,p)
=
p\cdot(x-X)
+
\int_0^z
\left(
P_s\cdot\dot Q_s-H(Q_s,P_s)
\right)\,ds,
\]
where
\(
x=x_z(X,p)
\) and \(
(Q_s,P_s)=\varphi_H^s(x,p).
\)

Then $F$ is smooth, $
F(0,X,p)=0$, and for
\(
x=x_z(X,p)
\) and \(
P=P_z(x_z(X,p),p),
\) 
the differential of $F$ is
\[
dF
=
(P-p)\cdot dX
+
(x-X)\cdot dp
-
H(X,P)\,dz.
\]

In particular, $F$ satisfies the Hamilton--Jacobi equation 
$
\partial_zF
+
H\bigl(X,p+\partial_XF\bigr)
=
0.
$
\end{enumerate}

Consequently, 
\[
\Lambda_H
=
\left\{
\left(
z,X,p;
-H(X,P),P-p,x-X
\right)
\ \middle|\
(X,P)=\varphi_H^z(q,p), |z|<\varepsilon
\right\}\subset T^*\bigl(
(-\varepsilon,\varepsilon)
\times \mathbb{R}^n_X
\times \mathbb{R}^n_p
\bigr)
\]
is the graph of $dF$.
\end{Prop}

\begin{proof}
(1)~Since
$
DX_H=JD^2H,
$
the Hamiltonian vector field $X_H$ has uniformly bounded derivative by $
\|D^2H\|_\infty\leq L
$. The linearized Hamiltonian equation gives
\[
\frac{d}{dz}D\varphi_H^z
=
DX_H(\varphi_H^z)D\varphi_H^z,
\qquad
D\varphi_H^0=\operatorname{id}.
\]
By Gronwall's inequality,
\[
\|D\varphi_H^z-\operatorname{id}\|_\infty
\leq
e^{L|z|}-1.
\]
Choose $\varepsilon>0$ such that
$
e^{L\varepsilon}-1<1.
$
Then, for every $|z|<\varepsilon$,
\[
\left\|
\frac{\partial Q_z}{\partial q}(q,p)
-\operatorname{id}
\right\|_\infty
<1.
\]

Therefore, by the Banach fixed point theorem, we have
\[
\widehat{\Psi}(z,q,p)
=
\bigl(z,Q_z(q,p),p\bigr).
\]
is a diffeomorphism, whose inverse has the form
\[
\widehat{\Psi}^{-1}(z,X,p)
=
\bigl(z,x_z(X,p),p\bigr),
\]
where $x_z(X,p)$ is smooth.

(2)~To see that $F$ is smooth at $z=0$, one can rewrite the integral in the definition of $F$, via the change of variable $s=zu$, as
\[
\int_0^z\left(P_s \cdot \dot{Q}_s - H(Q_s, P_s)\right) ds = z\int_0^1
\left(
p_{uz}\cdot
\partial_pH(x_{uz},p_{uz})
-
H(x_{uz},p_{uz})
\right)\,du.
\]

At $z=0$ one has $q=X$, and hence
$
F(0,X,p)=0.
$

Let
\[
A(z,X,x)
=
\int_0^z
\left(
P_s\cdot\dot Q_s-H(Q_s,P_s)
\right)\,ds
\]
denote the Hamiltonian action along the corresponding trajectory. Then
\[
dA
=
P\cdot dX
-
p\cdot dx
-
H(X,P)\,dz.
\]
Since
\(
F=A+p\cdot(x-X),
\) 
we obtain
\begin{align*}
dF
&=
P\cdot dX
-
p\cdot dx
-
H(X,P)\,dz
+
p\cdot dx
+
(x-X)\cdot dp
-
p\cdot dX \\
&=
(P-p)\cdot dX
+
(x-X)\cdot dp
-
H(X,P)\,dz.
\end{align*}

The description of $\Lambda_H$ as $\operatorname{Graph}(dF)$ follows
immediately.
\end{proof}

Now, we start our construction on sheaf quantization. Let us be careful about the variable. What we really need is a graph of $\varphi_H^z(q,p)$ over $q,Q$ variables. Now we consider
\[V(z,x,X,p)=F(z,X,p)+p\cdot (X-x).\]

Let 
\(\pi_{p}: (-\varepsilon,\varepsilon)_z\times \bR^n_x\times \bR^n_X\times \bR^n_p \to (-\varepsilon,\varepsilon)_z\times \bR^n_x\times \bR^n_X.\)
We define
\[\cK_1=\pi_{p!}1_{\{(z,x,X,p,t)\mid t+V(z,x,X,p)\geq 0\}}[n]\in \Sh((-\varepsilon,\varepsilon)_z\times \bR_x^n\times \bR^n_X,\sT).\]

\begin{Lemma}We have $\cK_1|_{z=0}=1_{\{x,X,t\}\mid x=X,t\geq 0\}}$ and \[\underline{ss}(\cK_1)\subset \{(x,-p,\varphi_H^z(x,p),z,-H \circ \varphi_H^z(z,p))\mid (z,x,p)\in (-\varepsilon,\varepsilon)_z\times T^*\bR^n_x\}.\] 
In other words, $\cK_1$ is a sheaf quantization of $\varphi_H$ and $\cK_1 \simeq \cK(H)|_{(-\varepsilon, \varepsilon)}$.
\end{Lemma}
\begin{proof}
By definition, we have
\[
V(z,x,X,p)
=
F(z,X,p)+p\cdot(X-x).
\]
Hence by the computation of $dF$ we have
\[
\partial_pV
=
\partial_pF+X-x
=
x_z(X,p)-x.
\]
We know that the graph $\mathrm{Graph}(dF)$ represents the graph of the Hamiltonian flow $\Lambda_H$. Therefore,
\[
\partial_p V=0
\quad\Longleftrightarrow\quad
(X,P)=\varphi_H^z(x,p),
\]
where
\[
P=\partial_XV =p+\partial_XF(z,X,p).
\]

By the standard microsupport estimate \cite[Proposition 5.4.4]{KS90} we obtain
\begin{align*}
    \underline{ss}(\cK_1)
&\subset
\left\{
(x,\partial_x V, X,\partial_X V,z,-H(X,P))
\ \mid \partial_p V=0
\right\}\\
&= \left\{
(x,-p,\varphi_z(x,p),z,-H \circ \varphi_H^z(z,p))
\ \middle|\
(z,x,p)\in
(-\varepsilon,\varepsilon)\times T^*\bR^n
\right\}.
\end{align*}

It remains to compute the restriction at $z=0$. We have
\[
F(0,X,p)=0,\quad V(0,x,X,p)
=
p\cdot (X-x ).
\]
Therefore
\[
\cK_1|_{z=0}
=
\pi_{p!}
1_{\{(x,X,p,t)\mid
t+ p\cdot(X-x)\geq0\}}[n]\simeq
1_{\{(x,X,t)\mid x=X,\, t\geq0\}}.\qedhere
\]
The fact that $\cK_1 \simeq \cK(H)|_{(-\varepsilon,\varepsilon)}$ follows from the uniqueness of sheaf quantization.
\end{proof}

On $(-N\varepsilon,N\varepsilon)$, we set $i_N:(-N\varepsilon,N\varepsilon) \to (-\varepsilon,\varepsilon), i_N(z)=z/N$, and we define the sheaf 
\[\cK_N\coloneqq i_{N}^*(\cK_1\circ_{z} \cdots \circ_z\cK_1 ) \in \Sh(\bR^n_x\times \bR^n_X\times (-N\varepsilon,N\varepsilon);\sT),\]
where $\circ_z$ is the relative convolution over $z\in (-\varepsilon,\varepsilon)$.
By the usual composition microsupport estimation and the uniqueness of sheaf quantization, we have:

\begin{Coro}
We have that $\cK_N$ is a sheaf quantization of $\varphi_H$ and
\[\cK_N\simeq \cK(H)|_{(-N\varepsilon,N\varepsilon)}.\]
\end{Coro}

At the end, we relate this $\cK_N$ to the broken trajectory generating function.

Setting
$x_0=x$, $x_N=X$, and
\[
V_N(z,x,X;\mathbf{x},\mathbf{p})
=
\sum_{j=1}^N
V\left(\frac{z}{N},x_{j-1},x_j,p_j\right)=\sum_{j=1}^N
\left[
F\left(\frac{z}{N},x_j,p_j\right)
+
 p_j\cdot (x_j-x_{j-1})
\right].
\]
This is exactly \cite[Equation (9.2.4)]{mcduff2017introST} (up to suitable sign convention), the broken trajectory generating function of Chaperon.

On the other hand, the standard 6-functor computation gives
\[
\cK_N
\simeq
\pi_{\mathbf{x},\mathbf{p}!}
1_{\{t+V_N\geq 0\}}[nN].
\]

\begin{eg}\label{eg: rotation quantization using V-type GF}Here, we consider the $2n$-dimensional rotation $\varphi_z(u)=\exp{(-2i\pi z)}u$ on $\bC^n_u$, which is the Hamiltonian flow of the Hamiltonian function $H(u)=\pi|u|^2 $ with constant Hessian. Under the real coordinate, we have
\[\varphi_z(x,p)=(x\cos(2\pi z)+p\sin(2\pi z),p\cos(2\pi z)-x\sin(2\pi z)).\]

Then we can apply \autoref{prop: existence of function F} to see that one can take $\varepsilon=1/4$ and the function $F$ to be 
\[F(z,X,p)=\Big(\frac{1}{\cos(2\pi z)}-1\Big) X\cdot p-\frac{\tan(2\pi z)}{2} (|X|^2+|p|^2).\]

Then we have
\[V(z,x,X,p)=F(z,X,p)+p\cdot(X-x)=\frac{1}{\cos(2\pi z)} X\cdot p - x\cdot p-\frac{\tan(2\pi z)}{2} (|X|^2+|p|^2).\]

In particular, if $z\neq 0$, we have
 \[V(z,x,X,p)=S_H(z,x,X)-\frac{\tan(2\pi z)}{2}\left| p- p^*\right|^2,\]  
for $S_H$ in \autoref{equation: formula of S_H}, and $p^*=\frac{X-x\cos(2\pi z)}{\sin(2\pi z)}$.

For this function $V$, we set, for $|z|<1/4$,  
\[\cS_1=\pi_{p!}1_{\{(z,x,X,p,t)\mid t+V(z,x,X,p)\geq 0\}}[n].\]

A direct 6-functor computation gives \autoref{equation: formula of GF sheaf quantization} on $0<|z|<1/4$: by the formula $V=S_H-\frac{\tan(2\pi z)}{2}\left| p- p^*\right|^2$, one can see that $\cS_1$ is supported on $\{t+S_H\geq 0\}$ and has expected stalks. 

To obtain \autoref{equation: formula of GF sheaf quantization} on $0<|z|<N/4$, one can either use $\cS_N\coloneqq i_{N}^*(\cS_1\circ_{z} \cdots \circ_z\cS_1 )$ by $\cS_1=1_{\{ t+S_H\geq 0\}}$ (with a suitable shift depending on $z$), or use 
\[
\cS_N
\simeq
\pi_{\mathbf{x},\mathbf{p}!}
1_{\{t+V_N\geq 0\}}[nN].
\]
with
\[V_N(z,x,X;\mathbf{x},\mathbf{p})=\sum _j S_H(z/N,x_j,x_{j+1})-\frac{\tan(2\pi z/N)}{2}\sum_j\left| p_j- p^*_j\right|^2\]
with $p_j^*=\frac{x_{j+1}-x_j\cos(2\pi z/N)}{\sin(2\pi z/N)}$ by a similar argument.
\end{eg}

\begin{eg}For compactly supported Hamiltonians, the existence of a $1$-parameter family version of generating function of $\varphi^H_z$ can be obtained by \cite{Sikorav1987}, see also \cite[Lemma 2.2-(i)]{large_non_equeezing_GF}. By the uniqueness of sheaf quantization, we can also define $\cK(H)|_{[0,1]}$ using this abstract existence result. 
\end{eg}

\bibliographystyle{bingyu}
\clearpage
\phantomsection
\bibliography{bibtexarxiv}

@book{mcduff2017introST,
	doi = {10.1093/oso/9780198794899.001.0001},
	year = {2017},
	month = {jun},
	publisher = {Oxford University Press},
	author = {Dusa McDuff and Dietmar Salamon},
	title = {Introduction to Symplectic Topology}
}

@InProceedings{tamarkin2013,
author="Tamarkin, Dmitry",
editor="Hitrik, Michael
and Tamarkin, Dmitry
and Tsygan, Boris
and Zelditch, Steve",
title="Microlocal condition for non-displaceability",
booktitle="Algebraic and Analytic Microlocal Analysis",
year="2018",
publisher="Springer",
pages="99--223",
doi={10.1007/978-3-030-01588-6\_3}
}

@article{GKS2012,
  title={Sheaf quantization of Hamiltonian isotopies and applications to nondisplaceability problems},
  author={{Guillermou}, St{\'e}phane and {Kashiwara}, Masaki and {Schapira}, Pierre},
  journal={Duke Mathematical Journal},
  volume={161},
  number={2},
  pages={201--245},
  year={2012},
  publisher={Duke University Press},
  doi={10.1215/00127094-1507367}
}

@incollection{GS2014,
	doi = {10.1007/978-3-319-06514-4\_3},
	year = {2014},
	publisher = {Springer},
	pages = {43--85},
    author={{Guillermou}, St{\'e}phane and {Schapira}, Pierre},
	title = {Microlocal theory of sheaves and Tamarkin's non Displaceability theorem},
	booktitle = {Homological mirror symmetry and tropical geometry}
}

@book{guillermou2019sheaves,
    author={ St{\'e}phane {Guillermou}},
    title={Sheaves and symplectic geometry of cotangent bundles},
    series = {Ast\'erisque},
    collection={Ast\'erisque},
    publisher = {Soci\'et\'e math\'ematique de France},
    doi = {10.24033/ast.1199},
    year = {2023},
    volume = {440},
}

@book{KS90,
  title={Sheaves on manifolds: With a short history.   Les d{\'e}buts de la th{\'e}orie des faisceaux. By Christian Houzel},
  author={Masaki {Kashiwara} and Pierre {Schapira}},
  volume={292},
  year={1990},
  publisher={Springer},
  doi={10.1007/978-3-662-02661-8}
}

@unpublished{chen1,
      title={On operadic open-closed maps in characteristic $p$}, 
      author={Zihong Chen},
      year={2024},
 note={\href{https://arxiv.org/abs/2402.06183}{arXiv:2402.06183}},
}

@unpublished{mc-Tamarkin,
  title={Microlocal Category},
  author={{Tamarkin}, Dmitry},
  note={\href{https://arxiv.org/abs/1511.08961}{arXiv:1511.08961}},
  year={2015}
}

@article{Capacities2021,
    title={Capacities from the Chiu-Tamarkin complex},
    author={Bingyu {Zhang}},
    DOI={10.4310/jsg.241001211759}, 
number={3}, 
journal={Journal of Symplectic Geometry}, 
publisher={International Press of Boston}, 
volume={22},
year={2024}, 
pages={441–524}
}

@unpublished{CyclicZHANG,
  title={Idempotence of microlocal kernels and the $S^1$-equivariant Chiu-Tamarkin invariant},
  author={ Bingyu {Zhang}},
    year={2023},
 note={\href{https://arxiv.org/abs/2306.12316}{arXiv:2306.12316}},
}

@unpublished{APT-Tatsuki-Zhang,
  title={Almost mathematics, Persistence module, and {T}amarkin category},
  author={ Tatsuki {Kuwagaki} and Bingyu {Zhang}},
    year={2025},
note={\href{https://arxiv.org/abs/2503.15933}{arXiv:2503.15933}},
}

@unpublished{Hochschild-Kuo-Shende-Zhang,
  title={On the {Hochschild} cohomology of {Tamarkin} categories},
  author={Christopher {Kuo} and Vivek {Shende} and  Bingyu {Zhang}},
    year={2023},
 note={\href{https://arxiv.org/abs/2312.11447}{arXiv:2312.11447}},
}

@article{guillermouviterbo_gammasupport,
author={St{\'e}phane {Guillermou} and Claude {Viterbo}},
title={The singular support of sheaves is {$\gamma$}-coisotropic},
  journal={Geometric and Functional Analysis},
  volume={34},
  number={4},
  pages={1052--1113},
  year={2024},
  publisher={Springer},
doi={10.1007/s00039-024-00682-x}
}

@article{asanoike2022complete,
  title={Completeness of derived interleaving distances and sheaf quantization of non-smooth objects},
  author={Tomohiro Asano and Yuichi Ike},
  journal={Mathematische Annalen},
  volume={390},
  number={2},
  pages={2991--3037},
  year={2024},
  doi={10.1007/s00208-024-02815-x},
}

@article{chiu2017,
  title={Non-squeezing property of contact balls},
  author={  Sheng-Fu {Chiu}},
  journal={Duke Mathematical Journal},
  volume={166},
  number={4},
  pages={605--655},
  year={2017},
  publisher={Duke University Press},
  doi={10.1215/00127094-3715517}
}

@article{gromov1985,
  title={Pseudo holomorphic curves in symplectic manifolds},
  author={ Mikhael {Gromov}},
  journal={Inventiones Mathematicae},
  volume={82},
  number={2},
  pages={307--347},
  year={1985},
  publisher={Springer},
  doi={10.1007/bf01388806 },
}

@article{viterbo1992symplectic,
  title={Symplectic topology as the geometry of generating functions},
  author={ Claude {Viterbo}},
  journal={Mathematische Annalen},
  volume={292},
  number={1},
  pages={685--710},
  year={1992},
  publisher={Springer},
  doi={10.1007/bf01444643 },
}

@article{viterbo-functorsI,
  title={Functors and computations in Floer homology with applications, I},
  author={Viterbo, Claude},
  journal={Geometric and Functional Analysis},
  volume={9},
  number={5},
  pages={985--1033},
  year={1999},
  publisher={Springer},
  doi={10.1007/s000390050106},
}

@article{kuo2021wrapped,
title = {Wrapped sheaves},
journal = {Advances in Mathematics},
volume = {415},
pages = {108882},
year = {2023},
doi = {10.1016/j.aim.2023.108882},
author = {Christopher Kuo},
}

@article{EkelandHoferII,
author = {  Ivar {Ekeland}  and Helmut {Hofer} },
journal = {Mathematische Zeitschrift},
number = {4},
pages = {553-568},
title = {Symplectic Topology and Hamiltonian Dynamics II.},
volume = {203},
year = {1990},
doi={10.1007/bf01221255},
publisher={Springer},
}

@unpublished{Chiu-kernel2021,
  title={Microlocal projector for complete flow},
  author={ Sheng-Fu  {Chiu}},
  note={\href{https://arxiv.org/abs/2112.00483}{arXiv:2112.00483}},
  year={2021}
}

@article{gutthutchings2018capacities,
  title={Symplectic capacities from positive $S^1$--equivariant symplectic homology},
  author={ Jean {Gutt} and Michael {Hutchings}},
  journal={Algebraic \& Geometric Topology},
  volume={18},
  number={6},
  pages={3537--3600},
  year={2018},
  publisher={Mathematical Sciences Publishers},
  doi={10.2140/agt.2018.18.3537}
}

@article{hutchings2011,
author = "Michael {Hutchings}",
doi = "10.4310/jdg/1320067647",
fjournal = "Journal of Differential Geometry",
journal = {Journal of Differential Geometry},
month = "06",
number = "2",
pages = "231--266",
publisher = "Lehigh University",
title = "Quantitative Embedded Contact Homology",
volume = "88",
year = "2011"
}

@article{Cieliebak-Frauenfelder-Oancea2010,
	doi = {10.24033/asens.2137},
	year = 2010,
	publisher = {Societe Mathematique de France},
	volume = {43},
	number = {6},
	pages = {957--1015},
	author = {Kai Cieliebak and Urs Frauenfelder and Alexandru Oancea},
	title = {Rabinowitz Floer homology and symplectic homology},
	journal = {Annales Scientifiques de l'École Normale Supérieure}
}

@unpublished{SiegelHigherCapacities,
  title={Higher symplectic capacities},
  author={Kyler {Siegel}},
  note={\href{https://arxiv.org/abs/1902.01490}{arXiv:1902.01490}},
  year={2019}
}

@article{CompSiegelHigherCapacities,
    author = {Siegel, Kyler},
    title = "{Computing Higher Symplectic Capacities I}",
    journal = {International Mathematics Research Notices},
    year = {2021},
    month = {04},
    issn = {1073-7928},
    doi = {10.1093/imrn/rnaa334},
}

@book{LodayCyclic,
	doi = {10.1007/978-3-662-11389-9},
	year = 1998,
	publisher = {Springer},
	author = {Jean-Louis Loday},
	title = {Cyclic homology},
}

@book{BredonSheaf,
	doi = {10.1007/978-1-4612-0647-7},
	year = 1997,
	publisher = {Springer},
	author = {Glen E. Bredon},
	title = {Sheaf Theory}
}

@article{GanatraS1CY, 
title={Cyclic homology, {$S^1$}–equivariant {F}loer
cohomology and {C}alabi–{Y}au structures}, 
volume={27}, 
ISSN={1465-3060}, 
DOI={10.2140/gt.2023.27.3461}, 
number={9}, 
journal={Geometry \& Topology}, 
publisher={Mathematical Sciences Publishers}, 
author={Ganatra, Sheel}, 
year={2023}, 
month=dec, 
pages={3461–3584} 
}

@article{Ike_2019,
	doi = {10.4171/prims/55-4-3},
	year = 2019,
	month = {oct},
	publisher = {European Mathematical Society},
	volume = {55},
	number = {4},
	pages = {737--778},
	author = {Yuichi Ike},
	title = {Compact Exact Lagrangian Intersections in Cotangent Bundles via Sheaf Quantization},
	journal = {Publications of the Research Institute for Mathematical Sciences}
}

@article{large_non_equeezing_GF,
	doi = {10.1007/s11784-025-01188-1},
	year = 2025,
	month = {April},
	publisher = {Springer Science and Business Media {LLC}},
	volume = {27},
  author={Maia {Fraser} and Sheila {Sandon} and  Bingyu {Zhang}},
  title={Contact non-squeezing at large scale via generating functions},
	journal = {Journal of Fixed Point Theory and Applications}
}

@unpublished{Kuo-Li-Duality,
 title={Duality and spherical adjunction from microlocalization}, 
      author={Christopher Kuo and Wenyuan Li},
  note={\href{https://arxiv.org/abs/2210.06643}{arXiv:2210.06643}},
  year={2022}
}

@article{Efimov-K-theory,
  title={K-theory and localizing invariants of large categories},
  author={Alexander I. Efimov},
  journal={Publications mathématiques de l'IHÉS},
  note={To appear; \href{https://arxiv.org/abs/2405.12169}{arXiv:2405.12169}},
  year={2024},
}

@article{zhaoperiodic,
doi = {10.4310/jsg.2019.v17.n5.a9},
year = 2019,
publisher = {International Press of Boston},
volume = {17},
number = {5},
pages = {1513--1578},
author = {Jingyu Zhao},
title = {Periodic symplectic cohomologies},
journal = {Journal of Symplectic Geometry},
}

@phdthesis{zhangthesis, 
 title={Sheaf techniques and symplectic geometry},
 author={Zhang, Bingyu}, 
 school = {Université Grenoble Alpes},
 url={https://byzhang.cn/Files/PhD_Thesis.pdf},
 year={2022}
 }

@unpublished{HA,
        note={Preprint, available at \href{https://www.math.ias.edu/~lurie/papers/HA.pdf}{https://www.math.ias.edu/~lurie/papers/HA.pdf}},
	year = 2017,
	author = {Jacob Lurie},
	title = {Higher Algebra}
}

@unpublished{Biolley,
  title={Floer homology, symplectic and complex hyperbolicities},
  author={Biolley, Anne-Laure},
    note={\href{https://arxiv.org/abs/math/0404551}{arXiv:math/0404551}},
  year={2004}
}

@unpublished{Asano-Ike-peg,
      title={The rectifiable rectangular peg problem}, 
      author={Tomohiro Asano and Yuichi Ike},
      year={2024},
note={\href{https://arxiv.org/abs/2412.21057}{arXiv:2412.21057}},
}

@article{Zhang_K_Remark, 
 title={A remark on continuous K-theory and Fourier-Sato transformation}, 
 volume={154}, 
 ISSN={0002-9939}, 
DOI={10.1090/proc/17589}, 
number={6}, 
journal={Proceedings of the American Mathematical Society}, publisher={American Mathematical Society (AMS)}, 
author={Zhang, Bingyu}, 
year={2026}, 
month=Mar, 
pages={2357–2368} }

@article{ToricMSrevisted-Shende, 
title={Toric mirror symmetry revisited}, 
volume={360}, 
ISSN={1778-3569}, 
DOI={10.5802/crmath.304}, 
number={G7}, 
journal={Comptes Rendus. Mathématique}, 
publisher={Cellule MathDoc/Centre Mersenne}, 
author={Shende, Vivek}, 
year={2022}, 
month=jun, 
pages={751–759} }

@unpublished{StacksProject,
    AUTHOR = " {Stacks Project}, Authors of",
    TITLE = "Stacks Project",
    note={\href{https://stacks.math.columbia.edu/}{https://stacks.math.columbia.edu/}},
}

@unpublished{SunFunctorialityToricCCC,
note={\href{https://arxiv.org/abs/2504.21483}{arXiv:2504.21483}},
	year = 2025,
  author = {Yuze Sun},
  title = {Functoriality of toric coherent-constructible correspondence},
}

@unpublished{Irie_ChainModel,
    title={A conjectural chain model for positive $S^1$-equivariant symplectic homology of star-shaped toric domains in $\mathbb{C}^2$},
    author = {Kei Irie},
    note={\href{https://arxiv.org/abs/2304.08720}{arXiv:2304.08720}},
    year={2023}
}

@incollection{quantitivesymplectic2010,
  author    = {Kai {Cieliebak} and Helmut {Hofer} and Janko {Latschev} and Felix {Schlenk}},
  title     = {Quantitative symplectic geometry},
  booktitle = {Dynamics, Ergodic Theory, and Geometry},
  editor    = {Boris Hasselblatt},
  doi       = {10.1017/cbo9780511755187.002},
  year      = {2010},
  publisher ={Cambridge University Press},
}

@unpublished{Kuo-Li-CY,
  title={Relative {Calabi--Yau} structure on microlocalization},
  author={Kuo, Christopher and Li, Wenyuan},
  year={2024},
  note={\href{https://arxiv.org/abs/2408.04085}{arXiv:2408.04085}}
}

@article{Humiliere2008, 
title={On some completions of the space of hamiltonian maps}, volume={136}, 
ISSN={2102-622X},  
DOI={10.24033/bsmf.2560}, 
number={3}, 
journal={Bulletin de la Société mathématique de France}, 
publisher={Societe Mathematique de France}, 
author={Humili\'ere, Vincent},
year={2008}, 
pages={373–404} }

@article{zhang-cutoff, title={Non-linear microlocal cut-off functors}, ISSN={2240-2926}, DOI={10.4171/rsmup/174}, journal={Rendiconti del Seminario Matematico della Università di Padova}, publisher={European Mathematical Society - EMS - Publishing House GmbH}, author={Zhang, Bingyu}, year={2025}, month=jan }

@article {Bourgeois-Oancea-automomous,
    AUTHOR = {Bourgeois, Fr\'{e}d\'{e}ric and Oancea, Alexandru},
     TITLE = {Symplectic homology, autonomous {H}amiltonians, and
              {M}orse-{B}ott moduli spaces},
   JOURNAL = {Duke Mathematical Journal},
  FJOURNAL = {Duke Mathematical Journal},
    VOLUME = {146},
      YEAR = {2009},
    NUMBER = {1},
     PAGES = {71--174},
      ISSN = {0012-7094},
   MRCLASS = {53D40},
  MRNUMBER = {2475400},
MRREVIEWER = {Tobias Ekholm},
       DOI = {10.1215/00127094-2008-062},
}

@article {Oancea-kunneth,
    AUTHOR = {Oancea, Alexandru},
     TITLE = {The {K}\"{u}nneth formula in {F}loer homology for manifolds
              with restricted contact type boundary},
   JOURNAL = {Mathematische Annalen},
  FJOURNAL = {Mathematische Annalen},
    VOLUME = {334},
      YEAR = {2006},
    NUMBER = {1},
     PAGES = {65--89},
      ISSN = {0025-5831,1432-1807},
   MRCLASS = {53D40},
  MRNUMBER = {2208949},
MRREVIEWER = {Timothy\ Perutz},
       DOI = {10.1007/s00208-005-0700-0},

}

@article {CieliebakMohnke,
    AUTHOR = {Cieliebak, Kai and Mohnke, Klaus},
     TITLE = {Punctured holomorphic curves and {L}agrangian embeddings},
   JOURNAL = {Inventiones Mathematicae},
  FJOURNAL = {Inventiones Mathematicae},
    VOLUME = {212},
      YEAR = {2018},
    NUMBER = {1},
     PAGES = {213--295},
      ISSN = {0020-9910,1432-1297},
   MRCLASS = {53D12 (53D05 53D35)},
  MRNUMBER = {3773793},
MRREVIEWER = {Luc\ Vrancken},
       DOI = {10.1007/s00222-017-0767-8},

}

@incollection{Viterbo_Seminaire_Bourbaki,
     author = {Viterbo, Claude},
     title = {Capacit\'es symplectiques et applications},
     booktitle = {S\'eminaire Bourbaki : volume 1988/89, expos\'es 700-714},
     series = {Ast\'erisque},
     note = {talk:714},
     pages = {345--362},
     year = {1989},
     publisher = {Soci\'et\'e math\'ematique de France},
     number = {177-178},
     mrnumber = {1040580},
     zbl = {0698.53019},
     language = {fr},
}

@book{Floer_book,
  author    = {Hofer, Helmut and Roberts, Siobhan},
  title     = {The {Floer} Jungle: Charting the Development of a Theory},
  series    = {History of Mathematics},
  volume    = {51},
  publisher = {American Mathematical Society},
  address   = {Providence, RI},
  year      = {2026},
  isbn      = {978-1-4704-8514-6},
  url       = {https://bookstore.ams.org/hmath-51}
}

@unpublished{Floer_notes,
  author={Floer, Andreas},
  title={Capacities, etc.},
  note={Notes taken by Weinstein}
}

@unpublished{Chekanov_notes,
  author={Chekanov, Yuri},
  title={Equivariant Morse-Floer theory},
  note={Notes taken by Schlenk}
}

@article {AntieauMathewNikolaus,
    AUTHOR = {Antieau, Benjamin and Mathew, Akhil and Nikolaus, Thomas},
     TITLE = {On the {B}lumberg-{M}andell {K}\"{u}nneth theorem for {TP}},
   JOURNAL = {Selecta Mathematica},
  FJOURNAL = {Selecta Mathematica. New Series},
    VOLUME = {24},
      YEAR = {2018},
    NUMBER = {5},
     PAGES = {4555--4576},
      ISSN = {1022-1824,1420-9020},
   MRCLASS = {14F30 (16E40 19D55)},
  MRNUMBER = {3874698},
MRREVIEWER = {Andrea\ Solotar},
       DOI = {10.1007/s00029-018-0427-x},
 
}

@article{MNN_nilpotence_equivariant_homotopy,
  author  = {Akhil Mathew and Niko Naumann and Justin Noel},
  title   = {Nilpotence and descent in equivariant stable homotopy theory},
  journal = {Advances in Mathematics},
  volume  = {305},
  year    = {2017},
  pages   = {994--1084},
  doi     = {10.1016/j.aim.2016.09.027}
}

@book {EKMM,
    AUTHOR = {Elmendorf, A. D. and Kriz, I. and Mandell, M. A. and May, J.
              P.},
     TITLE = {Rings, modules, and algebras in stable homotopy theory},
    SERIES = {Mathematical Surveys and Monographs},
    VOLUME = {47},
      NOTE = {With an appendix by M. Cole},
 PUBLISHER = {American Mathematical Society, Providence, RI},
      YEAR = {1997},
     PAGES = {xii+249},
      ISBN = {0-8218-0638-6},
   MRCLASS = {55N20 (19D10 19D55 55P42 55T25)},
  MRNUMBER = {1417719},
MRREVIEWER = {Donald\ M.\ Davis},
       DOI = {10.1090/surv/047},
     
}

@article {Goodwillie85,
    AUTHOR = {Goodwillie, Thomas G.},
     TITLE = {Cyclic homology, derivations, and the free loopspace},
   JOURNAL = {Topology},
  FJOURNAL = {Topology. An International Journal of Mathematics},
    VOLUME = {24},
      YEAR = {1985},
    NUMBER = {2},
     PAGES = {187--215},
      ISSN = {0040-9383},
   MRCLASS = {18F25 (18G15 19D55 55P35)},
  MRNUMBER = {793184},
MRREVIEWER = {Ross\ Staffeldt},
       DOI = {10.1016/0040-9383(85)90055-2},

}

@article{Pereira,
  title={Equivariant symplectic homology, linearized contact homology and the {L}agrangian capacity},
  author={Pereira, Miguel},
  journal={\href{https://arxiv.org/abs/2205.13381}{arXiv:2205.13381}},
  year={2022}
}

@article{Faisal-Li,
  title={Lagrangian capacity and chain level string topology},
  author={Shah Faisal and Yin Li},
  journal={\href{https://arxiv.org/abs/2606.20051}{arXiv:2606.20051}},
  year={2026},
}

@book{Viterbo_Princeton_Lecture,
title="An Introduction to Symplectic Topology through Sheaf theory",
  author={ Claude {Viterbo}},
publisher ="Lectures at Princeton Fall and New York Spring 2011",
year="2011"
}

@article {Bourgeois-Oancea-17,
    AUTHOR = {Bourgeois, Fr\'{e}d\'{e}ric and Oancea, Alexandru},
     TITLE = {{$S^1$}-equivariant symplectic homology and linearized contact homology},
   JOURNAL = {International Mathematics Research Notices},
  FJOURNAL = {International Mathematics Research Notices. IMRN},
      YEAR = {2017},
    NUMBER = {13},
     PAGES = {3849--3937},
      ISSN = {1073-7928,1687-0247},
   MRCLASS = {53D42 (53D40 55N91 57R17 57R58)},
  MRNUMBER = {3671507},
MRREVIEWER = {Daniel\ V.\ Mathews},
       DOI = {10.1093/imrn/rnw029},
          volume = {2017}, 
     
}

@article {Seidel-Equivariant,
    AUTHOR = {Seidel, Paul},
     TITLE = {The equivariant pair-of-pants product in fixed point {F}loer
              cohomology},
   JOURNAL = {Geometric and Functional Analysis},
  FJOURNAL = {Geometric and Functional Analysis},
    VOLUME = {25},
      YEAR = {2015},
    NUMBER = {3},
     PAGES = {942--1007},
      ISSN = {1016-443X,1420-8970},
   MRCLASS = {53D40 (55N91 57R58)},
  MRNUMBER = {3361776},
MRREVIEWER = {Darko\ Milinkovi\'{c}},
       DOI = {10.1007/s00039-015-0331-x},
  
}

@unpublished{BeckertGroth_Cubical_homotopy,
  title={Abstract cubical homotopy theory},
  author={Falk Beckert and Moritz Groth},
  year={2018},
  note={\href{https://arxiv.org/abs/1803.06022}{arXiv:1803.06022}},
}

@book{Weibel1994,
  author    = {Weibel, Charles A.},
  title     = {An Introduction to Homological Algebra},
  series    = {Cambridge Studies in Advanced Mathematics},
  volume    = {38},
  publisher = {Cambridge University Press},
  address   = {Cambridge},
  year      = {1994},
  doi       = {10.1017/CBO9781139644136}
}

@misc{Li-Zhang26-1,
    author = {Li, Wenyuan and Zhang, Bingyu},
    title = {Cyclic homology of sheaf categories},
    year = {2026},
    note = {in preparation}
}

@article{JonesPetrack1990,
  author  = {Jones, J. D. S. and Petrack, S. B.},
  title   = {The fixed point theorem in equivariant cohomology},
  journal = {Transactions of the American Mathematical Society},
  volume  = {322},
  number  = {1},
  pages   = {35--49},
  year    = {1990},
  doi     = {10.1090/S0002-9947-1990-1010411-X}
}

@article{HWZ1998,
  author  = {Hofer, Helmut and Wysocki, Krzysztof and Zehnder, Eduard},
  title   = {The dynamics on three-dimensional strictly convex energy surfaces},
  journal = {Annals of Mathematics},
  series  = {2},
  volume  = {148},
  number  = {1},
  pages   = {197--289},
  year    = {1998},
  doi     = {10.2307/120994}
}

@article{Zhou2022kdilations,
  author  = {Zhou, Zhengyi},
  title   = {Symplectic fillings of asymptotically dynamically convex
             manifolds {II}--{$k$}-dilations},
  journal = {Advances in Mathematics},
  volume  = {406},
  pages   = {108522},
  year    = {2022},
  doi     = {10.1016/j.aim.2022.108522}
}

@unpublished{Homological-perturbation,
    author = {Crainic, Marius},
    title = {On the perturbation lemma, and deformations},
    year={2004},
    note = {\href{https://arxiv.org/abs/math/0403266}{arXiv:math/0403266}}
}

@article{Franz_E2formality, 
 title={Homotopy Gerstenhaber formality of Davis–Januszkiewicz spaces}, 
 volume={23}, 
 ISSN={1532-0081}, 
 DOI={10.4310/hha.2021.v23.n2.a17}, 
 number={2}, 
 journal={Homology, Homotopy and Applications}, 
 publisher={International Press of Boston}, 
 author={Franz, Matthias}, 
 year={2021}, 
 pages={325–347} }

@article{Horel_E2formality,
     author = {Horel, Geoffroy},
     title = {$E_2$-formality via obstruction theory},
     journal = {Annales Henri Lebesgue},
     pages = {699--720},
     year = {2025},
     publisher = {\'ENS Rennes},
     volume = {8},
     doi = {10.5802/ahl.246},
}

@article{AlbersCieliebakFrauenfelder2016,
  author  = {Peter Albers and Kai Cieliebak and Urs Frauenfelder},
  title   = {Symplectic Tate homology},
  journal = {Proceedings of the London Mathematical Society},
  volume  = {112},
  number  = {1},
  pages   = {169--205},
  year    = {2016},
  doi     = {10.1112/plms/pdv065},
}

@article {Seidel-graded,
    AUTHOR = {Seidel, Paul},
     TITLE = {Graded {L}agrangian submanifolds},
   JOURNAL = {Bull. Soc. Math. France},
  FJOURNAL = {Bulletin de la Soci\'{e}t\'{e} Math\'{e}matique de France},
    VOLUME = {128},
      YEAR = {2000},
    NUMBER = {1},
     PAGES = {103--149},
      ISSN = {0037-9484,2102-622X},
   MRCLASS = {53D40 (53D12 57R58)},
  MRNUMBER = {1765826},
MRREVIEWER = {David\ E.\ Hurtubise},
       URL = {http://www.numdam.org/item?id=BSMF_2000__128_1_103_0},
}

@article {RobbinSalamon,
    AUTHOR = {Robbin, Joel and Salamon, Dietmar},
     TITLE = {The {M}aslov index for paths},
   JOURNAL = {Topology},
  FJOURNAL = {Topology. An International Journal of Mathematics},
    VOLUME = {32},
      YEAR = {1993},
    NUMBER = {4},
     PAGES = {827--844},
      ISSN = {0040-9383},
   MRCLASS = {58F05 (58E05)},
  MRNUMBER = {1241874},
MRREVIEWER = {J.\ S.\ Joel},
       DOI = {10.1016/0040-9383(93)90052-W},
}

@article {CieliebakFloerHoferWysocki-application,
    AUTHOR = {Cieliebak, K. and Floer, A. and Hofer, H. and Wysocki, K.},
     TITLE = {Applications of symplectic homology. {II}. {S}tability of the
              action spectrum},
   JOURNAL = {Math. Z.},
  FJOURNAL = {Mathematische Zeitschrift},
    VOLUME = {223},
      YEAR = {1996},
    NUMBER = {1},
     PAGES = {27--45},
      ISSN = {0025-5874,1432-1823},
   MRCLASS = {58F05 (57R99 58E99)},
  MRNUMBER = {1408861},
MRREVIEWER = {Karl\ Friedrich\ Siburg},
       DOI = {10.1007/PL00004267},
}

@article{Seidel2018Connections,
  author  = {Seidel, Paul},
  title   = {Connections on equivariant {Hamiltonian} {Floer} cohomology},
  journal = {Commentarii Mathematici Helvetici},
  volume  = {93},
  number  = {3},
  year    = {2018},
  pages   = {587--644},
  doi     = {10.4171/CMH/445}
}

@article{Hermann2000,
  author  = {Hermann, David},
  title   = {Holomorphic curves and {Hamiltonian} systems in an open set
             with restricted contact-type boundary},
  journal = {Duke Mathematical Journal},
  volume  = {103},
  number  = {2},
  year    = {2000},
  pages   = {335--374},
  doi     = {10.1215/S0012-7094-00-10327-4}
}

@article{Hermann2004,
  author  = {Hermann, David},
  title   = {Inner and outer {Hamiltonian} capacities},
  journal = {Bulletin de la Soci\'et\'e Math\'ematique de France},
  volume  = {132},
  number  = {4},
  year    = {2004},
  pages   = {509--541},
  doi     = {10.24033/bsmf.2472}
}

@article{Connes1983,
  author  = {Connes, Alain},
  title   = {Cohomologie cyclique et foncteurs {$\mathrm{Ext}^n$}},
  journal = {C. R. Acad. Sci. Paris S{\'e}r. I Math.},
  volume  = {296},
  number  = {23},
  year    = {1983},
  pages   = {953--958}
}

@incollection{GreenleesMay1995,
  author    = {Greenlees, J. P. C. and May, J. P.},
  title     = {Completions in algebra and topology},
  booktitle = {Handbook of Algebraic Topology},
  editor    = {James, I. M.},
  publisher = {North-Holland},
  address   = {Amsterdam},
  year      = {1995},
  pages     = {255--276}
}

@phdthesis{Sen2023,
  author = {Sen, Debtanu},
  title  = {The {$\mathbb{Z}/p$} {Gysin} Sequence in Symplectic Cohomology},
  school = {University of Southern California},
  year   = {2023}
}

@article {DwyerGreenlees,
    AUTHOR = {Dwyer, W. G. and Greenlees, J. P. C.},
     TITLE = {Complete modules and torsion modules},
   JOURNAL = {Amer. J. Math.},
  FJOURNAL = {American Journal of Mathematics},
    VOLUME = {124},
      YEAR = {2002},
    NUMBER = {1},
     PAGES = {199--220},
      ISSN = {0002-9327,1080-6377},
   MRCLASS = {16E30 (16D90 18E30)},
  MRNUMBER = {1879003},
MRREVIEWER = {Henning\ Krause},
       URL =
              {http://muse.jhu.edu/journals/american_journal_of_mathematics/v124/124.1dwyer.pdf},
}

@article{CoteKartal24,
    author = {C\^{o}t\'{e}, Laurent and Kartal, Yusuf Bari\c{s}},
    title = {Recovering generalized homology from {F}loer homology: the complex oriented case},
  journal={Ann. Sc. norm. super. Pisa},
  note={To appear; \href{https://arxiv.org/abs/2404.02776}{arXiv:2404.02776}},
  year={2024},
}

@article {GrothPontoShulman,
    AUTHOR = {Groth, Moritz and Ponto, Kate and Shulman, Michael},
     TITLE = {The additivity of traces in monoidal derivators},
   JOURNAL = {J. K-Theory},
  FJOURNAL = {Journal of K-Theory. K-Theory and its Applications in Algebra,
              Geometry, Analysis \& Topology},
    VOLUME = {14},
      YEAR = {2014},
    NUMBER = {3},
     PAGES = {422--494},
      ISSN = {1865-2433,1865-5394},
   MRCLASS = {18D10 (18E30 18G55 55U35)},
  MRNUMBER = {3349323},
MRREVIEWER = {J\v{i}r\'{\i}\ Rosick\'{y}},
       DOI = {10.1017/is014005011jkt262},
}

@article {HainePortaTeyssier-conicality,
    AUTHOR = {Haine, Peter J. and Porta, Mauro and Teyssier, Jean-Baptiste},
     TITLE = {Exodromy beyond conicality},
   JOURNAL = {J. Topol.},
  FJOURNAL = {Journal of Topology},
    VOLUME = {19},
      YEAR = {2026},
    NUMBER = {3},
     PAGES = {Paper No. e70091, 98},
      ISSN = {1753-8416,1753-8424},
   MRCLASS = {18N60 (32B20 32S60)},
  MRNUMBER = {5106933},
MRREVIEWER = {Mee\ Seong\ Im},
       DOI = {10.1112/topo.70091},
}

@unpublished{GPR22,
    author = {Gutt, Jean and Pereira, Miguel and Ramos, Vinicius G. B.},
    title = {Cube normalized symplectic capacities},
    note = {\href{https://arxiv.org/abs/2208.13666}{arXiv:2208.13666}},
    year={2022},
}

@article {HoferWysockiZehnder,
    AUTHOR = {Hofer, H. and Wysocki, K. and Zehnder, E.},
     TITLE = {A characterization of the tight {$3$}-sphere. {II}},
   JOURNAL = {Comm. Pure Appl. Math.},
  FJOURNAL = {Communications on Pure and Applied Mathematics},
    VOLUME = {52},
      YEAR = {1999},
    NUMBER = {9},
     PAGES = {1139--1177},
      ISSN = {0010-3640,1097-0312},
   MRCLASS = {53D35 (57R17)},
  MRNUMBER = {1692144},
MRREVIEWER = {Karl\ Friedrich\ Siburg},
       DOI = {10.1002/(sici)1097-0312(199909)52:9<1139::aid-cpa5>3.0.co;2-l},
}

@article {Cieliebak-Oancea-ES,
    AUTHOR = {Cieliebak, Kai and Oancea, Alexandru},
     TITLE = {Symplectic homology and the {E}ilenberg-{S}teenrod axioms},
      NOTE = {Appendix written jointly with Peter Albers},
   JOURNAL = {Algebr. Geom. Topol.},
  FJOURNAL = {Algebraic \& Geometric Topology},
    VOLUME = {18},
      YEAR = {2018},
    NUMBER = {4},
     PAGES = {1953--2130},
      ISSN = {1472-2747,1472-2739},
   MRCLASS = {53D40 (55N40 57R17 57R90)},
  MRNUMBER = {3797062},
MRREVIEWER = {Sonja\ Hohloch},
       DOI = {10.2140/agt.2018.18.1953},
}

@article {HutchingsNelson,
    AUTHOR = {Hutchings, Michael and Nelson, Jo},
     TITLE = {Cylindrical contact homology for dynamically convex contact forms in three dimensions},
   JOURNAL = {J. Symplectic Geom.},
  FJOURNAL = {The Journal of Symplectic Geometry},
    VOLUME = {14},
      YEAR = {2016},
    NUMBER = {4},
     PAGES = {983--1012},
      ISSN = {1527-5256,1540-2347},
   MRCLASS = {53D42 (53D35 57R17)},
  MRNUMBER = {3601881},
MRREVIEWER = {David\ E.\ Hurtubise},
       DOI = {10.4310/JSG.2016.v14.n4.a1},
}

@incollection {GuttHutchingsRamos,
    AUTHOR = {Gutt, Jean and Hutchings, Michael and Ramos,
              Vinicius G. B.},
     TITLE = {Examples around the strong {V}iterbo conjecture},
 BOOKTITLE = {Symplectic geometry---a {F}estschrift in honour of {C}laude
              {V}iterbo's 60th birthday},
     PAGES = {677--698},
 PUBLISHER = {Birkh\"{a}user{/}Springer, Cham},
      YEAR = {2022},
      ISBN = {978-3-031-19110-7},
       DOI={10.1007/s11784-022-00949-6},
}

@article{Sch18,
  author  = {Schlenk, Felix},
  title   = {Symplectic embedding problems, old and new},
  journal = {Bulletin of the American Mathematical Society},
  volume  = {55},
  number  = {2},
  pages   = {139--182},
  year    = {2018},
  doi     = {10.1090/bull/1587},
}

@article{HainePortaTeyssier2023,
  author        = {Haine, Peter J. and Porta, Mauro and Teyssier, Jean-Baptiste},
  title         = {The homotopy-invariance of constructible sheaves},
  journal       = {Homology, Homotopy and Applications},
  volume        = {25},
  number        = {2},
  pages         = {97--128},
  year          = {2023},
  doi           = {10.4310/HHA.2023.v25.n2.a6},
  eprint        = {2010.06473},
  archivePrefix = {arXiv},
  primaryClass  = {math.AT}
}

@article{Chaperon1984,
  author  = {Chaperon, Marc},
  title   = {Une id{\'e}e du type ``g{\'e}od{\'e}siques bris{\'e}es'' pour les syst{\`e}mes hamiltoniens},
  journal = {C. R. Acad. Sci. Paris S{\'e}r. I Math.},
  volume  = {298},
  number  = {13},
  pages   = {293--296},
  year    = {1984}
}

@book{RockafellarWets1998,
  author    = {Rockafellar, R. Tyrrell and Wets, Roger J.-B.},
  title     = {Variational Analysis},
  series    = {Grundlehren der mathematischen Wissenschaften},
  volume    = {317},
  publisher = {Springer},
  address   = {Berlin},
  year      = {1998},
  doi       = {10.1007/978-3-642-02431-3},
  isbn      = {978-3-540-62772-2}
}

@unpublished{liang2026monotone,
title={Symplectic capacities of 4-dimensional monotone toric domains}, 
      author={Yuanpu Liang},
      year={2026},
           note={\href{https://arxiv.org/abs/2611.14691}{arXiv:2611.14691}},
}

@article{Sikorav1987,
  author  = {Sikorav, Jean-Claude},
  title   = {Probl{\`e}mes d'intersections et de points fixes
             en g{\'e}om{\'e}trie hamiltonienne},
  journal = {Commentarii Mathematici Helvetici},
  volume  = {62},
  number  = {1},
  pages   = {62--73},
  year    = {1987}
}

@article{HaimKislevOstrover2023,
  author  = {Haim-Kislev, Pazit and Ostrover, Yaron},
  title   = {Remarks on symplectic capacities of {$p$}-products},
  journal = {International Journal of Mathematics},
  year    = {2023},
  volume  = {34},
  number  = {4},
  pages   = {2350021},
  doi     = {10.1142/S0129167X23500210}
}

@article{McDuffSchlenk2012,
  author  = {McDuff, Dusa and Schlenk, Felix},
  title   = {The embedding capacity of 4-dimensional symplectic ellipsoids},
  journal = {Annals of Mathematics},
  year    = {2012},
  volume  = {175},
  number  = {3},
  pages   = {1191--1282},
  doi     = {10.4007/annals.2012.175.3.5}
}

@article{Hausmann2022,
  author  = {Hausmann, Markus},
  title   = {Global group laws and equivariant bordism rings},
  journal = {Annals of Mathematics},
  year    = {2022},
  volume  = {195},
  number  = {3},
  pages   = {841--910},
  doi     = {10.4007/annals.2022.195.3.2}
}

@article{Quillen1969FormalGroupLaws,
  author  = {Quillen, Daniel},
  title   = {On the formal group laws of unoriented and complex cobordism theory},
  journal = {Bulletin of the American Mathematical Society},
  year    = {1969},
  volume  = {75},
  number  = {6},
  pages   = {1293--1298},
  doi     = {10.1090/S0002-9904-1969-12401-8},
}

@misc{HutchingsConvexYoutube,
  author = "Michael {Hutchings}",
  title= {Classification of Some Open Toric Domains},
  howpublished ={\url{https://www.youtube.com/watch?v=JGF20oz2CsQ}},
  year =     {2024},
}

@article{Schlenk2003,
  author  = {Schlenk, Felix},
  title   = {Symplectic embeddings of ellipsoids},
  journal = {Israel Journal of Mathematics},
  volume  = {138},
  year    = {2003},
  pages   = {215--252},
  doi     = {10.1007/BF02783427}
}

\end{document}